\documentclass[10pt,leqno]{article}

\usepackage[francais,english]{babel}
\usepackage{amsmath,amssymb,amscd}
\usepackage{a4wide}
\usepackage[T1]{fontenc}
\usepackage[latin1]{inputenc}
\usepackage[all]{xy} 
\usepackage{fancyhdr}
\usepackage{hyperref}
\usepackage{chngcntr}
\title{A spectral expansion for the symmetric space $\GL_n(E)/\GL_n(F)$} 
\author{Pierre-Henri Chaudouard}
\date{}

\newenvironment{paragr}[1][]{\refstepcounter{subsubsection} \noindent \textbf{\thesubsubsection . \ #1}}{\medskip}

\newenvironment{theoreme}{ \medskip\refstepcounter{theo}  \noindent\textbf{Theorem \thetheo}. ---\em}{\em \medskip}
\newenvironment{proposition}{\medskip\refstepcounter{theo}   \noindent\textbf{Proposition \thetheo}. ---\em}{\em\medskip}
\newenvironment{corollaire}{\medskip\refstepcounter{theo}  \noindent\textbf{Corollary \thetheo}. ---\em}{\em\medskip}

\newenvironment{lemme}{\medskip\refstepcounter{theo}   \noindent\textbf{Lemma \thetheo}. ---\em}{\em\medskip}

\newenvironment{preuve}[1][]{\noindent \textbf{Proof.} #1 --- }{\hfill
  \ensuremath{\square} \medskip}

\newenvironment{remarque}{\medskip\refstepcounter{theo}  \noindent\textbf{Remark \thetheo}. ---}{\medskip}

\DeclareMathOperator{\nd}{nd}

\DeclareMathOperator{\vol}{vol}

\DeclareMathOperator{\Ad}{Ad}

\DeclareMathOperator{\supp}{supp}

\DeclareMathOperator{\Hom}{Hom}

\DeclareMathOperator{\Int}{Int}
\DeclareMathOperator{\Ind}{Ind}

\DeclareMathOperator{\Id}{Id}

\DeclareMathOperator{\Res}{Res}
\DeclareMathOperator{\cusp}{cusp}
\DeclareMathOperator{\disc}{disc}
\DeclareMathOperator{\Ker}{Ker}

\DeclareMathOperator{\trace}{trace}

\newcommand{\ZZ}{\mathbb{Z}}

\newcommand{\SG}{\mathfrak{S}}
\newcommand{\NN}{\mathbb{N}}
\newcommand{\RR}{\mathbb{R}}
\newcommand{\AAA}{\mathbb{A}}
\newcommand{\CC}{\mathbb{C}}

\newcommand{\QQ}{\mathbb{Q}}

\newcommand{\ga}{\gamma}

\newcommand{\rc}{\mathcal{R}}
\newcommand{\Sc}{\mathcal{S}}
\newcommand{\tc}{\mathcal{T}}
\newcommand{\uc}{\mathcal{U}}

\newcommand{\ec}{\mathcal{E}}

\newcommand{\mc}{\mathcal{M}}

\newcommand{\lc}{\mathcal{L}}

\newcommand{\fc}{\mathcal{F}}
\newcommand{\pc}{\mathcal{P}}

\newcommand{\Ac}{\mathcal{A}}
\newcommand{\bc}{\mathcal{B}}

\newcommand{\Ic}{\mathcal{I}}
\newcommand{\Jc}{\mathcal{J}}

\newcommand{\zc}{\mathcal{Z}}

\newcommand{\ggo}{\mathfrak{g}}

\newcommand{\mgo}{\mathfrak{m}}

\newcommand{\ago}{\mathfrak{a}}
\newcommand{\cgo}{\mathfrak{c}}

\newcommand{\ugo}{\mathfrak{u}}

\newcommand{\Xgo}{\mathfrak{X}}

\newcommand{\GL}{\mathrm{GL}}

\newcommand{\al}{\alpha}
\newcommand{\be}{\beta}

\newcommand{\om}{\omega}
\newcommand{\Om}{\Omega}

\newcommand{\la}{\lambda}

\newcommand{\Ga}{\Gamma}
\newcommand{\back}{\backslash}

\newcommand{\Cc}{C_c^\infty}

\newcommand{\bg}{\langle}
\newcommand{\bd}{\rangle}

\newcommand{\eps}{\varepsilon}

\renewcommand{\leq}{\leqslant}
\renewcommand{\geq}{\geqslant}

\newcommand{\La}{\Lambda}

\newcommand{\bpar}{\begin{paragr}}
\newcommand{\epar}{\end{paragr}}

\begin{document}
\selectlanguage{english}
\maketitle
\counterwithin{equation}{subsubsection}
\begin{abstract}
  In this  article we  state and prove the spectral expansion of theta series attached to the symmetric space $\GL_n(E)/\GL_n(F)$ where $n\geq 1$ and $E/F$ is a quadratic extension of number fields. This is an important step towards the fine spectral expansion of relative trace formulas based on this symmetric space such as   the Jacquet-Rallis trace formula for general linear groups.  To obtain our result, we extend the work of Jacquet-Lapid-Rogawski on  intertwining periods to the case of discrete automorphic representations. The expansion we get is an absolutely convergent integral of relative characters built upon Eisenstein series and   intertwining periods.   We also establish a crucial but technical ingredient whose interest lies beyond  the focus of the article: we prove bounds for discrete Eisenstein series of $\GL_n$ on a neighborhood of the imaginary axis extending previous works of Lapid on cuspidal Eisenstein series. We even need a variant of such bounds on some shifts of the imaginary axis.
\end{abstract}



\tableofcontents

\section{Introduction}

\subsection{Motivations}

\begin{paragr}[The problem.] ---  Let $n\geq 1$ be an integer and let $E/F$ be a quadratic extension of number fields. Let $G=G_n$ be the group $\GL_E(n)$ viewed by restriction of scalars as an $F$-group. Let $\iota$ be the Galois involution and $G'=G_n'=\GL_F(n)$ be its fixed point subgroup. Let  $S$ be the symmetric space  $G/ G'$. Note that there is a left action of $G$ on $S$ by left translation. Let $\AAA$ be the ring of adèles of $F$. For any Schwartz function $\Phi$ on $S(\AAA)$ and $g\in G(\AAA)$ we can form the theta series
  $$\theta_\Phi(g)=\sum_{\sigma \in S(F)} \Phi(g^{-1}\sigma ).$$
  This gives a smooth function on $[G]=G(F)\back G(\AAA)$ and the purpose of the article is to provide a spectral decomposition of $\theta_\Phi(g)$ in terms of objects attached to the $L^2$-automorphic spectrum of $G$.  This problem is raised by Jacquet in \cite {Jac-Edin}. We can slightly restate the problem noting that we have $S(\AAA)=G(\AAA)/ G'(\AAA)$.   Any Schwartz function $\Phi$ can be obtained from $f$ in the space $\Sc(G(\AAA))$ of Schwartz functions  on $G(\AAA)$ in the following way:
  $$\Phi(\sigma)=\int_{G'(\AAA)} f(gh)\, dh$$
  where $dh$ is some Haar measure on $G'(\AAA)$ and $\sigma \in S(\AAA)$ is represented by the class $gG'(\AAA)$ for some  $g\in G(\AAA)$. Then we introduce the  automorphic kernel 
  \begin{align*}
    K_f(x, y)=\sum_{\gamma\in G(F)} f(x^{-1}\gamma y), \ \ x,y\in G(\AAA)
  \end{align*}
  and the so-called  Flicker-Rallis period of the automorphic kernel defined by the convergent integral
\begin{align*}
    J^G(g, f)=\int_{[G']} K_f(g, h)\, dh, \ \ g\in G(\AAA)
  \end{align*}
  where $[G']=G'(F)\back G'(\AAA)$ and $dh$ is the quotient measure.   By the equality $S(F)=G(F)/G'(F)$ we get $\theta_\Phi(g)=J^G(g,f)$ and the problem becomes to decompose spectrally the distribution $f\mapsto J^G(g,f)$.
\end{paragr}

\begin{paragr}[Jacquet-Rallis trace formula.] --- Let $\det:G'\to \mathbb{G}_{m,F}$ be the determinant and let $\eta$ be the quadratic character of $\AAA^\times$ attached to $E/F$ by class field theory. We can  consider the twisted  Flicker-Rallis period of the automorphic kernel defined by
   \begin{align*}
J^{G,\eta}(g, f)=\int_{[G']} K_f(g, h)\eta(\det(h))^{n+1}\, dh.
  \end{align*}
Note that our results and methods give also the spectral decomposition of this distribution. We shall leave the reader make the obvious modifications. The Jacquet-Rallis trace formula for linear groups as stated in \cite{Z3} is the regularized version of the (in general divergent) integral
  \begin{align*}
    \int_{[G_n]} J^{G_n,\eta}(g,f_n)J^{G_{n+1},\eta}(g,{f_{n+1}})\,dg
  \end{align*}
  for  Schwartz functions $f_n$ and $f_{n+1}$ respectively on $G_n(\AAA)$ and  $G_{n+1}(\AAA)$. Here, denoting by $(e_1,\ldots,e_{n+1})$ the canonical basis of $F^{n+1}$, we have identified $G_n$ with the subgroup of $G_{n+1}$  preserving the $E$-subspace of $E^{n+1}$ generated $(e_1,\ldots,e_n)$ and fixing the last vector $e_{n+1}$. Recall that the  Jacquet-Rallis trace formula plays  a central role in the recent proof of Gan-Gross-Prasad and Ichino-Ikeda conjectures for unitary groups, see among others \cite{Zhang2,BPLZZ,BCZ,BPC}. The more spectral contributions we understand, the more results we can extract from the (comparison of) relative trace formulae.  To any cuspidal datum $\chi$ of $G$ we can attach the $\chi$-component $K_\chi$ of the kernel. By integrating it over $[G']$ we get the distribution  $J^{G,\eta}_\chi(g, f)$. A key ingredient  in \cite{BCZ} is the spectral decomposition of $J^{G,\eta}_\chi(g, f)$ for some specific cuspidal data $\chi$ called $*$-regular, see \cite[theorem 4.3.3.1]{BCZ}. We shall provide in the article an extension of this decomposition to any cuspidal datum. Thus the article  solves an important  step towards the fine spectral decomposition of the Jacquet-Rallis trace formula from which we expect new applications; for a recent application of some of the results of the article see  \cite{BPC}.  More generally, the results presented here should be useful for the spectral expansion of any trace formula based on the symmetric space $S$. For example,  from our main result  theorem \ref{thm-intro:spec-kernel} below, it is possible not only to get a proof of  the fine spectral expansion of Jacquet's relative trace formula which is different from the proof of Lapid given in \cite[theorem 10.1]{LapFRTF} but also to get an explicit computation of the constants $\cgo(M,\pi)$ in \cite[theorem 10.1]{LapFRTF}.
\end{paragr}

\subsection{Statement of the main result}

\begin{paragr} In the following we fix $g\in G(\AAA)$. As usual, all parabolic subgroups of $G$ are assumed to be defined over $F$. They are standard, resp. semi-standard,  if they contain the standard Borel of $G=\GL_E(n)$, resp. the group of diagonal matrices, viewed as an $F$-group. The Levi subgroups of $G$ are the Levi factors defined over $F$ of the parabolic subgroups of $G$.

  Let $P=MN_P$ be a standard parabolic subgroup of $G$ with its standard Levi decomposition, $N_P$ being its unipotent radical. Let $\Pi_{\disc}(M)$ be the  set of discrete automorphic representations of $M(\AAA)$ with central character trivial on the central subgroup $A_M^\infty$, see §§ \ref{S:Haar-meas} and \ref{S:disc-aut}. Let $\pi\in \Pi_{\disc}(M)$ and let $\Ac_{P,\pi}(G)$ be the space of (smooth) automorphic functions on the quotient $A_M^\infty M(F)N_P(\AAA)\back G(\AAA) $ that belong to the ``$\pi$-component'', see § \ref{S:disc-aut}. Let $\la$ be an element in the $\CC$-vector space $\ago_{P,\CC}^{G,*}$  of complex unramified characters of $P(\AAA)$ trivial on $A_G^\infty$. We have a map $\varphi \mapsto \varphi_\la$ that identifies $\Ac_{P,\pi}(G)$ with the induced representation $I_P^G(\pi\otimes \la)$. Then $f\in \Sc(G(\AAA))$ acts on $I_P^G(\pi\otimes \la)$ and thus on $\Ac_{P,\pi}(G)$  by transport of structure. We denote by  $I(\la,f)$ the action we get, see  § \ref{S:disc-aut}. 
\end{paragr}

\begin{paragr}[Intertwining periods.]  --- We define a subset denoted by $\lc_2(M)$ of the set of  Levi subgroups of $G$ that contain $M$, see § \ref{S:fc2}. On the space $\Ac_{P,\pi}(G)$ we shall consider several linear forms $J_L(\varphi,\la)$ attached to Levi subgroups $L\in \lc_2(M)$ and $\la\in \ago_{M,\CC}^{L,*}$.   These are the so-called  intertwining periods $J_L(\varphi,\la)$ introduced in \cite[section VII]{JLR} when $\pi$ is cuspidal. We need to extend their definition to the discrete case. For the introduction we shall give the definition only when $L$ is moreover standard, for the general case see section \ref{sec:inter-period}.    One can attach to $L$  a permutation matrix $\xi$ of order $2$. Recall $\iota$ is the Galois involution of $G$. By the choice of an element $\tilde\xi\in G(F)$ such that  $\tilde\xi \iota(\tilde\xi)^{-1}=\xi$, one can identify the subgroup of $P$ fixed by the involution $\Int(\xi)\circ\iota$ with a subgroup $P_\xi$ of $G'$ which comes with a natural Levi decomposition $P_\xi=M_\xi N_\xi$, see § \ref{S:Jxitilde} (note that we slightly simplify the notation in the introduction). Then the intertwining period is given by
  $$J_L(\varphi,\la)=\int_{  N_{\xi} (\AAA) A_{M_{\xi}}^\infty  M_{\xi}(F) \back G'(\AAA)} \varphi_\la(\tilde\xi h)\,dh.$$
Here  the integral is convergent for $\la$ in some cone, see proposition \ref{prop:convergence}, and admits a meromorphic continuation in general, see corollary \ref{cor:period-entrelac}.
  \end{paragr}

  \begin{paragr}[Relative characters.]  --- The relative characters we consider are distributions on  $\Sc(G(\AAA))$ attached to $\pi$,  $L\in \lc_2(L)$ and $\la\in i\ago_{M}^{L,*}$ that is $\la$ is a unitary character. They are given by:
    \begin{align}
      \label{eq:intro-relchar}
      \Jc_{L,\pi}(g,f,\la)= \sum_{\varphi}   E(g,I(\la,f)\varphi,\la) \overline{ J_L(\varphi,\la)} 
    \end{align}
    where the sum is over some Hilbert basis of $\Ac_{P,\pi}(G)$  for the Petersson inner product (see subsection \ref{ssec:comput-limit}) and where $E(g,\varphi,\la)$ denotes an  Eisenstein series, see § \ref{S:Eis-series}. Note that whereas $E(g,I(\la,f)\varphi,\la)$ is holomorphic on  $i\ago_{M}^{L,*}$ the intertwining period $J_L(\varphi,\la)$ may have singularities which are described through theorems \ref{thm:period-entrelac} and \ref{thm:singularities}. However it turns out that the product $E(g,I(\la,f)\varphi,\la) \overline{ J_L(\varphi,\la)}$ is smooth on $ i\ago_{M}^{L,*}$. Moreover the relative character $\Jc_{L,\pi}(g,f,\la)$ is smooth in $\la\in  i\ago_{M}^{L,*}$ and continuous on  $\Sc(G(\AAA))$, see proposition \ref{prop:relch-JL}.
    \end{paragr}

  \begin{paragr}[Main theorem.] --- The main result of the paper is the following theorem. 

      \begin{theoreme} \label{thm-intro:spec-kernel} (see theorem \ref{thm:spec-kernel})
 For any $f\in \Sc(G(\AAA))$, we have
 \begin{align}\label{eq:intro-thm}
      J^G(g,f)= \sum_{(M,L,\pi)}c_M^L \int_{i\ago_M^{L,*}} \Jc_{L,\pi}(g,f,\la)\, d\la.
    \end{align}
The sum above is absolutely convergent. 
  \end{theoreme}

  In  \eqref{eq:intro-thm} the sum is  over  triples $(M,L,\pi)$ where $M$ is a standard  Levi  subgroup,  $L\in \lc_2(M)$ and $\pi\in \Pi_{\disc}(M)$. We also set
$$c_M^L= |\pc(M)|^{-1}   2^{-\dim(\ago_L)} $$
where $\pc(M)$ is the finite set of parabolic subgroups admitting $M$ as a Levi factor and $|\pc(M)|$ is its cardinality. Let $\chi$ be a cuspidal datum of $G$. One obtains the spectral decomposition of $J^G_\chi(g,f)$ by restricting  the sum to triples $(M,L,\pi)$ such that $\pi$ belongs to the subset $\Pi_{\chi}(M)\subset \Pi_{\disc}(M)$ defined in § \ref{S:Pichi}.

The discrete part of the spectral contribution is given by
 \begin{align*}
      J^G_{\disc}(g,f)= \sum_{(M,\pi)}c_M^M  \Jc_{M,\pi}(g,f,0)
    \end{align*}
    where the sum is over pairs $(M,\pi)$ where $M$ is a standard  Levi  subgroup and $\pi\in \Pi_{\disc}(M)$. The permutation matrix is trivial in this case and the relative character takes the simple form:
     \begin{align*}
      \Jc_{M,\pi}(g,f,0)= \sum_{\varphi}   E(g,I(\la,f)\varphi,0) \int_{  N' (\AAA) A_{M'}^\infty  M'(F) \back G'(\AAA)} \overline{\varphi( h)}\,dh
     \end{align*}
     where $P'=P\cap G'$ and $M'N'$ is its Levi decomposition. In particular, the contribution of $(M,\pi)$ vanishes unless $\pi$ is $M'$-distinguished in the sense of § \ref{S:distinction} below. The classification of $M'$-distinguished discrete representations in terms of  distinguished cuspidal  representations  is known by the work of Yamana, see  \cite[theorem 1.2]{Yquad}. We will not explicitly deal here with the question of the non-vanishing of the relative characters $\Jc_{L,\pi}(g,f,\la)$. Let us mention that the intertwining periods $J_L(\varphi,\la)$ are related to some intertwining periods of cuspidal representations by proposition \ref{prop:etape2}. The latter are related to  integrals of Whittaker functions, see  \cite[theorem 5.5.1.1]{BPC}, and this gives a way to study the non-vanishing of  the intertwining periods.
\end{paragr}

\subsection{Bounds for discrete Eisenstein series}

\begin{paragr}
  A technical but crucial point in the proof of theorem \ref{thm-intro:spec-kernel} is the  majorization of Eisenstein series for the group $G$ in the neighborhood of the imaginary axis. For cuspidal Eisenstein series the bounds are those introduced by Lapid in \cite{LapHC} and \cite{LapFRTF}. Here we generalize the bounds to the case of  discrete Eisenstein series, namely Eisenstein series built from a discrete automorphic representation.

  \begin{theoreme} \label{thm-intro:maj-Eis}(for a stronger statement see theorem \ref{thm:maj-Eis})
   There exists  $N>0$ such that for all $q>0$ there is a continuous semi-norm $\|\cdot\|$  on $\Sc(G(\AAA))$ such that for all $f\in \Sc(G(\AAA))$, all standard parabolic subgroups $P=MN_P$,  $\pi\in \Pi_{\disc}(M)$, $\la\in i\ago_{P}^{G,*} $ and $x\in G(\AAA)^1$  we have
    \begin{align*}
     \sum_{\varphi}  | E(x,I(\la,f)\varphi,\la) |^{2}  \leq \frac{  \|x\|_{G}^N\|f\|^2    }{(1+\|\la\|^2)^q(1+\|\chi_\pi\|^2)^{q}}
    \end{align*}
   where the sum is over some Hilbert basis of $\Ac_{P,\pi}(G)$.
  \end{theoreme}

  Here $\|\chi_\pi\|$ denotes the norm of the infinitesimal character $\chi_\pi$ of the Archimedean component of $\pi$. As usual, we view $\chi_\pi$ as an orbit under the Weyl group in the dual of a Cartan subalgebra of the complexified Lie algebra of $G(F_\infty)$ where $F_\infty$ is the product of the Archimedean completions of $F$. The norm is Euclidean and invariant under the Weyl group. The other notations are explained in subsections \ref{ssec:general-notations} and \ref{ssec:autom}. For simplicity we have stated the theorem  for $\la$ on the imaginary axis. However we shall need and prove  a stronger version where the bound holds on a neighborhood of the imaginary axis depending on $\pi$. Note also we use in the statement of theorem \ref{thm:maj-Eis} a slightly different invariant $\La_\pi$. The link  with the formulation above is given in remark \ref{rq:La invariant}.

  We mention that we also  need a variant of theorem \ref{thm-intro:maj-Eis} where we bound a sum of Eisenstein series on a shift  of the space $i\ago_{P}^{G,*}$, see theorem \ref{thm:maj-EisR}.
\end{paragr}

\begin{paragr}
  To get theorem  \ref{thm-intro:maj-Eis} we proceed as in  \cite{LapHC} and \cite{LapFRTF}, namely we majorize the Eisenstein series by the Petersson norm of a truncated Eisenstein series. Then we need a new ingredient namely the explicit computation of the truncated product of two Eisenstein series. This is provided by theorem \ref{thm:maj-Eis}. By explicit, we mean a combinatorial expression that involves only intertwining operators and Petersson products of discrete automorphic representations. Thus theorem \ref{thm:maj-Eis} both generalizes the classical statement for cuspidal Eisenstein series and the well-known Arthur's asymptotic formula of \cite{ar-truncated}.  The proof of the theorem is closely related to the works \cite{JLR} of  Jacquet-Lapid-Rogawski  and \cite{Linner} of Lapid. More precisely we use the fact that several intricate expressions provide  families of meromorphic invariant bilinear forms on some induced representations. By an observation due to Bernstein they must vanish. We also rely on the precise computation of the exponents of the discrete automorphic representations based on their description by M\oe glin-Waldspurger in \cite{MW} from which we extract some geometric properties, see lemma \ref{lem:negativity}. To conclude we need to establish  bounds on intertwining operators in some neighborhood of the imaginary axis shifted by such exponents, see proposition \ref{prop:bound-intertwining}.
\end{paragr}

\subsection{About the proof of the main theorem}

\begin{paragr}
  The starting point to get the spectral expansion of $J^G(g,f)$ is to approximate it by
  \begin{align}\label{eq-intro:JGT}
    J^{G,T}(g, f)=\int_{[G']} (K_f\La^T_m) (g, h)\, dh \ \ g\in G(\AAA).
  \end{align}
  Here $T$ is a parameter and $\La^T_m$ is the ``mixed'' (as opposed to Arthur's one) truncation operator introduced by  Jacquet-Lapid-Rogawski  in \cite{JLR}. When applied to the map $x\in G(\AAA) \mapsto K_f(g,x)$ we get the expression  $(K_f\La^T_m) (g, \cdot)$. On the one hand the limit  of $ J^{G,T}(g, f)$ when $T$ goes to infinity  is $J^G(g,f)$. On the other hand, we can  easily get the spectral expansion of  $J^{G,T}(g, f)$ using the Langlands spectral decomposition of the kernel $K_f$, see proposition \ref{prop:cv-spectral-T}:
 \begin{align*}
         J^{G,T}(g,f)= \frac12 \sum_{ (M,\pi)} |\pc(M)|^{-1} \int_{i\ago_M^{G,*}}    \ec_{\pi}^T(g,f,\la) \, d\la.
       \end{align*}
       The sum above is over pairs $(M,\pi)$ where $M$ is a standard Levi subgroup and $\pi\in \Pi_{\disc}(M)$. The relative character here is  defined by:
       \begin{align*}
        \ec_{\pi}^T(g,f,\la)= \sum_{\varphi}   E(g,I(\la,f)\varphi,\la) \overline{ \int_{[G']_0}     \La^T_m E(\varphi,\la) }.
       \end{align*}
       The main difference with \eqref{eq:intro-relchar} is that we have  replaced the intertwining period by the   Flicker-Rallis period of the mixed truncated discrete Eisenstein series.   The study of such periods is undertaken is section \ref{sec:FR}. The starting point is that the Flicker-Rallis period of a  truncated discrete Eisenstein series can be expressed in terms of Jacquet-Lapid-Rogawski regularized periods of Eisenstein series, see proposition \ref{prop:MS-periods}. Let us mention that the proposition is not completely formal: once again we rely on the explicit description of exponents of discrete automorphic representations. As in Arthur's asymptotic formula for the scalar product of two Eisenstein series, some of the contributions are negligible when $T$ goes to infinity. Among the main terms, many regularized periods of Eisenstein series vanish, see proposition \ref{prop:vanishing}. When they do not obviously vanish, they are related to the intertwining periods, see theorem \ref{thm:period-entrelac}.   At this point, we can show that the Flicker-Rallis period of a  truncated discrete Eisenstein series is asymptotic to an explicit combination of intertwining periods. Using the stronger version of theorem \ref{thm-intro:maj-Eis} from which we borrow the notations one can show:

       \begin{proposition}(see  proposition \ref{prop:maj-rel-char3} \label{prop-intro:maj-rel}for a stronger and more precise statement.)  There exists  $\eps>0$ such that for all $q>0$  there exists  a continuous semi-norm $\|\cdot\|$  on $\Sc(G(\AAA))$ such that for all pairs $(M,\pi)$ as above, all $f\in \Sc(G(\AAA))$, all $\la\in i\ago_{M}^{ G,*}$ 
    \begin{align*}
      |\ec_{\pi}^T(g,f,\la)- \sum_{Q \in \fc_2(M) } 2^{-\dim(\ago_Q^G)} \Jc_{Q,\pi}(g,\varphi,\la) \frac{\exp(-\bg \la, T_Q^G\bd)}{\theta_Q(-\la)}|  \leq \frac{  \exp(-\eps \| T\| )\|f\|  }{(1+\|\la\|^2)^q(1+\|\chi_\pi\|^2)^{q}}
          \end{align*}
for all $T$ suitably regular.         
       \end{proposition}

       Here the set  $\fc_2(M)$ is the set of semi-standard parabolic subgroups whose semi-standard Levi factor belongs to $\lc_2(M)$ and $ \Jc_{Q,\pi}(g,\varphi,\la) $ is yet another relative character whose definition is given in \eqref{eq:JcQpi}. Moreover $\theta_Q$ is the familiar polynomial from Arthur's theory, see § \ref{S:Haarmeasures}.

       For any $L\in \lc_2(M)$ the family $ (\Jc_{Q,\pi}(g,\varphi,\la))_{Q\in \pc(L)}$ indexed by the set of parabolic subgroups of Levi $L$ is a $(G,L)$-family (in the sense of Arthur) of Schwartz functions on $i\ago_M^{G,*}$. This has several consequences. First  the value $\Jc_{Q,\pi}(g,\varphi,\la)$ on $i\ago_M^{L,*}$ does not depend on $Q\in \pc(L)$ and is  in fact equal to  $\Jc_{L,\pi}(g,\varphi,\la)$. Second the expression 
\begin{align*}
  \sum_{Q \in \pc(L) } \Jc_{Q,\pi}(g,\varphi,\la) \frac{\exp(-\bg \la, T_Q^G\bd)}{\theta_Q(-\la)}
\end{align*}
 defines a Schwartz function on $i\ago_M^{G,*}$  and we have
  \begin{align*}
    \lim_{T\to +\infty} \int_{i\ago_M^{G,*}}    \sum_{Q \in \pc(L) } \Jc_{Q,\pi}(g,\varphi,\la) \frac{\exp(-\bg \la, T_Q^G\bd)}{\theta_Q(-\la)} \, d\la=\int_{i\ago_M^{L,*}} \Jc_{L,\pi}(g,\varphi,\la)\,d\la.
  \end{align*}
This gives theorem \ref{thm-intro:spec-kernel}.
  \end{paragr}

  \subsection{Organization of the paper}

  \begin{paragr}
    In section \ref{sec:prelim}, we collect the  notations and some elementary lemmas about sets of inversions of Weyl group elements, polynomial exponential maps and $(G,M)$-families.
  \end{paragr}

  \begin{paragr}
    The section 3 is devoted to the proof of  theorem \ref{thm-intro:maj-Eis} above. First in subsection \ref{ssec:truncated-scalar} we obtain an explicit expression for the scalar product of two truncated Eisenstein series, see theorem \ref{thm:inversion-calculee}. 
    In subsection \ref{ssec:bound-intert}, a bound for intertwining operators is stated, see proposition \ref{prop:bound-intertwining}. The next four subsections are devoted to the proof of  proposition \ref{prop:bound-intertwining}. The proof is divided into several tasks: after the introduction of normalization factors in subsection \ref{ssec:normalization inter} we establish a bound for them in subsection \ref{ssec:bd normalization}. Then we  show in subsection \ref{ssec:holom norm} that  the normalized  intertwining operators are holomorphic in some region of interest. Last we get a bound for these operators in subsection \ref{ssec:bd norm inter}. 
    In subsection \ref{ssec:bound-scalar} following a work of Lapid, we use the expression of the scalar product of two truncated Eisenstein series to get a bound on this scalar product, see proposition \ref{prop:bound-scalaire}. In subsection \ref{ssec:bds-herm}, the next step is to deduce a  bound for some relative characters, see proposition \ref{prop:maj-trace-rel2}. Finally in subsection \ref{ssec:bds-Eis} we state and prove theorems \ref{thm:maj-Eis} and \ref{thm:maj-EisR} which are stronger versions of   theorem \ref{thm-intro:maj-Eis} above.
  \end{paragr}

  \begin{paragr}
    The section \ref{sec:FR} is devoted to the study of truncated Flicker-Rallis periods of discrete Eisenstein series and their regularized versions. The definition of the mixed truncation operator of Jacquet-Lapid-Rogawski is recalled in subsection \ref{ssec:mixed}. The subsection \ref{ssec:reg-period-Eis} gives the definition (after Jacquet-Lapid-Rogawski) and the first properties of regularized periods of discrete Eisenstein series which are meromorphic functions defined in terms of truncated Flicker-Rallis periods. In subsection \ref{ssec:asym-period}, the process is reversed and the truncated Flicker-Rallis period of a discrete Eisenstein series is expressed in terms of regularized periods of Eisenstein series, see proposition \ref{prop:MS-periods} and corollary \ref{cor:MS-periods}. In this expansion, we can distinguish the ``main part'' \eqref{eq:PTRm} and we establish smoothness properties of it in proposition \ref{prop:type1}. This will be the main part when the parameter goes to infinity. In subsection \ref{ssec:sing-period} we show that most of the regularized periods of discrete Eisenstein series vanish. Then we determine their possible singularities on the imaginary axis, see theorem \ref{thm:singularities}.
  \end{paragr}

  \begin{paragr}
    The section \ref{sec:inter-period} extends the definition of intertwining periods, also due to Jacquet-Lapid-Rogawski, to the case of discrete automorphic representation.  The main results are stated in subsection \ref{ssec:an-cont}. For  some cones,  intertwining periods are given by a quite simple convergent integral see \eqref{eq:defJQ} and proposition \ref{prop:convergence}. In theorem \ref{thm:period-entrelac} we show that on the convergence domain they coincide with regularized periods of discrete Eisenstein series. As a consequence they admit a meromorphic continuation, see corollary \ref{cor:period-entrelac}. The bulk of the  proof of theorem  \ref{thm:period-entrelac}  is given in subsection \ref{ssec:comp-regper}. Note that it relies on  Jacquet-Lapid-Rogawski results on  intertwining periods of cuspidal representation. A crucial step is in fact to relate the  intertwining period of a discrete representation to that of the cuspidal representation that is part of M\oe glin-Waldspurger classification, see proposition \ref{prop:etape2}.
  \end{paragr}

  \begin{paragr}
    The goal of section \ref{sec:rel-char} is to introduce majorizations of various relative characters that depend on a truncation parameter and that will play a role in spectral decomposition. The first example given in subsection \ref{ssec:EX1} is built upon Flicker-Rallis periods of truncated discrete Eisenstein series. In subsection \ref{ssec:2nd-ex}, a second example is introduced built upon the ``main part''  defined in  subsection \ref{ssec:asym-period}. Moreover proposition \ref{prop:maj-rel-char3} shows that the two examples coincide asymptotically, see also proposition \ref{prop-intro:maj-rel} above.
  \end{paragr}

  \begin{paragr}
    The final section \ref{sec:Spec-exp} states and proves the spectral expansion namely  theorem \ref{thm:spec-kernel} (and theorem \ref{thm-intro:spec-kernel} above). The starting point given in subsection \ref{ssec:asalim} is the spectral decomposition \eqref{eq-intro:JGT} given by  proposition \ref{prop:cv-spectral-T} and proposition \ref{prop:limT} which reduces the problem to the computation of the limit of integrals of the relative characters associated to the main part. The computation of the limit is achieved in subsection \ref{ssec:comput-limit} and theorem \ref{thm:Jchi} which gives the spectral expansion of $J_\chi^G(g,f)$ for a cuspidal datum $\chi$. Finally in subsection \ref{ssec:spec-FR} the spectral decomposition is stated and an argument based on Müller's work is given for the absolute convergence.
  \end{paragr}

  \begin{paragr}[Acknowledgement.] ---  I thank the Institut universitaire de France (IUF) for its support. The impetus of this work was given by some joint projects we had with Micha\l\    Zydor. I warmly thank him  for his insights  and the countless discussions on regularized periods we had in past years.
    I am greatly indebted to the anonymous referee  who observed that some shift in the expression of the truncated scalar product and  the definition of regularized  Flicker-Rallis period  was overlooked in a previous version of the paper. I would like to thank him for his careful reading and his many useful comments that helped me to improve the exposition. I thank Raphaël Beuzart-Plessis and Paul Boisseau who also pointed out to me the missing shift. Finally I thank François Digne and Jean Michel for a useful conversation.
  \end{paragr}
  
\section{Preliminaries}\label{sec:prelim}

\subsection{General notations}\label{ssec:general-notations}

\begin{paragr}
Let $ F $ be a field of characteristic $0$. 
\end{paragr}

\begin{paragr}
Let $ G $ be a linear algebraic group defined over $ F $. Let  $N_G$ be the unipotent radical of $G$  and $ X^*(G) $ be the group of  algebraic morphisms  $ G \to \mathbb{G}_{m,F}$ defined over $ F $. Let $\ago_G^*=X^*(G)\otimes_\ZZ\RR$ and $\ago_G=\Hom_\ZZ(X^*(G),\RR)$. We denote by 
\begin{align}\label{eq:pairing}
\bg \cdot,\cdot\bd : \ago_G^*\times \ago_G\to \RR
\end{align}
the canonical pairing.
\end{paragr}

\begin{paragr} Let us  assume that $ G $ is moreover reductive. We shall recall  Arthur's notations. Let $ P_0 $ be a parabolic subgroup of $ G $ defined over $ F $  and minimal for these properties. Let $ M_0 $ be a Levi factor of $ P_0$ defined over  $ F $. A parabolic  (resp. and semi-standard, resp. and standard) subgroup of $G$ is a parabolic subgroup of $G$ defined over $ F $ (resp.  which contains $ M_0 $, resp. which contains $ P_0 $). For any semi-standard parabolic subgroup $P$, we have a Levi decomposition $ P = M_P N_P $ where $ M_P $ is the unique Levi factor that contains $ M_0 $.  A Levi subgroup of $ G $ (resp. semi-standard,  resp. standard) is a  Levi factor defined over $ F $ of a parabolic subgroup of $G$ (resp. semi-standard, resp.  standard).
\end{paragr}

\begin{paragr}\label{S:AP}
  Let  $A_G$ be the maximal central $F$-split torus of $G$.  For any semi-standard parabolic subgroup $P$ of $G$, we set  $A_P=A_{M_P}$. The restrictions maps  $X^*(P)\to X^*(M_P)\to X^*(A_P)$ induce isomorphisms $\ago_P^*\simeq \ago_{M_P}^* \simeq \ago_{A_P}^*$.
  We set $\ago_0^*=\ago_{P_0}^* $,  $\ago_0=\ago_{P_0}$  and $A_0=A_{P_0}$.
\end{paragr}

\begin{paragr}\label{S:proj} For any semi-standard parabolic subgroups $P\subset Q$ of $G$, the restriction map $X^*(Q)\to X^*(P)$ induces maps  $\ago_Q^*\to \ago_P^*$ and $\ago_P\to \ago_Q$. The first one is injective whereas the kernel of the second one is denoted by $\ago_P^Q$. The restriction map $X^*(A_P)\to X^*(A_Q)$ gives a surjective map  $\ago_P^*\to \ago_Q^*$ whose kernel is denoted by  $\ago_P^{Q,*}$. We get also an  injective map $\ago_Q\to \ago_P$. In this way, we get dual decompositions $\ago_P=\ago_P^Q\oplus \ago_Q$ and $\ago_P^*=\ago_P^{Q,*}\oplus \ago_Q^*$. Thus we have  projections $\ago_0\to \ago_P^Q$  (resp. $\ago_P$) and $\ago_0^*\to \ago_P^{Q,*}$  (resp. $\ago_P^*$) denoted by $X\mapsto X_P^Q$ (resp. $X_P$). They depend only on the Levi factors $M_P$ and $M_Q$ and they will be also denoted by  $X\mapsto X_{M_P}^{M_Q}$, resp. $X_{M_P}$.

We denote by $\ago_{P,\CC}^{Q,*}$ and $\ago_{P,\CC}^Q$ the $\CC$-vector spaces obtained by extension of scalars from $\ago_{P}^{Q,*}$ and  $\ago_{P}^{Q}$. We still denote by $\bg\cdot,\cdot\bd$ the pairing we get from \eqref{eq:pairing}  by extension of  scalars to $\CC$. We have a decomposition
$$\ago_{P,\CC}^{Q,*}=\ago_{P}^{Q,*}\oplus i\ago_{P}^{Q,*}$$
where $i^2=-1$.  We shall denote by $\Re$ and $\Im$ the real and imaginary  parts associated to this decomposition and by $ \bar{\la}$ the  complex conjugate   of $\la\in \ago_{P,\CC}^{Q,*}$.
\end{paragr}

\begin{paragr} \label{S:root-coroot} For any parabolic subgroup $P$ of $G$ containing $P_0$, we denote by  $\Delta_0^P=\Delta_{P_0}^P$  the set of simple roots of  $A_{0}$ in  $M_P\cap P_0$.   It is a subset of the set $\Phi_{0}^P=\Phi_{P_0}^P$ of roots of $A_{0}$ in  $M_P\cap P_0$. If moreover $Q$ is a parabolic subgroup containing $P$ we denote by $\Delta_P^Q$ be the image of $\Delta_{0}^Q\setminus \Delta_{0}^P$ (viewed as a subset of $\ago_0^*$) by the projection $\ago_{0}^*\to \ago_P^*$.  Note that $\Delta_P^Q$ is a subset of the set $\Phi_P^Q$ of roots of $A_P$ in $M_Q\cap P$. In the same way  one defines the sets $\Delta_P^{Q,\vee}\subset \Phi_P^{Q,\vee}\subset \ago_P^{Q}$ of simple coroots and coroots. The sets $\Delta_P^Q$ and $\Delta_P^{Q,\vee}$ are respectively bases of $\ago_P^{Q,*}$ and $\ago_P^{Q}$. Taking the respective dual bases,  we get the  sets  $\hat{\Delta}_P^{Q,\vee}$ and  $\hat{\Delta}_P^Q$ respectively of simple coweights and simple weights.   The sets $\Delta_P^Q$ and  $\hat{\Delta}_P^Q$ determine open  cones  in  $\ago_{0}$ whose characteristic functions are denoted respectively by $\tau_P^Q$ and  $\hat{\tau}_P^Q$. As usual we omit the exponent $Q$ if $Q=G$.
\end{paragr}

\begin{paragr}[Euclidean norms.] ---   Let $W$ be the Weyl group of $(G,A_0)$ that is the  quotient by $M_0(F)$ of the normalizer of  $A_0$ in $G(F)$. The group acts on $\ago_0$ and its dual $\ago_0^*$. We fix once and for all an invariant inner product on $\ago_0^*$. Then we can identify $\ago_0^*$ with its dual $\ago_0$ and all the decompositions we introduced are orthogonal for the inner product. We denote by $\|\cdot\|$ the Euclidean norm.
\end{paragr}

\begin{paragr}\label{S:Haarmeasures} Let $P$ be a semi-standard parabolic subgroup of $G$. We equip $\ago_P$ with the Haar measure that gives a covolume $1$ to the lattice $\Hom(X^*(P),\ZZ)$. The space  $i\ago_{P}^{*}$ is then equipped with  the dual Haar measure so that we have 
\begin{align*}
\int_{i\ago_{P}^*} \int_{\ago_{P}}  \phi(H) \exp (-\bg \la ,H\bd)\,dHd\la=\phi(0) 
\end{align*}
for all $\phi\in \Cc(\ago_{P})$.

For any basis $B$ of $\ago_P^Q$ let us denote by  $\ZZ(B)$ the lattice generated by $B$ and by  $\vol(\ago_P^Q/\ZZ(B))$ the  covolume of this lattice where  $\ago_P^Q\simeq \ago_P/\ago_Q$ is  provided with the quotient Haar measure. We have on $\ago_0^*$ the  polynomial functions:
$$\hat{\theta}_P^Q(\la)= \vol(\ago_P^Q/\ZZ(\hat{\Delta}_P^{Q,\vee}))^{-1} \prod_{\varpi^\vee \in \hat{\Delta}_P^{Q,\vee}}\bg \la,\varpi^\vee\bd
$$
and
$$\theta_P^Q(\la)= \vol(\ago_P^Q/\ZZ(\Delta_P^{Q,\vee}))^{-1} \prod_{\al \in \Delta_P^Q} \bg \la,\al^\vee\bd.
$$
\end{paragr}

\begin{paragr}[Truncation parameter.] ---\label{S:trunc-param}
  For any $T\in \ago_P^Q$, we set 
  \begin{align*}
    d_P^Q(T)=\inf_{\al\in \Delta_P^Q} \bg \al,T\bd.
  \end{align*}
Let $\ago_P^{Q+}=\{T\in \ago_P^Q \mid \   d_P^Q(T) \geq 0\}$.
If $Q=G$, the exponent $G$ is omitted. We set $d(T)=d_{P_0}^G(T)$.

Let $T\in \ago_0$ such that $d(T)\geq 0$. We shall say that  $T$ is enough positive (sufficiently regular in the terminology of \cite[p. 937]{ar1}) if $d(T)$ is large enough. Often we simply say $T$ is a truncation parameter. Since most of the constructions we are interested in do not depend on the choice of $T$, the  precise lower bound is irrelevant. For the constructions that really depend on $T$, we will be in fact interested in their asymptotic behaviour when $d(T)\to +\infty$.

For any semi-standard parabolic subgroup $P$ and any point $T\in \ago_0$,   we define a point $T_P\in \ago_{P}$ in the following way : we choose $w\in W$ such that $wP_0w^{-1}\subset P$,  the point $T_P$ is defined to be the orthogonal projection of $w\cdot T$ on $ \ago_{P}$ that is $T_P=(w\cdot T)_{M_P}$ (with the notations of § \ref{S:proj}).  One can check that this does not depend on the choice of $w$. If $P$ is standard, then we have $T_P=T_{M_P}$.
\end{paragr}

\subsection{Weyl groups and sets of Levi subgroups}

\begin{paragr}[Double cosets.] --- \label{S:rep} 
  
  Let  $P=MN_P$ and  $Q=LN_Q$ be standard  parabolic subgroups of $ G $ with standard Levi decompositions. Let $W^Q=W^{M_Q}$, resp. $W^{P}=W^{M_P}$, be  the Weyl group of $(M_Q,A_0)$, resp. $(M_P,A_0)$. Let $_QW_P$ be the set of  $w\in W$ such that 
\begin{itemize}
\item $M\cap w^{-1}P_0w=M \cap P_0$ ;
    \item $L\cap wP_0w^{-1}=L\cap P_0$.
\end{itemize}
This is a set of representatives of the double quotient $W^Q\back W/W^P$. Moreover, $M \cap w^{-1}L w$ is the Levi factor of  the  standard  parabolic subgroup $P_w=(M \cap w^{-1}Qw)N_P $ included in $P$. In the same way, $L \cap w M w^{-1}$ is the Levi factor of  the  standard  parabolic subgroup  $Q_w=(L \cap w Pw^{-1})N_Q$   with $Q_w\subset Q$. We introduce:
\begin{align*}
&  W(P;Q)=\{w\in \!_QW_P \mid P_w=P\}= \{w\in \!_QW_P \mid M\subset w^{-1}Lw\}\\
 & W(P,Q)=\{w\in \!_QW_P \mid M= w^{-1}Lw\}.
\end{align*}

When $P=Q$, the group $W(P,P)$ is denoted by  $W(P)$. For any  $w\in W(P,Q)$ one has  $w\Delta^P_0=\Delta_0^Q$. Note also that $w\in W(P_w,Q_w)$ for all $w\in \, _QW_P$. We have $W(P,Q)\subset W(P;Q)$ and  we can also write the set $W(P;Q)$ as the disjoint union
\begin{align}
  \label{eq:WPQ}
W(P;Q)=\bigcup_{R\subset Q} \{w\in W(P,R) \mid w^{-1}\Delta_0^Q>0\}
\end{align}
where $R$ runs over the set of standard parabolic subgroups of $Q$ and $>0$ means positive relatively to $P_0$. The contribution of $R$ is empty unless $R$ is associated to $P$.

Let
\begin{align*}
W_2(P)=\{w\in W(P)\mid w^2=1\}.
\end{align*}
 We view it as a subgroup of $W$. If $P\subset Q$, we set $W_2^Q(P)=W_2(P)\cap W^Q$. More generally, for any standard parabolic subgroup $R$ of $G$, we will denote by an upperscript $R$ an  object relative to its Levi factor $M_R$ equipped with the minimal pair $(P_0\cap M_R, M_0)$. The notations above hold if we replace the standard parabolic subgroups by their standard Levi components, e.g. $W(M;L)$ and $W_2(M)$ respectively mean $W(P;Q)$ and $W_2(P)$ if $M=M_P$ and $L=M_Q$.

\begin{lemme}\label{lem:ww'}
  Let $P,Q,R$ be standard parabolic subgroups of $G$. We assume $Q\subset R$.
  \begin{enumerate}
  \item For any $w\in \,_RW_P$, we have  $\, _QW^R_{R_w} w\subset  \,_QW_P$. 
  \item   For any $w_2\in \,_QW_P$, there is a unique  decomposition $w_2=w_1w$ with  $w\in \,_RW_P$ and $w_1\in  \, _QW^R_{R_w}$. Moreover, $w_2\in W(P;Q)$ if and only if $w_1\in W^R(R_w;Q)$ and $w\in W(P;R)$.
  \end{enumerate}
  
\end{lemme}

\begin{preuve} Let $w\in \,_RW_P$ and $w_1\in  \, _QW^R_{R_w}$. First, we show that  $w_2=w_1w$ belongs to $\,_QW_P$.
  
One has $M_{R_w}\cap w_1^{-1}P_0 w_1=M_{R_w}\cap P_0$. Thus $M_{P_w}\cap w_2^{-1}P_0 w_2=M_{P_w}\cap w^{-1} P_0 w=M_{P_w}\cap P_0$. Then one gets $M_P\cap w_2^{-1}P_0 w_2\subset M_P\cap w^{-1}R w=M_P \cap P_w\subset P_w$ (indeed one observes that $M_P\cap w^{-1}R w $ contains $M_P\cap w^{-1}P_0 w =M_P\cap P_0$), hence one has  $P_w\cap  w_2^{-1}P_0 w_2 \subset (M_{P_w}\cap w_2^{-1}P_0 w_2)N_{P_0}\subset P_0$ and  $M_P\cap w_2^{-1}P_0 w_2=M_P\cap P_0$.

  One also has
\begin{align*}
  M_Q\cap w_2P_0 w_2^{-1}&\subset w_1( M_R\cap wP_0w^{-1})w_1^{-1}\\
                         &\subset w_1( M_R\cap P_0)w_1^{-1}\subset w_1P_0 w_1^{-1}.
\end{align*}
Hence we have 
\begin{align*}
  M_Q\cap w_2P_0 w_2^{-1}=M_Q\cap w_1P_0 w_1^{-1}=M_Q\cap P_0.
\end{align*}
  Conversely, any $w_2\in  \,_QW_P$ can be uniquely written $w_2=w_1 w$ with $w_1\in W^R$ and $w\in   \,_RW_P$. We want to show $w_1\in  \, _QW^R_{R_w}$. On the one hand, we have
  \begin{align*}
    M_{R_w}\cap w_1^{-1}P_0 w_1&= M_R \cap w M_P w^{-1} \cap  w_1^{-1}P_0 w_1\\
                               &=M_R \cap w (M_P  \cap  w_2^{-1}P_0 w_2) w^{-1}\\
                               &=M_R \cap w (M_P  \cap  P_0 ) w^{-1}\\
                               &= M_R \cap w (M_P  \cap  w^{-1}P_0 w) w^{-1}\\
    &\subset P_0.
  \end{align*}

    Hence $M_{R_w}\cap w_1^{-1}P_0 w_1=M_{R_w}\cap P_0$. On the other hand, 
\begin{align*}
  M_Q\cap w_1P_0 w_1^{-1}&\subset M_R\cap w_1P_0w^{-1}_1&= & w_1(M_R\cap P_0)w^{-1}_1\\
    &                  &= &w_1(M_R \cap wP_0 w^{-1}) w^{-1}_1\\
     &                    &\subset& w_2 P_0 w_2^{-1}.
                           \end{align*}
Hence we have $  M_Q\cap w_1P_0 w_1^{-1}=M_Q\cap w_2P_0 w_2^{-1}=M_Q\cap P_0.$

Let us prove the last claim. Let $w_1\in W^R(R_w;Q)$ and $w\in W(P;R)$. Set $w_2=w_1w$. Then
$$M_P = w^{-1}M_{R_w} w\subset  w^{-1}w_1^{-1}M_{Q}w_1 w=w_2^{-1}M_Q w_2$$
hence $w_2\in  W(P;Q)$.

Conversely let $w\in \,_RW_P$ and $w_1\in  \, _QW^R_{R_w}$ such that $w_2=w_1w \in W(P;Q)$. We have
$$M_P\subset w_2^{-1}M_Q w_2 \subset w_2^{-1}M_R w_2 = w^{-1}M_R w$$
hence $w\in W(P;R)$. Then 
$$M_{R_w}=M_R\cap w M_P w^{-1}=w M_P w^{-1}\subset ww_2^{-1}M_Q w_2w^{-1}=w_1^{-1}M_Q w_1$$
hence $w_1\in W^R(R_w;Q)$.
\end{preuve}
\end{paragr}

\begin{paragr}\label{S:pcM-W}
  Let $M$ a semi-standard Levi subgroup of $G$. We denote by $\lc(M)$ (resp. $\fc(M)$)  the set of Levi subgroups (resp. parabolic subgroups) of $G$ that contain $M$. Let $\pc(M)$ be the set of minimal elements of $\fc(M)$ that is the elements $P\in \fc(M)$ such that $M_P=M$. We have a disjoint union
  $$\fc(M)=\bigcup_{L\in \lc(M)} \pc(L).$$
  If $M$ is moreover  standard, it defines a standard parabolic subgroup  $P=MN_P$. Then the map
  \begin{align}
    \label{eq:wQ}
(Q,w)\mapsto w^{-1}\cdot Q,
  \end{align}
  where we set  $w^{-1}\cdot Q=w^{-1}Qw$, induces a  bijection from the disjoint union $\bigcup_Q W(P;Q)$ where $Q$ runs over the set of standard parabolic subgroups of $G$ onto $\fc(M)$.  It also induces  a  bijection from the disjoint union $\bigcup_Q W(P,Q)$ onto $\pc(M)$ (in this case only parabolic subgroups $Q$ that are associated to $P$ do contribute to the source). 
  We denote the sets above $\lc^G(M)$, $\pc^G(M)$ and $\fc^G(M)$ if we want to emphasize the dependence on $G$.
\end{paragr}

\begin{paragr}\label{S:fc2}
  Let $P=MN_P$ be a standard parabolic subgroup of $G$ (with its standard decomposition).  Any $\xi\in W_2(M)$  induces an involution of $\ago_M$ and $\ago_M^*$. We have  a decomposition $\ago_M=\ago_M^{\xi}\oplus \ago_M^{-\xi}$ where $\ago_M^{\pm\xi}$ denotes the $(\pm1)$-eigenspace spaces of $\xi$. The same notation holds for the dual space $\ago_M^*$. There exists a unique $L_\xi\in \lc(M)$ such that $\ago_M^{L_\xi}=\ago_M^{-\xi}$. We get a bijective map $\xi \mapsto L_\xi$ from $W_2(M)$ onto a subset of $\lc(M)$ denoted by $\lc_2(M)$. We denote the inverse map $\lc_2(M)\to  W_2(M)$ by $L\mapsto \xi_M^L$. We say that the involution $\xi$ is standard  if $L_\xi$ is standard (such a $\xi$ is called minimal  in \cite{JLR} section VII.16). Let $W_2^{\mathrm{st}}(M)\subset W_2(M)$ be the subset of standard $\xi$'s. Let $\fc_2(M)\subset \fc(M)$ be the subset of $Q\in \fc(M)$ such that $M_Q\in \lc_2(M)$. We get a map 
$$Q\mapsto \xi_M^Q=\xi_M^{M_Q}$$
from $\fc_2(M)$ to $W_2(M)$.
\end{paragr}

\subsection{Sets of inversions}\label{ssec:set inversion}

\begin{paragr}[Set of inversions.]--- \label{S:elem sym}Let $P_1$ be a standard parabolic subgroup of $G$ and let $M_1$ be its standard Levi factor. Let
    \begin{align*}
      W_{M_1}&=\cup_{Q}W(P_1,Q)\\
      &=\{w\in W\mid w\Delta_0^{M_1}\subset \Delta_0\}
    \end{align*}
    where the union is taken over the set of standard parabolic subgroups $Q$. To shorten the notations, we set  $\Delta_1=\Delta_{P_1}$ and $\Phi_1=\Phi_{P_1}$. Let $\al\in \Delta_1$ and let $M_\al$ be the standard Levi subgroup defined by $\Delta_{M_1}^{M_\al}=\{\al\}$. The set $W_{M_1}^{M_\al}$ (defined relatively to the ambient group $M_\al$) is included in $W_{M_1}$ and is reduced to two elements. We denote by $s_\al$ the non trivial element: it is called an elementary symmetry. Let $w\in W_{M_1}$. For a root $\be$ of $wA_{M_1}w^{-1}$ in $G$ we write $\be>_{wM_1w^{-1}}0$ or simply $\be>0$ if $\be $ is  a root of $wA_{M_1}w^{-1}$ in the standard parabolic subgroup of Levi factor $wM_1w^{-1}$. We write $\be<0$ if $-\be>0$.

Let $\Phi_1^{\nd}\subset   \Phi_1$ be the subset of non divisible roots.   We define the set of inversions of $w$ by
    \begin{align*}
      \Ga_w=\{\al\in \Phi_1^{\nd} \mid w\al<0\}.
    \end{align*}
  \end{paragr}

  \begin{paragr}[Length.] --- \label{S:length} The cardinality of $\Ga_w$ is denoted by $\ell_{M_1}(w)$, or simply $\ell(w)$ if the context is clear, and is called the length of $w$. Note that $\Ga_{s_\al}=\{\al\}$.
    
    \begin{lemme}\label{lem:additivite}
  Let $w_1\in W_{M_1}$ be such that $\Ga_{w_1}\subset \Ga_w$. We have $ww_1^{-1}\in W_{M'}$ with $M'=w_1Mw_1^{-1}$ and $\ell_{M_1}(w)=\ell_{M'}(ww_1^{-1})+\ell_{M_1}(w_1)$. 
\end{lemme}

\begin{preuve}
  The set  $\Ga_w$ is the disjoint union of $\Ga'$ and $\Ga_{w_1}$ defined by:
    \begin{align*}
      \Ga'=\{ \al  \in \Phi_1^{\nd} \mid w_1\al>0 \text{ and } w\al<0   \}\\
      \Ga_{w_1}=\{ \al\in \Phi_1^{\nd}\mid w_1\al<0 \text{ and } w\al<0 \}.
    \end{align*}
   Set $w_2=ww_1^{-1}$. Let $P'$ be the standard parabolic subgroup of Levi factor $M'$. The map $\al\mapsto \be=w_1\al$  induces a  bijection from   $\Ga'$ onto the set
    \begin{align*}
      \{ \be\in \Phi_{P'}^{\nd} \mid w_1^{-1}\be>0 \text{ and } w_2\be<0   \}.
    \end{align*}
On the one hand the  cardinality of this set is $|\Ga_w|-|\Ga_{w_1}|=\ell_{M_1}(w)-\ell_{M_1}(w_1)$.  On the other hand it is equal to  $\ell_{M'}(w_2)$ if and only if the two conditions  $\be>_{M'}0$ and $w_2\be<0 $ imply $ w_1^{-1}\be>0 $ that is if and only the two conditions   $w_1\al>0$ and $w\al<0$ imply $\al>0$. Let us check the latter. Let  $\al<0$ be such that $w_1\al>0$. We have  $-\al\in \Ga_{w_1}\subset \Ga_w$ so $w\al>0$.    
\end{preuve}
    
  \end{paragr}

 \begin{paragr}[Closed and coclosed subsets.] ---     A subset $\Ga\subset \Phi_1^{\nd}$ is said to be closed if for all $\al,\be\in  \Ga$ such that $\al+\be\in \Phi_1^{\nd}$ we have $\al+\be\in \Ga$. It is called  coclosed if its complement $\Phi_1^{\nd}\setminus \Ga$ is closed. 

   \begin{lemme}
     \label{lem:set of inv}
     A subset $\Ga\subset \Phi_1^{\nd}$ is the set of inversions $\Ga_w$ for some element $w\in W_{M_1}$ if and only if it is closed and coclosed.
   \end{lemme}

   \begin{preuve}
   The condition is obviously necessary.  Let us start from a closed and coclosed subset  $\Ga\subset \Phi_1^{\nd}$ and let us prove the converse. Let $\tilde\Ga$ be the set of $\al\in \Phi_{P_0}^{\nd}\setminus \Phi^{P_1,{\nd}}_{P_0}$ such that the restriction of $\al$ to $A_{M_1}$ belongs to $\Ga$. It is clear that  $\tilde\Ga$ is a closed and coclosed subset of $\Phi_0=\Phi_{P_0}^{\nd}$. Following \cite[VI § 1 exercise 16]{Bki-Lie}, we introduce the set $\Phi_0'=(\Phi_{0}\setminus \tilde\Ga)\cup (-\tilde\Ga)$. By  \cite[VI § 1.7 corollaire 1]{Bki-Lie}, there exists a unique $w\in W$ such that $w^{-1}\Phi_0'=\Phi_0$. It is clear that $\tilde\Ga$ is exactly the set of $\al\in \Phi_0$ such that $w\al<0$. We claim that $w\in W_{M_1}$ that is $w\Delta_0^{M_1}\subset \Delta_0$. It is easy then to show that we have $\Ga_{w}=\Ga$. Let us prove the claim: let $\al\in \Delta_0^{M_1}$. Since $\al\notin\tilde\Ga$, we have $w\al\in \Phi_0$. Assume that $w\al=\be_1+\be_2$ with $\be_i\in \Phi_0$. Then we may and shall assume $w^{-1}\be_1>0$ and $w^{-1}\be_2 <0$. In particular $-w^{-1}\be_2 \in \tilde\Ga$ and $w^{-1}\be_1\in \Phi_0\setminus \tilde\Ga$. However, $w^{-1}\be_1=\al+ (-w^{-1}\be_2 )$ has the same restriction to $A_{M_1}$ as   $-w^{-1}\be_2 $ and thus belongs to $\tilde\Ga$. This is a contradiction.
 \end{preuve}

 For all  $\be\in \Delta_1$ let $\Phi_1^{\be}\subset \Phi_1^{\nd}$  be the subset of non divisible positive roots that are   sums of elements in $\Delta_1 \setminus\{\be\}$. 
  
  \begin{lemme}\label{lem:coclos} The set
    \begin{align}
      \label{eq:Ga be w}
      \Ga^{\be}_w=\Ga_w\cap \Phi_1^{\be}
    \end{align}
     is closed and coclosed.
      \end{lemme}

      \begin{preuve} It is clear that $\Phi_1^{\be}$ is closed and thus $\Ga^{\be}_w$ is  closed. Let $\ga,\ga'\in \Phi_1^{\nd}\setminus \Ga^{\be}_w$ be such that  $\ga+\ga'\in\Phi_1^{\nd}$. We have to show that $\ga+\ga'\notin  \Ga^{\be}_w$. The result is clear if   both $\ga$ and $\ga'$ do not belong to $\Ga_w$, resp. $\Phi_1^\be$. Without loss of generality we may assume that $\gamma\in \Phi_1^\be$ and $\ga'\in \Ga_w$. But then  $\ga'\notin \Phi_1^\be$ and $\ga+\ga'\notin \Phi_1^\be $.
      \end{preuve}
  
      The following lemma will also be useful.

       \begin{lemme}\label{lem:simplicite}
      Let $\al\in \Delta_1$ be such that $w\al>0$. The following assertions are equivalent:
  \begin{enumerate}
  \item The root $w\al$ is simple.
  \item For all $\gamma\in \Ga_w$ such that $\al+\gamma\in \Phi_1^{\nd}$, we have $\al+\gamma\in \Ga_w$.    
  \end{enumerate}
\end{lemme}

\begin{preuve}
  1 implies 2. Let  $\gamma\in \Ga_w$ be such that $\al+\gamma\in \Phi_1^{\nd}\setminus \Ga_w$. Then $w\al= w(\al+\gamma)-w\gamma$ is a sum of two positive roots.
  
    2 implies 1. Assume that $w\al$ is the sum of non divisible positive roots $\be_1$ and $\be_2$. We have  $\al=\al_1+\al_2$ with  $\al_i=w^{-1}\be_i$. Since $\al$ is simple we may assume without loss of generality that we have $\al_1>0$ and $\al_2<0$. Then we have  $-\al_2\in \Ga_w$ and $\al-\al_2=\al_1\in  \Phi_1^{\nd}\setminus \Ga_w$.    

  \end{preuve}
\end{paragr}

\begin{paragr}[Decomposition into elementary symmetries.] --- \label{S:dec elem sym} There exist standard  parabolic subgroups $P_2,\ldots, P_{\ell+1}$ with standard Levi subgroups $M_2,\ldots,M_{\ell+1}$ with $\ell=\ell(w)$, $\al_i\in \Delta_{M_i}$, $M_{i+1}=s_{\al_i}M_i s_{\al_i}^{-1}$ such that we have the ``reduced'' decomposition
   \begin{align}\label{eq:decomp}
     w=s_{\al_\ell}\cdots s_{\al_1}.
   \end{align}

   Moreover for such a decomposition we have, see \cite[lemme I.1.8]{MWlivre}, 
   \begin{align*}
     \Ga_w=\{\theta_i \mid i=1,\ldots,\ell\}\text{ with }  \theta_i= s_{\al_1}^{-1} s_{\al_2}^{-1}\cdots   s_{\al_{i-1}}^{-1}(\al_i)
   \end{align*}
   We define also
   \begin{align}
     \label{eq:thetavee}
      \theta_i^\vee= s_{\al_1}^{-1} s_{\al_2}^{-1}\cdots   s_{\al_{i-1}}^{-1}(\al_i^\vee).
   \end{align}
 \end{paragr}

\begin{paragr}[$w$-tableau.]  ---  This §  is inspired by \cite[section 2]{Kraski}. Each reduced decomposition \eqref{eq:decomp} defines a map $T:\Phi_{1}^{\nd}\to \{1,\ldots,\ell\}\cup \{+\infty\}$ such that  $T(\al)=+\infty$ for all $\al\in \Phi_1^{\nd}\setminus \Ga_w$ and $T(\al)=i$ if $\al=\theta_i$. We shall say that a map  $T:\Phi_{1}^{\nd}\to \{1,\ldots,\ell(w)\} \cup \{+\infty\}$ is  a $w$-tableau  if and only if it satisfies the following three conditions:
   \begin{itemize}
   \item For all $\al\in\Phi_{1}^{\nd}$ we have $T(\al)=+\infty$ if and only if $\al\notin\Ga_w$;
   \item $T$ induces a bijection between the sets $\Ga_w$ and $\{1,2,\ldots,\ell(w)\}$;
   \item for all $\al,\be\in  \Phi_1^{\nd}$ such that $\al+\be\in \Phi_1^{\nd}$ we have:
     \begin{align}\label{eq:couple ineg}
 \text{either } T(\al)\leq T(\al+\be)\leq T(\be)\text{ or }  T(\be)\leq T(\al+\be)\leq T(\al).
     \end{align}
   \end{itemize}

   \begin{proposition}\label{prop:w-tableau}
     Any reduced decomposition defines a map $T$ which is a $w$-tableau and conversely any $w$-tableau comes from a reduced decomposition.
   \end{proposition}

   \begin{preuve}
     We proceed by induction on the length $\ell$ of $w$. We start from the reduced decomposition \eqref{eq:decomp}.   Let $T$ and $T'$ be the maps defined respectively by $w$ and $w'=s_{\al_\ell}^{-1}w$. By induction, we assume that $T'$ is a $w'$-tableau. Note that $T$ and $T'$ coincides on $\Phi_1^{\nd}\setminus\{\theta_\ell\}$.   Let $\al,\be\in \Phi_1^{\nd}$ be such that $\al+\be\in \Phi_1^{\nd}$. We have to prove that $T(\al+\be)$ satisfies one of the two inequalities \eqref{eq:couple ineg}.  If none of the elements $\al,\be$ and $\al+\be$ is equal to $\theta_\ell$  the conclusion is clear. First assume that $\al=\theta_\ell$ so $\be$ and $\al+\be$ are not equal to $\al$. We have $T'(\al)=+\infty$. So $T'(\be)\leq T'(\al+\be)$ and so $T(\be)\leq T(\al+\be)$. If $T'(\be)=+\infty$ we have $T(\al+\be)=T'(\al+\be)=+\infty$ and  the conclusion is  clear. So we assume $T(\be)=T'(\be)<+\infty$. Since $\Ga_w$ is closed by lemma \ref{lem:set of inv}, we have $T(\al+\be)<+\infty$ and so $T(\al+\be)< T(\al)=\ell$. Assume now that we have $\al+\be=\theta_{\ell}$. Since $\Ga_w$ and $\Ga_{w'}$ are both  closed and coclosed exactly one  of $T(\al)=T'(\al)$ and $T(\be)=T'(\be)$ is equal to $+\infty$.  Without loss of generality we shall assume $T(\al)=+\infty$. So we  have  $T(\be)<T(\al+\be)=\ell<T(\al)=+\infty$.

   Conversely, let us start with a $w$-tableau $T$. Let $\theta_\ell=T^{-1}(\ell)$. Let $T'$ be the map which coincides on $\Phi_1^{\nd}\setminus\{\theta_\ell\}$ with $T$ and such that $T'(\theta_\ell)=+\infty$. Using lemma \ref{lem:set of inv}, we see that there exists $w'\in W_{M_1}$ such that $T'$ is a $w'$-tableau with $\Ga_{w'}=(T')^{-1}(\{1,\ldots,\ell-1\})$. By induction, $T'$ comes from a reduced decomposition $w'=s_{\al_{\ell-1}}\cdots s_{\al_1}$. Let $s_\ell=w(w')^{-1}$. By lemma \ref{lem:additivite}, we see that $\ell(s_\ell)=1$ that is $s_{\ell}$ is an elementary symmetry attached to some simple root $\al_{\ell}$. So $ s_{\al_1}^{-1} s_{\al_2}^{-1}\cdots   s_{\al_{\ell-1}}^{-1}(\al_\ell)$ is the unique element in $\Ga_w\setminus\Ga_{w'}$ namely $\theta_{\ell}$. In this way $T$ comes from the  reduced decomposition $w=s_{\al_{\ell}}\cdots s_{\al_1}$.
\end{preuve}
\end{paragr}

\subsection{Polynomial exponential}

\begin{paragr}[Polynomial exponential.] --- \label{S:py-exp}Let $V$ be a  real  vector space of finite dimension. Let $V^*$ be its dual. Let  $V_\CC^*$ and $V_\CC$ be the $\CC$-vector space obtained from $V^*$ and $V$ by extension of scalars to $\CC$. We denote by $\bg \cdot, \cdot\bd: V_\CC^*\times V_\CC\to \CC$ the canonical pairing.

  Let  $\Om\subset V$ be a non-empty open subset.   By a polynomial exponential on  $\Om$, we mean a map $E: \Om \to \CC$ which belongs to the $\CC$-vector space generated by $T\in \Om \mapsto p(T) \exp(\bg \la,T\bd)$ where $\la\in V_\CC^*$ and $p\in \CC[V]$ is a polynomial. Any  polynomial exponential $E$ on $\Om$ can be uniquely  written $\sum_{\la \in V_\CC^*} p_\la(T) \exp(\bg \la,T\bd)$ for a unique map $\la\in V_\CC^*\mapsto  p_\la \in \CC[V]$ with finite support. The finite set of $\la\in V_\CC^*$ such that $p_\la\not=0$ is the set of exponents of $E$. We shall say that $p_\la$ is the polynomial part of $E$ of exponent $\la$.  If $\la=0$, we shall simply say that $p_0$ is the polynomial part.

  Let $\om$ be a non-empty open subset of a $\CC$-vector space $Z$ of finite dimension.  By a meromorphic function on $\om$ with hyperplane singularities we mean a meromorphic  function $f$ on $\om$  such that for any $z_0\in \om$, there exists a finite family of linear forms $(L_i)_{i\in I }$ (not necessarily two by two distinct) such that the product $\prod_{i\in I} L_i(z-z_0) f(z)$ is holomorphic in a neighborhood of $z_0$ in $\om$. Let $\mc(\om)$ be the $\CC$-algebra of  meromorphic functions on $\om$ with hyperplane singularities.
  
The following lemma is a slight variant of a lemma used by Offen (see \cite{Offen-symplectic} lemma 8.1)

\begin{lemme}\label{lem:py-exp-lisse}
  Let $I$ be a finite set. For each $i\in I$, we consider the following objects:
  \begin{itemize}
  \item an affine map $q_i:Z\to V_\CC^*$;
  \item  a polynomial $p_i\in \CC[V]$;
  \item a meromorphic function $f_i\in \mc(\om)$.
  \end{itemize}
 Let $z_0\in \om$.  We assume that for each $T\in \Om$, the map
  \begin{align*}
    z\in \om \mapsto E(z,T)=\sum_{i\in I} p_i(T) f_i(z) \exp(\bg q_i(z),T\bd)
  \end{align*}
  is holomorphic in a neighborhood of  $z_0$ in $\om$.

  \begin{enumerate}
  \item The map $T\in \Om \mapsto E(z_0,T)$ is a polynomial exponential whose set of exponents is  included in $\{ q_i(z_0), i\in I\}$.
  \item Let  $j\in I$ and  $I_j=\{i\in I \mid q_i(z_0)=q_j(z_0)\}$. The sum
      \begin{align*}
    z\in \om \mapsto E_j(z,T)=\sum_{i\in I_j} p_i(T) f_i(z) \exp(\bg q_i(z)-q_j(z_0),T\bd)
      \end{align*}
      is also holomorphic  in a neighborhood of  $z_0$ in $\om$ and $E_j(z_0,T)$ is the  polynomial part of   $E(z_0,T)$ of exponent $q_j(z_0)$.
    \end{enumerate}
  \end{lemme}

  \begin{preuve}
    Clearly assertion 1 follows from assertion 2. Let us prove assertion 2. Without loss of generality we may and shall assume $z_0=0$. We may also shrink $\om$ if necessary. Let $(L_k)_{1\leq k \leq r}$  be a finite family of non-zero linear forms on $Z$. Let $H_k=\Ker(L_k)$. We assume that  hyperplanes  $H_k$ are two by two distinct for $1\leq k \leq r$. We assume also that for each $k$ there is an integer $n_k\geq 1$  such that  for  $L=\prod_{1\leq k \leq r} L_k^{n_k}$ the map $z\mapsto L(z) f_i(z)$ is holomorphic  in a neighborhood of $0$ in $\om$ for all $i\in I$. Let $z_k\in Z$ such that $L_k(z_k)=1$. We denote by $\partial_k$ the complex derivative along the vector $z_k$. We shall use the following elementary result: for any $f\in \mc(\om)$ such that the map $z\mapsto L(z)f(z)$ is holomorphic  on $\om$, the map $f$ is holomorphic on $\om$ if and only if for each  $1\leq k \leq r$ and $0\leq m \leq n_k-1$ the map $z\mapsto \partial_k^{m}(L f)(z)$ vanishes on  $\om \cap H_k$.   

By assumption, we have for each  $1\leq k \leq r$ and $0\leq m \leq n_k-1$
\begin{align}\label{eq:vanish-LE}
  \partial_k^{m} (L (z) E(z,T)) =\sum_{i\in I} p_i(T)  \partial_k^{m} (L(z)f_i(z) \exp(\bg q_i(z),T\bd))
\end{align}
vanishes for all $T\in \Om$ and $z\in \om \cap H_k$.  For all $z\in \om \cap H_k$, the map $T\mapsto \partial_k^{m} (L(z)f_i(z) \exp(\bg q_i(z),T\bd))$ is a polynomial exponential with a unique exponent namely $q_i(z)$. We may and shall shrink $\om$ such that that for all $i,j\in I$ such $q_i(z)=q_{j}(z)$ for some $z\in \om\cap H_k$ then $q_i(0)=q_{j}(0)$ and thus $i\in I_j$. Let $j\in I$. Identifying the summand of exponents $q_i(z)$ for $i\in I_j$ in \eqref{eq:vanish-LE} we get that
\begin{align}\label{eq:vanish-LE2}
  \partial_k^{m} (L (z) E_j(z,T)) \exp(\bg q_j(0),T\bd) =\sum_{i\in I_j} p_i(T)  \partial_k^{m} (L(z)f_i(z) \exp(\bg q_i(z),T\bd))
\end{align}
vanishes for all $T\in  \Om$ and $z\in \om \cap H_k$. As a consequence $E_j(z,T)$ is holomorphic on $\om$. Let us prove that  $E_j(0,T)$ is a polynomial in $T$. Without loss of generality, we may and shall assume  that $q_j(0)=0$. Then $q_i$ is linear for $i\in I_j$. Let $z\in \om$ outside the hyperplanes $H_k$ for  $1\leq k \leq r$. Then $E_j(0,T)$  is the degree $0$ contribution  in the Taylor expansion at $t=0$ of 
$$t\in \RR\mapsto \sum_{i\in I_j} p_i(T) f_i(t z) \exp(t \bg q_i(z),T\bd)$$
which is clearly a polynomial in $T$.
\end{preuve}
\end{paragr}

\begin{paragr} Let $n\geq 1$ an integer.   For any $z\in \CC^n$, we set $\|(z_1,\ldots,z_n)\|^2=\sum_{1\leq i\leq n}|z_i|^2$. For any $0<r\leq +\infty$ we denote by $B_r$ the set $\{z\in i\RR^n \mid \|z\|\leq r\}$. The following lemma is a variant of \cite[lemme 13.2.2]{LabWal}. For $k\geq 1$ and $c,\La>0$ we also set
  \begin{align*}
    \rc_{\La,c,k}=\{\la\in \CC^n\mid \|\Re(z)\|    (1+\La+\|\Im(z)\|)^{k} < c\}.
  \end{align*}
  
\begin{lemme}\label{lem:maj-f-F}  
  Let $(L_i)_{i\in I}$ be a finite family of non-zero real linear forms on $\CC^n$, resp. and let $k\geq 1$ and $\La>0$.   For any holomorphic differential operator  $D$ with polynomial coefficients on $\CC^n$ and any  $0<r\leq +\infty$, resp. any $c>0$, there exist differential operators $D_1,\ldots, D_\ell$ of the same nature and  $1>\al>0$ such that for any holomorphic function $f$ on a neighborhood of $i \RR^n$, resp. on   $\rc_{\La,c,k}$, we have:
  \begin{align}
    \sup_{z\in B_{\al r}} |Df(z)|\leq \sum_{i=1}^\ell   \sup_{z\in B_{r}} |D_iF(z)|
  \end{align}
  resp.
   \begin{align}
    \sup_{z\in \rc_{\La,\al c ,k}} |Df(z)|\leq \sum_{i=1}^\ell  \sup_{z\in \rc_{\La,c,k}} |D_iF(z)|
  \end{align}
where $F= f \prod_{i\in I} L_i$.
\end{lemme}

\begin{preuve}   For the convenience of the reader, we sketch a proof which is a variant of that of \cite[lemme 13.2.2]{LabWal}.  First, by recursion,  we may assume that there is only one linear  form $L$. There exists a real invertible matrix $A$ of size $n$ such that $L (Az)=z_1$ where $z=(z_1, \ldots,z_n)$. Since there exist $1>\al,\be>0$ such that $B_{\al r} \subset A (B_{\be r}) \subset B_r$, resp. $\rc_{\La,\al c,k}\subset A(\rc_{\La,\beta c,k})\subset \rc_{\La,c,k}$, it suffices to bound $f\circ A$ on $B_{\be r}$, resp. $\rc_{\La,\beta c,k}$. Replacing $f$ by $f\circ A$,  $L$ by $L\circ A$ and $r$ by $\be r$, resp. $c$ by $\be c$ we are reduced to bound $Df$ on $X$ in terms of $F=z_1f$ where $X=B_r$, resp. $X= \rc_{\La,c,k}$.  Then we are reduced to the case $D=D' D''$ where $D'$ depends only on $z_1$  and $D''$ depends only on the other coordinates $z_2,\ldots,z_n$. Since $D''(F)=LD''(f)$ we may assume that $D$ itself depends only on $z_1$. In this way, we are reduced to the case where $Df=\frac{\partial^i}{\partial z_1^i} (z_1^j f(z))$. The case $j\geq 1$ is obvious. So we assume that $D=\frac{\partial^i}{\partial z_1^i} $. There exists a holomorphic differential operator  $D_i$ with polynomial coefficients such that $Df= z_1^{-(i+1)} D_iF$. For any $z\in X$ let $z'=(0,z_2,\ldots,z_n)$. We have $z'\in X$. Since $D_iF$ vanishes at $L=0$ at the order $i+1$,  the mean value theorem implies that for all $z\in X$
  \begin{align*}
| D_iF (z)| \leq |z_1|^{i+1} \sup_{u\in [z;z']} |   \frac{\partial^{i+1}}{\partial z_1^{i+1}} D_i F (u)|
  \end{align*}
  hence
  \begin{align*}
    | Df(z)| \leq \sup_{u\in [z;z']} |   \frac{\partial^{i+1}}{\partial z^{i+1}} D_i F (u)        |.
  \end{align*}
and
 \begin{align*}
    \sup_{z\in  X}| Df(z)| \leq \sup_{u\in X} |   \frac{\partial^{i+1}}{\partial z^{i+1}} D_i F (u)        |.
 \end{align*}
\end{preuve}
\end{paragr}

\subsection{A result on $(G,M)$-families}

\begin{paragr} Let $M$ be a semi-standard Levi subgroup of $G$. A family $(c_P)_{P\in \pc(M)}$ of functions on $\ago_{P,\CC}^{G,*}$ is a Arthur's $(G,M)$-family (cf. \cite{ar-inv}) if for each $P\in \pc(M)$ the map $c_P$ is smooth on $i\ago_{P}^{G,*}$ and   for each pair $(P,P')$ of adjacent parabolic subgroups in $\pc(M)$ we have $c_P(\la)=c_{P'}(\la)$ for $\la\in i\ago_{P}^{G,*}$ such that  $\bg \la,\al^\vee\bd=0$ where $\{\al\}=\Delta_P\cap (-\Delta_{P'})$.

  \begin{proposition}\label{prop:GMfam}
    Let $M$ be a semi-standard Levi subgroup of $G$ and $\om$ be a neighborhood of  $i\ago_{M}^{G,*}$ in $\ago_{M,\CC}^{G,*}$. Let  $(c_P)_{P\in \fc(M)}$ be a family of    holomorphic  functions on $\om$. Assume that there exists a non-empty subset $\Om$ of $\ago_0$ such that for all $T\in \Om$ the map
    \begin{align*}
    \la\mapsto E(\la,T)=  \sum_{P\in \fc(M)} c_P(\la) \frac{\exp(\bg \la,T_P\bd)}{\theta_P(\la)}
    \end{align*}
    is holomorphic on $\om$. Then for any $L\in \lc(M)$ and $\la_0\in i\ago_{M}^{L,*}$ the family $(c_P(\la_0+\la))_{P\in \pc(L)}$ is a $(G,L)$-family of functions of the variable $\la\in \ago_{L,\CC}^{G,*}$.
  \end{proposition}

   \begin{remarque}\label{rq:GMrenforce}
We shall use the following reinforcement of  proposition \ref{prop:GMfam} when $\om$ is a domain. For each pair $(P,P')$ of adjacent parabolic subgroups in $\pc(L)$ we have $c_P(\la)=c_{P'}(\la)$ for $\la\in\om$ such that  $\bg \la,\al^\vee\bd=0$ where $\{\al\}=\Delta_P\cap (-\Delta_{P'})$.
   \end{remarque}

  \begin{preuve} Let $L\in \lc(M)$. We shall assume $L\subsetneq G$ otherwise there is nothing to prove. Let $(P, P')$ be a pair   of  adjacent parabolic subgroups in $\pc(L)$. We shall assume that $P$ is standard. There is no loss of generality because the choice of $P_0$ intervenes only in the definition of $T_P$. In  general there exists $w_0\in W$ such that $w_0 P_0 w_0^{-1}\subset P$ and $T_P$ is the projection of $w_0 Tw_0^{-1}$ on $\ago_P$. Then we can replace $P_0$ by $w_0 P_0 w_0^{-1}$  and accordingly $\Om$ by $w_0^{-1}\Om w_0$.

    The set  $\Delta_{P}\cap (-\Delta_{P'})$ is a singleton  whose unique element is denoted by $\al$. Let $P\subset P_\al$ be defined by $\Delta_P^{P_\al}=\{\al\}$. Let $L_\al=M_{P_\al}$. We have 
\begin{align*}
  \Ker(\al^\vee)= \{\la\in \ago_M^{G,*} \mid \bg \la,\al^\vee\bd=0\}=\ago_M^{L,*}\oplus \ago_{L_\al}^{G,*}.
\end{align*}
We need  to show the equality
    \begin{align}
      \label{eq:cPcP'}
c_P(\mu)=c_{P'}(\mu).
    \end{align}
for all $\mu\in  i  \Ker(\al^\vee)$. Because both $c_P$ and $c_{P'}$ are holomorphic on $\om$, it suffices to prove the equality for generic elements of $ i  \Ker(\al^\vee)$ in the sense that they do not belong to the  finite set of proper subspaces of $ i  \Ker(\al^\vee)$ we will encounter in our discussion. From now on, $\mu$ is a generic element of $ i  \Ker(\al^\vee)$.

Let $\fc_\al(M)$ be the set of $Q\in \fc(M)$ such that $\theta_Q(\mu)=0$ that is $\theta_Q$ vanishes identically on $\Ker(\al^\vee)$. Of course we have $P,P'\in \fc_\al(M)$. More generally we have   $Q\in \fc_\al(M)$ if and only there is $\be^\vee\in \Delta_Q^\vee$ such that  $ \Ker(\al^\vee)\subset \Ker(\be^\vee)$ that is $\be^\vee\in \RR\al^\vee$.  So for such a $Q$ we have

\begin{align}
  \label{eq:inclusion-hyperplan}
  \ago_M^{Q,*}\subset \Ker(\al^\vee).
\end{align}
 We introduce:
\begin{align*}
  E_\al(\la,T)=  \sum_{Q\in \fc_\al(M)} c_Q(\la) \frac{\exp(\bg \la,T_Q\bd)}{\theta_Q(\la)}.
\end{align*}
Since $E(\la,T)$ is holomorphic on $\om$ and since $E(\la,T)-E_\al(\la,T)$ is clearly holomorphic at $\la=\mu$, we deduce that  $E_\al(\la,T)$ is also holomorphic at $\la=\mu$.  First we want to determine the  summand of   $E_\al(\mu,T)$ of exponent $T\mapsto \bg \mu,T_P\bd= \bg \mu_L,T\bd$ (recall that $P$ is assumed to be standard).

Let $Q\in \fc_\al(M)$ be such that $\bg \mu,T_Q\bd= \bg \mu_L,T\bd$ for all $T\in \Om$. Let $w\in W$ be such that $wP_0 w^{-1}\subset Q$. Then $\bg \mu,T_Q\bd=\bg w^{-1}\cdot(\mu_{M_Q}),T\bd$. Thus we must have $\mu_L=w^{-1}\cdot(\mu_{M_Q})$. By genericity, we have  $\la_L=w^{-1}\cdot(\la_{M_Q})$ for all $\la\in i \Ker(\al^\vee)$. Since $\ago_M^{L,*}\subset \Ker(\al^\vee)$ we deduce $\ago_M^{L,*}\subset \ago_M^{Q,*}$. But by \eqref{eq:inclusion-hyperplan} we have also $\ago_M^{Q,*}\subset \Ker(\al^\vee)$ and we deduce  $\ago_M^{Q,*}\subset \ago_M^{L,*}$. As a consequence we have $\ago_M^{Q,*}=\ago_M^{L,*}$ and $L=M_Q$ that is $Q\in \pc(L)$. Once again by genericity of $\mu$, we deduce that $w$ stabilizes $\ago_{L_\al}^{G,*}$ and in fact acts trivially on it. As a consequence $w\in W^{L_\al}$. There are only two possibilities: either  $Q=P$ or  $Q=P'$. By lemma \ref{lem:py-exp-lisse}, the summand of    $E_\al(\mu,T)$ of exponent $\bg \mu_L,T\bd$ is given by the value at $\la=\mu$ of the expression
\begin{align*}
   D_{P,P'}(\la)=c_P(\la) \frac{\exp(\bg \la,T_P\bd)}{\theta_P(\la)}+c_{P'}(\la) \frac{\exp(\bg \la,T_{P'}\bd)}{\theta_{P'}(\la)}
\end{align*}
which is holomorphic at $\la=\mu$. In particular, $\bg \la,\al^\vee\bd   D_{P,P'}(\la)$ vanishes at $\la=\mu$. Let $d_P$ be the value at $\la=\mu$ of 
$\bg \la,\al^\vee\bd   \frac{\exp(\bg \la,T_P\bd)}{\theta_P(\la)}$. Let $d_{P'}$ be the value we get when $P$ is replaced by $P'$. One can check that $d_{P'}=-d_P\not=0$. Thus we must have $d_P(c_P(\mu)-c_{P'}(\mu))=0$ and thus $c_P(\mu)=c_{P'}(\mu)$.

\end{preuve}
  \end{paragr}

  \subsection{Automorphic forms}\label{ssec:autom}

\begin{paragr}
  We continue with the notations of the previous sections. From now on we assume that the base field $F$ is a number field. Let  $\AAA $ be its adèle ring. Let $ \AAA_f $ be  the ring of finite adèles and $F_\infty=F\otimes_{\QQ} \RR$. Let $ V_F $ be the set of places of $F$ and $V_{F,\infty}\subset V_F$ be the subset of Archimedean places.  For every $ v \in V_F $, we let $ F_v $ be the local field obtained by completion of  $ F $ at $ v $. We denote by  $ | \cdot | $  the morphism $ \AAA^\times \to \RR_+^\times $ given by the product of normalized absolute values $|\cdot|_v$ on each $ F_v $. 
\end{paragr}

\begin{paragr}
  Let $ K=\prod_{v\in V_F} K_v\subset G (\AAA) $ be a ``good'' maximal compact subgroup in good position relative to $M_0$ (that is it is admissible in the sense of \cite[p.9]{ar-inv}). Let $K_f=\prod_{v\in V_F\setminus V_{F,\infty}} K_v\subset G(\AAA_f)$ and $K_\infty=\prod_{v\in V_{F,\infty}} K_v\subset G(F_\infty)$.
\end{paragr}

\begin{paragr} Let $ P $ be a semi-standard parabolic subgroup. We have a canonical homomorphism  
\begin{align}\label{eq:HG}
H_P:P(\AAA)\to \ago_P
\end{align}
characterized by $ \bg \chi, H_P(g) \bd = \log | \chi (g) | $ for any $g\in P(\AAA)$ and $\chi\in X^*(P)$. Its kernel is denoted by  $ P(\AAA) ^1 $. We extend it to the Harish-Chandra map 
$$ H_P: G (\AAA) \to \ago_ {P} $$
 that satisfies: for every $g \in G(\AAA)$ we have $H_P(g)=H_P(p)$ whenever $ g \in pK $ with $ p \in P (\AAA) $. Let $H(g)=H_0(g)=H_{P_0}(g)$.
\end{paragr}

\begin{paragr}[Haar measures.] --- \label{S:Haar-meas} For any parabolic subgroups $P\subset Q\subset G$, let  $\rho_P^Q$ be the unique element of $\ago_{P}^{Q,*}$ such that for every $m\in M_P(\AAA)$
$$|\det(\Ad_P^Q(m))|=\exp(\bg 2\rho_P^Q,H_P(m)\bd).
$$
where $\Ad_{P}^Q$ is  the adjoint action of $M_P$ on the Lie algebra of $M_Q\cap N_P$.  For $Q=G$, the exponent $G$ is omitted.

Let $A_P^\infty$ be the neutral component of $A_{P,\QQ}(\RR)$ (for the Archimedean topology) where $A_{P,\QQ}\subset \Res_{F/\QQ}(A_P)$ is the maximal $\QQ$-split subtorus. If $P\subset Q$ we denote by  $A_P^{Q,\infty}=A_P^\infty\cap \Ker(H_Q)$.   The group $ A_P ^ \infty $  is equipped with the Haar measure compatible with the isomorphism $ A_P ^ \infty \simeq  \ago_ {P} $ induced by the homomorphism \eqref{eq:HG}. Its subgroup $ A_P ^{Q, \infty} $ is in the same way identified with $\ago_ {P}^Q $  and thus inherits the Haar measure of $\ago_P^Q$.

If $N$ is a unipotent group, $N(\AAA)$ is provided with the Haar measure whose quotient by the counting measure on $N(F)$ gives the total volume $1$ to $[N]=N(F)\back N(\AAA)$. We fix Haar measures on $G(\AAA)$ and $K$. We assume that for any semi-standard parabolic groups $P$ of $G$ the Haar measure on $M_P(\AAA)$ is such that
  \begin{align*}
    \int_{G(\AAA)} f(g)\,dg=\int_{N(\AAA)\times M_P(\AAA) \times K} f(n mk) \exp(-\bg 2\rho_P,H_P(m)\bd)\, dn dmdk
  \end{align*}
for all continuous and compactly supported function $f$ on $G(\AAA)$.

We set
 \begin{align*}
   [G]_{P,0}=A_P^\infty M_P(F)N_P(\AAA)\back G(\AAA) \text{ and   } [G]_{P}=M_P(F)N_P(\AAA)\back G(\AAA).
 \end{align*}

 We equip  $[G]_{P,0}$ with the ``quotient measure'': it is the right-invariant functional on the space of continuous functions $f$ on $ G(\AAA)$ such that $f(pg)=\exp(\bg 2\rho_P,H_P(p)\bd)f(g)$ for all  $p\in A_P^\infty M_P(F)N_P(\AAA)$ that satisfies the following property:  for all continuous and compactly supported functions $f$ on $G(\AAA)$ we have
 \begin{align*}
    \int_{[G]_{P,0}} \int_{ A_P^\infty}  \sum_{\gamma\in M_P(F)}\int_{N_P(\AAA)  }  f(a\gamma ng)\, dn da \,dg=\int_{G(\AAA)} f(g)\,dg.
  \end{align*}
  If $P=G$ we set $[G]_0=[G]_{G,0}$ and $[G]=[G]_{G}$. We put on $G(\AAA)^1$ the Haar measure such that the natural isomorphism $ G(\AAA)^1\times A_G^\infty\to G(\AAA)$ is compatible with the  Haar measures (on the source we put the product of the Haar measures).   We also set
  \begin{align*}
       [G]^1=G(F)\back G(\AAA)^1.
  \end{align*}
  Both  $[G]$ and  $[G]^1$ are  equipped with the quotient measure by the counting measure. More generally we introduce  $[G]^1_P= M_P(F)N_P(\AAA)\back P(\AAA)^1K$. We have an obvious surjective map $[M_P]^1\times K\to [G]^1_P$ and we use it to push-forward the product measure on  $[M_P]^1\times K$.
\end{paragr}

\begin{paragr}   Let $\ggo_\infty$ be the Lie algebra of $G(F_\infty)$ and $\uc(\ggo_\infty)$ be the enveloping algebra of its complexification and $\zc(\ggo_\infty)\subset \uc(\ggo_\infty)$ be its center. 
\end{paragr}

\begin{paragr}  We fix a height $\|\cdot\|:G(\AAA)\to \RR_+$ as in \cite[§  I.2.2]{MW}. A function $f:G(\AAA)\to \CC$  is said to be smooth if it is right invariant by a compact-open subgroup $J$ of $G(\AAA_{f})$ and for every $g_f\in G(\AAA_{f})$ the function $g_\infty\in G(F_\infty)\mapsto f(g_fg_\infty)$ is smooth in the usual sense. 
\end{paragr}

\begin{paragr}[Level.] ---   By a level $J$ (of $G$) we mean  a normal  open compact subgroup of $K_f$.
\end{paragr}

\begin{paragr}[Schwartz algebra.] ---   This is the algebra (for the convolution product) denoted by  $\Sc(G(\AAA))$ and defined, as a topological vector space, as the locally convex topological direct limit  of the spaces  $\Sc(G(\AAA),C,J)$ over all pairs $(C,J)$ consisting of a compact subset of $G(\AAA_f)$ and a level $J$  where $\Sc(G(\AAA),C,J)$ is  the space of smooth functions $f:G(\AAA)\to \CC$ which are
\begin{itemize}
\item biinvariant by $J$ and  supported in the subset $G(F_\infty)\times C$ 
\item such that 
\begin{align*}
\|f\|_{r,X,Y}=\sup_{g\in G(\AAA)} \|g\|^r  |(R(X)L(Y)f)(g)|<\infty
\end{align*}
for every integer $r\geq 1$ and $X,Y\in \uc(\ggo_\infty)$. 
\end{itemize}
The family of semi-norms  $\|\cdot\|_{r,X,Y}$ defines the topology on $\Sc(G(\AAA),C,J)$. We denote by $\Sc(G(\AAA))^J$ the subalgebra of  $J$-biinvariant function.
\end{paragr}

\begin{paragr} Let $P$ be a standard parabolic subgroup of $G$. For any $g\in G(\AAA)$, we define
  $$\|g\|_P=\inf_{\delta\in M_P(F)N_P(\AAA)} \|\delta g\|.$$

  For all $N\in \ZZ$,  $X\in \uc(\ggo_\infty)$ and  any smooth function $\varphi: [G]_P\to \CC$ we define
$$ \| \varphi\|_{N,X}=\sup_{x\in [G]_P}\|x\|_P^{N}\lvert (R(X)\varphi)(x)\rvert \in \RR_+ \cup\{+\infty\}.$$
\end{paragr}

\begin{paragr}[Functions of uniform moderate growth.] --- For any $N\geq 0$ and any  level $J$   let  $\tc_N([G]_P)^J$ be the space of smooth functions $\varphi: [G]_P\to \CC$ that are right-invariant by $J$ and such that  that for every $X\in \uc(\ggo_\infty)$ we have $ \| \varphi\|_{-N,X}<\infty$. Then $\tc_N([G]_P)^J$ is a Fréchet space equipped with these  semi-norms. We defined  $\tc_N([G]_P)$ as the  locally convex topological direct limit over $J$ of the spaces $\tc_N([G]_P)^J$. The space $\tc([G]_P)$ of functions of uniform moderate growth on $[G]_P$ is the  locally convex topological direct limit of the spaces $\tc_N([G]_P)$ over the integers $N\geq 0$.
\end{paragr}

\begin{paragr}[Rapidly decreasing functions.] ---  For any  level $J$ let $\Sc([G]_P)^J$ be  the space of smooth functions $\varphi:[G]_P\to \CC$ that are right-invariant by $J$ and such that for every $N>0$ and $X\in \uc(\ggo_\infty)$ we have  $ \| \varphi\|_{N,X}<\infty$. Then $\Sc([G]_P)^J$ is a Fréchet space equipped with these  semi-norms. We define the Schwartz space  $\Sc([G]_P)$ as  the  locally convex topological direct limit of the spaces $\Sc([G]_P)^J$ over the subgroups $J$.
\end{paragr}

\begin{paragr}[Automorphic forms.] ---  The space $\Ac_P(G)$ of {\em automorphic forms} on $[G]_P$ is the subspace of $\zc(\ggo_\infty)$-finite functions in $\tc([G]_P)$. Let $\Ac_{P}^0(G)$ be the subspace of $\varphi\in \Ac_P(G)$ such that 
\begin{align*}
  \forall g\in G(\AAA), \forall a\in A_P^\infty, \ \ \varphi(ag)=\exp(\bg \rho_P,H_P(a)\bd) \varphi(g).
\end{align*}
Let   $\Ac_{P,\disc}(G)$ be the subspace of $\varphi\in \Ac_P(G)$ such that  
\begin{align*}
  \forall g\in G(\AAA), \forall a\in A_P^\infty, \ \ |\varphi(ag)|=\exp(\bg \rho_P,H_P(a)\bd) |\varphi(g)|
\end{align*}
and  the Petersson norm defined by
 \begin{align*}
   \|\varphi\|_{P}^2=\bg \varphi,\varphi\bd_P= \int_{[G]_{P,0}} |\varphi(g)|^2 \, dg 
 \end{align*}
 is finite. Let $\Ac_{P,\disc}^0(G)=\Ac_{P,\disc}(G)\cap \Ac_{P}^0(G)$. For any ideal  $\Jc\subset \zc(\ggo_\infty)$  of finite codimension, let $\Ac_{P, \disc,\Jc}(G)$ be the subspace of $\varphi\in \Ac_{P,\disc}(G)$ killed by $\Jc$. There exists $N\geqslant 1$ such that $\Ac_{P, \disc,\Jc}(G)$ is a closed subspace of $\tc_N([G]_P)$. Then $\Ac_{P, \disc,\Jc}(G)$  is equipped  with the induced topology from $\tc_N([G]_P)$. This topology does not depend on the choice of $N$. Then $\Ac_{P,\disc}(G)$ is provided with the  locally convex direct limit topology.  The group $G(\AAA)$ acts on  $\Ac_{P,\disc}(G)$ by right translation. When $P=G$, we omit the subscript $P$. 
\end{paragr}

\begin{paragr}[Discrete automorphic representations.] ---    \label{S:disc-aut}A discrete automorphic representation  $\pi$ of $G(\AAA)$ is a topologically irreducible subrepresentation of $\Ac_{\disc}(G)$.  Let   $\Pi_{\disc}(G)$  be the set of  discrete  automorphic representations $\pi$ of $G(\AAA)$ that are subrepresentations of $\Ac_{\disc}^0(G)$. For such a $\pi$, we denote by $\Ac_{\pi}(G)$ the $\pi$-isotypic component of  $\Ac_{\disc}^0(G)$. Let $\pi\in \Pi_{\disc}(M_P)$ and let  $\Ac_{P,\pi}(G)$ be the subspace of $\varphi\in \Ac_{P,\disc}^0(G)$ such that the map $m\in [M_P]\mapsto  \exp(-\bg \rho_P,H_P(m)\bd) \varphi(m g)$ belongs to $\Ac_{\pi}^0(M)$ for all $g\in G(\AAA)$. Then $\Ac_{P,\pi}(G)$ is a closed subspace of $\Ac_{P,\disc}^0(G)$ equipped with the induced topology. For any $\lambda\in \ago_{P,\CC}^*$, we define $\pi_\lambda=\pi\otimes \exp(\bg\lambda,H_M(\cdot)\bd)$. The map $\varphi\mapsto\varphi_\la$, defined by 
$$\forall g\in G(\AAA)\ \ \varphi_\la(g)= \exp(\bg\lambda,H_P(g)\bd) \varphi(g),$$
identifies $\Ac_{P,\pi}(G)$ with a subspace of  $\Ac_{P}(G)$. By transport of structure, we denote by $I_{P,\pi}(\la)$, or simply $I_P(\la)$ if the context is clear, the action of $G(\AAA)$ on $\Ac_{P,\pi}(G)$ we get from that on $\Ac_{P}(G)$ by right translation. The right convolution gives an action of $\Sc(G(\AAA))$ denoted by $I_{P,\pi}(\la,f)$, or simply $I_P(\la,f)$ for $f\in \Sc(G(\AAA))$.

 For any $\varphi\in \Ac_P(G)$ and any standard parabolic subgroup $Q\subset P$ we define the constant term along $Q$ by
\begin{align*}
  \varphi_Q(g)=\int_{[N_Q]}\varphi(ng)\, dn, \ g\in G(\AAA)
\end{align*}
Let $\Ac_{\cusp}^0(G)\subset \Ac^0(G) $ be the subspace of functions whose constant terms vanish for all proper parabolic subgroups. Let  $\Pi_{\cusp}(G)$ be the set  of topologically irreducible subrepresentations of $\Ac_{\cusp}^0(G)$.

Let $J$ be a level and let $\Ac_{P,\pi}(G)^J\subset \Ac_{P,\pi}(G)$ be the subspace of right-$J$-invariant functions. We define the subset  $\Pi_{\disc}(M_P)^J\subset \Pi_{\disc}(M_P)$  by the condition 
$\Ac_{P,\pi}(G)^J\not=\{0\}$.

Let $\hat K_\infty$ be the set of isomorphism classes of irreducible unitary representations of $K_\infty$. For any $\tau\in \hat K_\infty$, let $\Ac_{P,\pi}(G)^\tau$ be the subspace of functions whose right-$K_\infty$-translate belongs to the $\tau$-isotypic component. If $J$ is a level we set
$$\Ac_{P,\pi}(G)^{\tau,J}=\Ac_{P,\pi}(G)^{\tau}\cap\Ac_{P,\pi}(G)^{J}.$$
\end{paragr}

\begin{paragr}[Intertwining operators.] ---   Let $P$ and $Q$ be standard parabolic subgroups of $G$. Let $w\in W(P,Q)$. For any $\la\in  \ago_{P,\CC}^{*}$, we have the intertwining operator
  \begin{align}\label{eq:NEW integ inter}
    M(w,\la):\Ac_P(G)\to \Ac_{Q}(G).
  \end{align}
 For any $\varphi\in \Ac_P(G)$ and any $\la\in \ago_{P,\CC}^{*}$ such that $\Re(\bg \la,\al^\vee\bd)$ is large enough for every $\al\in \Delta_P$ such that $w\al$ is negative, it is defined by the integral
\begin{align}\label{eq:integ inter}
   (M(w,\la)\varphi)_{w\la }(g)=   \int_{  (N_Q\cap w N_Pw^{-1})(\AAA) \back N_Q(\AAA) } \varphi_\la(w^{-1}ng)\, dn.
\end{align}
In general, it is given on $K$-finite functions by meromorphic continuation from the previous integral. By \cite[corollary 2.4]{Lap-remark}, it extends to a continuous operator on $\Ac_P(G)$. Note that  $M(w,\la)$ maps $\Ac_P^0(G)$ into $\Ac_{Q}^0(G)$.

We shall use the following simple formula for $\la,\mu\in \ago_{P,\CC}^*$
\begin{align}\label{eq:torsion mu}
  M(w,\la)(\varphi_\mu)=(M(w,\la+\mu)\varphi)_{w\mu}.
\end{align}

One has  $w\Delta_0^{P}=\Delta_0^Q$. Let $Q'\subset Q$ and $P'\subset P$ be  standard parabolic subgroups such that $w\Delta_0^{P'}=\Delta_0^{Q'}$. Using the subscript $P'$ or $Q'$ to denote the constant term we have
\begin{align}\label{eq:cst intert}
  (M(w,\la)\varphi)_{Q'}=M(w,\la)(\varphi_{P'}).
\end{align}

\end{paragr}

\begin{paragr}[Eisenstein series.] ---\label{S:Eis-series}
  Let $P\subset Q$ be   standard parabolic subgroups of $G$. For any $\varphi\in \Ac_{P,\disc}(G)$, $g\in G(\AAA)$ and $\la\in \ago_{P,\CC}^{*}$, we denote by
  \begin{align*}
    E^Q(g,\varphi,\la)=\sum_{\delta\in P(F)\back Q(F)} \varphi_\la(\delta g)
  \end{align*}
the Eisenstein series where the right-hand side is convergent for $\Re(\la)$ in a suitable cone. For $K$-finite functions it admits a meromorphic  continuation. This extends to $\Ac_{P,\disc}(G)$  and moreover $E^Q(\varphi,\la)$ has only hyperplane singularities,  see § \ref{S:py-exp}. Moreover, for  $\pi$  a discrete automorphic representation, the map $\varphi\in \Ac_{P,\pi}(G)\to E^Q(\varphi,\la)_{|M_Q(\AAA)}  \in \tc([M_Q])$ is continuous, see \cite[theorem 2.2]{Lap-remark}.

Let $E_Q^G$ be the constant term of $E^G$  along $Q$. We have for all $g\in G(\AAA)$
  \begin{align}
\nonumber    E_Q^G(g,\varphi,\la)&=\int_{[N_Q]}  E^G(ng,\varphi,\la)\,dn\\
\label{eq:cst-term-Eis}&=\sum_{w\in \, _QW_P} E^Q(g,M(w,\la)\varphi_{P_w},w\la). 
  \end{align}
where $M(w,\la)\varphi_{P_w}\in \Ac_{Q_w}$, see e.g. \cite[lemma 6.10]{BL-mero}. In general, we shall omit the upper script $G$ if the context is clear.
\end{paragr}

\section{On discrete Eisenstein series}\label{sec:disc}

\subsection{Scalar product of two truncated discrete Eisenstein series}\label{ssec:truncated-scalar}

\begin{paragr}In this section, the group $G$ is $G_n=\GL_F(n)$ for some $n\geq 1$. We denote by  $P_0$  the standard Borel subgroup and $T_0$ the diagonal maximal torus. The group $K$ is the standard maximal compact subgroup  of $G_n(\AAA)=\GL_F(n,\AAA)$. We identify naturally $\ago_0$ with $\RR^n$ equipped with the canonical scalar product. Thus $\ago_{0,\CC}^*$ and $\ago_{0,\CC}$ are identified with $\CC^n$ equipped with the canonical definite positive hermitian form whose associated norm is denoted by $\|\cdot\|$. We shall freely use the truncation operator $\La^{T,Q}$ attached to a truncation parameter $T$ and a parabolic subgroup $Q$ introduced by Arthur in \cite{ar2}.
  
\end{paragr}

\begin{paragr}\label{S:exposant-disc}
  Let $P=MN_P$ be a standard parabolic subgroup of $G$. Let $\pi\in \Pi_{\disc}(M)$. We write $M=G_{n_1}\times \ldots \times G_{n_k}$ with $n_1+\ldots +n_k=n$ and $\pi=\pi_1\boxtimes \ldots \boxtimes \pi_k$ accordingly. By the classification of discrete automorphic representations of general linear groups  (see \cite{MW}),  there exist integers $r_i,d_i\geq 1$ such that $n_i=r_id_i$ and $\sigma_i\in \Pi_{\cusp}(G_{r_i})$ such that elements of $\Ac_{\pi_i}(G_{n_i})$ are obtained as residues of Eisenstein series built from elements of $\Ac_{P_{\pi_i},\sigma_i^{\otimes d_i}}$ where $P_{\pi_i}\subset G_{n_i}$ is the  standard parabolic subgroup of Levi factor $G_{r_i}^{d_i}$. Let $P_\pi=P_{\pi_1}\times\ldots \times P_{\pi_k}$.

  Let $P_\pi\subset R \subset Q\subset P$ be  parabolic subgroups. We write $Q\cap M=Q_1\times \ldots \times Q_k$ and $R\cap M=R_1\times \ldots \times R_k$ with $R_i\subset Q_i\subset G_{n_i}$.  We set
  $$\nu_{R,\pi}^Q=(-\rho_{R_i}^{Q_i}/r_i)_{1\leq i \leq r}$$
  written relative to $\ago_{R}^{Q,*}=\oplus_{i=1}^k \ago_{R_i}^{Q_i,*}$. It will be useful later to have an explicit expression for $\nu_{R,\pi}^P$. Let $1\leq i\leq k$. Since we have $P_{\pi_i}\subset R_i$, the standard Levi factor of $R_i$ may be  identified with $G_{e_1 r_i}\times \ldots \times G_{e_k r_i}$ for some integers $(e_j)_{1\leq j \leq k}$ such that $e_1+\ldots+e_k=d_i$. Using the basis given by the determinants of  the blocks  $G_{e_j r_i}$ for $1\leq j \leq k$ we identified $\ago_{R_i}^* $ with $\RR^k$ and we get:
  \begin{align}
    \label{eq:calcul de nu}
    -\rho_{R_i}^{  G_{n_i}}/r_i=(\frac{e_1-d_i}{2}, \frac{(2e_1+e_2)-d_i}{2},\frac{(2e_1+2e_2+e_3)-d_i}{2},\ldots, \frac{d_i-e_{k}}{2}).
  \end{align}

  For convenience, we also set  $\nu_{R,\pi}^Q=0$ for standard parabolic subgroups $R\subset Q$ such that $P_\pi\subsetneq R$.  If $Q=P$ we omit the upperscript $P$ that is $\nu_{R,\pi}=\nu_{R,\pi}^P$. Abusing the notation, we shall also omit the subscript $\pi$ if the context is clear.

  Let $\varphi\in \Ac_{P,\pi}(G)$. The reason to introduce $\nu_R$ is that $\varphi_{R,-\nu_R}$ belongs to $ \Ac_R^{0}(G)$  as one can check. Note that $ \varphi_R(g)=0$   unless $P_\pi\subset R$.  Moreover $\varphi_{R,-\nu_R}$ is also square-integrable as it is shown in the next lemma.

  \begin{lemme}\label{lem:comparison-norm}
    Let $J$ be a level. For any large enough $N>0$ there exists a finite family $(X_i)_{i\in I}$ of elements of $\uc(\ggo_\infty)$ such that for all $(P,\pi)$ as above, all $P_\pi\subset Q\subset P$ and all $\varphi\in \Ac_{P,\pi}(G)^J$ we have:
    \begin{align*}
      \|\varphi_{Q,-\nu_Q}\|_Q \leq \sum_{i\in I}  \|\varphi\|_{-N,X_i}.
    \end{align*}
  \end{lemme}

  \begin{preuve}
We have to consider
    \begin{align}
      \nonumber \|\varphi_{Q,-\nu_Q}\|_Q^2=\int_{[G]_{Q,0}}     |\varphi_Q(g) |^2   \exp(-\bg 2\nu_Q,H_Q(g)\bd) \, dg\\
\label{eq:norm-phiQ}      = \int_{[M_Q]_0} \int_K       | \varphi_Q(mk) |^2   \exp(-\bg 2\nu_Q+2\rho_Q,H_Q(m)\bd) \, dk dm
    \end{align}
    We fix a truncation parameter $T$ and we use the inversion formula, see \cite[lemma 1.5]{ar2}:
    \begin{align*}
      \varphi_Q(g)=\sum_{P_\pi  \subset R\subset Q}\sum_{\delta\in R(F) \back Q(F)} \tau_R^Q(H_R(\delta g)-T) \La^{T,R}\varphi_R (\delta  g).
    \end{align*}
   Let $R$ be as in the sum above. We consider the contribution of $R$ in \eqref{eq:norm-phiQ} namely:
    \begin{align*}
      \int_{ (R\cap M_Q)(F) A_Q^\infty\back M_Q(\AAA)} \int_K     \tau_R^Q(H_R(m)-T) \La^{T,R}\varphi_R (mk) \overline{\varphi_Q(mk)}  \exp(-\bg 2\nu_Q+2\rho_Q,H_Q(m)\bd) \, dk dm\\
      = \int_{[M_R]^1    }  \int_{A_R^{Q,\infty}} \tau_R^Q(H_R(a)-T)  \exp(-\bg 2\nu_Q+2\rho_R,H_R(a)\bd)    \int_K\La^{T,R}\varphi_R (amk) \overline{\varphi_R(amk)}\, dk da dm.
    \end{align*}
    This last expression is the product of

    \begin{align}
\label{eq:interg-nu}      &\int_{A_R^{Q,\infty}} \tau_R^Q(H_Q(a)-T)  \exp(\bg 2\nu_R^Q,H_Q(a)\bd)da \\
\nonumber&      \text{and}\\
  \label{eq:produit-scal} &\int_{[M_R]^1    }      \int_K\La^{T,R}\varphi_R (mk) \overline{\varphi_R(mk)}\, dk dm.
    \end{align}
    The integral \eqref{eq:interg-nu} is convergent. When the data $\pi,Q,R$ vary,  its values belong to a finite set and  thus there is an absolute bound for \eqref{eq:interg-nu}.
    Let $r\in \NN$ large enough such that $\int_{[M_R]^1    }   \|m\|^{-r} dm<\infty$. Let $N>0$.
    By \cite[lemma 1.4]{ar2},  there exists a finite family $(Y_i)_{i\in I'}$ of element of $\uc(\mgo_{R,\infty})$  such that  for any smooth function $\psi$ on  $[M_R]$ that is right-invariant by $J\cap M_R(\AAA)$
    $$\forall m\in [M_R]^1 \ \ \ |\La^{T,R}\psi (m) |\leq   \left(\sum_{i\in I'} \| \psi\|_{-N,Y_i} \right) \|m\|_R^{-N-r}.$$
    We deduce that there exists      a finite family $(Y_i)_{i\in I'}$ as above (maybe larger) such that for all $\psi$ as above
    $$    \int_{[M_R]^1    }     |\La^{T,R}\psi (m) \overline{\psi(m)}|\, dm \leq   (\sum_{i\in I'} \| \psi\|_{-N,Y_i})^2
    $$
    To conclude it suffices to observe that there exists $N' >0$ and a finite family $(X_i)_{i\in I}$ of elements of $\uc(\ggo_\infty)$ such that for all $k\in K$
    $$\sum_{i\in I'} \| \varphi_R(\cdot k)\|_{-N,Y_i} \leq \sum_{i\in I}  \|\varphi\|_{-N',X_i}.$$
Indeed the elements $kY_i k^{-1}$ stay in a finite dimensional space with bounded coefficients in a fixed basis.
      \end{preuve}
  
\end{paragr}

\begin{paragr}
  Let $\varphi, \psi\in \Ac_{P,\pi}(G)$. The goal of this subsection is to prove the following theorem which gives an exact expression in terms of intertwining operators of the following integral:
  \begin{align}\label{eq:scalar-Eis}
\bg   \Lambda^{T} E(\varphi,\la), \overline{  E(\psi,\la')}\bd_G=   \int_{[G]_0} \Lambda^{T} E(g,  \varphi,\la)    E(g,  \psi,\la')\, dg
  \end{align}
  for all $\la,\la'\in \ago_{P,\CC}^{G,*}$ that are not singular for the corresponding Eisenstein series. Note that for such $\la,\la'$ the integral is absolutely convergent by the basic properties of Arthur's truncation operator.

\begin{theoreme} For any $\la,\la'$  in general position in $\ago_{P,\CC}^{G,*}$, we have:
  \label{thm:inversion-calculee}
  \begin{align*}
    &\bg   \Lambda^{T} E(\varphi,\la), \overline{  E(\psi,\la')}\bd_G\\
    &=  \sum_{R\subset G} \sum_{ w,w'} \bg M(w,\la+\nu_{P_w})\varphi_{w} ,  \overline{M(w',\la'+\nu_{P_{w'}})\psi_{w'}}\bd_{R} \frac{\exp(\bg w\la+w'\la'+w\nu_{P_w} +w'\nu_{P_{w'}} ,T_R \bd)}{\theta_R(w\la+w'\la'+w\nu_{P_w} +w'\nu_{P_{w'}})} 
  \end{align*}
  where the sum is over $w,w'\in \, _RW_{P}$ such that $  P_\pi\subset P_w\cap P_{w'} $ and $R_w=R$ and $R_{w'}=R$.

  For any standard parabolic subgroup $Q$ of $G$ and any  $w\in \, _QW_{P}$ we have set   $\varphi_{w}=    \varphi_{P_w,-\nu_{P_w}}$ and $\psi_{w'}=    \psi_{P_{w'},-\nu_{P_{w'}}}$.
\end{theoreme}

\begin{remarque}
  In the inner sum, the pairing is vanishing unless $P_\pi\subset P_w\cap P_{w'}$. In this case the corresponding denominators are non-vanishing for $\la,\la'$ in general position, see lemma \ref{lem:non-vanish-theta}.
\end{remarque}

The proof of theorem \ref{thm:inversion-calculee} is given in § \ref{S:preuve-inv} and  is based on lemmas given in §  \ref{S:lemmas}. With the notations of theorem \ref{thm:inversion-calculee}, using \eqref{eq:torsion mu}, we see that  the formula \eqref{eq:cst-term-Eis} takes the following shape:
\begin{align}\label{eq:cst-term-Eis NEW}
   E_Q^G(g,\varphi,\la)=\sum_{w\in \, _QW_P} E^Q(g,M(w,\la+\nu_{P_w})\varphi_{w},w(\la+\nu_{P_w})). 
  \end{align}

\end{paragr}

\begin{paragr}[Proof of theorem \ref{thm:inversion-calculee}.] ---  \label{S:preuve-inv}We consider $T$ a truncation parameter and $T'\in \ago_0^+$. We start from the following formula for any function $\phi$  on $[G]_0$ and $g\in G(\AAA)$:
  \begin{align*}
    (\La^{T+T'}\phi)(g)=\sum_{Q} \sum_{\delta\in Q(F)\back G(F)} \Ga_Q(H_Q(\delta g)-T,T')(\La^{T,Q}\phi)(\delta g)
  \end{align*}
  where the sum is over standard parabolic subgroups $Q$ and $\Ga_Q(H,X)$ satisfies for any parabolic subgroup $P\subset G$, see \cite[p.13]{ar-inv},
\begin{align*}
  \hat\tau_P(H-X)=\sum_{ P\subset Q} (-1)^{\dim(\ago_Q^G)} \hat\tau_P^Q(H) \Ga_Q(H,X).
\end{align*}
The formula is  easily deduced from the equality above and the definition of the truncation operator in \cite{ar2}. We apply the formula to the truncated Eisenstein series in the bracket. Using the computation \eqref{eq:cst-term-Eis NEW} of the constant term of the Eisenstein series we find:
\begin{align*}
&  \bg   \Lambda^{T+T'} E(\varphi,\la), \overline{  E(\psi,\la')}\bd_G=\\
&\sum_Q \sum_{w,w'\in \, _QW_{P}   } \bg  \Lambda^{T,Q} E^Q( M(w,\la+\nu_{P_w})\varphi_{w},(w(\la+\nu_{P_w}))^Q) , \overline{  E^Q(  M(w',\la'+\nu_{P_{w'}})\psi_{w'},(w'(\la'+\nu_{P_{w'}}))^Q)} \bd_Q \times\\
  &\int_{\ago_Q^G}  \Ga_Q(H-T,T') \exp(\bg w(\la +\nu_{P_w})  +w'(\la'+\nu_{P_{w'}}),H \bd)\, dH.
\end{align*}

By \cite[lemma 2.2]{ar-inv}, the last line above is equal to:
\begin{align*}
  \exp(\bg w(\la +\nu_{P_w})  +w'(\la'+\nu_{P_{w'}} ),T_Q \bd) \cdot \sum_{Q\subset R} (-1)^{\dim(\ago_Q^R)} \frac{\exp(\bg w(\la+\nu_{P_w} )+w'(\la'+\nu_{P_{w'}}) ,T'_R \bd)}{(\hat\theta_Q^R\theta_R)(w(\la+\nu_{P_w})+w'(\la' +\nu_{P_{w'}}))}.
\end{align*}
Note that for $\la,\la'$ in general position the right-hand side is well-defined, see lemma \ref{lem:non-vanish-theta} below, provided $P_\pi\subset P_w\cap P_{w'}$. We shall use this expression and invert the sums over $Q$ and over $R$. To do this we use lemma \ref{lem:ww'}. In this way we get:

   \begin{align*}
 &\bg   \Lambda^{T+T'} E(\varphi,\la), \overline{  E(\psi,\la')}\bd_G\\
&=\sum_R\sum_{w,w'\in \, _RW_{P},  P_\pi\subset P_w\cap P_{w'}  } \frac{\exp(\bg w(\la+\nu_{P_w} )+w'(\la'+\nu_{P_{w'}}) ,(T+T')_R \bd)}{\theta_R(w(\la+\nu_{P_w} )+w'(\la'+\nu_{P_{w'}}))}  \pc^{T,R}(\varphi,\psi,\la,\la',w,w')
   \end{align*}
where we set:
\begin{align*}
&  \pc^{T,R}(\varphi,\psi,\la,\la',w,w')=\sum_{Q\subset R}    (-1)^{\dim(\ago_Q^R)} \sum_{w_1,w_1'} \\
&   \bg  \Lambda^{T,Q} E^Q( M(w_1,\la+\nu_{P_{w_1}})\varphi_{w_1},(w_1(\la+\nu_{P_{w_1}}))^Q) , \overline{  E^Q(  M(w'_1,\la'+\nu_{P_{w'_1}})\psi_{w'_1},(w'_1(\la'+\nu_{P_{w_1'}}))^Q)} \bd_Q  \\
&\times\frac{\exp(\bg w_1(\la+  \nu_{w_1}^w ) + w_1'   (\la'+\nu_{w_1'}^{w'}),T^R\bd)}{\hat\theta_Q^R( w_1(\la+  \nu_{w_1}^w ) + w_1'   (\la'+\nu_{w_1'}^{w'}))}
\end{align*}
 where the inner sum is over $w_1\in \, _QW^R_{R_w}w $ and $w_1'\in \, _QW^R_{R_{w'}}w'$  and we set  $\nu_{w_1}^w =\nu_{P_{w_1}}^{P_w}$. More generally we shall use  the expression $\pc^{T,R}(\varphi,\psi,\la,\la',w,w')$ for elements  $\la\in \ago_{P_w,\CC}^{G,*}$ and $\la'\in \ago_{P_{w'},\CC}^{G,*}$ in general position. In fact it does not depend on  $T$: indeed this is nothing else but the regularization in the sense of Jacquet-Lapid-Rogawski, see \cite[section 12]{JLR} and  also \cite[§  5.2]{Lap-asym}, of the (in general non-convergent) pairing:
 \begin{align}\label{eq:inner-product}
 \bg E^R(M(w,\la+\nu_{P_{w}})\varphi_{w},(w(\la+\nu_{P_{w}}))^R), \overline{E^R( M(w',\la'+\nu_{P_{w'}})\psi_{w'},(w'(\la'+\nu_{P_{w'}}))^R)} \bd_R.
 \end{align}
The functions $\varphi_{w}$ and $\psi_{w'}$ belong respectively to subspaces  for some discrete representations $\sigma\in \Pi_{\disc}(M_{P_w})$ and $\sigma'\in \Pi_{\disc}(M_{P_{w'}})$. Then \eqref{eq:inner-product} is closely related to the regularized version of the pairing
 \begin{align*}
 (\varphi',\psi')\in \Ac_{R_w, w\sigma}(G)\times \Ac_{R_{w'},w' \sigma'}\mapsto \bg E^R(\varphi',\la^R), \overline{E^R( \psi',\la'^R)} \bd_R
 \end{align*}
 which is well-defined for $\la\in \ago_{R_w,\CC}^{G,*}$ and $\la'\in \ago_{R_{w'},\CC}^{G,*}$ in general position: moreover it gives a  pairing on the product of the spaces $\Ac_{R_w,w\sigma}(G)\times \Ac_{R_{w'},w'\sigma'}(G)$ which is invariant for the diagonal restriction to $G(\AAA)$ of the action $I_{R_w,w\sigma}(\la^R)\times I_{R_{w'},w'\sigma'}(\la'^R)$. By Bernstein argument, see \cite[p. 208]{JLR}, such a regularized pairing must vanish unless $R_w=R$ and $R_{w'}=R$. We deduce that the regularized version of \eqref{eq:inner-product} must also vanish unless the same condition is satisfied: then it reduces to the (convergent) pairing:
\begin{align}\label{eq:inner-product2}
 \bg M(w,\la+\nu_{P_{w}})\varphi_{w}, \overline{M(w',\la'+\nu_{P_{w'}})\psi_{w'}} \bd_R.
 \end{align}
This finishes the proof.
\end{paragr}

\begin{paragr}[Some lemmas.] ---\label{S:lemmas}
  
  \begin{lemme}\label{lem:non-vanish-theta}
    Let $Q\subset R$ be standard parabolic subgroups and $w,w'\in \, _QW_{P}$ such that $P_\pi\subset P_w\cap P_{w'}$. Then  the map
    \begin{align*}
      (\la,\la')\mapsto (\hat\theta_Q^R\theta_R)(w\la+w'\la'+w\nu_{P_w} +w'\nu_{P_{w'}})
    \end{align*}
does not vanish identically on $i\ago_{P}^{G,*}\times i\ago_{P}^{G,*}$.
  \end{lemme}

  \begin{preuve}
    The statement is obvious if the map $(\la,\la')\mapsto (\hat\theta_Q^R\theta_R)(w\la+w'\la')$ does not vanish identically on $i\ago_{P}^{G,*}\times i\ago_{P}^{G,*}$. Otherwise there exists $\gamma^\vee\in \hat\Delta_Q^{R,\vee} \cup \Delta_R^\vee$ such that 
    \begin{align*}
      \gamma^\vee\in w\ago_0^{P}\cap w'\ago_0^{P}.
    \end{align*}
    One has to consider two cases:
    \begin{itemize}
    \item $\gamma^\vee\in \hat\Delta_Q^{R,\vee} $. In this case, there exists a maximal parabolic subgroup $Q\subset S\subsetneq  R$ such that $\hat\Delta_S^{R,\vee}=\{\gamma^\vee\}$. Since $\dim(\ago_S^R)=1$  there exists $c>0$ such that  $c \gamma^\vee\in \Delta_S^{R,\vee}\subset  \Delta_S^{\vee} $.
    \item  $\gamma^\vee\in \Delta_R^{\vee} $.
    \end{itemize}
   In both cases we can apply lemma \ref{lem:negativity} below to conclude that 
\begin{align*}
  \bg w\nu_{P_w} +w'\nu_{P_{w'}},\gamma^\vee \bd <0.
\end{align*}

  \end{preuve}

  \begin{lemme} \label{lem:negativity}Let $Q\subset R$ be standard parabolic subgroups of $G$   and let $w\in \, _QW_{P}$ such that $P_\pi\subset P_w$. For any  $\gamma^\vee\in \Delta_R^\vee\cap  w\ago_0^{P}$ we have 
\begin{align*}
  \bg  w\nu_{P_w} ,\gamma^\vee\bd <0.
\end{align*}
  \end{lemme}

  \begin{preuve} First we observe that $w^{-1}\gamma^\vee\in \ago_{P_w}^{P}$. Indeed, for $\al\in \Delta_0^{P_w}$ we have
    \begin{align*}
      \bg \al, w^{-1}\gamma^\vee\bd=\bg w\al, \gamma^\vee\bd=0
    \end{align*}
since $w\al\in w\Delta_0^{P_w}=\Delta_0^{Q_w}\subset \Delta_0^Q\subset \ago_0^{Q,*}\subset \ago_0^{R,*}$.  In particular, we may and shall assume $P_w\subsetneq P$. 

Let $R\subsetneq S \subset G$ be the standard parabolic subgroup defined by  $\Delta_R^{S,\vee}=\{\gamma^\vee\}$.  The Levi factor $M_R$ is identified with a product $G_{m_1}\times \ldots \times G_{m_r}$ of general linear groups with $m_1+\ldots+m_r=n$. There exists $1\leq i<r$ such that the Levi factor $M_S$  is identified with the  product $G_{m_1}\times \ldots \times G_{m_{i-1}} \times  G_{m_i+m_{i+1}} \times G_{m_{i+2}}\times \ldots \times G_{m_r}$. In the usual way, we identify $\ago_0$ with $\RR^n$. Then, we have
\begin{align}\label{eq:gamma}
  \gamma^\vee= (0,\ldots, 0, 1/m_i, \ldots, 1/m_{i}, -1/m_{i+1},\ldots, -1/m_{i+1}, 0, \ldots, 0)
\end{align}
where the entries are repeated respectively $m_1+\ldots+m_{i-1}$, $m_i$, $m_{i+1}$  and $m_{i+2}+\ldots+m_{r}$ times.

We identify  $M_P$ with the product $G_{n_1}\times \ldots \times G_{n_k}$ with $n_1+\ldots +n_k=n$. Then $P_w\cap M_P$ is written accordingly $S_1\times \ldots \times S_k$ where
$S_j\subset G_{n_j}$ is a standard parabolic subgroup. In this way,  we have
$$\ago_{P_w}^{P,*}=\oplus_{j=1}^k \ago_{S_j}^{G_{n_j},*}.$$
Then according to the computation of  § \ref{S:exposant-disc}, the component of $\nu_{P_w}$ on $\ago_{S_j}^{G_{n_j},*}$ is $-\rho^{G_{n_j}}_{S_j}/r_j$ for some divisor $r_j$ of $n_j$. Dually we have $\ago_{P_w}^{P}=\oplus_{j=1}^k \ago_{S_j}^{G_{n_j}}$. For $1\leq j \leq k$, let $\beta_j^\vee\in\ago_{S_j}^{G_{n_j}}$ be the component of $w^{-1}\gamma^\vee\in  \ago_{P_w}^{P}$ on $\ago_{S_j}^{G_{n_j}}$. Then it is enough to show the following claim: for any $\varpi\in \hat\Delta_0^{G_{n_j}}$ we have $\bg \varpi,\be_j^\vee\bd \geq 0$ (the subscript $0$ refers to the standard Borel subgroup of $G_{n_j}$). Indeed with the claim we get:
 \begin{align*}
   -\bg  w\nu_{P_w},\gamma^\vee\bd &=\sum_{j=1}^k \bg \rho^{G_{n_j}}_{S_j}/r_j, \be_j^\vee\bd\\
&=\sum_{j=1}^k \bg \rho^{G_{n_j}}_{0}, \be_j^\vee\bd/r_j\\
&=\sum_{j=1}^k \sum_{  \varpi\in \hat\Delta_0^{G_{n_j}} } \bg \varpi,\be_j^\vee\bd /r_j \geq 0.
 \end{align*}
 Moreover the inequality is strict because otherwise we would have $\bg \varpi,\be_j^\vee\bd =0$ for all $j$ and all $\varpi\in \hat\Delta_0^{G_{n_j}}$ and thus $\gamma^\vee=0$.

 Let us prove the claim.  We identify naturally  $\ago_{S_j}^{G_{n_j},*}$ with a subspace of the hyperplane of $\RR^{n_j}$ of tuples whose sum of coefficients vanishes. So $\beta_j^\vee$ is identified with an element of this hyperplane. According to \eqref{eq:gamma}, the possible entries of $\beta_j^\vee$ are $0,  1/m_i, -1/m_{i+1}$ with multiplicities denoted by $N_{j,0}, N_{j,+}$ and $N_{j,-}$. We have $N_{j,0}+N_{j,+}+N_{j,-}=n_j$. Because the sum of the entries is $0$, we have  $N_{j,+}/m_i=N_{j,-}/m_{i+1}$. Let $\varpi\in \hat\Delta_0^{G_{n_j}}$. We observe that  the pairing  $\bg \varpi,\be_j^\vee\bd$ is equal, up to a positive multiplicative constant, to the sum of the first $l$ entries of $\be_j^\vee$ for some $1\leq l \leq n_j-1$ determined by $\varpi$. To conclude it suffices to show that  the entries $1/m_i$ of $\be_j^\vee$ are ``before'' the entries  $-1/m_{i+1}$. If this were not true, we could find a positive root  $\al$  inside $G_{n_j}$ such that 
 \begin{align*}
   \bg \al, \be_j^\vee\bd = -( 1/m_i +1/m_{i+1}).
 \end{align*}
 We can also view $\al$ as a positive root   inside $M_P$ and we observe that $w\al$ is also positive since $w\in \, _QW_{P}$. Thus we have:
 \begin{align*}
   \bg \al, \be_j^\vee\bd =  \bg \al, w^{-1}\gamma^\vee\bd=  \bg w\al, \gamma^\vee\bd.
 \end{align*}
 However, for a positive root $\be$, the values of the pairing  $\bg \be, \gamma^\vee\bd$ belong to
 $$\{0, \pm  1/m_i ,  1/m_i +1/m_{i+1}, \pm 1/m_{i+1}\}$$
which does not contain $-( 1/m_i +1/m_{i+1})$. This is the  contradiction we were looking for.
\end{preuve}

\end{paragr}

\subsection{Bounds for intertwining operators}\label{ssec:bound-intert}

\begin{paragr} The goal  of this subsection is to state proposition \ref{prop:bound-intertwining} below which provides useful bounds on intertwining operators.  The proof of this proposition will be dealt with in subsections \ref{ssec:normalization inter} to \ref{ssec:bd norm inter}. We continue with the notations of subsection \ref{ssec:truncated-scalar}. The other notations are borrowed from section \ref{sec:prelim}.
\end{paragr}

\begin{paragr}[Numerical invariants.] --- \label{S:numerical}We denote by $\Om_G$ and $\Om_{K_\infty}$ the Casimir operators of $G$ and $K_\infty$ respectively associated to the standard bilinear form on $\ggo_\infty$ associated to the trace, see e.g. \cite[p.323]{Muller98}. We set
  \begin{align*}
    \Delta=1-\Om_G+2\Om_{K_\infty}\in \uc(\ggo_\infty).
  \end{align*}
For any $\tau\in \hat K_\infty$, let $\la_\tau$ be the eigenvalue of $\Om_{K_\infty}$. Let $P=MN_P$ be a standard Levi subgroup.  We define similarly the Casimir operator $\Om_M$. Let $\pi_\infty$ be an irreducible unitary representation of $M(F_\infty)$.  We shall attach two invariants to $\pi_\infty$. Let $\la_{\pi_\infty}$ be the eigenvalue of $\Om_M$ on $\pi_\infty$. We set:
  \begin{align}
    \label{eq:LaG}
\La_{\pi_\infty}^M= \sqrt{ \la_{\pi_\infty}^2+\la_\tau^2 }
  \end{align}
  where $\tau$ is a minimal $K_\infty\cap M(F_\infty)$-type of $\pi_\infty$ and
\begin{align}
    \label{eq:LaGP}
\La_{\pi_\infty}^G= \min_\tau (\sqrt{ \la_{\pi_\infty}^2+\la_\tau^2 })
  \end{align}
where the minimum is taken over the finite set of minimal $K_\infty$-types of the induced representation $\Ind_{P(F_\infty)}^{G(F_\infty)}(\pi_\infty)$. This invariant was introduced by Müller in \cite[p. 695]{Muller02}. If $M=G$, the invariants \eqref{eq:LaG} and \eqref{eq:LaGP} are the same as the notation suggests.
Let $\pi\in \Pi_{\disc}(M)$ and let $\pi_\infty$ be the Archimedean component of $\pi$. We set 
\begin{align}
\la_\pi=\la_{\pi_\infty},  \ \La_\pi^M=\La_{\pi_\infty}^M ,  \ \La_\pi=\La_{\pi_\infty}^G.
\end{align}

\begin{remarque}
  \label{rq:La invariant} The referee asked whether it would be possible to replace the invariant $\La_{\pi}$ by the norm of the infinitesimal character $\chi_\pi$ of $\pi_\infty$. It follows from \cite[p.611]{FLM-limit-mult} and \cite[Remark 2.8]{FL-Iran} that there exist constants $c_1,c_2>0$ depending only on $G$ such that for all $\pi$ as above we have
  \begin{align*}
    c_1(1+\|\chi_\pi\|^2)^2 \leq 1+\La_\pi^2 \leq   c_2(1+\|\chi_\pi\|^2)^2.
  \end{align*}
  In the statement of our main results,  most of the time, we could have used $\|\chi_\pi\|$ instead of $\La_\pi$ (or  $\La_\pi^M$) as we did in the introduction. However,  for convenience of references, we prefer to stick to Müller's invariant.
\end{remarque}

\end{paragr}

\begin{paragr}[Pairs and triples.] ---  \label{S:J-triple}For any level $J$, we consider a $J$-pair that is a pair $(P,\pi)$ such that  $P=MN_P$ is a standard parabolic subgroup and $\pi\in \Pi_{\disc}(M)^{J\cap M(\AAA_f)}$. A $J$-triple, or simply a triple if the context is clear, is a triple $(P,\pi,\tau)$ such that  $(P,\pi)$ is a $J$-pair,  $\tau\in \hat K_\infty$  and  $\Ac_{P,\pi}(G)^{\tau,J}\not=\{0\}$. Let  $e_\tau$ be the measure supported on  $K_\infty$  given by $\deg(\tau) \trace(\tau(k)) dk$ where $dk$ is the Haar measure on $K_\infty$ giving the total volume $1$. We have $e_\tau*e_\tau=e_\tau$. In particular, for $\varphi\in \Ac_{P,\pi}(G)^{J}$ we have $\varphi*e_\tau=\varphi$ if and only if $\varphi\in \Ac_{P,\pi}(G)^{\tau,J}$.
For any complex function $f$ on a group, we define  $f^\vee(x)=\overline{f(x^{-1})}$. 
\end{paragr}

\begin{paragr}[Bounds for intertwinings operators.] ---  For  any standard parabolic subgroup   $P=MN_P$ of $G$, any discrete automorphic representation $\pi$ of $M(\AAA)$, any $c>0$ and  $l>0$ we define, following \cite[§ 3.3]{LapHC},
  \begin{align*}
    \rc_{\pi,c,l}=\{\la\in \ago_{P,\CC}^{G,*}\mid \|\Re(\la)\| < c(1+\La_{\pi}^M+\|\Im(\la)\|)^{-l}\}.
  \end{align*}

 For $w\in \, _RW_{P}$, we set $\nu_w=\nu_{P_w}$, see § \ref{S:exposant-disc}. Recall that for $\varphi\in \Ac_{P,\pi}(G)$ we set $\varphi_{w}=\varphi_{P_w,-\nu_w}$.

  \begin{proposition}\label{prop:bound-intertwining} There exist $k,l>0$ and for each standard parabolic subgroup $P$ of $G$ there exists $\ell_P$, a product of non-trivial real linear forms on $\ago_{P}^{G,*}$,   such that for any level  $J$ there exist $c>0$ and $C>0$ such that for all $J$-triples $(P,\pi,\tau)$ and all  standard parabolic subgroups $R$ of $G$,  all $w\in \, _RW_{P}$  such that  $  P_\pi\subset P_w$ and $R_w=R$ we have:
\begin{align*}
  \| \ell_P(\la) M(w,\la+\nu_w)\varphi_{w}\|_R \leq C (1+\|\la\|^2+\la_\tau^2+(\La_{\pi}^M)^2)^k\|\varphi_w\|_{P_w}
\end{align*}
for all $\varphi\in \Ac_{P,\pi}(G)^{\tau,J}$ and for all $\la\in \rc_{\pi,c,l}$.
\end{proposition}

\begin{remarque}\label{rq:holom op entr}
  We even show that $\ell_P(\la) M(w,\la+\nu_w)\varphi_{w}$ is holomorphic for $\la\in \rc_{\pi,c,l}$.  Moreover  we may take $\ell_P=1$ if we restrict ourselves to elements $w\in W(P;R)$ (see remark \ref{rq:bd norm factor}).
\end{remarque}
\end{paragr}

\begin{paragr}[Proof of  proposition \ref{prop:bound-intertwining}.] --- To prove  proposition \ref{prop:bound-intertwining}, we separately bound the normalized intertwining operators and the normalization factors. In the latter case, the bounds are given below in proposition \ref{prop:bounds norm fact} and are based on standard arguments. However, to bound the normalized intertwining operators, in the range we consider, requires much more efforts.   First in subsection \ref{ssec:holom norm} we show the holomorphy of  normalized intertwining operators in some interesting region, see proposition \ref{prop:le but}. Once the holomorphy has been proved, we manage to get the bounds we need  in  subsection \ref{ssec:bd norm inter}. All in all proposition \ref{prop:bound-intertwining} is a straightforward consequence of proposition \ref{prop:bounds norm fact} and proposition \ref{prop:NEWbd-norm-inter}.
\end{paragr}

\subsection{Normalization of intertwining operators}\label{ssec:normalization inter}

\begin{paragr}\label{S:L de paire}   For $i=1,2$ let $n_i\geq 1$ be an integer and let $\pi_i$ be a discrete  automorphic representation of $G_{n_i}$ with trivial  central character on $A_{G_{n_i}}^\infty$. Let $\pi_i^\vee$ be the contragredient of $\pi_i$. Let $L(s,\pi_1\times \pi_{2}^\vee)$ and $\eps(s,\pi_1\times \pi_{2}^\vee)$ be the  completed Rankin-Selberg $L$-function and the epsilon factor defined in \cite{JPSS}. We have $\eps(s,\pi_1\times \pi_{2}^\vee)=\eps_0 q^{\frac12-s}$ for some  $\eps_0\in \CC$  of modulus $1$ and an integer $q\geq 1$ called the arithmetic conductor. We write
  \begin{align*}
          L(s,\pi_1\times \pi_{2}^\vee)=L_\infty(s,\pi_1\times \pi_{2}^\vee)L^\infty(s,\pi_1\times \pi_{2}^\vee)
         \end{align*}
         as a product of the Archimedean contribution and the non-Archimedean one. We have:
         \begin{align*}
         L_\infty(s,\pi_1\times \pi_{2}^\vee)=\prod_{i=1}^m \Ga_\RR(s-\al_i)
         \end{align*}
         where $\Ga_\RR(s)=\pi^{-s/2}\Ga(s/2)$ (with the usual $\Ga$-function), $m=n_1n_2[F:\QQ]$ and $\al_i\in \CC$.  We define the analytic conductor $\cgo(\pi_1\times \pi_{2}^\vee,s)=q\prod_{j=1}^m (1+|s-\al_j|)$. We set  $\cgo(\pi_1\times \pi_{2}^\vee)=\cgo(\pi_1\times \pi_{2}^\vee,0)$. By specializing this to the case $G_{n_i}\times G_1$ and the trivial representation of $G_1$, one gets $\cgo(\pi_i)$. We have  $\cgo(\pi_1\times \pi_{2}^\vee)\leq \cgo(\pi_1)^{n_2}\cgo(\pi_2)^{n_1}$, see \cite[(24) and comments below]{LapHC}. For the external tensor product, we set $\cgo(\pi_1\boxtimes\pi_2)=\cgo(\pi_1) \cgo(\pi_2)$.  When  $\pi_1$ and $\pi_2$ are furthermore cuspidal, the product $(s(s-1))^{\delta} L(s,\pi_1\times \pi_{2}^\vee)$, where $\delta=\delta_{\pi_1,\pi_2}\in \{0,1\}$ is the Kronecker symbol, is an entire function of order $1$ (see \cite[§ 2.4]{RuSar}). 
         
         We set
         \begin{align}\label{eq:xipi1pi2}
           \xi_{\pi_1,\pi_2}(s)=\frac{L(s,\pi_1\times \pi_{2}^\vee)}{ \eps(s,\pi_1\times \pi_{2}^\vee) L(s+1,\pi_1\times \pi_{2}^\vee)  } .
         \end{align}
       \end{paragr}

       \begin{paragr}[Decomposition of  intertwining operators.] --- \label{S:decomp inter} 
        Let $P$ be a standard parabolic subgroup of $G$ with Levi decomposition $P=MN_P$. Let $\pi\in \Pi_{\disc}(M)$. Let $P_1$ be a parabolic subgroup of $G$ such that  $P_\pi\subset P_1\subset P$, see § \ref{S:exposant-disc}. We shall use the notations of subsection \ref{ssec:set inversion}. In particular $M_1$ is the standard Levi factor of $P_1$. Let $w\in W_{M_1}$ be such that $w\Delta_0^P>0$. Note that this condition is equivalent to  $\Delta_{P_1}^P\cap \Ga_w=\emptyset$. For any parabolic subgroup $P_1\subset Q\subset P$, let $\pi_Q$ be the  discrete automorphic representation of $M_{Q}$ such that the functions $\varphi_{Q,-\nu_Q}$ belong to  $\Ac_{Q,\pi_Q}(G)$ when $\varphi$ describes $\Ac_{P,\pi}(G)$. Set $\pi_1=\pi_{P_1}$. Using the decomposition \eqref{eq:decomp} and the notations therein, we get  a discrete automorphic representation $\pi_i=(s_{\al_{i-1}}\cdots s_{\al_1})\pi_1$ of $M_{i}$ for all $2\leq i \leq \ell+1$. For all $\la\in \ago_{P_1,\CC}^*$, let  $M_{\pi_1}(w,\la)$ be the restriction of the intertwining operator  $M(w,\la)$  attached to $w\in W(P_1,P_{\ell+1})$  to the space $\Ac_{P_1,\pi_{1}}(G)$. We have the decomposition:
   \begin{align}\label{eq:dec-Mw}
      M_{\pi_1}(w,\la) =                M_{\pi_{\ell}}(s_{\al_\ell}, (s_{\al_{\ell-1}}\cdots s_{\al_1})\la)  M_{\pi_{\ell-1}}(s_{\al_{\ell-1}}, (s_{\al_{\ell-2}}\cdots s_{\al_1})\la)  \cdots M_{\pi_{1}}(s_{\al_1},\la)
\end{align}
where  $M_{\pi_{i}}(s_{\al_{i}}, \mu)$ for $\mu\in \ago_{P_i,\CC}^*$ is the restriction of the  intertwining operator $M(s_{\al_{i}}, \mu)$ to  $\Ac_{P_i,\pi_{i}}(G)$.

\end{paragr}

\begin{paragr}[Normalization.] --- \label{S:Norm factors}Let us identify $M_{1}$ with a product $G_{n_1}\times \ldots \times G_{n_r}$ and $\pi_1$ with $\sigma_{1}\boxtimes \ldots \boxtimes \sigma_r$ accordingly where $\sigma_i$ is a discrete automorphic representation of $G_{n_i}$.  Let $\be$ be a positive root of $A_{P_1}$: it is associated to the factors, say,  $G_{n_i}$ and $G_{n_j}$ with $i<j$.  For any $\la\in \ago_{P_1,\CC}^*$, we set
   \begin{align}\label{eq:NEW coef npi1}
      n_{\pi_1}(\be,\la)=\xi_{\sigma_i,\sigma_j}(\bg \la,\be^\vee\bd).
   \end{align}
   Then we introduce the  normalization factor
\begin{align}\label{eq:global norm fact}
      n_{\pi_1}(w,\la)=\prod_{\be\in \Ga_w}   n_{\pi_1}(\be,\la)
    \end{align}
and we get the  normalized  intertwining operator  by setting:
    \begin{align*}
      N_{\pi_1}(w,\la)=    n_{\pi_1}(w,\la)^{-1}M_{\pi_1}(w,\la).
    \end{align*}
  \end{paragr}

  \begin{paragr}[Local decomposition.] --- For all places $v\in V_F$, let $\pi_{1,v}$ be the local component of $\pi_1$ such that $\pi_1$ is identified with the restricted tensor product $\otimes_{v\in V_F} \pi_{1,v}$. The global intertwining  operator  $M_{\pi_1}(w,\la) $ which intertwines the actions $I_{P_1,\pi_1}(\la)$ on $\Ac_{P_1,\pi_1}(G)$ and  $I_{P_{\ell+1},\pi_{\ell+1}}(w\la)$ on $\Ac_{P_{\ell+1},\pi_{\ell+1}}(G)$ has a local  counterpart $M_{\pi_1,v}(w,\la) $  which  intertwines the normalized induced representations $I_{P_1}^{G}(\pi_{1,v}\otimes \exp(\bg \la,H_{P_1}(\cdot)\bd))$ and $I_{P_{\ell+1}}^{G}(\pi_{\ell+1,v}\otimes \exp(\bg w\la,H_{P_{\ell+1}}(\cdot)\bd))$.
    In some half-planes,  the $L$-functions and the epsilon factors of § \ref{S:L de paire}  are given by an Euler product. Using the local factors, we define for each place $v$ the local normalization factor  $n_{\pi_1,v}(w,\la)$ so that the global  normalization factor $n_{\pi_1}(w,\la)$ is given by an Euler product in some convergence region
\begin{align*}
  n_{\pi_1}(w,\la)=\prod_{v\in V_F}   n_{\pi_1,v}(w,\la).
\end{align*}
        We get the  normalized  local intertwining operator $N_v(w,\la)$ by setting:
    \begin{align*}
      N_{\pi_1,v}(w,\la)=    n_{\pi_1,v}(w,\la)^{-1}M_{\pi_1,v}(w,\la).
    \end{align*}
    We also have:
     \begin{align}\label{eq:dec-Nw}
      N_{\pi_1,v}(w,\la) =                N_{\pi_{\ell},v}(s_{\al_\ell}, (s_{\al_{\ell-1}}\cdots s_{\al_1})\la)  N_{\pi_{\ell-1},v}(s_{\al_{\ell-1}}, (s_{\al_{\ell-2}}\cdots s_{\al_1})\la)  \cdots N_{\pi_1,v}(s_{\al_1},\la).
\end{align}
and 
\begin{align}\label{eq:local decomposition}
 N_{\pi_1}(w,\la)\simeq \otimes_v N_{\pi_1,v}(w,\la).
\end{align}

\end{paragr}

\subsection{Bounds for the normalization factors}\label{ssec:bd normalization}

\begin{paragr}In this subsection we give bounds for the normalization factors of § \ref{S:Norm factors}. These are given by the following proposition. The notations are those of § \ref{S:decomp inter}. 

  \begin{proposition} \label{prop:bounds norm fact}There exist $k,l>0$ and for each standard parabolic subgroup $P$ of $G$ there exists $\ell_P$, a product of non-trivial real linear forms on $\ago_{P,\CC}^{G,*}$,   such that for any level  $J$ there exist $c>0$ and $C>0$ such that for all $J$-pair $(P,\pi)$ and all  standard parabolic subgroups $P_1$ of $G$ with standard Levi factor $M_1$ such that $P_\pi\subset P_1\subset P$,  all $w\in W_{M_1}$  such that  $w\Delta_0^P>0$  we have:
\begin{align*}
  | \ell_P(\la)   n_{\pi_1}(w,\la+\nu_1)| \leq C (1+\|\la\|^2+(\La_{\pi}^M)^2)^k
\end{align*}
for  all $\la\in \rc_{\pi,c,l}$.
  \end{proposition}

  \begin{remarque}\label{rq:bd norm factor}
  As it follows from the proof, see remark \eqref{rq:l=1}, one may take  $ \ell_P=1$ if $\nu_1=0$ that is if $P_1=P$.
\end{remarque}

The proof of proposition \ref{prop:bounds norm fact} occupies the remainder of this subsection. 
\end{paragr}

\begin{paragr}\label{S:NEW calculs expl} Let $P$ be a standard parabolic subgroup of $G$ with Levi decomposition $MN_P$.  Let $\pi\in \Pi_{\disc}(M)$. Let $P_1$ be a parabolic subgroup of $G$ such that  $P_\pi\subset P_1\subset P$, see § \ref{S:exposant-disc}. Then as in  § \ref{S:decomp inter} we get $\pi_1\in \Pi_{\disc}(M_1)$ where $M_1$ is the standard Levi factor of $P_1$. We have an identification $M=G_{N_1}\times \ldots \times G_{N_r}$ with $N_1+\ldots +N_r=n$ (recall that $G=G_n$). According to this decomposition we write $P_1\cap M=P_1^1\times\ldots \times P_1^r$ and  $P_\pi\cap M=P_\pi^1\times\ldots \times P_\pi^r$. We also have decompositions  $\pi=\pi^1\boxtimes \cdots \boxtimes \pi^r$ and   $\pi_1=\pi^1_1\boxtimes \cdots \boxtimes \pi^r_1$ where $\pi^i$ and $\pi_1^i$ are discrete automorphic representations respectively of $G_{N_i}$ and the standard Levi factor $M_1^i$ of $P_1^i$. Let $1\leq i\leq r$. The standard Levi factor of $P_\pi^i$ is identified with $G_{r_i}^{d_i}$ for some integers $r_i,d_i\geq 1$ such that $N_i=r_id_i$. Then there exists a cuspidal representation $\tau_i$ of $G_{r_i}$ such that   $\pi^i$ is the discrete automorphic representation of  $G_{N_i}$ associated to $\tau_i$ in the M{\oe}glin-Waldspurger classification.

  Let $w \in W_{M_1}$  such that  $w\Delta_0^P>0$. For $\mu\in \ago_{P_1,\CC}^*$ the normalization factor  $n_{\pi_1}(w,\mu)$ defined in \eqref{eq:global norm fact} is given by a product over the subset $\Ga_w\subset \Phi_{P_1}$. We shall fix a root $\al\in\Phi_{P}$ such that $\Ga_w(\al)$ is non-empty and we shall restrict ourselves to the subset $\Ga_w(\al)\subset \Ga_w$ of roots whose restriction to $A_P$ is $\al$. So we have $1\leq i_0<j_0\leq r$ such that the root $\al$   is associated to the blocks  $G_{N_{i_0}}$ and $G_{N_{j_0}}$ in the decomposition of $M$. The standard Levi factor $M_1^{i_0}$ of $P_1^{i_0}$ is identified with $G_{e_1r_{i_0}}\times \cdots \times G_{e_kr_{i_0}}$  for positive integers $k$ and $(e_i)_{1\leq i\leq k}$ such that $e_1+\ldots +e_k=d_{i_0}$. In the same way, we attach to $M_1^{j_0}$ positive integers $l$ and $(f_i)_{1\leq i\leq l}$ such that $f_1+\ldots +f_l=d_{j_0}$. We have decompositions $\pi_1^{i_0}=\sigma_1\boxtimes \cdots\boxtimes \sigma_k$ where for $1\leq j\leq k$ the discrete automorphic representation  $\sigma_j$ of  $G_{e_j r_{i_0}}$ is associated to $\tau_{i_0}$ by  M{\oe}glin-Waldspurger classification. Similarly we have $\pi_1^{j_0}=\sigma_1'\boxtimes \cdots\boxtimes \sigma_l'$ where $\sigma_j'$ is the discrete automorphic representation of  $G_{f_j r_{j_0}}$ associated to $\tau_{j_0}$. The roots in $\Phi_{P_1}$ whose restriction to $A_P$ is $\al$ are the  roots $\be_{i,j}$ where $1\leq i\leq k$, $1\leq j\leq l$ and $\be_{i,j}$ relates the block $G_{e_ir_{i_0}}$ of $M_1^{i_0}$ to the block $G_{f_j r_{j_0}}$  of $M_1^{j_0}$. By proposition \ref{prop:w-tableau} there exist $1\leq k'\leq k$ and  a non-decreasing sequence $(l_i)_{k'\leq i\leq k}$ of integers $1\leq l_i\leq l$ such that the elements of $\Ga_w(\al)$ are the roots $\be_{i,j}$ such that $k'\leq i\leq k$ and $1\leq j\leq l_i$. Alternatively there exist $1\leq l'\leq l$ and a  non-decreasing sequence $(k_j)_{1\leq j\leq l'}$ of integers  $1\leq k_j\leq k$  such that the elements of $\Ga_w(\al)$ are the roots $\be_{i,j}$ such that   $1\leq j\leq l'$ and $k_j\leq i\leq k$.

  Let $\la\in \ago_{P,\CC}^*$. The factor we are considering is the following:
  \begin{align*}
    \prod_{\be\in \Ga_w(\al)}  n_{\pi_1}(\be,\la+\nu_1)=E(s) \prod_{i=k'}^k \prod_{j=1}^{l_i} \frac{L(s+ \bg \nu_1,\be_{i,j}^\vee\bd,\sigma_i\times {\sigma_{j}'}^{\vee})}{ L(s+ 1+\bg \nu_1,\be_{i,j}^\vee\bd,\sigma_i\times {\sigma_{j}'}^{\vee})}
  \end{align*}
  where $s=\bg \la,\be_{i,j}^\vee\bd$ does not depend on $i,j$ and we set
  \begin{align*}
    E(s)= \prod_{i=k'}^k \prod_{j=1}^{l_i}\eps(s+ \bg \nu_1,\be_{i,j}^\vee\bd,\sigma_i\times {\sigma_{j}'}^{\vee})^{-1}.
  \end{align*}
Note that on any vertical strip the factor $E(s)$ is bounded by a constant that depends on the level only. So we will focus on the quotient of $L$-functions considered in next lemma. 
  
  \begin{lemme}
    We have:
    \begin{align}
\nonumber&      \prod_{i=k'}^k \prod_{j=1}^{l_i} \frac{L(s+ \bg \nu_1,\be_{i,j}^\vee\bd,\sigma_i\times {\sigma_{j}'}^{\vee})}{ L(s+ 1+\bg \nu_1,\be_{i,j}^\vee\bd,\sigma_i\times {\sigma_{j}'}^{\vee})}\\
       \label{eq:rapp 1}  &    =\prod_{i=k'}^k \prod_{t=0}^{e_i-1}  \frac{  L(s+\frac{d_{j_0}-d_{i_0}}2 +(e_1+\ldots+e_i)-(f_1+\ldots+f_{l_i})     -t,\tau_{i_0}\times \tau_{j_0}^\vee)}{L(s+\frac{d_{j_0}-d_{i_0}}2 +(e_1+\ldots+e_i)-t,\tau_{i_0}\times \tau_{j_0}^\vee)}\\
    \label{eq:rapp 2}     &= \prod_{j=1}^{l'} \prod_{t=0}^{f_j-1}  \frac{  L(s+\frac{d_{j_0}-d_{i_0}}2 +(e_1+\ldots+e_{k_j-1})-(f_1+\ldots+f_{j-1})     -t,\tau_{i_0}\times \tau_{j_0}^\vee)}{L(s+\frac{d_{j_0}+d_{i_0}}2 -(f_1+\ldots+f_{j-1})    -t,\tau_{i_0}\times \tau_{j_0}^\vee)}.
    \end{align}
  \end{lemme}

\begin{preuve}
  For all  $s\in \CC$  and all integers $1\leq i\leq k$ and $1\leq j\leq l$ we have:
   \begin{align*}
             L(s,\sigma_i\times {\sigma_{j}'}^{\vee})
             &=\prod_{t=0}^{e_i-1}  \prod_{t'=0}^{f_j-1} L(s+\frac{e_i+f_j}2-1-t-t',\tau_{i_0}\times \tau_{j_0}^\vee).
           \end{align*}

In particular
\begin{align*}
  \frac{  L(s,\sigma_i\times {\sigma_{j}'}^{\vee}) }{  L(s+1,\sigma_i\times {\sigma_{j}'}^{\vee})}&= \prod_{t=0}^{f_j-1} \frac{  L(s+\frac{f_j-e_i}2-t,\tau_{i_0}\times \tau_{j_0}^\vee)}{L(s+\frac{f_j+e_i}2-t,\tau_{i_0}\times \tau_{j_0}^\vee)}\\
&  = \prod_{t=0}^{e_i-1} \frac{  L(s+\frac{e_i-f_j}2-t,\tau_{i_0}\times \tau_{j_0}^\vee)}{L(s+\frac{e_i+f_j}2-t,\tau_{i_0}\times \tau_{j_0}^\vee)}
\end{align*}
Using the explicit computation of $\nu_1$, see \eqref{eq:calcul de nu}, we get
\begin{align*}
&  \prod_{i=k'}^k \prod_{j=1}^{l_i} \frac{L(s+ \bg \nu_1,\be_{i,j}^\vee\bd,\sigma_i\times {\sigma_{j}'}^{\vee})}{ L(s+ 1+\bg \nu_1,\be_{i,j}^\vee\bd,\sigma_i\times {\sigma_{j}'}^{\vee})}\\
&=\prod_{i=k'}^k \prod_{j=1}^{l_i}  \prod_{t=0}^{e_i-1}  \frac{  L(s+\frac{d_{j_0}-d_{i_0}}2 +(e_1+\ldots+e_i)-(f_1+\ldots+f_{j})     -t,\tau_{i_0}\times \tau_{j_0}^\vee)}{L(s+\frac{d_{j_0}-d_{i_0}}2 +(e_1+\ldots+e_i)-(f_1+\ldots+f_{j-1})     -t,\tau_{i_0}\times \tau_{j_0}^\vee)}\\
  &=\prod_{i=k'}^k \prod_{t=0}^{e_i-1}  \frac{  L(s+\frac{d_{j_0}-d_{i_0}}2 +(e_1+\ldots+e_i)-(f_1+\ldots+f_{l_i})     -t,\tau_{i_0}\times \tau_{j_0}^\vee)}{L(s+\frac{d_{j_0}-d_{i_0}}2 +(e_1+\ldots+e_i)-t,\tau_{i_0}\times \tau_{j_0}^\vee)}.
\end{align*}
We get \eqref{eq:rapp 1}. A similar computation gives \eqref{eq:rapp 2}. 
\end{preuve}
\end{paragr}

\begin{paragr}\label{S:j0 ged  i0}
  In the following we assume $d_{j_0}\geq d_{i_0}$ and we shall bound expression \eqref{eq:rapp 1} on a neighborhood of the imaginary axis. To do this, we split  \eqref{eq:rapp 1} into three factors: the non-Archimedean component of the denominator which will be bound thanks to a result of Brumley, the non-Archimedean component of the numerator which will be bound by standard techniques based on the Phragmén-Lindelöf principle and the remaining Archimedean factors which is a ratio of $\Gamma$-functions.   The case $d_{j_0}< d_{i_0}$ is similar but requires to use expression \eqref{eq:rapp 2}: it is left to the reader.

Let $k'\leq i\leq k$ and $0\leq t\leq e_i-1$.    We observe that we have
         \begin{align}\label{eq:NEW di0}
           \frac{d_{j_0}-d_{i_0}}2 +(e_1+\ldots+e_i)-t\geq 1.
         \end{align}
and the inequality is even strict unless we have $d_{j_0}=d_{i_0}$, $i=1$, $t=e_1-1$ and thus $k'=1$ . 
Using Kronecker symbols, we set $\delta=\delta_{\tau_{i_0},\tau_{j_0}}$ and $\delta'=\delta \times \delta_{d_{i_0},d_{j_0}}\times\delta_{k',1}\in \{0,1\}$.

We write:
 \begin{align*}
         L_\infty(s,\tau_{i_0}\times \tau_{j_0}^\vee)=\prod_{i=1}^m \Ga_\RR(s-\al_i)
         \end{align*}
         where $m=r_{i_0}r_{j_0}[F:\QQ]$ and $\al_i\in \CC$ are such that $\Re(\al_i)<1$ (which is the Jacquet-Shalika bound). 
  We set
  \begin{align*}
&    A_\infty(s)=\prod_{j=1}^m \prod_{i=k'}^k \prod_{t=0}^{e_i-1}   \Ga_\RR(s+\frac{d_{j_0}-d_{i_0}}2 +(e_1+\ldots+e_i)-(f_1+\ldots+f_{l_i})     -t-\al_j),\\
&B_\infty(s)=    \prod_{j=1}^m\prod_{i=k'}^k \prod_{t=0}^{e_i-1}  \Ga_\RR(s+\frac{d_{j_0}-d_{i_0}}2 +(e_1+\ldots+e_i)-t-\al_j),
  \end{align*}
  so that $Q_\infty=A_\infty(s)B_\infty(s)^{-1}$ is the Archimedean component of \eqref{eq:rapp 1}.  We also define the non-Archimedean components:
  \begin{align*}
    &    A^\infty(s)=\prod_{i=k'}^k \prod_{t=0}^{e_i-1}  L^\infty(s+\frac{d_{j_0}-d_{i_0}}2 +(e_1+\ldots+e_i)-(f_1+\ldots+f_{l_i})-t,\tau_{i_0}\times \tau_{j_0}^\vee);\\
      &B^\infty(s)= \prod_{i=k'}^k \prod_{t=0}^{e_i-1}  L^\infty(s+\frac{d_{j_0}-d_{i_0}}2 +(e_1+\ldots+e_i)-t,\tau_{i_0}\times \tau_{j_0}^\vee).
  \end{align*}

\end{paragr}

\begin{paragr}[Bounds on the non-Archimedean denominator.]   

  \begin{lemme}
     There are $c_1,c_2,m_1,m_2>0$  independent of  $\tau_{i_0}$ and $\tau_{j_0}$ such that we have:
         \begin{align}\label{eq:maj-Bfin}
           (s^{\delta'}B^\infty(s)  )^{-1}\leq  c_1 (1+\cgo(\tau_{i_0})+\cgo( \tau_{j_0})+|s|)^{m_1}
         \end{align}
         for all $s\in \CC$ such that $\Re(s)\geq -c_2 \cgo(\tau_{i_0}\times\tau_{j_0}^\vee,s)^{-m_2}$.
       \end{lemme}

       \begin{preuve} The lemma is a straightforward consequence of  \eqref{eq:NEW di0} and the following result due to Brumley \cite[proposition 3.5]{LapHC}: there are $c_1,c_2,m_1,m_2>0$    independent of  $\tau_{i_0}$ and $\tau_{j_0}$ such that we have:
         \begin{align*}
           \left|s^{\delta} L^\infty(1+s,\tau_{i_0}\times \tau_{j_0}^\vee)\right|^{-1}\leq  c_1 (1+\cgo(\tau_{i_0})+\cgo( \tau_{j_0})+|s|)^{m_1}
         \end{align*}
         for all $s\in \CC$ such that $\Re(s)\geq -c_2 \cgo(\tau_{i_0}\times\tau_{j_0}^\vee,s)^{-m_2}$.
       \end{preuve}
\end{paragr}

\begin{paragr}[Bounds on the Archimedean component.] --- By elementary properties of the $\Gamma$ function, see for example the discussion in \cite[p. 33]{Ar-PW}, we can choose a monic polynomial $\ell'$ (independent on $\tau_{i_0}$ and $\tau_{j_0}$) and $N,c>0$ such that for all $s\in \CC$ such that $\Re(s) \geq -1$ the expression $\prod_{j=1}^m \ell'(s-\al_j)A_\infty(s)$ is holomorphic and non-vanishing and we have
  \begin{align}\label{eq:maj Q infin}
     |\prod_{j=1}^m  \ell'(s-\al_j) Q_\infty(s) |\leq c (\prod_{j=1}^m (1+|s-\al_i|)^N.
  \end{align}
  
\end{paragr}

\begin{paragr}[Bounds on the non-Archimedean numerator.] --- \label{S:NEW bd non-Archi}Set
  \begin{align*}
&    \ell^\flat(s)=\prod_{i=k'}^k \prod_{t=0}^{e_i-1}  (s+\frac{d_{j_0}-d_{i_0}}2 +(e_1+\ldots+e_i)-(f_1+\ldots+f_{l_i})-t)^\delta \\
  &  \ell(s)=\ell^\flat(s)\ell^\flat(s-1).
  \end{align*}
  Then the map $s\mapsto  \ell(s)  A_\infty(s)A^\infty(s)$    is holomorphic on $\CC$. We deduce that the map $s\mapsto           \frac{ s^{\delta '}\ell(s)}{\prod_{j=1}^m  \ell'(s-\al_j) } A^\infty(s)$    is holomorphic on the set of $s\in \CC$ such that $\Re(s)\geq -1$.    By an application of the Phragmén-Lindelöf principle (see e.g. \cite[section 2.1]{Michel}), one shows that there are $c',m'>0$  such that  we have on this subset
     \begin{align}\label{eq:maj-Afin}
    |  \frac{ s^{\delta'}\ell(s)}{ \prod_{j=1}^m  \ell'(s-\al_j) } A^\infty(s)      |\leq c' (1+\cgo(\tau_{i_0}\times \tau_{j_0}^\vee)+|s|)  ^{m'}.
     \end{align}

     \begin{remarque}\label{rq:l=1}
       If we assume $P=P_1$ we have $l=1$, $k=1$ and 
       \begin{align*}
   s^{\delta'} A_\infty(s) A^\infty(s)=       s^{\delta'}  \prod_{t=0}^{d_{i_0}-1}  L(s+\frac{d_{i_0}-d_{j_0}}2 -t,\tau_{i_0}\times \tau_{j_0}^\vee)
       \end{align*}
       is  holomorphic on the subset of $s\in \CC$ such that $|\Re(s)|\leq \frac14$.  Thus on this subset the map $s\mapsto           \frac{ s^{\delta '}}{\prod_{j=1}^m  \ell'(s-\al_j) } A^\infty(s)$    is  also holomorphic and a bound similar to \eqref{eq:maj-Afin} holds with $\ell=1$.
     \end{remarque}

     \begin{remarque}\label{rq:l=2}
If we remove from $\ell$ all the factors that do vanish at $s=0$, the bound \ref{eq:maj-Afin} still holds on the subset of $s\in \CC$ such that $|\Re(s)|\leq \frac14$ 
     \end{remarque}
   \end{paragr}

   \begin{paragr} For any level $J$, there exists $c_J>0$ such that for all $J$-pair $(P,\pi)$ we have  $\cgo(\pi)\leq c_J (1+\La^M_\pi)$ (see \cite[lemme 5.4]{Muller07}; we use also the fact that the arithmetic conductor is bounded in terms of the level $J$). Moreover, for $\tau_{i_0}$  as above, $\cgo(\tau_{i_0})$ is bounded in terms of a power of $\cgo(\pi)$. As a consequence, we get proposition \ref{prop:bounds norm fact} by gathering the bounds \eqref{eq:maj-Bfin}, \eqref{eq:maj Q infin} and \eqref{eq:maj-Afin}. \emph{A priori} the polynomials $\ell_P$ we get are more general than expected. But, taking into account remark \ref{rq:l=2}, we may and shall assume that the    polynomials $\ell$ in §~\ref{S:NEW bd non-Archi} are simply some powers of $s$.
  \end{paragr}

\subsection{Holomorphy of normalized intertwining operators}\label{ssec:holom norm}

\begin{paragr}
  The bulk of the subsection is devoted to the proof of  proposition \ref{prop:le but} which gives a crucial  domain of holomorphy of normalized intertwining operators. We continue with the notations of subsections \ref{ssec:set inversion} and  \ref{ssec:normalization inter}.
\end{paragr}

\begin{paragr}[Properties of normalized local intertwining operators.] ---  We fix a place $v\in V_F$.

  \begin{lemme}\label{lem:holom op simple}
    Let $1\leq i\leq \ell$. The normalized intertwining operator $\mu\in \ago_{P_i,\CC}\mapsto N_{\pi_i,v}(s_{\al_{i}}, \mu)$ is a meromorphic function of $\bg \mu, \al_i^\vee\bd$. Moreover it is holomorphic on the open subset of $\mu\in  \ago_{P_i,\CC}^*$ such that
  \begin{align*}
    \Re(\bg \mu, \al_i^\vee\bd)> -\frac{2}{n^2+1}.
  \end{align*}
  In particular, the operator $N_{\pi_i,v}(s_{\al_{i}}, (s_{\al_{i-1}}\cdots s_{\al_1})(\la+\nu_1))$  is holomorphic on the open subset of $\la\in  \ago_{P_1,\CC}^*$ such that 
\begin{align}\label{eq:positiv-n2}
  \Re(\bg \la+\nu_1, \theta_i^\vee\bd)> -\frac{2}{n^2+1}.
\end{align}
  \end{lemme}

\begin{remarque}\label{rq:linear form}
  The linear form $\la\mapsto \bg \la, \theta_i^\vee\bd$ is non trivial on $\ago_{P,\CC}^{G,*}$: otherwise we would have $\theta_i\in \Phi_{P_1}^P$ and thus $w\theta_i>0$, which is impossible since   $\theta_i\in \Ga_w$.
\end{remarque}

\begin{preuve}The first assertion follows from \cite[proposition I.10]{MW} and \cite[proposition 3.5]{MuSpeh} (based on the approximation of the  Ramanujan conjecture of \cite{LRS}). To get the second assertion, it suffices to observe that we have:
  \begin{align*}
    \bg s_{\al_{i-1}}\cdots s_{\al_1}(\la+\nu_1),\al_{i}^\vee  \bd=\bg \la+\nu_1, \theta_i^\vee\bd.
  \end{align*}
  
\end{preuve}
  
\end{paragr}

\begin{paragr}[A domain of holomorphy of normalized intertwining operators.] --- We fix a level   $J$ and a $J$-triple $(P,\pi,\tau)$. Let  $P_\pi\subset P_1\subset P$ and $w\in W_{M_1}$ be such that $w\Delta_0^P>0$.  Set $\Delta_1=\Delta_{P_1}$ and  $\Phi_1=\Phi_{P_1}$.  Let $r$ be the cardinality of $\Delta_1$. We fix a numbering on the roots $\be_1,\ldots,\be_r$ in  $\Delta_1$ in such a way that the elements 
 \begin{align*}
    \gamma_{k,l}=\be_k+\ldots+\be_{l-1}
  \end{align*}
  for $1\leq k< l\leq r+1$ are precisely the elements of $\Phi_1$.  Note the formula $\gamma_{k,l}+\gamma_{l,m}=\ga_{k,m}$ for $ 1\leq k< l<m\leq r+1$. Let $s$ be the cardinality of $\Ga_{w}\cap \Delta_1$: we have   $0\leq s\leq r$. Let  $1\leq i_1<\ldots <i_s\leq r$ be such that  $\Ga_{w}\cap \Delta_1=\{\be_{i_1}, \ldots \be_{i_s}\}$. Note that for $1\leq j\leq s$ we have $\be_{i_j}\notin \Delta_{1}^P$.  We equip  $\Phi_1$ with the partial order given by  $\al\prec  \be$ if $\be-\al$ is a sum of elements of  $\Delta_1$.
We set:
  \begin{align}\label{eq:def-Om-New}
    \Xi_P=\{\la\in  \ago_{P,\CC}^{G,*}\mid   \forall \al\in \Phi_P \ |\Re( \bg \la, \al^\vee\bd)|<\frac{2}{n^2+1}\}.
    \end{align}
    Let $\nu_1=\nu_{P_1}$ and $\varphi_1=\varphi_{P_1,-\nu_1}$ for $\varphi \in \Ac_{P,\pi}(G)$.  The goal of  the subsection  is to prove the following proposition.

    \begin{proposition}\label{prop:le but}
    For all $\varphi\in \Ac_{P,\pi}(G)$ the map  $\la\mapsto N_{\pi_1}(w,\la+\nu_1)\varphi_1$ is holomorphic on $\Xi_P$.
  \end{proposition}
  
The rest of the subsection is devoted to the proof of proposition \ref{prop:le but}.
\end{paragr}

\begin{paragr}[Auxiliary lemmas.] ---  In the proof of proposition \ref{prop:le but} we shall use the following lemmas.

   \begin{lemme}\label{lem:Delta1P}
  Let $\al\in \Delta_{P_1}$. We have $  \al\in \Delta_{P_1}^P $  if and only if  $\bg \nu_1,\al^\vee\bd < 0$.
  \end{lemme}

  \begin{preuve}
   Let $\tilde \al\in \Delta_0\setminus \Delta_0^{P_1}$ be the unique element whose restriction to $A_{P_1}$ is  $\al$. We have
    \begin{align*}
      \bg \nu_1,\al^\vee \bd =\bg \nu_1,\tilde \al^\vee \bd.
    \end{align*}
    We can identify  $M_P$ with  a product $\prod_{i} G_{n_i}$ of general linear groups and accordingly $P_1\cap M_P$ with a product $\prod_i R_i$ of standard parabolic subgroups $R_i\subset  G_{n_i}$ .

    First let us consider the case $\tilde \al\in \Delta_0^P$. Then there is an index $i$ such that $\tilde \al$ is identified with a simple root in $\Delta_0^{G_{n_i}}$ still denoted by $\tilde \al$.  Then $\bg \nu_1,\tilde \al^\vee \bd$ is equal to (up to a positive coefficient)
\begin{align*}
  \bg -\rho_{R_i}^{G_{n_i}},\tilde \al^\vee\bd&=\bg-\rho_{0}^{G_{n_i}}+  \rho_{0}^{R_i},\tilde \al^\vee\bd \\
  &\leq -1.
\end{align*}
Indeed we have $\bg\rho_{0}^{G_{n_i}}, \tilde \al^\vee\bd =1$ and  $\rho_{0}^{R_i}$ is a linear combination with non-negative coefficients of elements of  $\Delta_{0}^{R_i}$. But  for all $\be\in \Delta_0^{R_i}$ we have $\tilde \al\not=\be$ and thus $\bg \be,\tilde\al^\vee\bd\leq 0$.

Second we assume  $\tilde \al\in  \Delta_0\setminus \Delta_0^P$. Since  $\nu_1$ is a linear combination with non-positive coefficients of elements of  $\Delta_{0}^{P}$, the same argument gives $\bg \nu_1,\al^\vee\bd \geq 0$.
\end{preuve}

 \begin{lemme}\label{lem: Ga+al} Let $1\leq l\leq r$ be an integer such that $l>i_s$ or $l<i_1$. Let  $\Ga'$ be the union of  $\Ga_w$ and  the set of elements $\ga+\be_{l}$ for all  $\gamma\in \Ga_w$ such that $\ga+\be_{l}\in \Phi_1$. Then  $\Ga'$ is closed and coclosed in the sense of § \ref{S:elem sym}.
  \end{lemme}

  \begin{preuve} We shall just prove  the case $l>i_s$, the proof for  the case $l<i_1$ being similar.  Let $\gamma,\gamma'\in \Phi_1$ be such that  $\ga+\ga'\in \Phi_1$.

    First let us show that the set  $\Ga'$ is closed. Assume $\ga,\ga'\in \Ga'$. If moreover $\ga,\ga'\in \Ga_w$ then  $\ga+\ga'\in \Ga_w$. We cannot have both  $\ga$ and $\ga'$ in the complement $\Phi_1\setminus \Ga_w$: otherwise, when written in the basis of simple roots,  both would contain $\be_{l}$ as a summand and their sum could not be a root. So without loss of generality we may assume  $\ga\in \Ga_w$ and $\ga'=\ga''+\be_l\in \Phi_1\setminus \Ga_w$  with  $\ga''\in \Ga_w$.  We must have  $\ga''=\ga_{k,l}$ with  $1\leq k<l$ (the other possibility $\ga''=\ga_{l+1,k}$ with $l+1<k\leq r+1$ would imply  $\ga''\notin \Ga_w$). So we also have $\ga'=\ga_{k,l+1}$ and $\ga=\ga_{t,k}$ with  $1\leq t<k$. Thus $\ga+\ga''$ belongs to $\Phi_1$ and to $\Ga_w$ since $\Ga_w$ is closed.  As a consequence we get $\ga+\ga'=\ga+\ga''+\be_l\in \Ga'$.

Second let us show that the set $\Ga'$  is coclosed. Assume  $\ga,\ga'\in \Phi_1\setminus \Ga'$. It is clear that $\ga+\ga'\notin \Ga_w$. Assume there exists $\delta\in \Ga_w$ such that  $\ga+\ga'=\delta+\be_{l}$. We must have  $\delta=\gamma_{k,l}$ with $1\leq k<l$ (the other possibility  $\delta=\gamma_{l+1,k}$ with $l<k\leq r+1$ would imply $\delta\notin \Ga_w$). Without loss of generality we may assume that there is  $k<t\leq l$ such that  $\ga=\gamma_{k,t}$ and $\ga'=\gamma_{t,l+1}=\gamma_{t,l}+\be_l$. Then we have $\ga+\gamma_{t,l}=\delta \in \Ga_w$: this implies  $\gamma_{t,l}\in \Ga_w$ since $\ga\notin \Ga_w$. So $\ga'\in \Ga'$: contradiction.
  \end{preuve}

\end{paragr}

\begin{paragr}[Induction hypothesis.] ---
  Set $\ell=\ell(w)$, see § \ref{S:length}.  If $\ell=0$, proposition  \ref{prop:le but} is trivial. If  $\ell=1$,  the set $\Ga_w$ is reduced to a singleton $\{\al\}$ where $\al\in \Delta_1\setminus \Delta_1^P$. So proposition  \ref{prop:le but} follows from   lemmas \ref{lem:holom op simple} and \ref{lem:Delta1P}. Thus by induction we shall assume that proposition \ref{prop:le but} holds for all $w' \in W_{M_1}$ such that $w'\Delta_0^P>0$ and $\ell(w')<\ell$.  We shall also assume that proposition \ref{prop:le but} holds if one replaces  $P_1$ by a standard parabolic subgroup $Q$ such that $P_1\subsetneq Q\subset P$. Indeed if $Q=P$,  proposition \ref{prop:le but} is a straightforward consequence of  lemma \ref{lem:holom op simple}. 

  We shall use notations of § \ref{S:decomp inter}.  A first consequence of the induction hypothesis is the following lemma.

  \begin{lemme}\label{lem:rec terme cst}
    Assume that there exists $\be\in \Delta_1^P$ such that $w\be\in \Delta_{P_{\ell+1}}$. Then the induction hypothesis implies that proposition \ref{prop:le but} holds.
  \end{lemme}

  \begin{preuve}
    Let  $P_1\subsetneq Q\subset P$ be the  parabolic subgroup defined by  $\Delta_1^Q=\{\be\}$. Let $P_{\ell+1}\subset Q'$ be the   parabolic subgroup defined by  $\Delta_{P_{\ell+1}}^{Q'}=\{w\be\}$. We have $w\Delta_0^{Q}=\Delta_0^{Q'}$.    By \eqref{eq:torsion mu} and \eqref{eq:cst intert} one has for all $\la\in \ago_{P,\CC}^{G,*}$
    \begin{align*}
      M_{\pi_1}(w,    \la+\nu_1)\varphi_{1}&=  (M(w,    \la)\varphi_{P_1})_{-w\nu_1}\\
                                             &= (M(w,    \la)((\varphi_Q)_{P_1}))_{-w\nu_1}\\
      &=((M(w,    \la)(\varphi_Q))_{P_{\ell+1}})_{-w\nu_1}\\
      &=(M_{\pi_Q}(w,    \la+\nu_Q)(\varphi_{Q,-\nu_Q}))_{P_{\ell+1},-w\nu_1^Q}.
    \end{align*}
  Using lemma   \ref{lem:NEW pi1 piQ} below, we deduce that  we also have:
    \begin{align*}
       N_{\pi_1}(w,    \la+\nu_1)\varphi_{1}= (N_{\pi_Q}(w,    \la+\nu_Q)(\varphi_{Q,-\nu_Q}))_{P_{\ell+1},-w\nu_1^Q}.
    \end{align*}
    By the induction hypothesis,  the  right-hand side is holomorphic for $\la\in\Xi_P$. The conclusion is then obvious.
  \end{preuve}
  
  \begin{lemme}
    \label{lem:NEW pi1 piQ}
    With notations as in the proof of lemma \ref{lem:rec terme cst}, we have for all $\la\in \ago_{P,\CC}^{G,*}$ :
    \begin{align*}
      n_{\pi_Q}(w,\la+\nu_Q)=   n_{\pi_1}(w,\la+\nu_1).
    \end{align*}
  \end{lemme}

  \begin{preuve} We use the notations of the proof of lemma \ref{lem:rec terme cst}.  Let us add an upperscript $P_1$, resp. $Q$, to the sets $\Ga_w$ when we want to emphasize that they are relative to $P_1$, resp. $Q$.

    We have a surjective restriction map $\Phi_1\setminus\{\be\}\to \Phi_Q$. Let $\al\in \Phi_Q$. There is a dichotomy: either $\al$ has only one pre-image say $\al'$ or it has two pre-images say $\al',\al'+\be$. In the former case, it is clear that $\al\in \Ga^Q_w$ if and only if $\al'\in \Ga_w^{P_1}$. Moreover we readily check that we have $n_{\pi_1}(\al',\la+\nu_1)=n_{\pi_Q}(\al,\la+\nu_Q)$ where the coefficients are those defined in \eqref{eq:NEW coef npi1}. In the latter case, we see   that $\al\in \Ga^Q_w$ if and only if both $\al', \al'+\be\in \Ga_w^{P_1}$, see lemma \ref{lem:simplicite}. To conclude it suffices to check that we have
    \begin{align}\label{eq:NEW the result}
      n_{\pi_1}(\al',\la+\nu_1) n_{\pi_1}(\al'+\be,\la+\nu_1)=n_{\pi_Q}(\al,\la+\nu_Q).
    \end{align}
    We shall use the notations and the computations of  § \ref{S:NEW calculs expl}.  Recall that we have $1\leq i_0<j_0\leq r$. We assume that $\be$ is a root ``inside'' the factor $G_{N_{i_0}}$. The  case  where $\be$ is a root ``inside''  the factor $G_{N_{j_0}}$ is similar and left to the reader. So we assume that $\al'=\be_{i,j}$ with $2\leq i\leq k$ and $1\leq j\leq l$ and $\be= \be_{i-1,i}$ relates the block $G_{e_{i-1}r_{i_0}}$ to the block $G_{e_{i}r_{i_0}}$.  We set $s=\bg \la, \al^\vee\bd\in \CC$. We have  $\bg \la, (\al')^\vee\bd=  \bg \la, (\al')^\vee+\be^\vee\bd=s $. We also set
    \begin{align*}
      \delta=(d_{j_0}-d_{i_0})/2 -(f_1+\ldots+f_{j-1})-f_j/2 .
    \end{align*}
    We have
    \begin{align*}
 &     \bg \nu_1, (\al')^\vee\bd=e_1+\ldots+e_{i-1}+e_i/2+\delta.
   &   \bg \nu_1, (\al')^\vee+\be^\vee\bd=e_1+\ldots+e_{i-2}+e_{i-1}/2 +\delta.
    \end{align*}
    We get:
    \begin{align*}
      n_{\pi_1}(\al',\la+\nu_1)&= \prod_{t=0}^{e_i-1} \xi_{\tau_{i_0},{\sigma_j'}^\vee}(s+ \bg \nu_1, (\al')^\vee\bd+(e_i-1)/2-t)\\
                             &= \prod_{t=1}^{e_i} \xi_{\tau_{i_0},{\sigma_j'}^\vee}(s+e_1+\ldots+e_{i-1}+t+\delta-1/2)
    \end{align*}
    \begin{align*}
        n_{\pi_1}(\al'+\be,\la+\nu_1)&= \prod_{t=1}^{e_{i-1}}  \xi_{\tau_{i_0},{\sigma_j'}^\vee}(s+e_1+\ldots+e_{i-2}+t+\delta-1/2)
    \end{align*}
    So we get:
    \begin{align*}
      n_{\pi_1}(\al',\la+\nu_1)n_{\pi_1}(\al'+\be,\la+\nu_1)= \prod_{t=1}^{e_{i-1}+e_i}  \xi_{\tau_{i_0},{\sigma_j'}^\vee}(s+e_1+\ldots+e_{i-2}+t+\delta-1/2).
    \end{align*}
    We write $Q\cap M=Q^1\times \ldots\times Q^r$ according to the decomposition $M=G_{N_1}\times \ldots\times G_{N_r}$. Accordingly we write $\pi_Q=\pi_Q^{1}\boxtimes \ldots\boxtimes \pi_Q^{r}$ where $\pi_Q^i$ is a discrete automorphic representation of the standard Levi factor of $Q^i$. Note that we have $Q^i=P_1^i$ and  $\pi_Q^{i}=\pi_1^i$ if $i\not=i_0$. The standard Levi factor of $Q^{i_0}$ is naturally identified with $G_{e_1r_0}\times \ldots \times G_{e_{i-2}r_0}\times G_{(e_{i-1}+e_i)r_0}\times G_{e_{i+1}r_0}\times\ldots G_{e_kr_0}$ and  $\pi_Q^{i_0}$ is the external tensor product of the discrete automorphic representations associated to $\tau_{i_0}$ by M\oe glin-Waldspurger classification. So the same computation as above, applied to $n_{\pi_Q}(\al,\la+\nu_Q)$,  gives the equality \eqref{eq:NEW the result}.
  \end{preuve}
\end{paragr}

\begin{paragr}
  Locally the map  $\la\in \ago_{P,\CC}^*\mapsto N_{\pi_1}(w,\la+\nu_1)\varphi_1$ has only hyperplanes singularities along  hyperplanes $\bg \la,\theta^\vee\bd=c$ for some  $c\in \CC$ and $\theta\in \Ga_w$. Let $\theta \in \Ga_w$ and let $\theta_P^\vee$ be the projection of  $\theta^\vee$ on  $\ago_P$. We shall say that $\theta_P^\vee$  is the direction of the singular hyperplane   $\bg \la,\theta^\vee\bd=c$.
   In the following we fix $\theta\in \Ga_w$ and we will prove the following lemma \ref{lem: clef} \emph{under the induction hypothesis}. 

  \begin{lemme}\label{lem: clef}
    No singular hyperplane of direction $\theta_P^\vee$ of  $\la\in \ago_{P,\CC}^*\mapsto N_{\pi_1}(w,\la+\nu_1)\varphi_1$ meets $\Xi_P$.
      \end{lemme}

By varying the root $\theta$ we see that lemma \ref{lem: clef} implies  proposition \ref{prop:le but}.
\end{paragr}

\begin{paragr}[Start of the proof of lemma \ref{lem: clef}.] --- Recall that we fix  $\theta \in \Ga_w$.
    
  \begin{lemme} \label{lem:reduc theta} Assume  $\ga_{i_1,i_s+1}\not\prec \theta$. Then, under the induction hypothesis,  lemma \ref{lem: clef} holds.
  \end{lemme}

  \begin{preuve}
   We may assume that  $\theta$ belongs to  the set $\Ga_w^{\be_{i_s}}$ defined in \eqref{eq:Ga be w}  (the argument is similar if  $\theta \in \Ga_w^{\be_{i_1}}$). Since the set  $\Ga_w^{\be_{i_s}}$  is closed and coclosed by lemma \ref{lem:coclos} there exists $w_1\in W_{M_1}$ such that $\Ga_{w_1}=\Ga_w^{\be_{i_s}}$, see lemma \ref{lem:set of inv}. Set $w_2=ww_1^{-1}$. By  lemma \ref{lem:additivite}, we have $\ell(w)=\ell(w_2)+\ell(w_1)$ and 
    \begin{align*}
  N_{\pi_1}(w,\la+\nu_1)\varphi_1=    N_{w_2\pi_1}(w_2,w_1(\la+\nu_1)) N_{\pi_1}(w_1,\la+\nu_1)\varphi_1.
    \end{align*}
    The directions of the singular  hyperplanes of  $ \la\mapsto  N_{w_2\pi_1}(w_2,w_1(\la+\nu_1)) $ are distinct from $\theta_P^\vee$. Indeed for any element $\al\in w_1^{-1}\Ga_{w_2}$ we have 
  $\be_{i_s}\prec \al$. Since $\be_{i_s}\in \Delta_1\setminus \Delta_1^P$, the projections of $\al^\vee$ and $\theta^\vee$ on $\ago_P$ cannot be equal.  Since we have $\ell(w_1)<\ell(w)$, we can use the induction hypothesis to deduce the result.
\end{preuve}

By lemma \ref{lem:reduc theta}, we may and shall assume that  $\theta=\ga_{k_0,l_0}$ with $1\leq k_0\leq i_1$ and $ i_{s}+1\leq l_0\leq r+1$. Observe that if $\be_{k_0}\in \Delta_1^P$ then   $k_0< i_1$ and we see that $\ga_{k_0+1,l_0}$ belongs to  $\Ga_w$ and that  its projection onto $\ago_P$ is $\theta_P$. In the same way, if $\be_{l_0-1}\in \Delta_1^P$ then  $l_0-1>i_s$ and $\ga_{k_0,l_0-1}$ belongs to $\Ga_w$ with projection  $\theta_P$ onto $\ago_P$.  So we may and shall moreover assume that  $\theta$ is such that $\be_{k_0}\notin \Delta_1^P$ and $\be_{l_0-1}\notin \Delta_1^P$. Let $0\leq p\leq k_0-1$ and $0\leq q\leq r-l_0+1$ be the greatest integers such that the sets $\{\be_{k_0-p},\ldots, \be_{k_0-1}\}$ and $\{\be_{l_0},\ldots, \be_{l_0+q-1}\}$, of respective cardinality   $p$ and $q$, are  included in $ \Delta_1^P$.

Set
 \begin{align*}
   \Ga_\theta=\{\ga\in \Ga_w\mid \ga=\ga_{k,l} \text{ with }  k_0-p\leq k \leq k_0 \text{ and } l_0\leq l \leq l_0+q\}.
 \end{align*}
Note that   $\Ga_\theta$ is nothing else but the set of elements  in $ \Ga_w$ whose projection onto $\ago_P$ is  $\theta_P$.

\begin{lemme} \label{lem: translation} Let $1\leq l\leq r$ be an integer such that  $l>i_s$ or $l<i_1$ and  $\be_l\in\Delta_1^P$. Assume that for all $\gamma\in \Ga_\theta$ such that  $\ga+\be_l\in \Phi_1$ we have  $\ga+\be_l\in \Ga_\theta$.  Then, under the induction hypothesis,  lemma \ref{lem: clef} holds.
\end{lemme}

\begin{preuve}
By lemma \ref{lem: Ga+al} and lemma \ref{lem:set of inv}, there exists $w'\in W_{M_1}$ such that $ \Ga_{w'}$ is the union of  $\Ga_w$ and the set of  $\ga+\be_{l}$ for all  $\ga\in \Ga_w$ such that $\ga+\be_{l}\in \Phi_1$.We have
\begin{align*}
  N_{\pi_1}(w',\la+\nu_1)\varphi_1=N_{w\pi_1}(w'w^{-1},w(\la+\nu_1))N_{\pi_1}(w,\la+\nu_1)\varphi_1.
\end{align*}
By lemma \ref{lem:simplicite} and the very construction of $\Ga_{w'}$ we see that $w' \be_l$ is a simple root. Lemma \ref{lem:rec terme cst} and the induction hypothesis imply that the left-hand side is holomorphic on $\Xi_P$. Let $\gamma\in \Ga_{w'}$ be such that the projection of $\gamma$ onto $\ago_P$ is equal to $\theta_P$. We claim that  $\gamma\in\Ga_{\theta}$: indeed if $\gamma=\gamma'+\be_l$ with $\gamma'\in \Ga_w$ then we have $\gamma'\in \Ga_\theta$ and, by the main hypothesis of the lemma, we have $\gamma\in \Ga_\theta$ and thus $\gamma\in \Ga_w$. In particular, we see that no hyperplane of direction $\theta_P$ is singular for the operator $N_{w\pi_1}(w'w^{-1},w(\la+\nu_1))$ and its inverse. The lemma follows.
\end{preuve}

\begin{lemme}\label{lem:no N} Under the induction hypothesis,  lemma \ref{lem: clef} holds unless there exists an integer $\max(p-1,q-1)\leq N \leq \min(p,q)$ such that 
\begin{align*}
   \Ga_\theta=\{\ga_{k_0-i,l_0+j} \mid i+j\leq N, 0\leq i\leq p,0\leq j \leq q\}. 
\end{align*}
\end{lemme}

\begin{preuve}
  Recall that the roots in $ \Phi_1$ whose projection onto $\ago_P$ is  $\theta_P$ are the roots  $\ga_{k_0-i,l_0+j}$ with $0\leq i\leq p$ and $0\leq j \leq q$.   We shall implicitly use the following fact: if $\ga_{k_0-i,l_0+j}\in \Ga_\theta$ with $0\leq i\leq p$ and $0\leq j \leq q$ then $\ga_{k_0-i',l_0+j'}\in \Ga_\theta$  for  $ 0\leq i'\leq i $ and $0\leq j'\leq j$. This follows from the fact that $\Ga_w$ is coclosed and the roots $\be_{k_0-i'}$ and $\be_{l_0+j'}$ do not belong to $\Ga_w$ for $1\leq i'\leq p$ and $0\leq j'\leq q-1$.
  
  Let $N$ be the greatest integer  $i+j$ such that   $\ga_{k_0-i,l_0+j}\in \Ga_\theta $ with $0\leq i\leq p$ and $0\leq j \leq q$.  Note that $0\leq N\leq p+q$. Let $t$ be such that  $\max(0, N-q) \leq t \leq \min(p, N)$ and $\ga_{k_0-t,l_0+ N-t}\in \Ga_\theta $. So for all $0\leq i \leq t$ we have  $\ga_{k_0-i,l_0+ N-t}\in \Ga_\theta $ and for all $0\leq j \leq N-t$ we have  $\ga_{k_0-t,l_0+ j} \in \Ga_\theta $. 

 Observation (1): we assume $N-t>0$. Then we also have  $\ga_{k_0-i,l_0+ N-t-1}\in \Ga_\theta $   for all $0\leq i \leq t$. Assume moreover that either we have $t=p$ or we have $t<p$ and $\ga_{k_0-(t+1),l_0+ N-t-1}\notin \Ga_\theta $. Then we can see that the set of $\ga\in \Ga_\theta$ such that  $\ga+\be_{l_0+N-t-1}\in \Phi_1$ is exactly the set of roots  $\ga_{k_0-i,l_0+ N-t-1}$ for $0\leq i\leq t$. In particular, for all $\ga\in \Ga_\theta$ the condition  $\ga+\be_{l_0+N-t-1}\in \Phi_1$ implies  $\ga+\be_{l_0+N-t-1}\in \Ga_\theta$. Since we have $l_0+N-t-1 >l_0-1\geq i_s$, we conclude, by  lemma \ref{lem: translation}, that lemma \ref{lem: clef} holds under the induction hypothesis.

 Observation (2): we   assume $t>0$.   Then we also have  $\ga_{k_0-(t-1),l_0+ j}\in \Ga_\theta $   for all $0\leq j \leq N-t$. Assume moreover that either we have $N-t=q$ or we have $N-t<q$ and $\ga_{k_0-(t-1),l_0+ N-t+1}\notin \Ga_\theta $. Then the set of $\ga\in \Ga_\theta$ such that  $\ga+\be_{k_0-t}\in \Phi_1$ is exactly the set of roots  $\ga_{k_0-(t-1),l_0+ j}$ for $0\leq j\leq N-t$. So for all $\ga\in \Ga_\theta$ the condition  $\ga+\be_{k_0-t}\in \Phi_1$ implies  $\ga+\be_{k_0-t}\in \Ga_\theta$. Since $k_0-t< k_0\leq i_1$,  we conclude, by  lemma \ref{lem: translation},   that lemma \ref{lem: clef} holds under the induction hypothesis.

Observation (3): Assume there exists $t<i \leq \min(p,N)$ such that  $\ga_{k_0-i,l_0+ N-i}\notin \Ga_\theta $ and we assume $i$ is  the smallest such integer. So $\ga_{k_0-(i-1),l_0+ N-i+1}\in \Ga_\theta $ and $i-1<i\leq N$. By observation (1), we conclude that  lemma \ref{lem: clef} holds under the induction hypothesis.

Observation (4): Assume there exists $\max(0,N-q)\leq j<t$ such that  $\ga_{k_0-j,l_0+ N-j}\notin \Ga_\theta $ and we assume $j$ is  the greatest such integer. So $\ga_{k_0-(j+1),l_0+ N-j-1}\in \Ga_\theta $ and $0\leq j<j+1$. By observation (2), we conclude that  lemma \ref{lem: clef} holds under the induction hypothesis.

By observations (3) and (4), we see that we can conclude that  lemma \ref{lem: clef} holds  unless we have:
\begin{align*}
   \Ga_\theta=\{\ga_{k_0-i,l_0+j} \mid i+j\leq N, 0\leq i\leq p,0\leq j \leq q\}. 
\end{align*}
In the rest of the proof, we assume that the equality above holds. If $N>p$ then  $\ga_{k_0-p,l_0+N-p}\in \Ga_{\theta}$. By   observation (1), we conclude that  lemma \ref{lem: clef} holds. If $N>q$ then    $\ga_{k_0-(N-q),l_0+q}\in \Ga_{\theta}$. By   observation (2), we conclude that  lemma \ref{lem: clef} holds. So we assume $N\leq \min(p,q)$.  If $N\leq p-2$ then $\be_{k_0-(N+2)}\in \Delta_1^P$ with $k_0-(N+2)<i_1$ and there is no $\ga\in \Ga_\theta$ such that $\ga+\be_{k_0-(N+2)}\in\Phi_1$.  So by   lemma \ref{lem: translation}, one can conclude. If $N\leq q-2$ then $\be_{l_0+N+1}\in \Delta_1^P$ and $l_0+N+1>i_s$. Since  there is no $\ga\in \Ga_\theta$ such that $\ga+\be_{l_0+N+1}\in\Phi_1$, one can conclude by  lemma \ref{lem: translation}.
\end{preuve}

\end{paragr}

\begin{paragr}[End of the proof of lemma \ref{lem: clef}.] --- By lemma \ref{lem:no N} we shall assume  from now on  that there exists $\max(p-1,q-1)\leq N\leq \min(p,q)$ such that
\begin{align*}
  \Ga_\theta=\{\ga_{k_0-i,l_0+j} \mid i+j\leq N, 0\leq i\leq p,0\leq j \leq q\}.
\end{align*}
  We fix a reduced decomposition \eqref{eq:decomp} of  $w$  so that we have $\Ga_w=\{\theta_1,\ldots,\theta_\ell\}$. We have  $w=s_{\al_\ell} w_1$ for some $w_1\in W_{M_1}$ such that $w_1\Delta_0^P>0$. We have
  \begin{align*}
    N_{\pi_1}(w,\la+\nu_1)=N_{w_1\pi_1}(s_{\al_\ell}, w_1(\la+\nu_1)) N_{\pi_1}(w_1,\la+\nu_1).
  \end{align*}
If one of the following conditions is satisfied
  \begin{itemize}
  \item $\theta_\ell\notin \Ga_\theta$;
  \item $\bg \nu_1, \theta_{\ell}^\vee\bd \geq 0$,
  \end{itemize}
the operator $N_{w_1\pi_1}(s_{\al_\ell}, w_1(\la+\nu_1)) $ has no    singular hyperplanes of direction $\theta_P^\vee$  which meets $\Xi_P$ (see lemma \ref{lem:holom op simple}) and lemma \ref{lem: clef} follows from the induction hypothesis.  From now we assume  $\theta_\ell\in \Ga_\theta$ and  $\bg \nu_1, \theta_{\ell}^\vee\bd < 0$. By proposition \ref{prop:w-tableau}, the reduced decomposition of $w$ defines a $w$-tableau from which we see that we must have $\theta_\ell=\ga_{k_0-N+j,l_0+j}$ for some  $0\leq j \leq N$ (use inequalities \eqref{eq:couple ineg}). 

We have to distinguish several cases according to the values of $j$.

Let us first assume that we have $j\leq q-1$. Then  $\be_{l_0+j}\in \Delta_1^P$ and $\theta_\ell+\be_{l_0+j}\in \Phi_1$. By lemma \ref{lem: Ga+al}, the union of $\Ga_w$ and the set of  $\ga+\be_{l_0+j}$ for $\ga\in \Ga_w$ such that  $\ga+\be_{l_0+j}\in \Phi_1$ is closed and coclosed and thus is equal to  $ \Ga_{w'}$ for some element $w'\in W_{M_1}$ such that $w'\Delta_0^P>0$. We have
\begin{align*}
  N_{\pi_1}(w',\la+\nu_1)\varphi_1=N_{w \pi_1}(w'w^{-1},w(\la+\nu_1))N_{\pi_1}(w,\la+\nu_1)\varphi_1.
\end{align*}
By lemma \ref{lem:simplicite},  we see that $w' \be_{l_0+j}$ is a simple root. By  lemma \ref{lem:rec terme cst}, the left-hand side is holomorphic on  $\Xi_P$. Let $\ell=\ell(w)$ and $\ell'=\ell(w')$. By  lemma \ref{lem:additivite}, we have $\ell'=\ell(w'w^{-1})+\ell$. Let  $w'w^{-1}=s_{\al_{\ell'}} \cdots s_{\al_{\ell+1}}$  be a reduced decomposition. We have  
\begin{align*}
  N_{w \pi_1}(w'w^{-1},w(\la+\nu_1))^{-1}\\= N_{w\pi_1}( s_{\al_{\ell+1}},w(\la+\nu_1))^{-1} \cdots N_{ s_{\al_{\ell'}}^{-1}w'\pi_1 }( s_{\al_{\ell'}},    s_{\al_{\ell'}}^{-1}w'(\la+\nu_1))^{-1} .
\end{align*}
One can check that the set  of  $\ga\in \Ga_{w'}\setminus \Ga_w$ whose projection on $\ago_P$ is $\theta_P$ is the singleton $\{\theta_\ell+\be_{l_0+j}\}$. Thus there is in the decomposition above only one operator that may have singularities of direction $\theta_P^\vee$. It corresponds to the unique integer $k$ such that  $1\leq k\leq \ell'-\ell$  and  $\theta_\ell+\be_{l_0+j}=w^{-1}s_{\al_{\ell+1}}^{-1}\cdots s_{\al_{\ell+k-1}}^{-1}(\al_{\ell+k})$. We have
\begin{align*}
  N_{ (s_{\al_{\ell+k-1}}\cdots s_{\al_{\ell+1}}w)\pi_1  }(s_{\al_{\ell+k}}, (s_{\al_{\ell+k-1}}\cdots s_{\al_{\ell+1}}w)(\la+\nu_1))^{-1}\\=N_{(s_{\al_{\ell+k}}\cdots s_{\al_{\ell+1}}w)\pi_1  }(s_{\al_{\ell+k}}^{-1}, (s_{\al_{\ell+k}}\cdots s_{\al_{\ell+1}}w)(\la+\nu_1)).
\end{align*}
We claim that the  operator above is holomorphic on $\Xi_P$. Indeed by  lemma \ref{lem:holom op simple} and the fact that  $s_{\al_{\ell+k}}^{-1}=s_{-s_{\al_{\ell+k}}(\al_{\ell+k})}$,  it suffices to check that we have 
 \begin{align*}
   \Re(\bg     (s_{\al_{\ell+k}}\cdots s_{\al_\ell+1}w)(\la+\nu_1), -s_{\al_{\ell+k }}(\al_{\ell+k}^\vee)   \bd)> -\frac{2}{n^2+1}
 \end{align*}
 for all $\la\in \Xi_P$. But the left-hand side of the inequality is
 \begin{align*}
     -\Re(\bg     (s_{\al_{\ell+k-1}}\cdots s_{\al_\ell+1}w)(\la+\nu_1), \al_{\ell+k}^\vee   \bd)=- \Re(\bg \la+\nu_1,\theta_\ell^\vee+\be_{l_0+j}^\vee\bd)> -\frac{2}{n^2+1}
 \end{align*}
for all $\la\in \Xi_P$: this follows from  lemma \ref{lem:Delta1P} and our hypothesis $\bg \nu_1,\theta_\ell^\vee\bd<0$. We can conclude that lemma \ref{lem: clef} is true for
 \begin{align}\label{eq:dec N w'w}
   N_{\pi_1}(w,\la+\nu_1)\varphi_1=N_{w\pi_1}(w'w^{-1},w(\la+\nu_1))^{-1}N_{\pi_1}(w',\la+\nu_1)\varphi_1.
 \end{align}

 Then we assume  $j=q$. In particular we have $N=q\leq p$. If $p=0$ then $q=0$. In this case one can check that we have $\bg \nu_1, \theta_{\ell}^\vee\bd =0$ unlike our assumption  $\bg \nu_1, \theta_{\ell}^\vee\bd <0$. So we assume moreover that we have  $p\geq 1$ that is  $N\leq p+q-1$. Then  $\be_{k_0-N+q-1}\in \Delta_1^P$ and  $\theta_\ell+\be_{k_0-N+q-1}\in \Phi_1$.  We can conclude by the same argument as above by introducing  $w'\in W_{M_1}$ such that $ \Ga_{w'}$ is the union of $\Ga_w$ and the set of  $\ga+\be_{k_0-N+q-1}$ such that $\ga\in \Ga_w$ and $\ga+\be_{k_0-N+q-1}\in \Phi_1$.
\end{paragr}

\subsection{Bounds for normalized intertwining operators}\label{ssec:bd norm inter}

\begin{paragr} For all standard parabolic subgroups $P$ and all $\eps>0$ we set
 \begin{align*}
    \Xi_{P,\eps}=\{\la\in  \ago_{P,\CC}^{G,*}\mid   \|\Re(\la)\| < \eps\}.
  \end{align*}
In the following we assume that  $\eps$ is small enough so that we have   $\Xi_{P,\eps}\subset \Xi_P$ where $\Xi_P$ is defined in \eqref{eq:def-Om-New}. The bulk of this subsection is devoted to the proof of the following proposition.
  
  \begin{proposition}\label{prop:NEWbd-norm-inter}
     There exists $k>0$ such that for  any level   $J$  there exist $\eps>0$ and $c_J>0$ such that for all  $J$-triples $(P,\pi,\tau)$, all  standard parabolic subgroups $R$ of $G$,  all $w\in \, _RW_{P}$  such that  $  P_\pi\subset P_w$ and $R_w=R$ we have
\begin{align*}
  \| N_{\pi_{P_w}}(w,\la+\nu_w)\varphi_{w}\|_R \leq c_J (1+\la_\tau^2)^k\|\varphi_w\|_{P_w}
\end{align*}
for all $\varphi\in \Ac_{P,\pi}(G)^{\tau,J}$ and for all $\la\in  \Xi_{P,\eps}$.
  \end{proposition}

\end{paragr}

\begin{paragr}  \label{S:NEW level}Let $J$ be a level and $(P,\pi,\tau)$ be a $J$-triple. Without loss of generality we shall assume in the following that $J=\prod_{v\notin V_{F,\infty}}J_{v}$  where $J_v\subset K_v$ is a  normal  open compact subgroup of $K_v$.  Let $S'$ be the finite set of finite places $v$ such that $J_v\subsetneq K_v$. We set $S=S'\cup  V_{F,\infty}$.  Let $R$ be a parabolic subgroup  of $G$ and  $w\in \, _RW_{P}$  be such that  $  P_\pi\subset P_w$ and $R_w=R$. We set $P_1=P_w$. With  notations of § \ref{S:elem sym}, we have $w\in W_{M_1}$, we also fix a decomposition of $w$ into elementary symmetries as in \eqref{eq:decomp}: in particular  we have $R=P_{\ell+1}$ where $\ell=\ell(w)$. We use also the notations of § \ref{S:decomp inter}.
\end{paragr}

\begin{paragr} \label{S:NEW localnorm} We fix $1\leq i \leq \ell$.  We start by considering the normalized operator $N_{\pi_i}(s_{\al_i},\mu)$ on $\Ac_{P_i,\pi_i}(G)$ for $\mu\in \ago_{P_i,\CC}^*$. It decomposes as  a tensor product of local  intertwining operators $N_{\pi_i,v}(s_{\al_i},\mu)$  as in    \eqref{eq:local decomposition}. By  choosing  a decomposition of the Petersson norm on $\Ac_{P_i,\pi_i}(G)$, we  get a norm on each induced space $\Ind_{P_{i}}^G(\pi_{i,v})$. Let $v\in V_F\setminus V_{F,\infty}$ and $q_v$ be the order of the residue field at $v$.  Let $\Ind_{P_{i}}^G(\pi_{i,v})^{J_v}\subset \Ind_{P_{i}}^G(\pi_{i,v})$ be the     $J_v$-fixed subspace.
  Let   $N_{\pi_i,v}(s_{\al_i},\mu)^{J_v}$ be the operator from $\Ind_{P_{i}}^G(\pi_{i,v})^{J_v}$ to $\Ind_{P_{i+1}}^G(\pi_{i+1,v})^{J_v}$ we get by restriction of the operator  $N_{\pi_i,v}(s_{\al_i},\mu)$.  We denote by $\|N_{\pi_i,v}(s_{\al_i},\mu)^{J_v} \|$ the operator norm  of $N_{\pi_i,v}(s_{\al_i},\mu)^{J_v}$. We may and shall assume that, for $v\in V_F\setminus S$,  we have $\|N_{\pi_i,v}(s_{\al_i},\mu)^{J_v} \|=1$. So in the following we are mainly concerned with the case $v\in S$.
\end{paragr}

\begin{paragr}[Case of a finite place.] ---
  
  \begin{lemme} \label{lem:bd finite place}Let $v\in S'$.
    \begin{enumerate}
    \item There exist two polynomials  of the $z$-variable $Q_d(z)\in \CC[z]$ (of degree denoted by $d\geq 0$) and  $Q_n(z)$ with values in the space of linear maps   from $\Ind_{P_{i}}^G(\pi_{i,v})^{J_v}$ to $\Ind_{P_{i+1}}^G(\pi_{i+1,v})^{J_v}$ such that
      \begin{align*}
        N_{\pi_i,v}(s_{\al_i},\mu)^{J_v}=Q(q_v^{-\bg \mu,\al_i^\vee\bd}) \text{ with } Q=\frac{Q_n}{Q_d}
      \end{align*}
  for all $\mu\in \ago_{P_i,\CC}^*$.
\item The degree $d$ can be bounded independently of $J_v$ and $\pi$ and the degree  $\deg(Q_n)$ can be  bounded in terms of $J_v$ only.
  \item For all $z\in \CC$ such that $|z|\leq 1$ we have  $\|Q(z)\|\leq 1$. 
    \item     There exists a family of complex numbers $(z_i)_{1\leq i\leq d}$ of modulus $\leq q^{-2/(n^2+1)}$ such that we have 
      \begin{align*}
         \|\left(\prod_{j=1}^{d}(1-zz_j)\right)Q(z)\|\leq 2^{d} |z|^{\deg(Q_n)}
      \end{align*}
      for all $z\in \CC$ such that $|z|\geq 1$.
    \item Let $\frac{2}{n^2+1}>\eta>0$. There exists $C>1$ which depends only on $\eta$ and  $J_v$ such that for all $z\in \CC$ such that  $  |z|< q^{2/(n^2+1)-\eta}$ we have
      \begin{align*}
        \|Q(z)\|\leq C  .
      \end{align*}
    \end{enumerate}
  \end{lemme}

  \begin{preuve} This is a variation on  \cite[proofs of lemmas 3.10 and 3.11]{LapHC}. The first assertion is a basic property of normalized intertwining operators. The assertions on $\deg(Q_d)$ and $\deg(Q_n)$ follow respectively from  \cite[section IV.1]{Walds-Plan} and \cite[theorem 1 and proposition 5]{FLM-degree}.   Let $(z_j)_{1\leq j\leq d}$ be the roots of $Q_d$ counted with multiplicity.  Note  that the operator $N_{\pi_i,v}(s_{\al_i},\mu)^{J_v}$ is bounded for all $\mu\in \ago_{P_i,\CC}^*$ such that $\Re(\bg \mu,\al_i^\vee\bd)$ is sufficiently large: this follows from the proof of \cite[théorème  IV.1.1]{Walds-Plan} and the shape of the normalization factor. In particular, we get $Q_d(0)\not=0$.    By   lemma \ref{lem:holom op simple}, we may assume $|z_j|\geq q^{2/(n^2+1)}$ for $1\leq j\leq d$.  In particular $Q$ is holomorphic on the open subset of $z\in \CC$ such that $|z|< q^{2/(n^2+1)}$. Since  $N_{\pi_i,v}(s_{\al_i},\mu)^{J_v}$ is unitary for $\mu\in i\ago_{P_i}^*$, we get that $\|Q(z)\|=1$ for $|z|=1$. The assertion 3 follows from the  maximum modulus principle.   The operator 
\begin{align*}
  z^{\deg(Q)} \left(\prod_{j=1}^d \frac{z-z_j^{-1}}{1-\bar{z}_j^{-1}z}  \right)Q(z^{-1})= z^{\deg(Q_n)} \left(\prod_{j=1}^d  \frac{1-(z z_j)^{-1}}{1-\bar{z}_j^{-1}z}  \right)Q(z^{-1})
\end{align*}
is holomorphic for $|z|\leq 1$ and of norm  $1$ for $|z|=1$. Thus by the maximum modulus principle we get

\begin{align*}
  \| \left(\prod_{j=1}^d  (1-z z_j^{-1})\right)Q(z)\|\leq |z|^{\deg(Q_n)}  \prod_{j=1}^d  |1-\bar{z}_j^{-1}z^{-1}|\leq 2^{d} |z|^{\deg(Q_n)}  
\end{align*}
for all $z\in \CC$ such that $|z|\geq 1$. So we get assertion 4.

Finally from assertion 3 and 4,  for all $z\in \CC$ such that  $  |z|< q^{2/(n^2+1)-\eta}$ we have
  \begin{align*}
        \|Q(z)\|\leq  \max(1, 2^d q^{2\deg(Q_n) /(n^2+1)} (1-q^{-\eta})^{-d}).
      \end{align*}
Then  assertion 5 follows from assertion 2.
\end{preuve}
\end{paragr}

\begin{paragr}[Case of an Archimedean place.] --- We fix in this § an Archimedean place $v\in V_{F,\infty}$. Let $\tau_v$ be an irreducible unitary representation of $K_v$ and let $\la_{\tau_v}$ be the eigenvalue of the Casimir operator of $K_v$.   For all $\mu\in \ago_{P_i,\CC}^*$ we denote by  $N_{\pi_i,v}(s_{\al_i},\mu)^{\tau_v}$  the operator from  $\Ind_{P_{i}}^G(\pi_{i,v})^{\tau_v}$ to  $\Ind_{P_{i+1}}^G(\pi_{i+1,v})^{\tau_v}$ given by the restriction of $N_{\pi_i,v}(s_{\al_i},\mu)$. Here we have denoted by an upper script $\tau_v$  the $\tau_v$-isotypic components. As in § \ref{S:NEW localnorm} we  denote by $\|  N_{\pi_i,v}(s_{\al_i},\mu)^{\tau_v}  \|$ the operator norm of  $N_{\pi_i,v}(s_{\al_i},\mu)^{\tau_v}$. For all $s,a\in \CC$ we set:
  \begin{align*}
    B_a(s)=\frac{s-a}{s+\bar a}.
  \end{align*}
We shall use the following inequality which holds if moreover  $\Re(a)\Re(s)\geq 0$
    \begin{align}\label{eq:ineg reel}
     \left| B_a(s)\right|\geq \left|\frac{\Re(s-a)}{\Re(s+a)}  \right|.
    \end{align}

  \begin{lemme}\label{lem:norm inter archi}
    \begin{enumerate}
          \item For all $\mu\in \ago_{P_i,\CC}^*$ such that $\Re(\bg \mu,\al_i^\vee\bd)\geq 0$ the map $\mu \mapsto N_{\pi_i,v}(s_{\al_i},\mu)^{\tau_v}$ is holomorphic and we have
            \begin{align*}
              \| N_{\pi_i,v}(s_{\al_i},\mu)^{\tau_v}\|\leq 1.
            \end{align*}
          \item There exist an  integer  $r>0$ depending only on $G$, an integer $0\leq k\leq r$, a family $(s_l)_{1\leq l\leq k}$ of complex numbers of real parts  $\Re(s_l)\leq -2/(n^2+1)$ such that if  we set 
              \begin{align*}
              f(\mu)=\left(  \prod_{l=1}^k \prod_{j=0}^{m_l}   B_{s_l-j}(\bg \mu,\al^\vee_i\bd)  \right)N_{\pi_i,v}(s_{\al_i},\mu)^{\tau_v},
              \end{align*}
  where $m_l$ is the integer part of $\Re(s_l)+r(1+\la_{\tau_v}^2)$,   then for   all $\mu\in \ago_{P_i,\CC}^*$ such that $\Re(\bg \mu,\al_i^\vee\bd)\leq 0$ the operator   $f(\mu)$  is holomorphic and  we have  $              \|f(\mu)\|\leq 1$.
             \item We use the  notations of assertion 2 above. Let $s\in \CC$ be such that $\Re(s)<0$ and $\eta\in \RR$ be  such that  $0<\eta <\min(1/2(n^2+1),-\Re(s))$. Let  $I$ be the set of integers $1\leq l\leq k$ such that there exists  $0\leq j\leq m_l$ such that
  \begin{align}\label{eq:NEW ineg reelle II}
|   \Re(s-s_l+j)|\leq \frac1{2(n^2+1)}.
  \end{align}
  For $l\in I$, we set $s_l'=s_l-j$ where $j$ is the unique integer  $0\leq j\leq m_l$ that satisfies the inequality \eqref{eq:NEW ineg reelle II}.  There exist  $c,e>0$ (which depend only on $G$, $\eta$ and $s$) such that the operator
     \begin{align*}
              \tilde f(\mu)=\left(  \prod_{l\in I} B_{s_l'}(\bg \mu,\al^\vee_i\bd)     \right)N_{\pi_i,v}(s_{\al_i},\mu)^{\tau_v},
     \end{align*}
     is holomorphic and satisfies  $              \|\tilde f(\mu)\|\leq  c(1+\la_{\tau_v}^2)^{e}$     for all $\mu\in \ago_{P_i,\CC}^*$ such that $|\Re(\bg \mu,\al_i^\vee \bd -s)|<\eta $.  
        \item The results of assertion 3 hold for $s=0$,  $I=\emptyset$ and $0<\eta <1/2(n^2+1)$.
            \end{enumerate}
            
  \end{lemme}
  
  \begin{preuve} This is a variation on \cite[proof of lemma A.1]{MuSpeh}. For the reader's convenience, we recall the arguments. First we know that the intertwining operator $N_{\pi_i,v}(s_{\al_i},\mu)^{\tau_v}$ can be written  $Q(\bg \mu,\al_i^\vee\bd)$ for some ``rational fraction'' $Q(s)$ of the complex variable $s$. Since it is unitary on the imaginary axis we have $\| Q(s)\|=1$ for $s\in i\RR$.  By   lemma \ref{lem:holom op simple} we know that all the poles of $Q$ are of real part $\leq -2/(n^2+1)$. Then assertion 1 follows from a variant of the Phragmén-Lindelöf principle. Moreover the poles of $Q$ counted with multiplicities are included in $\cup_{i=1}^k (s_i-\NN)$ for some complex numbers $s_i$ of  real part  $\Re(s_i)\leq -2/(n^2+1)$ and the integer $k$ is bounded by an absolute constant $r$, see \cite[proposition A.2]{MuSpeh}. The same proposition shows that if $s_i-j$ with $j\in \NN$ is a pole then  we have $\Re(s_i-j)\geq - r(1+\la_{\tau_v}^2)$. So $j\leq  \Re(s_i)+ r(1+\la_{\tau_v}^2)$. The first assertion of 2  is then clear. Since $|B_{s_i-j}(s) |=1$ for $s\in i\RR$, the bound on the norm follows from an other application of a  variant of the Phragmén-Lindelöf principle.

    Let us prove assertion 3. Let $\mu\in \ago_{P_i,\CC}^*$ be such that $|\Re(\bg \mu,\al_i^\vee \bd -s)|<\eta $. Let $-N$ be the integer part of $2(\Re(s)-\eta)$. Note that $N\geq 1$.     In particular, we have
    \begin{align}\label{eq:NEW encadrement}
    -\Re(\bg \mu,\al_i^\vee \bd )-N<  \Re(\bg \mu,\al_i^\vee \bd )< \Re(s)+\eta<0.
    \end{align}
    
    Let $1\leq l\leq k$. By  \eqref{eq:NEW encadrement} we have  $\Re(\bg \mu,\al_i^\vee \bd )<0$. For all integers $0\leq j\leq m_l$, we have  $\Re(s_j-j) \leq \Re(s_j)<0$. Hence we can apply \eqref{eq:ineg reel} to get
    
 \begin{align}\label{eq:NEW minoration fonda}
     |     B_{s_l-j}(  \bg \mu,\al_i^\vee\bd)  |&\geq  \left|\frac{\Re(\bg \mu,\al_i^\vee\bd-s_l)+j)}{ \Re(\bg \mu,\al_i^\vee\bd+s_l)-j)}\right|.
 \end{align}

    Let $0\leq j_0\leq m_l$ be the least integer (if such an integer exists) such that
    \begin{align*}
        \Re(s-s_l)+j_0 > \frac1{2(n^2+1)}.
    \end{align*}
 For all integers $j \geq j_0$ we have
   \begin{align}
  \nonumber      \Re(\bg \mu,\al_i^\vee\bd-s_l)+j&\geq \Re(s-s_l)-\eta+j_0\\
   \label{eq:NEW minor2}    &\geq 1/2(n^2+1)-\eta>0.
   \end{align}

    We deduce that, for all integers $j_0\leq j\leq m_l$, the inequality \eqref{eq:NEW minoration fonda} implies that we have:
    \begin{align}
     \label{eq:NEW borner B 1}    |     B_{s_l-j}(  \bg \mu,\al_i^\vee\bd)  |             &\geq\frac{\Re(\bg \mu,\al_i^\vee\bd-s_l)+j}{ -\Re(\bg \mu,\al_i^\vee\bd+s_l)+j}>0.
    \end{align}
    Assume first  that we  have $j_0+N\leq m_l$. Then for any integer $j$ such that $j_0+N\leq j\leq m_l$  we shall use the following inequality which comes from  \eqref{eq:NEW encadrement}
    \begin{align}\label{eq:NEW borner B 2}
      |     B_{s_l-j}(  \bg \mu,\al_i^\vee\bd)  |&>\frac{-\Re(\bg \mu,\al_i^\vee\bd+s_l)+j-N}{ -\Re(\bg \mu,\al_i^\vee\bd+s_l)+j}>0.
    \end{align}
    Using the bound \eqref{eq:NEW borner B 1} for $j_0\leq j \leq j_0+N-1$ and the bound \eqref{eq:NEW borner B 2} for $j_0+N\leq j\leq m_l$  we get
   \begin{align*}
     \prod_{j=j_0}^{m_l}|     B_{s_l-j}(  \bg \mu,\al_i^\vee\bd)  |&\geq  \frac{  \prod_{j=j_0}^{j_0+N-1}(\Re(\bg \mu,\al_i^\vee\bd-s_l)+j)  \prod_{j=j_0}^{m_l-N}(-\Re(\bg \mu,\al_i^\vee\bd+s_l)+j)  }{\prod_{j=j_0}^{m_l}(-\Re(\bg \mu,\al_i^\vee\bd+s_l)+j)}\\
     &= \frac{  \prod_{j=j_0}^{j_0+N-1}(\Re(\bg \mu,\al_i^\vee\bd-s_l)+j)  }{\prod_{j=m_l-N+1}^{m_l}(-\Re(\bg \mu,\al_i^\vee\bd+s_l)+j)}.
   \end{align*}
   For all integers $j \leq m_l$, we have the following upper  bound for the denominator:
   \begin{align}
\nonumber     -\Re(\bg \mu,\al_i^\vee\bd+s_l)+j&\leq   -\Re(\bg \mu,\al_i^\vee\bd) +(-\Re(s_l)+m_l)\\
 \label{eq:NEW minor1}    &\leq -\Re(s)+\eta +r(1+\la_{\tau_v}^2).
   \end{align}
Using also the inequality \eqref{eq:NEW minor2} for the numerator,  we get:
  \begin{align}\label{eq:NEW puiss N}
     \prod_{j=j_0}^{m_l}|     B_{s_l-j}(  \bg \mu,\al_i^\vee\bd)  |&\geq   \left(\frac{ 1/2(n^2+1)-\eta}{-\Re(s)+\eta +r(1+\la_{\tau_v}^2)}\right)^N.
  \end{align}

  Now we assume   that we  have $j_0+N> m_l$. Using first the lower bound \eqref{eq:NEW borner B 1} and then the bounds \eqref{eq:NEW minor1}   and \eqref{eq:NEW minor2} we get for  $j_0\leq j\leq m_l$
      \begin{align*}
      |     B_{s_l-j}(  \bg \mu,\al_i^\vee\bd)  |&\geq\frac{  1/2(n^2+1)-\eta }{-\Re(s)+\eta +r(1+\la_{\tau_v}^2)}.
  \end{align*}
Since the right-hand side is $<1$ and $m_l-j_0+1\leq N$ we also get the inequality \eqref{eq:NEW puiss N} in this case.

  Let $0\leq j_1\leq m_l$ be the greatest integer (if such an integer exists) such that
    \begin{align*}
        \Re(s-s_l)+j_1 <- \frac1{2(n^2+1)}.
    \end{align*}
   
    Let $j$ be such that  $0\leq j\leq j_1$; we have:
     \begin{align*}
      \Re(\bg \mu,\al_i^\vee\bd-s_l)+j&\leq   \Re(\bg \mu,\al_i^\vee\bd-s)+\Re(s-s_l)+j_1 \\
  &    <\eta-\frac1{2(n^2+1)}<0.
          \end{align*}
Since we have also  $\Re(s_l-j)\leq \Re(s_l)<0$ and   $\Re(\bg \mu,\al_i^\vee \bd )<0$ by \eqref{eq:NEW encadrement}, we can use \eqref{eq:ineg reel} and \eqref{eq:NEW minor1} to get:
  \begin{align*}
      |     B_{s_l-j}(  \bg \mu,\al_i^\vee\bd)  |&\geq  \left|\frac{\Re(\bg \mu,\al_i^\vee\bd-s_l)+j)}{ \Re(\bg \mu,\al_i^\vee\bd+s_l)-j)}\right|\\
                                                 &=\frac{\Re(-\bg \mu,\al_i^\vee\bd+ s_l)-j}{ -\Re(\bg \mu,\al_i^\vee\bd+s_l)+j}\\
                                                     &\geq \frac{ 1/2(n^2+1)-\eta}{ -\Re(s)+\eta +r(1+\la_{\tau_v}^2)  }.
    \end{align*}
    Hence we have
      \begin{align}\label{eq:NEW3}
      |     \prod_{j=0}^{j_1}B_{s_l-j}(  \bg \mu,\al_i^\vee\bd)  |&\geq  \left(\frac{ 1/2(n^2+1)-\eta}{ -\Re(s)+\eta +r(1+\la_{\tau_v}^2)  }\right)^{-\Re(s)+1}
      \end{align}
         since we have  $j_1<-\Re(s)$.
Finally assertion 3 follows from  assertion 2, \eqref{eq:NEW puiss N} and  \eqref{eq:NEW3}. Assertion 4 is \cite[lemma A.1 and proposition A.2]{MuSpeh}.
  \end{preuve}

\end{paragr}

\begin{paragr}[End of the proof of proposition \ref{prop:NEWbd-norm-inter}.] --- We use the notations of § \ref{S:NEW level}. We shall prove by induction on the length $\ell=\ell(w)$ that   for any $w\in W_{M_1}$  there exists $e_w\geq 1$ and for  any level   $J$  as in § \ref{S:NEW level} there exist $\eps, c_J>0$ such that for all  $J$-triples $(P,\pi,\tau)$ such that  $P_\pi\subset P_1\subset P$ and  $w\Delta^{P}_0>0$ and all  $\la\in  \Xi_{P,\eps}$ we have
  \begin{align}\label{eq:ineq B}
  \| N_{\pi_1}(w,\la+\nu_1)\| \leq c_J (1+\la_\tau^2)^{e_w},
  \end{align}
  where $\|\cdot \|$ denotes the norm operator on $\Ac_{P,\pi}(G)^{\tau,J}$ and $\nu_1=\nu_{P_1,\pi}$, see § \ref{S:exposant-disc}.   Using the decomposition \eqref{eq:decomp} of $w$ into elementary symmetries, we write $w=s_{\al_\ell}w'$ and we have
  \begin{align*}
    N_{\pi_1}(w,\la+\nu_1)=    N_{w'\pi_1}(s_{\al_\ell},w'(\la+\nu_1))N_{\pi_1}(w',\la+\nu_1).
  \end{align*}
  By induction, the conclusion holds for the operator $N_{\pi_1}(w',\la+\nu_1)$.   Using the notations of §~\ref{S:elem sym}, we have to distinguish between different cases according to the sign of $\bg \nu_1,\theta_\ell^\vee\bd$. If it is positive then the result follows from lemma \ref{lem:bd finite place} assertion 3 and  lemma \ref{lem:norm inter archi} assertion 1. If $\bg \nu_1,\theta_\ell^\vee\bd=0$ then the  result follows from lemma \ref{lem:bd finite place} assertion 5 and  lemma \ref{lem:norm inter archi} assertion 4. So we may and shall assume $\bg \nu_1,\theta_\ell^\vee\bd<0$. Attached to the level $J$ we have defined  a finite set   $S\subset V_F$  in §~\ref{S:NEW level}. Then, according to   lemma \ref{lem:bd finite place} assertion 4, lemma \ref{lem:norm inter archi} assertion 3 and our induction hypothesis, there exists $\eps>0$ such that the bound we are looking for holds at least for the operator $B(\la+\nu_1)N_{\pi_1}(w,\la+\nu_1)$ on $\Ac_{P,\pi}(G)^{\tau,J}$ for all  $\la\in  \Xi_{P,\eps}$ where we have:
  \begin{itemize}
  \item $ B(\la)=\prod_{v\in S}B_v(\la)$;
  \item  If $v$ is finite, there is an integer $r$ independent of  $\pi$ and $J$ such that $B_v(\la)$ is a finite product of at most $r$  factors of the following shape $1-q_v^{-\bg \la,\theta_\ell^\vee\bd +s}$ with $\Re(s)\leq -2/(n^2+1)$;
  \item If $v$ is  an Archimedean place, there is an integer $r$ independent of  $\pi$ and $J$ such that $B_v(\la)$ is a finite product of at most $r$ factors  of the following shape  $B_{s}(\bg \la,\theta_\ell^\vee\bd )$ where $s\in \CC$ satisfies $-r(1+\la_{\tau_v}^2)\leq \Re(s)\leq -2/(n^2+1)$ and $| \bg \nu_1,\theta_\ell^\vee\bd- \Re(  s)|\leq 1/2(n^2+1)$.
  \end{itemize}

 Let us begin with an Archimedean factor  $B_{s}(\bg \la,\theta_\ell^\vee\bd )$ as above.  Let $\la\in \Xi_{P,\eps}$. If $| \bg \la+\nu_1,\theta_\ell^\vee\bd- s|\geq |\Re(s)| $ then 
  \begin{align} \label{eq:NEW min B par 3} |    B_{s}(\bg \la+\nu_1,\theta_\ell^\vee\bd )^{-1}|&= | 1+   \frac{2\Re(s)}{ \bg \la+\nu_1,\theta_\ell^\vee\bd -s }    |\leq 3.
  \end{align}
  If $| \bg \la+\nu_1,\theta_\ell^\vee\bd- s|\leq |\Re(s)| $ then
    \begin{align*}
      |\bg \la+\nu_1,\theta_\ell^\vee\bd +\bar s|&\leq  3|\Re(s)| \leq 3r(1+\la_{\tau_v}^2).
  \end{align*}

  Hence we have \begin{align*}
    |    B_{s}(\bg \la+\nu_1,\theta_\ell^\vee\bd )|&\geq   \frac{| \bg \la+\nu_1,\theta_\ell^\vee\bd- s|}{3r(1+\la_{\tau_v}^2)}.
  \end{align*}
  Let $B'$ be the function defined as $B$ except that  we replace the factor $B_{s}(\bg \la+\nu_1,\theta_\ell^\vee\bd )$ by the factor $ \bg \la+\nu_1,\theta_\ell^\vee\bd-  s$. The map $\la\mapsto
B'(\la+\nu_1)N_{\pi_1}(w,\la+\nu_1)$ is holomorphic for all $\la\in \Xi_{P,\eps}$ by proposition \ref{prop:le but}. The required bound holds for $B'(\la+\nu_1)N_{\pi_1}(w,\la+\nu_1)$ and  its derivatives (by Cauchy formula) for all $\la\in \Xi_{P,\eps}$ such that $| \bg \la+\nu_1,\theta_\ell^\vee\bd- s|\leq |\Re(s)| $. We can conclude by the mean value formula and \eqref{eq:NEW min B par 3} that we can simply remove the factor $B_{s}(\bg \la+\nu_1,\theta_\ell^\vee\bd )$ from $B$.   In this way, we can now assume that all the Archimedean factors of $B$ are equal to $1$. We are left with non-Archimedean factors. Let $v\in S$ be a finite place. The point is to bound from below the holomorphic function $z\mapsto 1-q^{-z}_v$ on a vertical strip $|\Re(s)|\leq C$. We can divide this strip into a union of rectangles $\rc_t$ indexed by integers $t$ and given by the additional condition $|\Im(z)- t\pi /\log(q_v)|\leq \pi/2\log(q_v)$. If $t$ is odd then $|1-q^{-z}_v|$  is bounded below on $\rc_t$ by a positive constant that does not depend on $t$. If $t$ is even then $|1-q^{-z}_v|$   is bounded below on $\rc_t$ by $|z-t\pi i/\log(q_v)|$ up to a positive constant which does  not depend on $t$. Then we can conclude as before.
  
\end{paragr}

\subsection{Bounds for the scalar product of truncated Eisenstein series}\label{ssec:bound-scalar}

\begin{paragr} The main result is the following:
  
\begin{proposition}\label{prop:bound-scalaire}
There exist  $k,l,r>0$ such that for any level  $J$  there exist $c>0$ and for any large enough $N>0$ there exists a finite family $(X_i)_{i\in I}$ of elements of $\uc(\ggo_\infty)$ such that for all  $J$-triples $(P,\pi,\tau)$ we have
\begin{align*}
  &\bg   \Lambda^{T} E(\varphi,\la),   \overline{E(\psi,\la')}\bd_G  \\
 & \leq  (1+\|\la\|^2+\|\la'\|^2+\la_\tau^2+(\La_{\pi}^M)^2)^k  \exp(r\|T\|) \left(\sum_{i\in I}  \|\varphi\|_{-N,X_i}\right)   \left(\sum_{i\in I}  \|\psi\|_{-N,X_i}\right) 
\end{align*}
for  all $\varphi,\psi \in \Ac_{P,\pi}(G)^{\tau,J}$, $T\in \ago_0^G$ sufficiently positive and all $\la,\la'\in \rc_{\pi,c,l}$.
  \end{proposition}
\end{paragr}

  \begin{paragr}[Proof of proposition \ref{prop:bound-scalaire}.] --- First we observe that the expression $\bg   \Lambda^{T} E(\varphi,\la),   \overline{E(\psi,\la')}\bd_G  $ is holomorphic for $\la,\la'$ in a neighborhood of  $i\ago_{P}^{G,*}$. Let   $(P,\pi,\tau)$ be a $J$-triple. We set
    \begin{align*}
&      L_{\pi}(\la,\la')=\prod_{(R, \al, w,w')}  \bg w\la+w'\la',\al^\vee\bd\\
  &    \tilde L_{\pi}(\la,\la')=\ell_P(\la) \ell_P(\la') L_{\pi}(\la,\la')
    \end{align*}
    where $\ell_P$ is given by proposition \ref{prop:bound-intertwining} and  the product is over the tuples $(R, \al, w,w')$ where:
    \begin{itemize}
    \item     $R$ is a standard parabolic subgroup of $G$ ;
    \item $w,w'\in \, _RW_{P}$ such that $  P_\pi\subset P_w\cap P_{w'} $ and $R_w=R$ and $R_{w'}=R$;
    \item  $\al\in \Delta_R$ such that $\bg \al^\vee,w\nu_{P_w}^P +w'\nu_{P_{w'}}^P\bd=0$.
    \end{itemize}
Note that $\tilde L_\pi$ is a product of non-zero linear forms, see lemma \ref{lem:negativity} and proposition \ref{prop:bound-intertwining}  . Observe also that the set of $\tilde L_\pi$ for different $(P,\pi,\tau)$ is in fact finite. Using theorem  \ref{thm:inversion-calculee} and proposition \ref{prop:bound-intertwining}, we see that there exist  $c,l>0$ such that the map
\begin{align*}
  (\la,\la')\mapsto \tilde L_\pi(\la,\la')\bg   \Lambda^{T} E(\varphi,\la),  \overline{E(\psi,\la')}\bd_G  )
\end{align*}
 is  holomorphic for  $\la,\la'\in \rc_{\pi,c,l}$. However none of the hyperplanes defined by the linear factors of $\tilde L_\pi$ can be singular for $\bg   \Lambda^{T} E(\varphi,\la),\overline{E(\psi,\la')}\bd_G  $ since the Eisenstein series hence the pairing are holomorphic on the imaginary axis. So $\bg   \Lambda^{T} E(\varphi,\la),  \overline{E(\psi,\la')}\bd_G  $ is even holomorphic on  $\rc_{\pi,c,l}$. By a variant of lemma \ref{lem:maj-f-F}, we are reduced to majorize $D( \tilde L_\pi(\la,\la')\bg   \Lambda^{T} E(\varphi,\la),   \overline{E(\psi,\la')}\bd_G)$  on $\rc_{\pi,c,l}\times\rc_{\pi,c,l} $ for some constants $c,l>0$  for a finite set, which does not depend on $\pi$, of holomorphic differential operators $D$. Using Cauchy formula, we are reduced to bound  $\tilde L_\pi(\la,\la')\bg   \Lambda^{T} E(\varphi,\la),   \overline{E(\psi,\la')}\bd_G  $ on $\rc_{\pi,c,l}\times\rc_{\pi,c,l} $ for some constants $c,l>0$. Using theorem \ref{thm:inversion-calculee}, we see that this latter is a finite sum, indexed by $R$ and $w,w'\in \, _RW_{P}$ such that $  P_\pi\subset P_w\cap P_{w'} $ and $R_w=R$ and $R_{w'}=R$ of the product of three factors
    \begin{align*}
      \bg \ell_P(\la)M(w,\la+\nu_{P_w})\varphi_{w} ,  \overline{ \ell_P(\la')M(w',\la'+\nu_{P_{w'}})\psi_{w'}}\bd_{R}
    \end{align*}
   
    \begin{align*}
      \exp(\bg w\la+w'\la'+w\nu_{P_w}+w'\nu_{P_{w'}} ,T_R \bd)
    \end{align*}
    and
    \begin{align*}
      \frac{L_\pi(\la,\la')}{\theta_R(w\la+w'\la'+w\nu_{P_w} +w'\nu_{P_{w'}})}.
    \end{align*}
    The first factor is  bounded by proposition \ref{prop:bound-intertwining} and lemma \ref{lem:comparison-norm}. The second factor, which is the only factor that depends on $T$,  is clearly bounded on $\rc_{\pi,c,l}\times\rc_{\pi,c,l} $ by some power of $\exp(\|T\|)$ . Finally the third factor is a rational function which, as we can assume, has no pole on  $\rc_{\pi,c,l}\times\rc_{\pi,c,l} $. So it is bounded by some power of $1+\|\la\|+\|\la'\|$. The conclusion is then clear.
    
  \end{paragr}

\subsection{Bounds for some Hermitian forms}\label{ssec:bds-herm}

\begin{paragr} We denote by $L^2([G]_{P,0})^{J,\infty}$ the subspace of $L^2([G]_{P,0})$ formed by right-$J$-invariant functions that are smooth and such that  $R(X)\varphi\in   L^2([G]_{P,0})$ for all $X\in  \uc(\ggo_\infty)$.
  
  \begin{lemme} \label{lem:maj-sobolev}Let $J$ be a level and $P$ be a parabolic subgroup of $G$.
    For any large enough $N>0$ and any $X\in \uc(\ggo_\infty)$, there exists $c>0$ and an integer $k\geq 0$ such that for all $\varphi\in L^2([G]_{P,0})^{J,\infty}$ we have
    \begin{align*}
      \|\varphi\|_{-N,X}^2\leq   c\sum_{i=0}^k \| R(\Delta^i)\varphi\|_P^2
    \end{align*}
  \end{lemme}

  \begin{preuve}
    First we can use Sobolev inequalities (see \cite[key lemma]{Ber}) to get that there exist  $c>0$ and $Y_1,\ldots , Y_r\in  \uc(\ggo_\infty)$ such that if $N$ is large enough we have:
    \begin{align*}
            \|\varphi\|_{-N,X}^2  \leq   c    \sum_{i=1}^r \| R(Y_i)\varphi\|_P^2,  \ \ \ \forall \varphi\in L^2([G]_{P,0})^{J,\infty}.
    \end{align*}
    But by \cite[proposition 3.5]{BK}, the topology on $L^2([G]_{P,0})^{J,\infty}$ is also given by the family of semi-norms $(\sum_{i=0}^k \| R(\Delta^i)\varphi\|_P^2)^{1/2})_{k\in \NN}.$ The conclusion is clear.
  \end{preuve}
\end{paragr}

\begin{paragr} Let $\tau_1,\tau_2\in \hat K_\infty$. For any $f\in \Sc(G(\AAA))$  we define 
$$f_{\tau_1,\tau_2}=\bar{e}^\vee_{\tau_1}*f*\bar{e}^\vee_{\tau_2}.$$

Let $J$ be a  level and $T$ be a truncation parameter. For any $J$-pair $(P,\pi)$ and $\la\in \ago_{P,\CC}^{G,*}$ we define a Hermitian form on $S(G(\AAA))^J$ by setting
  \begin{align}\label{eq:formeB}
    B_{ (P,\pi,\tau_1,\tau_2)}^T(\la, f)=\sum_{\varphi\in \bc_{P,\pi}(\tau_2,J)}  \bg   \Lambda^{T} E(I_{P,\pi}(\la,f_{\tau_1,\tau_2})\varphi,\la), E(I_{P,\pi}(\la,f_{\tau_1,\tau_2})\varphi,\la)\bd_G 
  \end{align}
  where $\bc_{P,\pi}(\tau_2,J)$ is an orthonormal basis for the Petersson norm  of the finite dimensional space $\Ac_{P,\pi}(G)^{\tau_2,J}$. It is well-defined if $\la$ is non-singular  for the Eisenstein series  $ E(\varphi,\la)$  for $\varphi \in \Ac_{P,\pi}(G)$. In this case, it does not depend on the choice of the orthonormal basis. Note that $B_{ (P,\pi,\tau_1,\tau_2)}^T=0$ unless both $(P,\pi,\tau_1)$ and $(P,\pi,\tau_2)$ are $J$-triples.

  \begin{remarque}\label{rq:formeB}
    The form $B$ does not depend on the choice of the level $J$, that is it induces a hermitian form on $S(G(\AAA))$. Indeed take an other level $J'\subset J$. We have  $S(G(\AAA))^J\subset S(G(\AAA))^{J'}$  and $\Ac_{P,\pi}(G)^{\tau_2,J}\subset \Ac_{P,\pi}(G)^{\tau_2,J'}$. The projector $p_J$ on $\Ac_{P,\pi}(G)^{\tau_2,J}$ is in fact an orthogonal projector.  Thus we can choose $\bc_{P,\pi}(\tau_2,J')$ to be the union of  an orthonormal basis of $\ker(p_J)$ and    $\bc_{P,\pi}(\tau_2,J)$. Since  the operator $I_{P,\pi}(\la,f_{\tau_1,\tau_2})$ factors through $p_J$ for $f\in S(G(\AAA))^J$ we see that we can replace in \eqref{eq:formeB} the basis $\bc_{P,\pi}(\tau_2,J)$ by $\bc_{P,\pi}(\tau_2,J')$. 
  \end{remarque}

\begin{remarque}\label{rq:positivite}
  Since $\La^T$ is (in some sense) a self-adjoint projector, see \cite[corollary 1.2 and lemma 1.3]{ar2} , we have:
\begin{align*}
  \int_{[G]_0}   | \La^TE(y,I_{P,\pi}(\la,f)\varphi,\la)|^2\, dy=  \bg \La^TE(I_{P,\pi}(\la,f)\varphi,\la),E(I_{P,\pi}(\la,f)\varphi,\la)\bd_G.
\end{align*}
As a consequence $B_{ (P,\pi,\tau_1,\tau_2)}^T(\la, f)$ is real and non-negative.
\end{remarque}

  \begin{proposition} \label{prop:maj-trace-rel} There exist $l,r>0$ such that for all $q>0$ and all levels $J$ there exist $c>0$ and  a continuous semi-norm $\|\cdot\|_{\Sc}$  on $\Sc(G(\AAA))^J$ such that for all $J$-pairs $(P,\pi)$, all $\tau_1,\tau_2\in\hat K_\infty$, all $f\in \Sc(G(\AAA))^J$, all enough positive $T\in \ago_0^G$  and all $\la\in \rc_{\pi,c,l}$ we have
    \begin{align*}
      |B_{ (P,\pi,\tau_1,\tau_2)}^T(f,\la)|\leq \frac{\|f_{\tau_1,\tau_2}\|^2_{\Sc}\exp(r\|T\|)    }{(1+\|\la\|^2)^q(1+\la_{\pi}^2+\la_{\tau_1}^2)^{q} (1+\la_{\pi}^2+\la_{\tau_2}^2)^{q}   }.
    \end{align*}
  \end{proposition}

  \begin{preuve} We may and shall assume that $f=f_{\tau_1,\tau_2}$ and that $(P,\pi,\tau_1)$ and $(P,\pi,\tau_2)$ are $J$-triples.

We start with the following observation. Let $(P,\pi,\tau)$ be a $J$-triple and  $K_{M,\infty}=K_\infty\cap M(F_\infty)$.  There exists  $\sigma\in \hat K_{M,\infty}$ an irreducible representation that appears in the decomposition of the restriction of $\tau$ to $K_{M,\infty}$  such that $\sigma$ is also a $K_{M,\infty}$-type of $\pi_\infty$, in particular $\la_\tau \geq \la_\sigma$, see \cite[proof of lemma 6.1]{Muller02}.  Hence  one can bound $(\La_{\pi}^M)^2$ by an absolute constant times $1+\la_{\pi}^2+\la_\tau^2$.

    Then, using proposition \ref{prop:bound-scalaire},  lemma \ref{lem:maj-sobolev} and the fact that $I_{P,\pi}(\la,f)\varphi\in \Ac_{P,\pi}(G)^{\tau_1,J}$ if $\varphi \in \Ac_{P,\pi}(G)^{J}$,  we get the existence of $k,l>0$ such that for any level $J$ there exist  $c,C>0$ and $N\in \NN$ such that for all  $J$-pairs $(P,\pi)$, all $\tau_1,\tau_2\in\hat K_\infty$,  all $f\in \bar{e}^\vee_{\tau_1}*\Sc(G(\AAA))^J*\bar{e}^\vee_{\tau_2}$ and all $\la\in  \rc_{\pi,c,l}$ we have
\begin{align*}
 |B_{ (P,\pi,\tau_1,\tau_2)}^T(f,\la)|  \leq   \\
C  (1+\|\la\|^2+\la_{\pi}^2+\la_{\tau_1}^2)^k  \exp(r\|T\|)  \sum_{ \varphi \in \bc_{P,\pi}(\tau,J)  }\left(\sum_{i=0}^{N}  \|R(\Delta^i)I_{P,\pi}(\la,f)\varphi\|_P^2\right).
\end{align*}

Before going further, we collect some general facts for a $J$-triple $(P,\pi,\tau)$.

First the operator $I_{P,\pi}(\la,\Delta)$ acts on  the subspace $\Ac_{P,\pi}(G)^{\tau,J}$ by the scalar $1-(\la,\la)-\la_{\pi}+2\la_\tau$, see \cite[eq. (6.7)]{Muller02} and \cite[proof of lemma 5.4] {Muller07} where $(\cdot,\cdot)$ is the standard quadratic form on $\ago_{0,\CC}^*$ (whose restriction to $\ago_{0}^*$ is the standard scalar product). Moreover, we have $\la_\tau\geq \la_{\pi}$, see \cite[lemma 6.1]{Muller02}, and $\la_\tau\geq 0$. So for $\varphi \in \Ac_{P,\pi}(G)^{\tau,J}$ we have:
\begin{align*}
    \|R(\Delta)\varphi\|_P=(1-\la_{\pi}+2\la_\tau) \|\varphi\|_P.
\end{align*}
There exists $C_0$ (for example one can take $C_0=3$) such that for all $J$-triples $(P,\pi,\tau)$ and all $\la\in\ago_{P,\CC}^{G,*}$ we have
$$(1-\la_{\pi}+2\la_\tau)\leq C_0 (1+\|\la\| ^2+\la_{\pi}^2+\la_\tau^2).$$
Hence there exists $C_0'$ such that for all $\varphi \in \Ac_{P,\pi}(G)^{\tau,J}$ 
\begin{align}\label{eq:point1}
  \sum_{i=0}^{N}  \|R(\Delta^i)\varphi\|_P^2 
  \leq C_0' (1+\|\la\| ^2+\la_{\pi}^2+\la_\tau^2)^{2N}  \|\varphi\|_P^2 .
\end{align}

On the other hand,  there exists $C_1>0$ such that  for any $\varphi\in \Ac_{P,\pi}(G)^{\tau,J}$, any $i\in \NN$ and $\la\in \rc_{\pi,c,l}$ (one can shrink $c$ if necessary)
  \begin{align}\label{eq:point2}
    (1+\|\la\|^2+\la_{\pi}^2+\la_\tau^2)^i \|\varphi\|_P \leq C_1^i \| I_{P,\pi}(\la,\Delta^{2i})\varphi\|_P.  
  \end{align}

Second, there exists $C_2>0$ and $k_1\in \NN$, see \cite[eq. (6.14)]{Muller02} such that
  \begin{align}\label{eq:point3}
    \dim(\Ac_{P,\pi}(G)^{\tau,J})\leq  C_2(1+\la_\tau^2+\la_{\pi}^2)^{k_1}\leq C_1(1+\|\la\|^2+\la_\tau^2+\la_{\pi}^2)^{k_1}
  \end{align}

 Using \eqref{eq:point1} and \eqref{eq:point3}, we get that there exists $C_3>0$ such that

\begin{align*}
|B_{ (P,\pi,\tau_1,\tau_2)}^T(f,\la)|  \leq   C_3(1+\|\la\|^2+\la_{\pi}^2+\la_{\tau_1}^2)^{k+2N}(1+\|\la\|^2+\la_{\pi}^2+\la_{\tau_2}^2)^{k_1} \times\\
  \exp(r\|T\|)\frac{\sum_{ \varphi \in \bc_{P,\pi}(\tau_2,J)  } \|I_{P,\pi}(\la,f)\varphi\|_P^2}{| \bc_{P,\pi}(\tau_2,J)   |  }.
\end{align*}  

Using \eqref{eq:point2}, for any $q>0$ we have: 
\begin{align*}
|B_{ (P,\pi,\tau_1,\tau_2)}^T(f,\la)|  \leq   C_3(1+\|\la\|^2+\la_{\pi}^2+\la_{\tau_1}^2)^{k+2N-2q}(1+\|\la\|^2+\la_{\pi}^2+\la_{\tau_2}^2)^{k_1-2q} \times\\
  \exp(r\|T\|)\frac{\sum_{ \varphi \in \bc_{P,\pi}(\tau_2,J)  } \|I_{P,\pi}(\la,L(\Delta^{2q}) R(\Delta^{2q})f)\varphi\|_P^2}{| \bc_{P,\pi}(\tau_2,J)   |  }.
\end{align*}  
Note that $I_{P,\pi}(\la,\Delta)I_{P,\pi}(\la,f_\tau)\varphi=I_{P,\pi}(\la,L(\Delta)f_\tau)\varphi$.  

It is easy to conclude since, for every $i\in \NN$ there is  a continuous semi-norm $\|\cdot\|_i$ on $\Sc(G(\AAA))^J$ such that for all pairs $(P,\pi)$ and all $\varphi\in \Ac_{P,\pi}(G)^{J}$, all $\la\in \rc_{\pi,c,l}$ we have
  $$\|I_{P,\pi}(\la,L(\Delta^{2i})R(\Delta^{2i})f)\varphi\|_P\leq      \|\varphi\|_P \|f\|_{i}.$$ 
\end{preuve}
\end{paragr}

\begin{paragr} Let $J$ be a level and let $(P,\pi)$ be a $J$-pair. We define for $f\in \Sc(G(\AAA))^J$
  \begin{align*}
    B_{ (P,\pi)}^T(\la, f)=\sum_{\tau_1,\tau_2\in  \hat K_\infty} B_{ (P,\pi,\tau_1,\tau_2)}^T(\la, f).
  \end{align*}
Recall that the terms in the sum above  are non-negative, see remark \ref{rq:positivite}. The sum is  convergent by the next proposition.

  \begin{proposition}
    \label{prop:maj-trace-rel2} There exist $l,r,q_0>0$ such that for all $q>q_0$ and all levels $J$ there exist $c>0$ and  a continuous semi-norm $\|\cdot\|_{\Sc}$  on $\Sc(G(\AAA))^J$ such that for  all $f\in \Sc(G(\AAA))^J$, all enough positive $T\in \ago_0^G$  and all $\la\in \rc_{\pi,c,l}$ we have
    \begin{align*}
      B_{ (P,\pi)}^T(f,\la)\leq \frac{\|f\|_{\Sc}^2\exp(r\|T\|) }{(1+\|\la\|^2)^q(1+\La_{\pi}^2)^{q}}.
    \end{align*}
  \end{proposition}

  \begin{preuve}
    Using proposition \eqref{prop:maj-trace-rel} and its notations, we have 
    \begin{align*}
\sum_{\tau_1,\tau_2\in  \hat K_\infty} B_{ (P,\pi,\tau_1,\tau_2)}^T(\la, f)\leq   \sum_{\tau_1,\tau_2 } \frac{\|f_{\tau_1,\tau_2}\|_{\Sc}^2\exp(r\|T\|)    }{(1+\|\la\|^2)^q(1+\la_{\pi}^2+\la_{\tau_1}^2)^{q}(1+\la_{\pi}^2+\la_{\tau_2}^2)^{q}}.
    \end{align*}
Note that in the right-hand side the sum is over $\tau_1,\tau_2\in \hat K_\infty$ such that $(P,\pi,\tau_1)$ and $(P,\pi,\tau_2)$ are $J$-triples. Using Cauchy-Schwartz inequality, we can bound the  right-hand side by
\begin{align*}
  \frac{\exp(r\|T\|)    }{(1+\|\la\|^2)^q }   (\sum_{\tau_1,\tau_2 \in \hat K_\infty}   \|f_{\tau_1,\tau_2}\|_{\Sc}^4)^{1/2}   \times \sum_{  \tau }(1+\la_\tau^2+\la_{\pi}^2)^{-2q}
\end{align*}
where the last sum is over $\tau\in \hat K_\infty$ such that $(P,\pi,\tau)$ is a $J$-triple.
We observe that $(\sum_{\tau_1,\tau_2 \in \hat K_\infty}   \|f_{\tau_1,\tau_2}\|_{\Sc}^4)^{1/4} $ is a continuous semi-norm on  $\Sc(G(\AAA))^J$. Moreover, there exists $C>0$ such that for all $P$ and $\pi\in \Pi_{\disc}(M)$
\begin{align*}
  \sum_{  \tau}(1+\la_\tau^2+\la_{\pi}^2)^{-2q}\leq  C  (1+\La_{\pi}^2)^{-q} ( \sum_{  \tau}(1+\la_\tau^2)^{-q})
\end{align*}
and $\sum_{  \tau\in \hat K_\infty}(1+\la_\tau^2)^{-q}<\infty$ if $q$ is large enough.
  \end{preuve}
\end{paragr}

\subsection{Bounds for Eisenstein series}\label{ssec:bds-Eis}

\begin{paragr} We start with a lemma.

  \begin{lemme}\label{lem:maj-Eis}
    Let $m\geq 0$ and $g\in  C_c^m(G(\AAA))$. There exists $l>0$ and $N>0$ such that for all $q>0$ and all levels $J$ there exist $c>0$ and  a continuous semi-norm $\|\cdot\|_{\Sc}$  on $\Sc(G(\AAA))^J$ such that for all $J$-pairs $(P,\pi)$, all $f\in \Sc(G(\AAA))^J$, all $\la\in \rc_{\pi,c,l}$ and  all $x\in G(\AAA)$ we have
    \begin{align*}
     \sum_{\varphi\in \bc_{P,\pi}(J)}  | E(x,I_{P,\pi}(\la,g*f)\varphi,\la) |^{2}  \leq \frac{  \|x\|_{G}^N\|f\|_{\Sc}^2    }{(1+\|\la\|^2)^q(1+\La_{\pi}^2)^{q}}
    \end{align*}
   where $\bc_{P,\pi}(J)$ is the union over $\tau\in \hat K_\infty$ of orthonormal bases of $\Ac_{P,\pi}(G)^{\tau,J}$.
  \end{lemme}

  \begin{preuve} Since the Eisenstein series is $A_G^\infty$-invariant, it suffices to prove the bound for $x\in G(\AAA)^1$. The main point is to express the square modulus of the Eisenstein series in terms of the truncated inner product and then to apply proposition \ref{prop:maj-trace-rel2}. To do this we follow \cite[beginning of the proof of proposition 6.1]{LapFRTF}.

    Let $g\in  C_c^m(G(\AAA))$. We may and shall assume that $x$ is in a fixed Siegel set of $G(\AAA)^1$.  According to a slight variant of \cite[lemma 6.2]{LapFRTF}, there exists an absolute constant $c_0>0$ such that we have
    \begin{align*}
  (\La^TE)(xy,I_{P,\pi}(\la,f)\varphi,\la)=E(xy,I_{P,\pi}(\la,f)\varphi,\la)
\end{align*}
for all truncation parameters $T$ and all  $y$ in the (compact) support  $\supp(g)$ of $g$  such that
\begin{align}\label{eq:condition-sur-T}
  \bg \varpi,T-H_0(xy)\bd >c_0 \text{ for  all  } \varpi\in \hat\Delta_0.
\end{align}
We have
\begin{align*}
  H_0(xy)=H_0(x)+H_0(k(x)y)
\end{align*}
where $k(x)\in K$ is such that $xk(x)^{-1}\in P_0(\AAA)$ (Iwasawa decomposition). In particular, there exists $c_1$ depending on  $\supp(g)$ of $g$ such that if   $T_1$ is a truncation parameter such that $ \bg \varpi,T_1\bd >c_1$ for  all $\varpi\in \hat\Delta_0$ then $T=T_1+H_0(x)$ is enough positive and satisfies  \eqref{eq:condition-sur-T} for all $y\in \supp(g)$. We fix  such a $T_1$ and we set $T=T_1+H_0(x)$. Then we get, see  \cite[p. 284]{LapFRTF}:
    \begin{align*}
      E(x,I_{P,\pi}(\la,g*f)\varphi,\la)=\int_{[G]_0} k_{g}(x,y)   (\La^TE)(y,I_{P,\pi}(\la,f)\varphi,\la)\, dy.
    \end{align*}
    where  we set
$$    k_{g}(x,y) =\sum_{\gamma\in G(F)} \int_{A_G^\infty} g(x^{-1} a\gamma y)\,da.$$
According to \cite[lemme I.2.4]{MWlivre} there exist $N>0$ and $c_2>0$ such that
\begin{align}\label{eq:majMW}
  |k_g(x,y)|\leq c_2 \|x\|_G^{N}, \ \ \ x\in G(\AAA)^1.
\end{align}
Hence we have by Cauchy-Schwartz inequality
\begin{align*}
  |E(x,I_{P,\pi}(\la,g*f)\varphi,\la)|^2
  &\leq \vol([G]_0) c_2^2\|x\|_G^{2N} \int_{[G]_0}   | (\La^TE)(y,I_{P,\pi}(\la,f)\varphi,\la)|^2\, dy.
\end{align*}
We deduce that (see remark \ref{rq:positivite}):

\begin{align*}
   \sum_{\varphi\in \bc_{P,\pi}(J)}  | E(x,I_{P,\pi}(\la,g*f)\varphi,\la) |^{2} \leq \vol([G]_0) c_2^2\|x\|_G^{2N} B^T_{P,\pi}(f,\la)
\end{align*}
Using the result and the notations of proposition \ref{prop:maj-trace-rel2}, we get the following bound:
\begin{align*}
   \sum_{\varphi\in \bc_{P,\pi}(J)}  | E(x,I_{P,\pi}(\la,g*f)\varphi,\la) |^{2} \leq \vol([G]_0) c_2^2\|x\|_G^{2N} \frac{\|f\|_{\Sc}^2\exp(r\|T\|) }{(1+\|\la\|^2)^q(1+\La_{\pi}^2)^{q}}.
\end{align*}
To conclude it suffices to observe that there exists $c_3$ and $N'$ such that $\exp(\|T\|)\leq \exp(\|T_1\| +\|H_0(x)\|)\leq c_3 \|x\|_G^{N'}$, see \cite[I.2.2]{MWlivre}.
  \end{preuve}

\end{paragr}

\begin{paragr} The following theorem is an extension of \cite[proposition 6.1]{LapFRTF} to the case of discrete Eisenstein series.  

  \begin{theoreme} \label{thm:maj-Eis}
   There exists $l>0$ and $N>0$ such that for all $q>0$ and all levels $J$ there exist $c>0$ and  a continuous semi-norm $\|\cdot\|_{\Sc}$  on $\Sc(G(\AAA))^J$ such that for all $J$-pairs $(P,\pi)$, all $f\in \Sc(G(\AAA))^J$, all $\la\in \rc_{\pi,c,l}$ and  all $x\in G(\AAA)$ we have
    \begin{align*}
     \sum_{\varphi\in \bc_{P,\pi}(J)}  | E(x,I_{P,\pi}(\la,f)\varphi,\la) |^{2}  \leq \frac{  \|x\|_{G}^N\|f\|_{\Sc}^2    }{(1+\|\la\|^2)^q(1+\La_{\pi}^2)^{q}}
    \end{align*}
   where $\bc_{P,\pi}(J)$ is the union over $\tau\in \hat K_\infty$ of orthonormal bases of $\Ac_{P,\pi}(G)^{\tau,J}$.
  \end{theoreme}

  \begin{preuve}    
    Following \cite[corollary 4.2]{ar1}, for a level $J$ and an integer  $m\geq 1$ large enough, we can find $Z\in \uc(\ggo_\infty)$, $g_1\in \Cc(G(\AAA))$ and $g_2\in C_c^m(G(\AAA))$ such that
  \begin{itemize}
  \item $Z$ is  invariant under $K_\infty$-conjugation ;
  \item $g_1$ and $g_2$ are  invariant under $K$-conjugation and $J$-biinvariant;
  \item  for any $f\in \Sc(G(\AAA))^{J}$  we have:
    \begin{align*}
      f=g_1*f+ g_2*(Z*f).
    \end{align*}
  \end{itemize}
Then the theorem is a straightforward consequence of lemma  \ref{lem:maj-Eis}.
  \end{preuve}

\end{paragr}

\begin{paragr}
  In the sequel we shall need the following extension of theorem \ref{thm:maj-Eis}.

  \begin{theoreme} \label{thm:maj-EisR}
    There exist $l,N>0$ and for each standard parabolic subgroup $P$ of $G$ there exists $\ell_P$, a product of non-trivial real linear forms on $\ago_{P}^{G,*}$,       such that for all $q>0$ and all levels $J$ there exist $c>0$ and  a continuous semi-norm $\|\cdot\|_{\Sc}$  on $\Sc(G(\AAA))^J$ such that for all standard parabolic subgroups $R$, for all $J$-pairs $(P,\pi)$, all $w\in \,_RW_P$, all $f\in \Sc(G(\AAA))^J$, all $\la\in \rc_{\pi,c,l}$ and  all $x\in G(\AAA)$ we have
    \begin{align*}
     \sum_{\varphi\in \bc_{P,\pi}(J)}  | \ell_P(\la)E^R(x,M(w,\la+\nu_w)(I_{P,\pi}(\la,f)\varphi)_w, (w(\la+\nu_w))^R) |^{2}  \leq \frac{  \|x\|_{R}^N\|f\|_{\Sc}^2    }{(1+\|\la\|^2)^q(1+\La_{\pi}^2)^{q}}
    \end{align*}
   where $\bc_{P,\pi}(J)$ is the union over $\tau\in \hat K_\infty$ of orthonormal bases of $\Ac_{P,\pi}(G)^{\tau,J}$. Moreover we can take $\ell_P=1$ if we restrict the statement to the elements $w\in W(P;R)$.
  \end{theoreme}

  \begin{preuve} It is straightforward to analyze the behaviour of the Eisenstein  series on $A_R^\infty$. In particular, we see that it suffices to prove the theorem for $x\in G(\AAA)$ such that $ H_R(x)=0$. Then the proof is very close to that of theorem \ref{thm:maj-Eis}, so we shall be brief.  For the discussion, we set $\psi=M(w,\la+\nu_w)\varphi_w$ and $\mu=w(\la+\nu_w)$. The starting point is that we have for $x$ in some Siegel set, $g\in   C_c^m(G(\AAA))$ and a suitable $T$ depending on $x$ 
\begin{align*}
      E^R(x,I_{R_w}(\mu,g*f)\psi,\mu^{R})=\int_{[G]_{R,0}} k_{R,g}(x,y)   \La^{T,R}E^R(y,I_{R_w}(\mu,f)\psi,\mu^R)\, dy
\end{align*}
    where  we set
$$    k_{R,g}(x,y) = \int_{[N_R]} \int_{A_R^\infty} \sum_{\gamma\in R(F)}g(x^{-1} a\gamma n y)\,da.$$
There exist $N>0$ and $c_2>0$ such that for all $x\in G(\AAA)$  such that $ H_R(x)=0$ we have
\begin{align*}
  |k_{R,g}(x,y)|\leq c_2 \|x\|_R^{N}.
\end{align*}
In this way, we are reduced to bound $\| \La^{T,R}E^R(I_{R_w}(\mu,f)\psi,\mu)\|_R$. If $w\in W(P;R)$ we have $\nu_w=0$ and the Eisenstein series $E^R(\psi,\mu)$ is holomorphic for $\la\in i\ago_{P}^{G,*}$. It is then easy to get a bound by a variant of proposition \ref{prop:bound-scalaire}   and the rest of the proof is nearly identical  to that of  theorem \ref{thm:maj-Eis}. In general, we have to introduce the polynomial  $\ell_P$ which appears in proposition \ref{prop:bound-intertwining}. We claim that the product $\ell_P(\la)E^R(\psi,\mu)$, hence also  $\ell_P(\la)\La^{T,R}E^R(\psi,\mu)$, is holomorphic for $\la\in i\ago_{P}^{G,*}$. To see this we first compute the constant term of $E^R(\psi,\mu)$  along a parabolic subgroup $Q\subset R$: it is given by
\begin{align*}
  \sum_{w'\in \, _QW^R_{R_w}w}  E^Q(M(w',\la+\nu_{w'})\varphi_{w'},w'(\la+\nu_{w'})).
\end{align*}
In particular, the cuspidal component along $Q$ is the sum 
\begin{align*}
  \sum_{w'}  (M(w',\la+\nu_{w'})\varphi_{w'})_{w'(\la+\nu_{w'})}
\end{align*}
over the set of elements $w'\in \, _QW^R_{R_w}w$ such that $Q_{w'}=Q$.
Using proposition \ref{prop:bound-intertwining} and remark \ref{rq:holom op entr} we infer that all the cuspidal components of $\ell_P(\la)E^R(\psi,\mu)$ are holomorphic  for  $\la\in i\ago_{P}^{G,*}$ and thus $\ell_P(\la)E^R(\psi,\mu)$ itself is holomorphic, see \cite[lemme I.4.10]{MWlivre}. So $\ell_P(\la)\La^{T,R}E^R(\psi,\mu)$ is also holomorphic for  $\la\in i\ago_{P}^{G,*}$. Then we can get easily variants of  proposition \ref{prop:bound-scalaire} and theorem \ref{thm:maj-Eis}.

  \end{preuve}
\end{paragr}

\section{Flicker-Rallis periods}\label{sec:FR}

\subsection{Notations}\label{ssec:nota-FR}

\begin{paragr}\label{S:notations}
  From now on   $E/F$ is a quadratic extension of number fields. Sometimes we will consider $\tau\in F^\times$ such that $E=F(\sqrt{\tau})$. 

  Let  $n\geq 1$ be an integer. Let  $G_n=\Res_{E/F} \GL_E(n)$ be the $F$-group obtained by restriction of scalars. Let  $\iota$ be the Galois involution of $G$ whose fixed point set is the subgroup $G'_n=\GL_F(n)$. The inclusion $G'_n\subset G_n$ gives an inclusion $A_{G'_n}\subset A_{G_n}$ which is in fact an equality. The  restriction map $X^*(G_n)\to X^*(G'_n)$ gives an isomorphism $\ago_{G_n}^*\simeq \ago_{G'_n}^*$.

  The minimal pair  $(P_{0,n}',M_{0,n}')$ for $G_n'$ is formed by  the Borel subgroup $P_{0,n}'$  of $G_n'$ of  upper triangular matrices and the diagonal maximal torus  $M_{0,n}'$ of $ G' _n$. Let $(P_{0,n},M_{0,n})$ be the minimal pair for $G_n$ deduced from $(P_{0,n}',M_{0,n}')$ by extension of scalars to $E$ and restriction to $F$. The words ``standard'' and ``semi-standard'' will refer to these  pairs. The map $ P '\mapsto P = \Res_{E / F} (P' \times_F E) $ induces a bijection between the sets of standard parabolic subgroups of $ G '_n$ and  $ G _n$ whose inverse bijection is given by
  $$ P \mapsto P '= P \cap G'_n. $$

Let $ P $ be a standard parabolic subgroup of $ G_n $. The restriction map $X^*(P)\to X^*(P')$ identifies $X^*(P)$ with a subgroup of $X^*(P')$ of index $2^{\dim(\ago_P)}$.  It induces an isomorphism $ \ago_{P'} \to   \ago_ {P }$ which does not preserve the Haar measures: the pullback to $\ago_{P'}$ of the Haar measure on $\ago_P$ is $2^{\dim(\ago_P)}$ times the  Haar measure on $\ago_{P'}$. In the same way,   the groups $A_P^\infty$ and $A_{P'}^\infty$ are  canonically identified but  the Haar  measure on $A_P^\infty$ is $2^{\dim(\ago_P)}$ times the Haar measure on $A_{P'}^\infty$.
For any standard parabolic subgroup $P\subset Q$, the restriction of the function $ \tau_P ^ Q $ to $\ago_{P'}$ coincides with the function $ \tau_ {P '} ^ {Q'} $.
\end{paragr}

\begin{paragr} Let $\AAA$ be the ring of adèles of $F$.  The groups $G_n(\AAA)$ and $G_n'(\AAA)$ come with their standard maximal compact subgroups respectively denoted by $K_n$ and $K_n'$. We have $ K '_n= K_n  \cap G'_n (\AAA)$.  Note that for all $x\in G'(\AAA)$
  \begin{align}
    \label{eq:rhoQ rhoQ'}
 \bg  \rho_P ^ Q, H_P (x) \bd = 2 \bg  \rho_ {P '} ^ {Q'}, H_ {P '} (x) \bd. 
  \end{align}
  
\end{paragr}

\begin{paragr}
  In most of the rest of the paper, the integer $n$ is fixed and will be omitted in the notation ($G=G_n$, $P_0=P_{0,n}$ etc.). As before, we identify $\ago_{0,\CC}$ and its dual with $\CC^n$ equipped with the usual non-degenerate positive definite Hermitian  form. We denote by  $\|\cdot\|$ the associated norm. The other notations are borrowed from the previous sections.
\end{paragr}

\begin{paragr}
  In this section we denote by $T$  a truncation parameter in $\ago_0$.
\end{paragr}

\subsection{Mixed truncation operator}\label{ssec:mixed}

\begin{paragr} Let $Q$ be a standard parabolic subgroup of $G$. Following Jacquet-Lapid-Rogawski, see \cite{JLR}, we define the (mixed) truncation operator $\Lambda^{T,Q}_m$  that associates to a function $\varphi$ on $[G]$ the following function of the variable $h  \in [G']_{Q'}$:
\begin{align}\label{eq:LaTm}
    (\Lambda^{T,Q}_m\varphi)(h)= \sum_{ P_0\subset P\subset Q} (-1)^{\dim(\ago_P^Q)}\sum_{\delta \in P'(F)\back Q'(F) }   \hat\tau_P^Q(H_P(\delta h)-T_P) \varphi_P(\delta h)
  \end{align}
where  $\varphi_P$ is the constant term along $P$. If $Q=G$, the exponent $G$ is omitted. We denote by $\Lambda^{T,M_Q}_m$ the operator on the space of functions on $[M_Q]$ given by the  formula \ref{eq:LaTm} where the parabolic subgroup $P$ is now interpreted as a standard parabolic subgroup of $M_Q$. By definition $\Lambda^{T,M_Q}_m$ and $\Lambda^{T,Q}_m$  depend only on the projection $T^Q$ of $T$ on $\ago_0^Q$.

One of the most important property of the truncation operator is given by  the following proposition whose proof is a variant of that of \cite[proposition 8]{JLR} and is omitted.

\begin{proposition}\label{prop:cont-LaTm}
  Let $J\subset K_f$ be a level and let $Q$ be a standard parabolic subgroup. For any $N,N'>0$ there exists a finite family $(X_i)_{i\in I}$ of elements of $\uc(\ggo_\infty)$ such that for any smooth and right-$J$-invariant function $\varphi$  on $[G]_Q$,  the function $\Lambda^{T,Q}_m\varphi$ is a  right-$J$-invariant  function on $[G']_{Q'}$, and we have
  \begin{align*}
    \sup_{g\in [G']_{Q'}^1}    \|g\|^{N'}_Q| \Lambda^{T,Q}_m\varphi(g)| \leq \sum_{i\in I } \|\varphi\|_{-N,X_i}.
  \end{align*}

\end{proposition}
\end{paragr}

\subsection{Regularized periods of discrete Eisenstein series}\label{ssec:reg-period-Eis}

\begin{paragr}
  Let $P=MN_P$ be a standard parabolic subgroup of $G$. Let $\pi\in \Pi_{\disc}(M)$.
\end{paragr}

\begin{paragr} \label{S:ITQ}  Let $\varphi\in \Ac_{P,\pi}(G)$.  Let $Q$ be a standard parabolic subgroup of $G$ and $w\in \, _QW_{P}$. For  $\la\in \ago_{P_w,\CC}^{G,*}$, we set
\begin{align*}
 \Ic^{T,Q}(\varphi,\la,w)= \int_{[G']_{Q',0}}  \Lambda_m^{T,Q} E^Q(g,  M(w,\la)\varphi_{w},(w\la)^Q)\, dg,
 \end{align*}
 where we set, see § \ref{S:exposant-disc}, $\varphi_{w}=    \varphi_{P_w,-\nu_{P_w}}$. This is what we call   the truncated Flicker-Rallis period of the Eisenstein series  $E^Q(g,  M(w,\la)\varphi_{w},(w\la)^Q)$.  Note that the integrand is left-equivariant  under the character
 $$x \in A_{Q'}^\infty M_{Q'}(F)N_{Q'}(\AAA) \mapsto \exp(\bg \rho_Q,H_Q(x)\bd)= \exp(\bg 2\rho_{Q'},H_{Q'}(x)\bd).$$
 So the ``integral'' makes formally sense. It is in fact convergent as soon as the Eisenstein series and the intertwining operator are well-defined: this follows from the moderate growth of Eisenstein series and proposition \ref{prop:cont-LaTm}. Note also that $\Ic^{T,Q}(\varphi,\la,w)$  depends only on the projection $T^Q$ of $T$ and it vanishes unless $P_\pi\subset P_w$. In particular, if $Q=P_0$ then the integral does not depend on $T$.  If $Q=G$, then  $w$ must be $1$ and  we get:

\begin{align*}
 \Ic^{T}(\varphi,\la)= \int_{[G']_0}  \Lambda_m^{T}E (g,  \varphi,\la)\, dg
 \end{align*}
where we omit $w$ and $G$ from the notation.

\begin{lemme}\label{lem:ITQ-analytic}
Let   $\om\subset \ago_{P_w,\CC}^{G,*}$ a compact subset. Let $r(\la)$  be a finite product of non-zero affine functions on $\ago_{P_w,\CC}^{G,*}$ such that $r(\la) E^Q(M(w,\la)\varphi_{w},(w\la)^Q)$ is holomorphic on $\om$.
\begin{enumerate}
\item For any $\varphi\in \Ac_{P,\pi}(G)$, the product $r(\la)  \Ic^{T,Q}(\varphi,\la,w)$ is also  holomorphic on $\om$.
\item There exists a continuous semi-norm $\|\cdot\|$ on $\Ac_{P,\pi}(G)$ such that 
\begin{align*}
  |r(\la)\Ic^{T,Q}(\varphi,\la,w)| \leq \|\varphi\|, \ \ \ \varphi\in \Ac_{P,\pi}(G), \ \la\in \om.
\end{align*}
\end{enumerate}
\end{lemme}

\begin{preuve}
  This is a straightforward consequence of the continuity of Eisenstein series of \cite[theorem 2.2]{Lap-remark} and the property of the truncation operator given in proposition \ref{prop:cont-LaTm}.
\end{preuve}
\end{paragr}

 \begin{paragr}[Regularized periods of  Eisenstein series.] --- \label{S:reg-period}   Let $R$ be a standard parabolic subgroup of $G$. Let $w\in \,_RW_{P}$. For any $w'\in  \,_QW^R_{R_w}  w $ we have $P_{w'}\subset P_w\subset  P$. If $P_\pi\subset P_{w'}$ we set
 \begin{align*}
\nu_{w'}^w=   \nu_{P_{w'}}^{P_w}.
 \end{align*}

For any $\la\in \ago_{P_w,\CC}^{G,*}$, we introduce the Jacquet-Lapid-Rogawski regularized period
\begin{align}\label{eq:PTRc}
  \pc^{T,R}(\varphi,\la,w)&=\sum_{ P_0\subset Q\subset R} (-2)^{-\dim(\ago_Q^R)}  \sum_{w ' \in \,_QW^R_{R_w}  w } \Ic^{T,Q}(\varphi,\la+\nu_{{w'}}^w,w')\cdot \frac{\exp(\bg w'(\la+\nu_{{w'}}^w),T_Q^R\bd)}{\hat{\theta}_Q^R(w'(\la+\nu_{{w'}}^w))}.
\end{align}
If $w=1$ we set  $\pc^{T,R}(\varphi,\la)= \pc^{T,R}(\varphi,\la,1)$.

\begin{remarque}\label{rq:interpretation}
  In the sum above, we tacitly  assume that a summand indexed by $w'$ such that $P_\pi\not\subset P_{w'}$ is understood to be $0$ whether the numerator $\hat{\theta}_Q^R(w'(\la+\nu_{{w'}}^w))$ vanishes or not. Recall that in this case we have $\Ic^{T,Q}(\varphi,\la,w')=0$. In particular, we have  $\pc^{T,R}(\varphi,\la,w)=0$ unless $P_\pi\subset P_w$.
\end{remarque}

\begin{remarque}
  \label{rq:descente-par}
  We can replace $P$ by $R_w$ and $\varphi$ by $M(w,\la)\varphi_{P_w}$. We get:
  \begin{align}\label{eq:chgt-w}
     \pc^{T,R}(\varphi,\la,w)&= \pc^{T,R}(M(w,\la)\varphi_{P_w},w\la,1).
  \end{align}

We can introduce analogous objects for the Levi factor $M_R$ in place of $G$ or $R$.  Then we have a parabolic descent (written for $w=1$)
   \begin{align}\label{eq:par-desc}
     \pc^{T,R}(\varphi,\la)&= \pc^{T,M_R}(\psi,\la^R).
   \end{align}
   where $\psi(m)=\int_{K'}    \varphi_{-\rho_R}(mk)dk$.
\end{remarque}

\begin{lemme}\label{lem:analyticite-m+c}
  The expression \eqref{eq:PTRc} is well-defined and holomorphic for $\la$ in a complement of hyperplanes  in $\ago_{P_w,\CC}^{G,*}$. 
\end{lemme}

\begin{remarque}\label{rq:analyticite-m+c}
  We will in fact later show that the statement is still  true if we replace $\ago_{P_w,\CC}^{G,*}$ by $\ago_{P,\CC}^{G,*}$, see remark \ref{rq:MS-periods}.
\end{remarque}

\begin{preuve}
  Let  $P_0\subset Q\subset R$ and $w'\in \,_QW^R_{P_w} w$.    Using lemma \ref{lem:ITQ-analytic} and remark \ref{rq:interpretation}, we are reduced to prove the statement for the rational maps $\la\mapsto \hat{\theta}_Q^R(w'(\la+\nu_{{w'}}^w))^{-1}$ where $w'$ is such that $P_\pi\subset P_{w'}$. Let us assume  that  $\hat{\theta}_Q^R$ vanishes identically on $w' \ago_{P_w,\CC}^{G,*}$ for such a $w'$.  More precisely let us assume that there is $\varpi^\vee\in \hat\Delta_Q^{R,\vee}$ such that $\varpi^\vee \in w'\ago_{0}^{P_w}$.  But then by  lemma \ref{lem:positivite-nu}  below we have $\bg w'\nu_{w'}^w, \varpi^\vee\bd\not= 0$. In this way, we see that $\la\mapsto \hat{\theta}_Q^R(w'(\la+\nu_{{w'}}^w))$ does not vanish identically on  $w' \ago_{P_w,\CC}^{G,*}$.
 \end{preuve}

 \begin{lemme} \label{lem:positivite-nu} Let  $P_0\subset Q\subset R$ and $w ' \in \,_QW^R_{R_w}  w$ be such that $P_\pi\subset P_{w'}$. 
   \begin{enumerate}
   \item For all $\varpi^\vee\in \hat\Delta_Q^{R,\vee}$, we have 
$$\bg w'\nu_{w'}^w, \varpi^\vee\bd\leq 0.$$
Moreover we have equality if and only if $\varpi^\vee\in w' \ago_{P_w}^G$.
\item The following conditions are equivalent:
  \begin{enumerate}
  \item $P_w=P_{w'}$
  \item  $(w'\nu_{w'}^w)_Q^R=0$;
  \item $w'\in W^R(R_w;Q)w$.
  \end{enumerate}
\item There exists $c_1>0$ such that for any $T\in \ago_0^+$, any   $P_0\subset Q\subset R$ and any $$w'\in \,_QW^R_{R_w}  w \setminus W^R(R_w;Q)w$$ such that $P_\pi\subset P_{w'}$ we have:
  \begin{align*}
       \bg (w' \nu_{{w'}}^w)^R_Q,T\bd \leq -c_1d(T).
  \end{align*}
\end{enumerate}
\end{lemme}

 \begin{preuve} The first two assertions are essentially in \cite[lemma 6]{Linner}. Since our setting is slightly different we give a proof for the reader's convenience.    

   1. Let $\varpi^\vee\in \hat\Delta_Q^{R,\vee}$. We have $\nu_{w'}^w=\sum_{\al \in \Delta_{w'}^w} c_\al \al$ where $\Delta_{w'}^w=\Delta_{P_{w'}}^{P_w}$ and $c_\al<0$. Observe that $w'\Delta_{w'}^w\subset \Delta_{Q_{w'}}^R$ and $\hat\Delta_Q\subset \hat\Delta_{Q_{w'}}$. Thus $\bg w'\al, \varpi^\vee\bd\geq 0 $ for all $\al\in \Delta_{w'}^w$. Hence the first assertion. The condition $\bg w'\nu_{w'}^w, \varpi^\vee\bd=0$ is equivalent to  $\bg w'\al, \varpi^\vee\bd= 0 $ for all $\al\in \Delta_{w'}^w$ that is $\varpi^\vee\in w'(\ago_0^{P_{w'}}\oplus \ago_{P_w}^G)$. But we have $w'\ago_{P_w}^G\subset w'\ago_{P_{w'}}=\ago_{Q_{w'}}$ and $\ago_Q\subset  \ago_{Q_{w'}}$. Thus $\varpi^\vee\in \ago_Q\cap w'(\ago_0^{P_{w'}}\oplus \ago_{P_w}^G)=  \ago_Q\cap w'\ago_{P_w}^G$.

2. Clearly (a) implies (b). Let us assume (b). Then by 1 we have $\ago_Q^R\subset w'\ago_{P_w}^G$. Note that $\ago_R \subset \ago_{R_w}=w\ago_{P_w}$ hence $\ago_R\subset w' \ago_{P_w}$ since  $w\in W^Rw'$. Thus we have $\ago_Q\subset w'\ago_{P_w}$ hence $w' M_{P_w}(w')^{-1}\subset M_Q$. So $M_{P_w}\subset M_P \cap (w')^{-1}M_Qw'=M_{P_{w'}}$. Hence $P_w=P_{w'}$ and  (b) implies (a).

Let us prove that (a) is equivalent to (c). Let $w_1\in \,_QW^R_{R_w}$ be such that $w'=w_1w$. We always have $w' M_{P_{w'}} (w')^{-1}\subset w' M_{P_w}(w')^{-1}$ hence $M_{Q_{w'}}\subset w_1 M_{R_w} w_1^{-1}$. If  $P_w=P_{w'}$ we have $M_{Q_{w'}}= w_1 M_{R_w} w_1^{-1}$ and $M_{R_w} \subset w_1^{-1}M_Q w_1$ so $w_1\in W^R(R_w;Q)$. Conversely if $M_{R_w} \subset w_1^{-1}M_Q w_1$ we have
\begin{align*}
  M_{P_{w'}}&=M_P\cap (w')^{-1}M_Q w'\\
            &=w^{-1}(wM_P w^{-1}\cap w_1^{-1}M_Q w_1)w\\
            &=w^{-1}(wM_P w^{-1}\cap M_R\cap w_1^{-1}M_Q w_1)w\\
            &=w^{-1}(M_{R_w}\cap w_1^{-1}M_Q w_1)w\\
              &=w^{-1}M_{R_w}w=M_{P_w}.
\end{align*}

3. By the equivalence of 2.(b) and 2.(c) we have  $(w' \nu_{{w'}}^w)^R_Q\not=0$.  Therefore there exists $\be\in \Delta_Q^R$ such that $\bg w' \nu_{{w'}}^w, \varpi^\vee_\be\bd \not=0$. We write $T_Q^R=\sum_{\al\in \Delta_Q^R} \bg \al, T\bd \varpi_\al^\vee $. For any $\al\in \Delta_Q^R$ and $T\in \ago_0^+$, we have $\bg \al, T\bd  \geq \bg \tilde{\al}, T\bd  \geq 0$ where $\tilde{\al}$ is the root in $ \Delta_0^R \setminus \Delta_0^Q$ that projects on $\al$.  We have also $\bg w' \nu_{{w'}}^w, \varpi^\vee_\al\bd \leq 0$ by assertion 1. Thus we get
\begin{align*}
  \bg (w' \nu_{{w'}}^w)^R_Q,T\bd = \sum_{\al\in \Delta_Q^R} \bg \al, T\bd  \bg w' \nu_{{w'}}^w, \varpi^\vee_\al\bd \leq \bg w' \nu_{{w'}}^w, \varpi^\vee_\be\bd  \bg \be, T\bd \leq   \bg w' \nu_{{w'}}^w, \varpi^\vee_\be\bd d(T).
\end{align*}
The result is clear.
 \end{preuve}

 The construction \eqref{eq:PTRc} is nothing else but an explicit version of the  (Flicker-Rallis) \emph{regularized period} that  we denote by
  \begin{align*}
     \pc^R(\varphi,\la,w)
   \end{align*}
 and that was introduced by Jacquet-Lapid-Rogawski in  \cite[section 7]{JLR} as a substitute for the (in general divergent) integral:
   \begin{align*}
     \int_{[G']_{R',0}}  E^{R}(g,  M(w,\la)\varphi_{w},(w\la)^R)\, dg\\
     =    \int_{[M_{R}']_0}  \int_{K'} \exp(-\bg 2\rho_{R'},H_{R'}(m)\bd)    E^{R}(mk,  M(w,\la)\varphi_{w},(w\la)^R)\, dm
   \end{align*}
   This is what we shall check  among other things in the next proposition. Note we are  in fact considering the obvious variant of their construction when one replaces $G'(F)\back G'(\AAA)^1$ by the quotient  $A_{R'}^\infty M_R'(F)\back M_R'(\AAA)$.

  \begin{proposition}    \label{prop:reg-period} 
\begin{enumerate}
    \item Let $\la\in \ago_{P_w,\CC}^{G,*}$. For any $\varphi\in \Ac_{P,\pi}(G)$ we have
      \begin{align*}
                \pc^R(\varphi,\la,w)=    \pc^{T,R}(\varphi,\la,w).
      \end{align*}
      In particular, the regularized period is well-defined, holomorphic in  a complement of hyperplanes  in $\ago_{P_w,\CC}^{G,*}$ and the right-hand side does not depend on the parameter $T$. 
    \item Let $\om\subset \ago_{P_w,\CC}^{G,*}$ a compact subset of non-empty interior. There exists $r(\la)$ a finite product of non-zero affine functions such that  $\la\in \om\mapsto r(\la) \pc^R(\varphi,\la,w)$ is regular for all $\varphi\in \Ac_{P,\pi}(G)$. Moreover for any holomorphic differential operator $D$ on  $\ago_{P_w,\CC}^{G,*}$, there  exists a continuous semi-norm $\|\cdot\|$ on $\Ac_{P,\pi}(G)$ such that 
\begin{align*}
  |D(r(\la) \pc^R(\varphi,\la,w))| \leq \|\varphi\|, \ \ \ \varphi\in \Ac_{P,\pi}(G), \ \la\in \om.
\end{align*}
\item Assume $P\subset R$ and $w=1$. The regularized period   $\varphi\mapsto \pc^R(\varphi,\la,w)$ gives, at each regular $\la$, a map $\Ac_{P,\pi}(G)\to \CC$ that is $G'(\AAA)$-invariant for the action $I_{P,\pi}(\la^R)$.
    \item If $R=P$ and $w\in W(P,P)$, the regularized period $\pc^P(\varphi,\la,w)$ reduces to the following  integral:
      \begin{align*}
        \int_{[G']_{P',0}} (M(w,\la)\varphi)(g)\,dg
      \end{align*}
      which is convergent outside the singularities of $M(w,\la)$.
    \end{enumerate}
      \end{proposition}

      \begin{preuve} Set $\psi= M(w,\la)\varphi_{w}$, $\la'=(w\la)^R$ and $S=R_w$.  Note that $\la'\in \ago_{S,\CC}^{R,*}$.    By definition,  $\pc^R(\varphi,\la,w)$ is the  sum  indexed by $P_0\subset Q\subset R$ of
        \begin{align*}
         & \int_{A_{R'}^\infty\back [G']_{Q'}}^*  \La^{T,Q}  E^{R}( g, \psi, \la') \tau_Q^R(H(g)-T)\, dg\\
        &=  \sum_{w'\in   \,_QW^R_{S}}     \int_{A_{R'}^\infty\back [G']_{Q'}}^*  \La^{T,Q}  E^{Q}( g, M(w',\la')\psi_{S_{w'}}, w'\la')  \tau_Q^R(H(g)-T)\, dg\\
          &=\sum_{w'\in   \,_QW^R_{S}}     (-2)^{-\dim( \ago^{R}_{Q})}  \int_{[G']_{Q',0}}  \La^{T,Q}  E^{Q}( g, M(w',\la'+\nu_{\psi_{S_{w'}}})\tilde{\psi}_{S_{w'}}, (w'(\la'+\nu_{\psi_{S_{w'}}}))^Q)  \, dg  \times \\
          &\frac{  \exp(\bg w'(\la'+\nu_{\psi_{S_{w'}}}), T_Q^R \bd }{ \hat\theta_Q^R(  w'(\la'+\nu_{\psi_{S_{w'}}}))}
                  \end{align*}
                  where  we defined $\nu_{\psi_{S_{w'}}}\in \ago_{S_w}^{S,*}$ and  $\tilde{\psi}_{S_{w'}}$  by the condition that $\tilde{\psi}_{S_{w'}}=\psi_{S_{w'},-\nu_{\psi_{S_{w'}}}} \in \Ac_{S_{w'}}^0(G)$. The reader  is advised to consult  \cite{JLR} or to take the first two lines just as a suggestive notation and the third one as a definition.  The power of $2$ is due to the discrepancy between the measure on $\ago_Q^R$ and that on $\ago_{Q'}^{R'}$ (see § \ref{S:notations}). 
             
                  Let $w'\in  \,_QW^R_{S}$. Set $w''=w'w$. Then $w''\in  \,_QW_{P}$ by lemma \ref{lem:ww'}. One has $P_{w''}\subset P_w\subset P$,
                  \begin{align*}
                                      M(w',\la'+\nu_{\psi_{S_{w'}}}) \tilde \psi_{S_{w'}}=M(w'',\la+\nu_{w''}^{w}) \varphi_{{w''}}
                  \end{align*}
                  and $\nu_{\psi_{S_{w'}}}=w\nu_{w''}^{w}$. 
Moreover we have $w'\la'=(w''\la)^R$,  $(w'\la')^Q=(w''\la)^Q$. Thus we have $ w'(\la'+\nu_{\psi_{S_{w'}}})=(w''\la)^R+ w''\nu_{w''}^{w}$.  Now the comparison with \eqref{eq:PTRc} is straightforward.

Then assertion 1 and assertion 3 come from lemma \ref{lem:analyticite-m+c} and \cite[theorem 9]{JLR}:  either the argument of the proof of \cite[theorem 9]{JLR} can be used in our context or one can use the parabolic descent of remark \ref{rq:descente-par} to reduce to the case of \cite[theorem 9]{JLR}.

        To prove assertion 2, the existence of $r(\la)$ follows from the definition of $\pc^R(\varphi,\la,w)$, properties of Eisenstein series and lemma \ref{lem:ITQ-analytic}. If we add some affine factors to $r(\la)$ we can even assume that each factor of each term in the definition of $\pc^R(\varphi,\la,w)$, see \eqref{eq:PTRc}, is  holomorphic. Then  for $D=1$ the bound follows from lemma \ref{lem:ITQ-analytic}. Using Cauchy's integral formula, we see that the assertion holds also for any $D$. Now we have to remove the extra factors we add to $r(\la)$. By recursion, we are reduced to the case where we add only one affine factor $l(\la)$. So the bound holds for $g(\la)=l(\la)r(\la)  \pc^R(\varphi,\la,w)$. But then  $r(\la)\pc^R(\varphi,\la,w)$ can be expressed as an integral over a compact subset of a derivative of $g$. The result follows easily.

For the assertion 4, we have
        \begin{align*}
          \int_{[G']_{P',0}} (M(w,\la)\varphi)(g)\,dg=  \int_{[M_{P'}]_{0}} \int_{K'}   \exp(-\bg 2\rho_{P'},H_{P'}(m)\bd) (M(w,\la)\varphi)(mk)\,dkdm.
      \end{align*}
Since $m\mapsto \exp(-\bg 2\rho_{P'},H_{P'}(m)\bd) (M(w,\la)\varphi )(m)$ is square-integrable on $[M_P]_0$, its integral over $[M_{P'}]_0$ is absolutely convergent (see \cite[lemma 3.1]{Yquad}). So assertion 4 comes also from \cite[theorem 9]{JLR}.
      \end{preuve}
    \end{paragr}   

\subsection{Truncated Flicker-Rallis  periods}\label{ssec:asym-period}

\begin{paragr} The following proposition gives an explicit expression for the truncated Flicker-Rallis period in terms of the regularized periods.

  \begin{proposition}\label{prop:MS-periods}
    For $\la\in \ago_{P,\CC}^{G,*}$ in general position, we have:
    \begin{align*}
      \Ic^{T}(\varphi,\la)=\sum_{P_0\subset R} 2^{-\dim(\ago_R^G)}  \sum_{w\in \,_RW_P} \pc^R(\varphi,\la+\nu_w,w) \cdot \frac{\exp(\bg w(\la+\nu_{w}),T_R^G\bd )}{\theta_R^G(w(\la+\nu_{w}))}.
    \end{align*}
  \end{proposition}

  \begin{remarque}\label{rq:MS-periods}
    As we shall see later many terms in fact vanish. By definition the term $\pc^R(\varphi,\la,w)$ is understood to be $0$ unless $P_\pi\subset P_w$. In this case the proof shows that, for  $\la\in \ago_{P,\CC}^{G,*}$ in general position,  not only  $\pc^R(\varphi,\la,w)$ is well-defined (we just knew \emph{a priori}  that $\pc^R(\varphi,\la,w)$ was defined for $\la$  in general position in the \emph{bigger} space $\ago_{P_w,\CC}^{G,*}$)   but also  the corresponding denominator is non-vanishing.
  \end{remarque}

  \begin{preuve}
      We proceed as in the proof of theorem \ref{thm:inversion-calculee}, see also \cite[proof of proposition 4.1]{LR}. With notations as in the proof of  theorem \ref{thm:inversion-calculee}, for $T'\in \ago_0^+$, we have:
\begin{align*}
&  \Ic^{T+T'}(\varphi,\la)\\
 & =\sum_{P_0\subset Q}   \int_{N_{Q'}(\AAA) M_{Q'}(F) A_{G'}^\infty \back G'(\AAA)} \Ga_Q(H_Q(g)-T,T')  \La_m^{T,Q} E(g,  \varphi,\la)\, dg.
 \end{align*}
 Let $Q$ be as in the sum above. The constant term $E_Q(g,  \varphi,\la)$ along $Q$ of the Eisenstein series $E(g,  \varphi,\la)$ is computed by the formula \eqref{eq:cst-term-Eis NEW}: the  term indexed by $w\in  \, _QW_P$  gives in the sum above the contribution:
 \begin{align*}
&   \int_{N_{Q'}(\AAA) M_{Q'}(F) A_{G'}^\infty \back G'(\AAA)} \Ga_Q(H_Q(g)-T,T')  \La_m^{T,Q} E^Q(g,M(w,\la+\nu_w)\varphi_{w},w(\la+\nu_w))\, dg\\
   &=  \int_{[G']_{Q',0}}   \int_{A_{G'}^\infty\back A_{Q'}^\infty } \exp(-2\rho_{Q'},H_{Q'}(ag)\bd ) \Ga_Q(H_Q(ag)-T,T')  \times
   \\ & \La_m^{T,Q} E^Q(ag,M(w,\la+\nu_w)\varphi_{w},w(\la+\nu_w))\, dadg\\
   &= \int_{[G']_{Q',0}} \left(  \int_{A_{G'}^\infty\back A_{Q'}^\infty }  \exp(w(\la+\nu_w)  ,H_{Q}(a)\bd ) \Ga_Q(H_Q(ag)-T,T')  \, da\right) \times\\
   &\La_m^{T,Q} E^Q(g,M(w,\la+\nu_w)\varphi_{w},w(\la+\nu_w))\,dg
 \end{align*}
 Here we have used the equality $2\rho_{Q'}=\rho_{Q}$, see \eqref{eq:rhoQ rhoQ'} and the fact that the projection of $\rho_{Q_w}$ on $\ago_Q^*$ is $\rho_Q$. Taking into account the different choices of measures on $\ago_{Q'}$ and $\ago_{Q}$ which gives a power of $2$, we see that the inner integral is:
 \begin{align*}
    2^{-\dim(\ago_Q^G)}\int_{\ago_Q^G } \Ga_Q(H+H_Q(g)-T,T')  \exp(\bg w(\la+\nu_w),H\bd)\, dH\\
   =\exp(\bg w(\la +\nu_{w})  ,T_Q -H_Q(g)\bd) \cdot \sum_{Q\subset R} (-1)^{\dim(\ago_Q^R)} \frac{\exp(\bg w(\la+\nu_{w} ) ,T'_R \bd)}{(\hat\theta_Q^R\theta_R)(w(\la+\nu_{w}))}.
 \end{align*}
 Thus the outer integral over $[G']_{Q',0}$ reduces to:
 \begin{align*}
   \int_{[G']_{Q',0}}  \La_m^{T,Q} E^Q(g,M(w,\la+\nu_w)\varphi_{w},(w(\la+\nu_w))^Q)\,dg \times \sum_{Q\subset R} (-1)^{\dim(\ago_Q^R)} \frac{\exp(\bg w(\la+\nu_{w} ) ,T'_R \bd)}{(\hat\theta_Q^R\theta_R)(w(\la+\nu_{w}))}\\
   = \Ic^{T,Q}(\varphi,\la+\nu_w,w)  \times \sum_{Q\subset R} (-1)^{\dim(\ago_Q^R)} \frac{\exp(\bg w(\la+\nu_{w} ) ,T'_R \bd)}{(\hat\theta_Q^R\theta_R)(w(\la+\nu_{w}))}
 \end{align*}
The expression $  \Ic^{T,Q}(\varphi,\la,w)$ vanishes unless $P_\pi\subset P_w$.  So you may and shall  assume $P_\pi\subset P_w$: in this case each denominator that appears is non-vanishing for $\la$ in general position by a variant of  lemma \ref{lem:non-vanish-theta}. Using this expression and inverting the sum over $R$ and $Q$ thanks to lemma \ref{lem:ww'} we find:
\begin{align*}
&  \Ic^{T+T'}(\varphi,\la)\\
  &=\sum_R 2^{-\dim(\ago_R^G)} \sum_{w\in\,_RW_P}    \frac{\exp(\bg w(\la +\nu_{w})  ,(T+T')_R \bd)}{\theta_R(w(\la+\nu_{w}))} \times\\
  &\left[ \sum_{Q\subset R} (-2)^{\dim(\ago_Q^R)}  \sum_{w'\in    \,_QW^R_{R_w}  w } \Ic^{T,Q}(\varphi,\la+\nu_{w'},w')\cdot \frac{\exp(\bg w'(\la+\nu_{{w'}}^w),T_Q^R\bd)}{\hat{\theta}_Q^R(w'(\la+\nu_{{w'}}^w))}. \right]
\end{align*}
However the bracket is nothing else but  $\pc^R(\varphi,\la+\nu_w,w)$, see \eqref{eq:PTRc} and proposition \ref{prop:reg-period}. It suffices to take $T'=0$ to conclude.
\end{preuve}

\end{paragr}

\begin{paragr}
  
\begin{corollaire}\label{cor:MS-periods}
   Let $R$ be a standard  parabolic subgroup and $w\in \,_RW_P$.  For $\la\in \ago_{P_w,\CC}^{G,*}$ in general position, we have:
    \begin{align}\label{eq:NEW-MS-periods}
      \Ic^{T,R}(\varphi,\la,w)=\sum_{P_0\subset Q\subset R} 2^{-\dim(\ago_Q^R)}  \sum_{w'\in \,_QW^R_{R_w}w} \pc^Q(\varphi,\la+\nu_{w'}^w,w') \cdot \frac{\exp(\bg w'(\la+\nu_{w'}^w),T_Q^R\bd )}{\theta_Q^R(w'(\la+\nu_{w'}^{w}))}.
    \end{align}
  \end{corollaire}

  \begin{preuve}
  The proof is similar to that of proposition \ref{prop:MS-periods}. Alternatively, by parabolic descent to $M_R$, see   remark \ref{rq:descente-par}, one is reduced to proposition \ref{prop:MS-periods}. Details are left to the reader.
  \end{preuve}

      \begin{corollaire}\label{cor:type-exponents}
We keep the notations of corollary \ref{cor:MS-periods}. The map
$$T\mapsto \Ic^{T,R}(\varphi,\la,w)$$ 
coincides with a polynomial exponential whose set of exponents  is included in the following set:
\begin{align}\label{eq:set-exponents}
S_\la= \{  (w'(\la+\nu_{w'}^w))_Q^R \mid  Q\subset R \ ; \ w ' \in \,_QW^R_{R_w}  w  \ ; \ P_\pi\subset P_{w'}  \}.
\end{align}
\end{corollaire}

\begin{preuve}
  It is a straightforward consequence of corollary \ref{cor:MS-periods}, the nature of regularized period (see proposition \ref{prop:reg-period}) and lemma \ref{lem:py-exp-lisse}.
\end{preuve}
  
\end{paragr}

    \begin{paragr}\label{S:type-exp}   Let $0<\eps<1$.  We fix $c,k>0$ and $\La\geq 0$ and we set
  \begin{align}\label{eq:set-rc}
    \rc=\rc_{\La,c,k}=\{\la\in \ago_{0,\CC}^* \mid \|\Re(\la)\| < c(1+\La+\|\Im(\la)\|)^{-k}\}.
  \end{align}
  As follows from the methods and results of section \ref{sec:disc}, see in particular (the proof of) theorem \ref{thm:maj-Eis}, we may and shall assume that $c,k$ and $\La$ are chosen so that for all standard parabolic subgroups $R$ and all  $w\in \, _RW_{P}$ the Eisenstein series $E^R(M(w,\la)\varphi_w,(w\la)^R)$ is holomorphic on $\ago_{P_w,\CC}^{G,*}\cap \rc$. In this case, by lemma \ref{lem:ITQ-analytic},  $\Ic^{T,R}(\varphi,\la,w)$  is also holomorphic on $\ago_{P_w,\CC}^{G,*}\cap \rc$.

\begin{lemme}\label{lem:types} 
  Let $Q\subset R$ be  standard parabolic subgroups, $w\in \, _RW_{P}$, $w ' \in \,_QW^R_{R_w}  w $ such that $ P_\pi\subset P_{w'}$. 
  Let   $\la\in \ago_{P_w,\CC}^{G,*} \cap \rc$ and  $\mu= (w'(\la+\nu_{{w'}}^w))_Q^R$.
  \begin{enumerate}
  \item  If $w'\in  W^R(R_w;Q)w$, we have $ \bg \Re(\mu),T\bd \geq -c \|T\|$ for all $T\in\ago_0$.
  \item    If $w'\notin  W^R(R_w;Q)w$, we have $ \bg\Re( \mu),T\bd \leq  (c- c_1\eps) \|T\|$  for all $T\in\ago_0$ such that $d(T)\geq  \eps \|T\|$ where $c_1>0$ appears in lemma \ref{lem:positivite-nu} assertion 3.
  \end{enumerate}
\end{lemme}

\begin{preuve}
  We have  $|\bg \Re(\mu),T\bd |\leq \|\Re(\mu)\| \|T\|$. In case $1$, we have  then $\Re(\mu)= (w'\Re(\la))_Q^R$ by lemma \ref{lem:positivite-nu} assertion 2. Thus $\|\Re(\mu)\| \leq \|\Re(\la)\| <c$.
  In case $2$, by  lemma \ref{lem:positivite-nu} assertion 3, for any $T\in \ago_0^+$ we have
    \begin{align*}
      \bg\Re( \mu),T\bd =  \bg (w'\Re( \la))_Q^R,T\bd+    \bg (w'\nu_{{w'}}^w)_Q^R,T\bd \leq  c \|T\|-c_1d(T).
    \end{align*}
    The result follows.
\end{preuve}

Let $\la\in \ago_{P_w,\CC}^{G,*}\cap \rc$. We shall say that an exponent $\mu$ in the set $S_\la$ defined in  \eqref{eq:set-exponents} is of   type 1, resp. of type 2,  if it satisfies the inequality 1, resp. 2, of lemma \ref{lem:types}.

\begin{remarque}\label{rq:types}
  We shall assume in the following that we have $c<c_1\eps/2$. Then  any point of $\ago_0^*$ that appears in  $S_\la$ for some $\la \in \ago_{P_w,\CC}^{G,*}\cap \rc$  is either of type 1 or of type 2 but cannot be  both.
\end{remarque}

We then define $S_\la^1\subset S_\la$ as the subset of exponents of type $1$. We define $S_\la^2$ as the complement of $S_\la^1$ in $S_\la$. As follows from remark \ref{rq:types}, we have $S_\la^1\cap S_{\la'}^2=\emptyset$ for $\la,\la'\in  \ago_{P_w,\CC}^{G,*}\cap \rc$.
\end{paragr}

\begin{paragr}  Let $w\in \,_RW_P$. We introduce the following expression (the superscript $m$ is for ``main''):
      \begin{align}\label{eq:PTRm}
        \pc^{T,R,m}(\varphi,\la,w)&=\sum_{ P_0\subset Q\subset R}   2^{-\dim(\ago_Q^R)}   \sum_{w'\in W^R(R_{w};Q)w} \pc^{Q}(\varphi,\la,w') \cdot \frac{\exp(\bg w'\la,T_Q^R\bd)}{\theta_Q^R(w'\la)}.                                    
      \end{align}
Observe that for all  $w'\in W^R(R_{w};Q)w$ we have $\nu_{w'}^w=0$: in particular the right-hand side of \eqref{eq:PTRm} is really a subsum of the right-hand side of the expression \eqref{eq:NEW-MS-periods}.  If $w=1$ we set $     \pc^{T,R,m}(\varphi,\la)=\pc^{T,R,m}(\varphi,\la,1)$. If $R=G$, we remove the exponent $G$ from the notation.

      \begin{proposition} \label{prop:type1}   Let $w\in \,_RW_P$. 
        \begin{enumerate}
        \item The expression $\pc^{T,R,m}(\varphi,\la,w)$ is well-defined and holomorphic  for $\la$ in a complement of hyperplanes  in $\ago_{P_w,\CC}^{G,*}$, resp. in $\ago_{P,\CC}^{G,*}$.
        \item For $\la\in \ago_{P_w,\CC}^{G,*}\cap \rc$, the map $T\mapsto \pc^{T,R,m}(\varphi,\la,w)$ coincides with the polynomial exponential given by the summand of exponents of type 1 of $T\mapsto \Ic^{T,R}(\varphi,\la,w)$ (see corollary \ref{cor:type-exponents}).
        \item The map $\la\mapsto \pc^{T,R,m}(\varphi,\la,w)$ is holomorphic on  $\ago_{P_w,\CC}^{G,*}\cap \rc$.
        \end{enumerate}
      \end{proposition}

      \begin{preuve}
        1. The case of  $\ago_{P_w,\CC}^{G,*}$ follows from (the proof of) proposition \ref{prop:reg-period}. The case  of  $\ago_{P,\CC}^{G,*}$ is observed in remark \ref{rq:MS-periods}.
        
        2. It is straightforward  to identify the summand of exponents of type 1 in corollary \ref{cor:MS-periods} thanks to lemma \ref{lem:types}.
        
        3. Once we have identified $\pc^{T,R,m}(\varphi,\la,w)$ with the summand of exponents of type 1 of $\Ic^{T,R}(\varphi,\la,w)$, the holomorphy on $\ago_{P_w,\CC}^{G,*}\cap \rc$ follows from the holomorphy of  $\Ic^{T,R}(\varphi,\la,w)$  on $\ago_{P_w,\CC}^{G,*}\cap \rc$  and lemma \ref{lem:py-exp-lisse}.
      \end{preuve}
    \end{paragr}

\subsection{Singularities of  regularized periods} \label{ssec:sing-period}

\begin{paragr}[Distinction.] ---\label{S:distinction}
  We follow the notations of § \ref{S:exposant-disc}. Recall that $P=MN_P$ is a standard parabolic subgroup, $\pi\in\Pi_{\disc}(M)$ and $\varphi\in \Ac_{P,\pi}(G)$. We shall say that $\pi$ is $M'$-distinguished if the map
  \begin{align*}
    \psi\in \Ac_{\pi}(M)\mapsto \int_{[M']_0} \psi 
  \end{align*}
does not vanish identically. Since $\psi$ is a discrete automorphic form, the integral is absolutely convergent  (see \cite[lemma 3.1]{Yquad}). We denote by $\pi^*$ the conjugate dual of $\pi$. If $\pi$ is distinguished then $\pi=\pi^*$: this is well known if $\pi $ is  a cuspidal representation, see the work of Flicker \cite{Flicker},  and follows from the work of Yamana \cite{Yquad} in general.
\end{paragr}

\begin{paragr}[Vanishing of regularized periods.] ---  As a generalization of  \cite[theorem 23]{JLR}, we show in the next proposition that the regularized period
  \begin{align*}
    \pc(\varphi,\la)=\pc^G(\varphi,\la,1)
  \end{align*}
often vanishes for $\la\in \ago_{P,\CC}^{G,*}$.

\begin{proposition}\label{prop:vanishing}
   Assume that the period $\pc(\varphi,\la)$ is well-defined at $\la\in \ago_{P,\CC}^{G,*}$ and non-zero. Then we are in one of the following situations:
\begin{enumerate}
\item $P=G$,  $\pi$ is $G'$-distinguished and $\pi=\pi^*$.
\item $n=2r$, $G=G_{2r}$, $M\simeq G_r\times G_r$ and there exists a discrete automorphic representation $\sigma$ of   $G_r$ such that $\pi=\sigma \boxtimes \sigma^*$.
\end{enumerate}
  \end{proposition}

  \begin{preuve}
1. It results from observations in § \ref{S:distinction}.

2. Assume $P\subsetneq G$. By assumption and proposition \ref{prop:reg-period}, the period $\pc(\varphi,\la)$ is  well-defined and non-zero on a non-empty open subset $\Om$ of $\ago_{P,\CC}^{G,*}$.  The representation $\pi$ can be written as a restricted tensor product $\otimes_{v\in V} \pi_v$ of   irreducible  representations of $M(F_v)$. Let  $v$ be a finite  place of  $F$ which is totally split in  $E$.  For such a place, $G'(F_v)$ is identified with the diagonal subgroup of $G(F_v)=G'(F_v)\times G'(F_v)$.   Locally we get a   non-zero and $G'(F_v)$-invariant linear  form on  the induced representation $\Ind_{P(F_v)}^{G(F_v)}(\pi_{\la})$. Here $\pi_\la(m)= \exp(\bg \la, H_P(m)\bd)\pi(m) $  for $m\in M(F_v)$. Let $\Om_v\subset \Om$ be a non-empty open subset such that for $\la\in \Om_v$ the representation  $\Ind_{P(F_v)}^{G(F_v)}(\pi_{\la})$ is irreducible.

According to the  identification $M(F_v)=M'(F_v)\times M'(F_v)$, one writes $\pi_v= \pi_{1}\boxtimes \pi_{2}$.  For any $\la\in \Om_v$,   the induced representation  $\Ind_{P'(F_v)}^{G'(F_v)}(\pi_{1,\la})$ is isomorphic to the contragredient representation of $\Ind_{P'(F_v)}^{G'(F_v)}(\pi_{2,\la})$ and thus there is $w\in W^{G'}(M')$ such that $\pi_{1,\la+w\la}\simeq w\pi_{2}^\vee$ where $\pi_2^\vee$ is the contragredient of $\pi_2$. Let $w\in W^{G'}(M')$. We view $w$ as an endomorphism of $\ago_{M'}^{G',*}$. If $\Ker(w+\Id)\subsetneq \ago_{M'}^{G',*}$, the set of $\la\in  \ago_{M'}^{G',*}$  such that  $\pi_{1,\la+w\la}\simeq w\pi_{2}^\vee$ is a closed subset with empty interior. So there are $\la\in \Om_v$ and $w\in W^{G'}(M')$ such that 
$\pi_{1,\la+w\la}\simeq w\pi_{2}^\vee$  and $\Ker(w+\Id)= \ago_{M'}^{G',*}$. The latter condition implies $n=2r$ and $w$ switches the two factors of $M'=G_r'\times G_r'$. The relation $\pi_{1,\la+w\la}\simeq w\pi_{2}^\vee$ boils down to
\begin{align}\label{eq:pi1wpi2}
  \pi_{1}\simeq w\pi_{2}^\vee. 
\end{align}

Globally we  have $M=G_r\times G_r$ and we can write $\pi=\sigma\boxtimes \tau$ accordingly. From \eqref{eq:pi1wpi2}, we infer that $\tau_v\simeq \sigma_v^*$ at any finite  place of  $F$ which is totally split in  $E$. Thus the automorphic  representations $\tau$ and $\sigma^*$ of $G_r$ are isobaric  (by the M{\oe}glin-Waldspurger classification of the discrete spectrum, cf. \cite{MW}) and  coincide on the set of  places  of $E$ of  degree $1$ over  $F$. It follows from Ramakrishnan's theorem (see \cite{Ram2} theorem A) that $\tau=\sigma^*$ (one can also consult \cite{Ram1} corollary B).
  \end{preuve}
\end{paragr}

\begin{paragr}[Singularities.] --- We consider a set $\rc$ as in § \ref{S:type-exp}. Let $R$ be a standard parabolic subgroup and let $w\in W(P;R)$.

  We assume that $\rc$ is chosen so that for all $\la\in \ago_{P,\CC}^*\cap\rc$  the intertwining operators $M(w',\la)$  induce a  holomorphic operator $\Ac_{P,\pi}(G)\to \Ac_{Q,w\pi}(G)$ for all $w'\in \cup_Q W(P,Q)$ where $Q$ runs over the standard parabolic subgroups (see remark \ref{rq:holom op entr}).

   We have  $M_{R_w}=wM_{P}w^{-1}$ and we write $M_{R_w}=G_{n_1}\times \ldots \times G_{n_r}$ and $w\pi=\sigma_1\boxtimes \ldots \boxtimes \sigma_r$ accordingly. We write $M_R=G_{m_1}\times \ldots \times G_{m_s}$. Any involution $\xi\in W_2(M_{R_w})$ permutes the blocks $G_{n_i}$. In this way, it is identified with a permutation of $\{1,\ldots,r\}$.

  \begin{theoreme}\label{thm:singularities}
    \begin{enumerate}
    \item The regularized period $\pc^{R}(\varphi,\la,w)$ vanishes identically on $\ago_{P,\CC}^{G,*}$ unless there exists an involution $\xi\in W_2^R(M_{R_w})$ such that $\ago_{R_w}^R=\ago_{R_w}^{-\xi}$  and $\xi w\pi=w\pi^*$. If $\xi(i)=i$ then  the representation $\sigma_i$ is moreover $G_{n_r}'$-distinguished.
    \item On $\rc\cap \ago_{P,\CC}^{G,*}$ the only possible poles of $\pc^{R}(\varphi,\la,w)$ are simple and along the hyperplanes $\bg w \la,\al^\vee\bd=0$ where $\al\in \Delta_{R_w}^R$.
    \item Let $\al\in \Delta_{R_w}^R$. If the hyperplane  $\bg w \la,\al^\vee\bd=0$ of $\ago_{P,\CC}^{G,*}$ is singular for  $\pc^{R}(\varphi,\la,w)$ then the following conditions are  satisfied:
      \begin{enumerate}
      \item $s_\al w\pi=w\pi$ where $s_\al$ is the elementary symmetry attached to $\al$ and $M(s_\al,0)$ acts by $-1$ on $\Ac_{R_w,w\pi}(G)$.
      \item If $s_\al(i)=i+1$ then we have $\sigma_i=\sigma_{i+1}$ and these representations are $G_{n_i}'$-distinguished.
      \end{enumerate}
    \end{enumerate}
  \end{theoreme}

  \begin{remarque}
    The possibility of singular hyperplanes in assertion 3 does not contradict remark \ref{rq:MS-periods}. Indeed $w^{-1}\al^\vee$ induces a non-zero linear form on $\ago_{P,\CC}^{G,*}$ for any $\al\in \Delta_{R_w}^R$ (otherwise we would have $\al^\vee\in \ago_{0}^{wM_Pw^{-1}}\cap\ago^{R}_{R_w}=(0)$).
  \end{remarque}

  \begin{preuve}
 By remark \ref{rq:descente-par},  we note that we can set $\varphi'=M(w,\la)\varphi $, $\la'=w\la$ and replace $P$ by $R_w$ and $\pi$ by $w\pi$. Without loss of generality, we see that we may and shall assume $w=1$ and $P\subset R$.

1.  By parabolic descent, see remark \ref{rq:descente-par}, one is reduced to consider regularized periods associated to $M_R$ and thus to proposition \ref{prop:vanishing}.

2. By assertion 1, we may and shall assume that there exists $\xi\in W_2(M_{})$ such that $\ago_{P}^R=\ago_{P}^{-\xi}$  and $\xi \pi=\pi^*$. We start from the fact that the expression
 \begin{align*}
        \pc^{T,R,m}(\varphi,\la)&=\sum_{ P_0\subset Q\subset R}   2^{-\dim(\ago_Q^R)}   \sum_{w\in W^R(P;Q)} \pc^{Q}(\varphi,\la,w) \cdot \frac{\exp(\bg w\la,T_Q^R\bd)}{\theta_Q^R(w\la)}                                 
      \end{align*}
is holomorphic on  $\rc\cap \ago_{P,\CC}^{G,*}$ (see proposition \ref{prop:type1}). For $w\in W^R(P;Q)$ the parabolic subgroups $P$ and $Q_w$ are $R$-associated. In our situation this implies $P=Q_w$ and $P\subset Q$. In this case we have $ W^R(P;Q)\subset W^R(M)$. Observe also that  we have here for $P\subset Q\subset R$ and $w\in W^R(M)$
\begin{align*}
  \theta_Q^R(w\la)=\eps_Q^R(w) \theta_Q^R(\la)  \ \text{ and } \ \theta_P^R(\la)= \theta_P^Q(\la) \theta_Q^R(\la)
\end{align*}
where $\eps_Q^R(w)\in \{\pm 1\}$.
So we have:
\begin{align*}
  \theta_P^R(\la) \pc^{R}(\varphi,\la) = \theta_P^R(\la) \pc^{T,R,m}(\varphi,\la) \\
-\sum_{ P \subset Q\subsetneq R}   2^{-\dim(\ago_Q^R)}   \sum_{w\in W^R(P;Q)} \eps_P^R(w)  \pc^{Q}(\varphi,\la,w) \theta_P^Q(w\la)  \exp(\bg w\la,T_Q^R\bd).
\end{align*}

By recursion, we may and shall assume that  $\pc^{Q}(\varphi,\la,w) \theta_P^Q(w\la) $ is holomorphic on $\rc\cap \ago_{P,\CC}^{G,*}$. It follows that $\theta_Q^R(\la) \pc^{G}(\varphi,\la) $ is also holomorphic on $\rc\cap\ago_{P,\CC}^{G,*}$. 

3. Let $\al\in \Delta_{P}^R$.   Let $P\subset P_\al, R_\al\subset R$ be defined by $\Delta_P^{P_\al}=\{\al\}$ and $\Delta_P^{R_\al}=\Delta_P^{R}\setminus\{\al\}$. Note that 
$$\ago_P^R=\ago_P^{P_\al}\oplus \ago_{P}^{R_\al}$$
and so $\ago_{P_\al}^R=\ago_{P}^{R_\al}$.
If one has $P\subset Q\subset R$ then either $Q\subset R_\al$ or $P_\al\subset Q$. In the former case we have
$$ W^R(P;Q)= W^{R_\al}(P;Q)W^{P_\al}(M)$$
and in the latter case $ W^R(P;Q)= W^{R_\al}(P;Q_\al)$ where $Q_\al=R_\al\cap Q$. The non trivial element of $W^{P_\al}(M)$ is $s_\al$.   Then $\pc^{T,R,m}(\varphi,\la)$  is the sum of the following three contributions:
\begin{align}\label{eq:contrib1}
  \sum_{ P_\al\subset Q\subset R}   2^{-\dim(\ago_Q^R)}   \sum_{w\in W^{R_\al}(P;Q_\al)} \pc^{Q}(\varphi,\la,w) \cdot \frac{\exp(\bg w\la,T_Q^R\bd)}{\theta_Q^R(w\la)}        
\end{align}
\begin{align}\label{eq:contrib2}
  & \pc^{T,R_\al,m}(\varphi,\la)\frac{\exp(\bg \la,T_P^{P_\al}\bd)}{\theta_P^{P_\al}(\la)} =\\
 \nonumber& \left(\sum_{ P \subset Q\subset R_\al}   2^{-\dim(\ago_Q^R)}   \sum_{w\in W^{R_\al}(P;Q)} \pc^{Q}(\varphi,\la,w) \cdot \frac{\exp(\bg w\la,T_Q^{R_\al}\bd)}{\theta_Q^{R_\al}(w\la)}        \right) \frac{\exp(\bg \la,T_P^{P_\al}\bd)}{\theta_P^{P_\al}(\la)} 
\end{align}

\begin{align}\label{eq:contrib3}
&\pc^{T,R_\al,m}(M(s_\al,\la)\varphi,s_\al\la)\frac{\exp(\bg s_\al \la,T_P^{P_\al}\bd)}{\theta_P^{P_\al}(s_\al \la)} \\
\nonumber&=-\pc^{T,R_\al,m}(M(s_\al,\la)\varphi,s_\al\la)\frac{\exp(-\bg  \la,T_P^{P_\al}\bd)}{\theta_P^{P_\al}( \la)} 
\end{align}

For  $P_\al\subset Q\subset R$ and $w\in W^{R_\al}(P;Q)$ we set  $\pc^{Q}_\al(\varphi,\la,w) =\bg \la,\al^\vee\bd \pc^{Q}(\varphi,\la,w) $.

Recall that $\pc^{T,R,m}(\varphi,\la)$ is holomorphic on $\rc\cap \ago_{P,\CC}^{G,*}$. So   $\bg \la,\al^\vee\bd \pc^{T,R,m}(\varphi,\la)$  vanishes for $\la\in \Ker(\al^\vee)\cap \rc\cap \ago_{P,\CC}^{G,*}$. Summing the three expressions above and multiplying them by $\bg \la,\al^\vee\bd $ we deduce that we have for such $\la$:

\begin{align*}
  \sum_{ P_\al\subset Q\subset R}   2^{-\dim(\ago_Q^R)}   \sum_{w\in W^{R_\al}(P;Q_\al)} \pc_\al^{Q}(\varphi,\la,w) \cdot \frac{\exp(\bg w\la,T_Q^R\bd)}{\theta_Q^R(w\la)}    \\
+    \pc^{T,R_\al,m}(\varphi,\la) - \pc^{T,R_\al,m}(M(s_\al,0)\varphi,\la)=0.
\end{align*}
This is  a polynomial exponential in $T$. Extracting the part of exponent $0$ (which corresponds in the double sum to $Q=R$ and thus $w=1$) we get:
\begin{align*}
  \pc_\al^{R}(\varphi,\la,1)  +\pc^{R_\al}(\varphi,\la) - \pc^{R_\al}(M(s_\al,0)\varphi,\la)=0.
\end{align*}

Assume $\bg \la,\al^\vee\bd$ is singular for  $\pc^{R}(\varphi,\la,1)$. Then  $\pc_\al^{R}(\varphi,\la,1) \not=0$ and
\begin{align}
  \label{eq:non-egal}
\pc^{R_\al}(\varphi,\la) \not=\pc^{R_\al}(M(s_\al,0)\varphi,\la).
\end{align}
By assertion 1, we have $\xi\pi=\pi^*$ for $\xi\in W_2^R(P)$ such that $\ago_{P}^R=\ago_{P}^{-\xi}$   . If  $\pc^{R_\al}(\varphi,\la)\not=0 $, resp. $\pc^{R_\al}(M(s_\al,0)\varphi,\la)\not=0$, we have $\xi' \pi=\pi^*$, resp.  $\xi' s_\al \pi=s_\al\pi^*$, for some $\xi'\in W^{R_\al}(P)$ such that $\ago_P^{-\xi'}=\ago_P^{R_\al}$. We have $s_\al=\xi'\xi$. Since $\pc^{R_\al}(\varphi,\la) $ and  $\pc^{R_\al}(M(s_\al,0)\varphi,\la)$ cannot both vanish, we have in any case $s_\al \pi=\xi'\xi \pi=\xi'\pi^*=\pi$ and the representations $\sigma_i=\sigma_{i+1}$ where $i$ is such that   $s_\al(i)=i+1$, are $G_{n_i}'$-distinguished by proposition \ref{prop:vanishing}. Since $s_\al\pi=\pi$ and $s_\al^2=1$, the intertwining operator $M(s_\al,0)$ acts by a scalar which must be $\pm 1$. By \eqref{eq:non-egal}, it is  $-1$.
\end{preuve}
\end{paragr}

\section{Intertwining periods}\label{sec:inter-period}

\subsection{Definition and analytic continuation}\label{ssec:an-cont}

\begin{paragr}
  Let $P=MN_P$ be a standard parabolic subgroup of $G$ and $\pi\in \Pi_{\disc}(M)$.
\end{paragr}

\begin{paragr} \label{S:Jxitilde} Let $\xi\in W_2^{\mathrm{st}}(M)$. We identify this element with the unique permutation matrix that shuffles  the diagonal blocks of $M$ without causing any internal change within each block. Let $\tilde{\xi}\in G(F)$ be such that $\tilde\xi \iota(\tilde\xi)^{-1}=\xi$ and 
\begin{align*}    P_{\tilde \xi}=G'\cap \tilde\xi^{-1}P \tilde\xi.
\end{align*}
This is an $F$-subgroup of $G'$ with a Levi decomposition $P_{\tilde \xi}=M_{\tilde \xi} N_{\tilde \xi}$ where $N_{\tilde \xi}=G'\cap \tilde\xi^{-1}N \tilde\xi$ is the unipotent radical of $P_{\tilde \xi}$ and the Levi factor is given by
\begin{align*}
  M_{\tilde \xi}=G'\cap \tilde\xi^{-1}M \tilde\xi.
\end{align*}
  
The group $M_{\tilde \xi}$ is reductive whereas $N_{\tilde \xi}$ is a unipotent subgroup. Moreover $P_{\tilde \xi}=M_{\tilde \xi} N_{\tilde \xi}$ is a Levi decomposition. 

For any $\la\in \ago_{P,\CC}^{*,-\xi}$, we define the \emph{intertwining period}
  \begin{align}\label{eq:inter-period}
    J(\xi,\varphi,\la)=\int_{  N_{\tilde\xi} (\AAA) A_{M_{\tilde \xi}}^\infty  M_{\tilde\xi}(F) \back G'(\AAA)} \varphi_\la(\tilde\xi h)\,dh.
  \end{align}
  We need some explanations. We equip the group $N_{\tilde\xi} (\AAA) A_{M_{\tilde \xi}}^\infty  M_{\tilde\xi}(F)$ with  the  right-invariant Haar measure we get by transport of the product measure from the composition map  
$$N_{\tilde\xi} (\AAA) \times A_{M_{\tilde \xi}}^\infty \times  M_{\tilde\xi}(F) \to N_{\tilde\xi} (\AAA) A_{M_{\tilde \xi}}^\infty  M_{\tilde\xi}(F).$$
This measure is not left-invariant and the modular character is given by 
$$\delta_{P,\tilde\xi}:x\mapsto \exp(\bg \rho_P,H_P(\tilde\xi x {\tilde \xi}^{-1})\bd),$$
 see \cite[VII p.221]{JLR}. Thus the integral in \eqref{eq:inter-period} has to be understood as a right-$G'(\AAA)$-invariant linear form on the space of $(N_{\tilde\xi} (\AAA) A_{M_{\tilde \xi}}^\infty  M_{\tilde\xi}(F),\delta_{P,\tilde\xi})$-equivariant functions and this space contains $\varphi_\la(\tilde\xi\cdot)$ for  $\la\in \ago_{P,\CC}^{*,-\xi}$. So the integral \eqref{eq:inter-period} makes sense at least formally. Note also that, for any $\xi$, the element $\tilde\xi$ is unique up to a right translation by an element of $G'(F)$. Such a change does not affect the definition of  $J(\xi,\varphi,\la)$, hence the notation. We show in the next proposition that the integral  \eqref{eq:inter-period} is absolutely convergent in some cone. 

\begin{proposition}\label{prop:convergence}
There exists $c\in \RR$ such that for any level $J\subset K_f$, any large $N>0$ there exists a finite family $(X_i)_{i\in I}$ of elements of $\uc(\ggo_\infty)$ such that for  all  $b>c$ there exists $C>0$ such that
\begin{align*}
  \int_{ A_{M_{\tilde \xi}}^\infty  M_{\tilde\xi}(F) N_{\tilde\xi} (\AAA)\back G'(\AAA)} |\varphi_\la|(\tilde\xi h)\,dh \leq C \sum_{i\in I}  \|\varphi\|_{-N,X_i}   
\end{align*}
for all  $\varphi \in\Ac_{P,\pi}(G)$ and  all $\la\in \ago_{P,\CC}^{*,-\xi}$ such that $b>\Re(\bg\la,\al^\vee\bd)>c$ for all $\al\in \Delta_P$ such that $\xi\al=-\al$.
\end{proposition}
\end{paragr}

\begin{paragr}[Proof of proposition \ref{prop:convergence}.] --- First we prove the following result.

  \begin{lemme}\label{lem:first-bound}
   For any level $J\subset K_f$, any large $N>0$ there exists a finite family $(X_i)_{i\in I}$ of elements of $\uc(\ggo_\infty)$ such that
\begin{align*}
  \int_{   [M_{\tilde\xi}]_0} |\varphi_{-\rho_P}|(\tilde\xi  m h)\,dm \leq \sum_{i\in I}  \|\varphi\|_{-N,X_i}     
\end{align*}
for all $ h\in G'(\AAA)$ and $\varphi \in\Ac_{P,\pi}(G)$.
  \end{lemme}

  \begin{preuve}
    Let $L$ be the standard Levi subgroup containing $M$ and defined by $\ago_M^L=\ago_M^{-\xi}$. We may and shall assume that $\tilde\xi\in L$. Then $M_{\tilde\xi}\subset L$. Let $Q$ be the standard parabolic subgroup of Levi $L$. Using  the Iwasawa decomposition $G'(\AAA)=N_{Q'}(\AAA)L'(\AAA)K'$, we may assume that  $h\in L'(\AAA)K'$. Note that $\sup_{k'\in K'}   \|\varphi(\cdot k')\|_{-N,X}$ is bounded by a finite sum of norms  $ \|\varphi\|_{-N_i,X_i}$ we may even assume that $h\in L'(\AAA)$. The group $L$ is a product of linear factors $G_r$. Since the discussion is the same for each factor, for simplicity we shall assume that $L=G$. Then there are two cases. First $\xi=1$ and $M=G$. In this case the proof is close to that of lemma \ref{lem:comparison-norm}: it first uses the inversion formula for the mixed truncation operator, see \cite[eq. (19) p.190]{JLR} and then proposition \ref{prop:cont-LaTm}. The details are left to the reader. We inspect the second case more carefully: we have then $G=G_n$ with $n=2r$, $M=G_r\times G_r$ and $\xi$ is the non-trivial element of $W_2(M)$. One can take (with $\tau$ as in § \ref{S:notations})
$$\tilde\xi=\begin{pmatrix}
  I_r & \sqrt{\tau} I_r\\ I_r& - \sqrt{\tau}I_r
\end{pmatrix}
$$
Then the map $m\mapsto \tilde\xi  m \tilde\xi^{-1}$ identifies $M_{\tilde\xi}$ with the  embedding $G_r\hookrightarrow G_r\times G_r\subset G$ given by $g\mapsto (g,\iota(g))$. Using the Iwasawa decomposition $G(\AAA)=N_P(\AAA)M(\AAA)K$ applied to the element $\tilde\xi h$, we see that it suffices to prove the existence of  $N>0$ and a family $(X_i)_{i\in I}$ of elements of  $\ugo_{\mgo_\infty}$ such that
\begin{align}\label{eq:bound-psi}
  \int_{   [G_r]_0} |\varphi_{}|(    g x,\iota(g)y)\,dg \leq \sum_I \|\varphi\|_{-N,X_i}
\end{align}
for all $x,y\in G_r(\AAA)$ and $\varphi\in \Ac_{\pi}(M)^J$ where $J$ is a fixed level in $K_f\cap M(\AAA)$. By a slight variation  on \cite[proof of lemme I.4.1]{MWlivre}, for all $\mu\in \ago_{0}^{P_\pi,*}$ there exist $N>0$ and $(X_i)_{i\in I}$ as above such that for all $\varphi\in \Ac_{\pi}(M)^J$ and  $m\in M(\AAA)$ we have
\begin{align}\label{eq:bound-varphi}
  |\varphi(m)| \leq   \left(\sum_I \|\varphi\|_{-N,X_i}  \right) \inf_{\delta} \exp(\bg \mu+\rho_{P_\pi} +\nu_{P_\pi},H_0(\delta m)\bd)
\end{align}
where the infimum is taken over the set of $\delta\in M(F)$ such that $\delta m$ belongs  to a fixed Siegel set $\SG^M$ of $M(\AAA)$. Note that if $m\in \SG^M$ then this set is contained in a finite set independent  of $m$.

We write $P_\pi=Q_1\times Q_2$ and $\nu_{P_\pi}=(\nu_{Q_1},\nu_{Q_2})$ accordingly. To get \eqref{eq:bound-psi}, we use \eqref{eq:bound-varphi} and Cauchy-Schwartz inequality so that we are reduced to find $\mu\in \ago_{P_{0,r}}^{Q_1,*}$ such that the integral 
\begin{align*}
  \int_{   [G_r]_0}  \left(\inf_{\delta} \exp(\bg \mu+\rho_{Q_1} +\nu_{Q_1},H_0(\delta (g)\bd)\right)^2 \,dg.
\end{align*}
is convergent where the infimum is taken over $\delta\in M(F)$ such that $\delta m$ belongs  to a fixed Siegel set $\SG^{G_r}$ of $G_r(\AAA)$. In particular we are reduced to show that the integral is finite:
\begin{align*}
  \int_{ \SG^{G_r}}   \exp(\bg 2\mu+2\rho_{Q_1} +\nu_{Q_1},H_0(g)\bd) \,dg.
\end{align*}
Hence it suffices to find $\mu$ as above such that the integral
\begin{align*}
  \int_{ \ago_{P_{0,r}} }\tau_{P_{0,r}}(H)   \exp(\bg 2\mu-2\rho^{Q_1}_{P_{0,r}}+2\nu_{Q_1},H\bd) \,dH
\end{align*}
is finite.  We must find $\mu\in \ago_{P_{0,r}}^{Q_1,*}$  so that for all $\al\in \Delta_{P_{0,r}}
^{G_r}$ we have $\bg \mu-\rho^{Q_1}_{ P_{0,r} } +\nu_{Q_1},\varpi_\al^\vee\bd <0$. If $\al\notin \Delta_0^{Q_1}$ we have $\varpi_\al^\vee\in \ago_{Q_1}$ and thus 
\begin{align*}
  \bg \mu-\rho^{Q_1}_{ P_{0,r} } +\nu_{Q_1},\varpi_\al^\vee\bd =\bg \nu_{Q_1},\varpi_\al^\vee\bd <0
\end{align*}
by the explicit formula for $\nu_{P_\pi}$, see § \ref{S:exposant-disc}. So it is enough  to take $\mu$ such that for  $\al\in \Delta_0^{Q_1}$ we have:
\begin{align*}
  \bg \mu,\varpi_\al^\vee\bd <\bg \rho^{Q_1}_{ P_{0,r} } -\nu_{Q_1},\varpi_\al^\vee\bd.
\end{align*}

\end{preuve}

Using lemma \ref{lem:first-bound}, we are reduced  to bound
\begin{align*}
  \int_{  P_{\tilde\xi} (\AAA)  \back G'(\AAA)} \exp(\bg \Re(\la)+\rho_P, H_P(\tilde\xi h)\bd)\,dh.
\end{align*}
We follow notations of the proof of lemma \ref{lem:first-bound}. We have a parabolic subgroup $Q$ containing $P_\xi$. By Iwasawa decomposition, we are reduced to bound
\begin{align*}
  \int_{  (P_{\tilde\xi}\cap M_{Q'}) (\AAA)  \back M_{Q'}(\AAA)} \exp(\bg \Re(\la)+\rho_P^Q, H_P(\tilde\xi h)\bd)\,dh.
\end{align*}
For $\la\in \ago_P^{G,*}$ such that for all $\al\in \Delta_P^Q$, the real part $\Re(\bg \la^Q,\al^\vee\bd )$  is large enough, the above integral is known to be convergent (we are easily  reduced to the case of  \cite[lemma 27 ]{JLR}).
\end{paragr}

\begin{paragr}[Intertwining periods.] --- \label{S:RwQ}   Let $Q\in \fc_2(M)$, see § \ref{S:fc2}. There exists a unique pair $(R,w)$ such that $w\cdot Q=R$ is a standard parabolic subgroup of $G$ and $w\in W(P;R)$. Let $L_w=wMw^{-1}$ be the standard Levi component of $R_w$, see  § \ref{S:rep}. Then $R\in \fc_2(L_w)$. Let $\xi=\xi_{L_w}^R\in W_2^{\mathrm{st}}(L_w)$, see § \ref{S:fc2}. Then for any $\varphi\in \Ac_{P,\pi}(G)$ and $\la\in \ago_{P,\CC}^{G,*}$ we set
  \begin{align}\label{eq:defJQ}
    J_Q(\varphi,\la)=J(\xi,M(w,\la)\varphi,(w\la)^R) 
  \end{align}
  where the right-hand side is the integral \eqref{eq:inter-period} written relatively to the standard Levi subgroup $L_w$ and the involution $\xi\in W_2^{\mathrm{st}}(L_w)$. The right-hand side makes sense a priori only for $\la$ in some cone given by proposition \ref{prop:convergence}. We shall show in corollary \ref{cor:period-entrelac} that the right-hand side in fact admits a meromorphic continuation to  $\ago_{P,\CC}^{G,*}$.

The intertwining period is closely related to the regularized period of some Eisenstein series. This is the content of the following theorem which generalizes to the case of discrete automorphic representations some of the fundamental results of Jacquet-Lapid-Rogawski (cf.  \cite{JLR}).

  \begin{theoreme}\label{thm:period-entrelac}
Let $R$ be a standard parabolic subgroup and  $w\in W(P;R)$ such that $Q=w^{-1}\cdot R\in \fc_2(M)$. For $\la$ in some open cone, we have
    \begin{align*}
      J_Q(\varphi,\la)=\pc^{R}(\varphi,\la,w)
    \end{align*}
where the left-hand side is defined by  \eqref{eq:defJQ} and the right-hand side is defined in § \ref{S:reg-period}.
  \end{theoreme}

  \begin{preuve}
    The case where $\pi$ is cuspidal is due to Jacquet-Lapid-Rogawski (see the  proof of lemma 32 in \cite{JLR}). By the definition \eqref{eq:defJQ} and the equality \eqref{eq:chgt-w} of remark   \ref{rq:descente-par}, we may and shall assume that $w=1$, the parabolic subgroup  $Q$ is equal to $R$ and thus  $J_Q(\varphi,\la)=J(\xi,\varphi,\la^R)$ with $\xi=\xi_M^R$. We shall use the construction given in § \ref{S:Jxitilde}. Here on can choose  $\tilde\xi\in M_R(F)$ so that $\tilde\xi\iota(\tilde\xi)^{-1}=\xi$. Then $N_{R'}=N_R\cap G'\subset N_{\tilde\xi}\subset R'$ and $M_{\tilde\xi}\subset R'$. By Iwasawa decomposition, one can deduce that in some cone we have:
    \begin{align*}
      J(\xi,\varphi,\la^R)=   \int_{   (N_{\tilde\xi}\cap M_{R'})(\AAA)  A_{M_{\tilde \xi}}^\infty  M_{\tilde\xi}(F) \back M_{R'}(\AAA)}\int_{K'}   \varphi_{\la-\rho_R}(\tilde\xi h k)\,dkdh.
    \end{align*}
    Using the parabolic descent \eqref{eq:par-desc}, we see that we are reduced to the case  $R=G$. If $M=G$ the theorem is an obvious consequence of the definitions, the right-hand side being explicitly described in proposition \ref{prop:reg-period} assertion 4. If $M\subset G$, we are in  the case $G=G_{2d}$ and $M=G_d\times G_d$ which is proved in the next section (see corollary \ref{cor:periode-carre-int}).
  \end{preuve}

\begin{corollaire}\label{cor:period-entrelac}
    The intertwining period  $J_Q(\varphi,\la)$  which was defined on some cone admits a meromorphic continuation to $\ago_{P,\CC}^{G,*}$ with hyperplane singularities.
  \end{corollaire}

  \begin{preuve}
    It is a direct consequence of theorem \ref{thm:period-entrelac} and  proposition \ref{prop:reg-period}.
  \end{preuve}

Let $R$ a standard parabolic subgroup and $w\in W(P;R)$. We have defined  $\pc^{T,R,m}(\varphi,\la,w)$ for $\la\in \ago_{P,\CC}^{G,*}$, see \eqref{eq:PTRm}. The following corollary gives an expression in terms of intertwining periods.

  \begin{corollaire}\label{cor:main-part}
Let $S=w^{-1}Rw$.    We have
    \begin{align*}
      \pc^{T,R,m}(\varphi,\la,w)=\sum_{Q\in \fc_2^S(M)}  2^{-\dim(\ago_Q^S)}  J_Q(\varphi,\la) \cdot \frac{\exp(\bg \la,T_Q^S\bd)}{\theta_Q^S(\la)}.                                    
    \end{align*}
  \end{corollaire}
  
  \begin{preuve}
    By definition we have
    \begin{align*}
        \pc^{T,R,m}(\varphi,\la,w)&=\sum_{ P_0\subset Q\subset R}   2^{-\dim(\ago_Q^R)}   \sum_{w'\in W^R(R_{w};Q)w} \pc^{Q}(\varphi,\la,w') \cdot \frac{\exp(\bg w'\la,T_Q^R\bd)}{\theta_Q^R(w'\la)}.      
    \end{align*}
    The map $(Q,w_1)\mapsto w_1^{-1}Qw_1$  induces a bijection from the disjoint union $\bigcup_{P_0\subset Q \subset R}W^R(R_w;Q)$ onto $\fc^{R}(wMw^{-1})$. Thus we get a bijection from the disjoint union $\bigcup_{P_0\subset Q \subset R}W^R(R_w;Q)w$ onto $\fc^{S}(M)$. Let $(Q,w')$  be in the source  and let $Q'\in \fc^S(M)$ be its  image. By proposition \ref{prop:vanishing}, the regularized period
 $\pc^{Q}(\varphi,\la,w')$ vanishes unless $Q\in \fc_2(wMw^{-1})$ and thus $Q'\in \fc_2(M)$. For such a parabolic subgroup we  have 
$$J_{Q'}(\varphi,\la)=\pc^{Q}(\varphi,\la,w').$$
The formula follows easily.
  \end{preuve}

\end{paragr}

\subsection{Computation of a regularized period}\label{ssec:comp-regper}

\begin{paragr} The main goal of this subsection is to prove theorem \ref{thm:periode-carre-int} and corollary \ref{cor:periode-carre-int} which was used in the proof of theorem \ref{thm:period-entrelac}. Thoughout this subsection, we set  $G=G_{2n}$. Let $P$ be a  standard  parabolic  subgroup of $G$ with standard Levi component $M=G_n\times G_n$. 
\end{paragr}

\begin{paragr} Let  $\pi\in\Pi_{\disc}(G_n)$. By M{\oe}glin-Waldspurger classification of the discrete spectrum reviewed in § \ref{S:exposant-disc}, one can attach to $\pi$  two integers $d,r  \geq 1$ with $dr=n$  and a cuspidal automorphic representation  $\sigma$ of $G_r$. 

Let $Q_L$ be the standard parabolic subgroup of $G$ of type the $2d$-tuple $(r,\ldots,r)$. Its standard Levi component, denoted by $L$, is naturally identified  with $G_r^{2d}$. The vector space $\ago_{L,\CC}^{*}$ is accordingly  identified with $\CC^{2d}$. Let
\begin{align*}
  \Lambda_L^M=(\frac{d-1}2, \frac{d-3}2,\ldots,-\frac{d-1}2,\frac{d-1}2, \ldots,-\frac{d-1}2)\in \ago_L^{P,*}
\end{align*}
and for any $\la=(\la_1,\ldots,\la_{2d})\in \ago_{L,\CC}^*$
\begin{align*}
  D_L^M(\la)=\prod_{\al\in \Delta_{Q_L}^P} (\bg \la,\al^\vee\bd-1).
\end{align*}

   Let $\varphi \in \Ac_{Q_L,\sigma^{\otimes d}\times \sigma^{*,\otimes d}}(G)$.  For any $\la\in \ago_{Q_L,\CC}^{G,*}$, one has the Eisenstein series $E^P(\varphi,\la)$ defined by analytic continuation from
  $$E^P(g,\varphi,\la)=\sum_{\delta\in Q_L(F)\back P(F)} \exp(\bg\la,H_P(\delta g)\bd) \varphi(\delta g)
  $$
for $g\in G(\AAA)$. Moreover, for any  $g\in G(\AAA)$, the map
\begin{align*}
  D_L^M(\la)\cdot E^P(g,\varphi,\la)
\end{align*}
is  holomorphic at  $\Lambda^M_L$. Let us denote by $E_{-1}^P(g,\varphi)$ its value at this point. Then  $E_{-1}^P(\varphi) \in \Ac_{P,\pi\otimes\pi^*}(G)$.
\end{paragr}

\begin{paragr}[Main results.] ---
    
  \begin{theoreme}
    \label{thm:periode-carre-int}
We have for all $\varphi \in \Ac_{Q_L,\sigma^{\otimes d}\times \sigma^{*,\otimes d}}(G)$ and all $\mu\in \ago_{P,\CC}^{G,*}$ in some cone
    \begin{align}\label{eq:periode-carre-int}
      J_G(E_{-1}^P(\varphi),\mu)=\pc^{G}(E_{-1}^P(\varphi),\mu).
    \end{align}
where the left-hand side is defined by \eqref{eq:defJQ} and the right-hand side in § \ref{S:reg-period}.
  \end{theoreme}

  \begin{preuve}
    The theorem is the combination of propositions \ref{prop:etape1} and \ref{prop:etape2} below whose proofs run over § \ref{S:preuve1}  to § \ref{S:preuvefinale}.
  \end{preuve}
  
\begin{corollaire} \label{cor:periode-carre-int}
    For all $\psi \in \Ac_{P,\pi\times \pi^{*}}(G)$ and all $\mu\in \ago_{P,\CC}^{G,*}$ in some cone we have
      \begin{align*}
        J_G(\psi,\mu)=\pc^G(\psi,\mu).
      \end{align*}
  \end{corollaire}

  \begin{preuve}
    This is an obvious consequence of theorem \ref{thm:periode-carre-int} : both members are continuous for  $\mu\in \ago_{P,\CC}^{G,*}$ in some cone (see proposition \ref{prop:convergence} and proposition \ref{prop:reg-period}) and coincides in the dense subspace generated by the functions $ E_{-1}^P(\varphi) $.
  \end{preuve}
\end{paragr}

\begin{paragr}\label{S:preuve1}
  By definition, the right-hand side in \eqref{eq:periode-carre-int} is the regularized period of the Eisenstein series $E(E_{-1}^P(\varphi),\mu)$ where  $\mu\in \ago_{P,\CC}^{G,*}$. We shall follow the same method as in \cite{Offen-symplectic} or \cite{Ysymp} to get a new expression for it.  It follows from proposition \ref{prop:MS-periods} and remark \ref{rq:MS-periods}, see also  \cite[proposition 8.4.1]{LR},  that, at least for generic $\mu\in \ago_{P,\CC}^{G,*}$, the integral
\begin{align}\label{eq:eis-tr}
  \int_{[G']_0} \Lambda^T_m     E(g,E_{-1}^P(\varphi),\mu)\,dg,
  \end{align}
is, as a function of $T$, a polynomial exponential whose purely polynomial term is in fact constant and equal to  $\pc^{G}(E_{-1}^P(\varphi),\mu)$. To compute this constant term, the starting point is the following lemma.

\begin{lemme}\label{lem:residu}
  \begin{align*}
    \int_{[G']_0} \Lambda^T_m     E(g,E_{-1}^P(\varphi),\mu)\,dg=\lim_{\la \to \Lambda^M_L}       D_L^M(\la)\cdot  \int_{[G']_0} \Lambda^T_m     E(g,\varphi,\la+\mu) \,dg.
  \end{align*}
\end{lemme}

\begin{preuve}
First one has
   \begin{align*}
    E(E_{-1}^P(\varphi),\mu)=\lim_{\la \to \Lambda^M_L}       D_L^M(\la)\cdot  E(\varphi,\la+\mu).
  \end{align*}
One has basically to invert the limit, the integration and the truncation. This can be easily justified by properties of the mixed truncation operator $\Lambda_m^T$ (see \cite{ar-truncated} proof of lemma 3.1 for a  similar result).
\end{preuve}

The second step is to use the following analogue of the ``Maa\ss-Selberg relations''.

\begin{proposition}(Jacquet-Lapid-Rogawski)
  We have for $\la\in \ago_{Q_L,\CC}^{G,*}$ in general position
  \begin{align}\label{eq:mscusp}
     \int_{[G']_0} \Lambda^T_m     E(g,\varphi,\la) \,dg= \sum_{Q \in \fc_2(L) } 2^{-\dim(\ago_Q^G)} J_Q(\varphi,\la) \frac{\exp(\bg \la, T_Q^G\bd)}{\theta_Q(\la)}.
  \end{align}
\end{proposition}

\begin{preuve}
  This is just a paraphrase of theorem 40 of \cite{JLR}. For the sake of clarity, let usw comment that the statement can be deduced from proposition \ref{prop:MS-periods}. Indeed since $\varphi$ is cuspidal the formula of proposition \ref{prop:MS-periods} reduces to the equality:
  $$\Ic^{T}(\varphi,\la)=\pc^{T,G,m}(\varphi,\la).$$
  Now by definition  the left-hand side of the formula above is the left-hand side of \eqref{eq:mscusp}. We can conclude with corollary \ref{cor:main-part} that $\pc^{T,G,m}(\varphi,\la)$ is the right-hand side of  formula \eqref{eq:mscusp}.  Of course, here we implicitly use theorem \ref{thm:period-entrelac} but only for cuspidal representations for which it is due to Jacquet-Lapid-Rogawski.
  \end{preuve}

We can compute the constant term  (as a function of $T$) in the limit  in lemma \ref{lem:residu} using the expression \eqref{eq:mscusp}. From lemma \ref{lem:py-exp-lisse},  this constant term  is obtained from the subsum restricted to the set of $Q\in\fc_2(L)$ such that $ (\Lambda_L^M+\mu)_Q^G=0$ for all $\mu\in \ago_P^{G,*}$. We shall provide a description of this. First we shall identify the group $W(L)$ with the symmetric group $\SG_{2d}$ (in $2d$ letters).   For any $\sigma\in \SG_d$, let $w_\sigma \in \SG_{2d}$ defined for any  $1\leq i\leq d$ by
\begin{align*}
  w_\sigma(i)=2\sigma(i)-1   \text{   and }w_\sigma(2d+1-i)=2\sigma(i).
\end{align*}
Let $R\subset G$ be the unique standard parabolic subgroup with standard Levi component isomorphic to $G_{2r}^{d}$. 

\begin{lemme}
  The map $\sigma\in\SG_d \mapsto w_\sigma^{-1}R$ induces a bijection from $\SG_d $ onto the subset of $Q\in \fc_2(M)$ such that $ (\Lambda_L^M+\mu)_Q^G=0$ for all $\mu\in \ago_P^{G,*}$.
\end{lemme}

\begin{preuve}
  Left to the reader.
\end{preuve}

Let $Q=w_\sigma^{-1}\cdot R$ with $\sigma\in\SG_d$. Note that $w_\sigma\in W(Q_L;R)$. Let $\xi=\xi_L^R\in W_2^{\mathrm{st}}(L)$. By the definition of $J_Q(\varphi,\la)$, see § \ref{S:RwQ}, we have:
\begin{align*}
  J_Q(\varphi,\la)=J(\xi, M(w_\sigma,\la)\varphi, (w_\sigma\la)^R)
\end{align*}
for $\la\in \ago_{L,\CC}^{G,*}$. Recall that the right-hand side is given by a convergent integral if $\Re(\bg \la,\al^\vee\bd)$ is large enough for all $\al\in \Delta_{Q_L}^R$. Thus  we get:

            \begin{lemme}\label{lem:lim-restr}
For $\mu\in \ago_{P,\CC}^{G,*}$, the regularized  period  $\pc^G(E_{-1}^P(\varphi),\mu)$ is equal to 
              \begin{align}\label{eq:lalim}
                2^{1-d} \lim_{\la \to \Lambda^M_L}     D_L^M(\la)  \sum_{\sigma\in \SG_d}  J(\xi, M(w_\sigma,\la+\mu)\varphi, (w_\sigma(\la+\mu))^R) \, \frac{1}{\theta_R ( w_\sigma(\la+\mu))}
              \end{align}
            \end{lemme}
          \end{paragr}

 \begin{paragr} For any $\sigma\in \SG_d$, we have 
$$D_{L}^{M}(\la)=D_{L,\sigma}^{M,+}(\la)D_{L,\sigma}^{M,-}(\la)$$
where we introduce
 $$D_{L,\sigma}^{M,\pm}(\la)= \prod_{\al\in \Delta_{Q_L}^P, \pm w_\sigma\al<0} (\bg \la, \al^\vee\bd -1).$$
Then the operator $D_{L,\sigma}^{M,+}(\la)M(w_\sigma,\la+\mu)$ has a limit at $\la=\Lambda^M_L$ denoted by
\begin{align*}
  M_{-1}(w_\sigma,\mu)= \lim_{\la \to \Lambda^M_L} D_{L,\sigma}^{M,+}(\la)M(w_\sigma,\la+\mu).
\end{align*}

However in the expression  \eqref{eq:lalim} of  lemma \ref{lem:lim-restr}, one cannot permute the limit and the sum. To get around this difficulty, we shall use a directional limit. To do this, we identify naturally $\ago_{L,\CC}^{*}$ with $\CC^{2d} $ and we take $\la=\Lambda_L^M + \eps x$ where $\eps\in \CC$ will go to  $0$ and $x=(x_1,\ldots,x_{2d})\in \CC^{2d}$ is a fixed point in general position. For any  $\sigma\in \SG_d$, we can write
\begin{align*}
              D_{L,\sigma}^{M,-}(\la)&  = \eps^{d-1}\prod_{i =1}^{d-1}   D^-_{\sigma,i}(x) 
            \end{align*}
           with
            \begin{align*}
               D_{\sigma,i}^-(x) =\left\lbrace
              \begin{array}{l}
                x_i-x_{i+1} \text{  if  } \sigma(i)<\sigma(i+1)\\
                x_{2d-i}-x_{2d+1-i} \text{ otherwise}.
              \end{array}
\right.
            \end{align*}
We also have
\begin{align*}
 \theta_R ( w_\sigma(\la+\mu)) = \vol(\ago_R^{G}/\ZZ(\Delta_R^\vee))^{-1} (\eps/2)^{d-1} \theta_{R,\sigma}(x)
\end{align*}
  with
  \begin{align*}
    \theta_{R,\sigma}(x)=\prod_{i=1}^{d-1}((x_{\sigma^{-1}(i)}-x_{\sigma^{-1}(i+1)})+(x_{2d+1-\sigma^{-1}(i)}-x_{2d+1-\sigma^{-1}(i+1)})).
  \end{align*}

  Using the continuity of $J(\xi,\varphi,\la)$ provided by proposition \ref{prop:convergence}, we get:

\begin{lemme} \label{lem:lim-restr2}
Let  $\al_0$ be the unique element of $\Delta_P$. For any $\mu\in \ago_P^{G,*}$  with  $\Re(\bg \mu,\al_0^\vee \bd)$ large, the  limit \eqref{eq:lalim} is equal to 
  \begin{align}\label{eq:lalim2} 
   \vol(\ago_R^{G}/\ZZ(\Delta_R^\vee)) \sum_{\sigma\in \SG_d}  J(\xi, M_{-1}(w_\sigma,\mu)\varphi, (w_\sigma(\La_L^M+\mu))^R) \, \frac{ \prod_{i =1}^{d-1}   D^-_{\sigma,i}(x)   }{\theta_{R,\sigma}(x)}.
  \end{align}
\end{lemme}

\begin{preuve}
  Let $\sigma\in \SG_d$ and let  $\tau=w_\sigma (\sigma^{\otimes d}\times \sigma^{*,\otimes d})$. By proposition \ref{prop:convergence}, for  $\Re(\bg \mu,\al_0^\vee \bd)$ is large enough, the integral that defines the map
  $$(\psi,\la)\in  \Ac_{Q_L}(G)\times \ago_{L,\CC}^{G,*} \mapsto J(\xi, \psi, (w_\sigma(\la+\mu))^R)$$
  is convergent and continuous  uniformly for   $\la$ in a compact neighborhood of $\Lambda_L^M$. One deduced that for $\la=\Lambda_L^M + \eps x$, one has:
  \begin{align*}
     \lim_{\eps\to 0} J(\xi, D_{L,\sigma}^{M,+}(\la) M(w_\sigma,\la+\mu)\varphi,  (w_\sigma(\la+\mu))^R) = J(\xi, M_{-1}(w_\sigma,\mu)\varphi, (w_\sigma(\La_L^M+\mu))^R).
  \end{align*}
  The lemma follows.
\end{preuve}
\end{paragr}

\begin{paragr}
  Of course, the expression \eqref{eq:lalim2} does not depend on a particular choice of $x$. It is also regular in a non-empty open subset of  the subspace defined by  $x_1=x_2=\ldots=x_d=0$.  Moreover on such a subset,  all the  terms in \eqref{eq:lalim2} vanish but the one corresponding  to the permutation $\sigma_0\in \SG_d$ that satisfies $\sigma_0(1)>\sigma_0(2)>\ldots>\sigma_0(d)$ that is  $\sigma_0(i)=d+1-i$ for $1\leq i \leq d$. We set
  $$w_0=w_{\sigma_0}.$$
Since we have 
$$\frac{ \prod_{i =1}^{d-1}   D^-_{\sigma_0,i}(x)   }{\theta_{R,\sigma_0}(x)}=1,$$
on such a subset. We observe that $w_0(\La_L^M+\mu)\in \ago_L^{R,*}$. We get:

\begin{proposition}\label{prop:etape1} For any $\mu\in \ago_P^{G,*}$  with large $\Re(\bg \mu,\al_0^\vee \bd)$, we have  
\begin{align*}
  \pc(E_{-1}^P(\varphi),\mu)  =   \vol(\ago_R^{G}/\ZZ(\Delta_R^\vee))    J(\xi, M_{-1}(w_0,\mu)\varphi, w_0(\La_L^M+\mu)).
\end{align*}
 for any $\mu\in \ago_P^{G,*}$ such that both sides are  regular at $\mu$.
\end{proposition}

\end{paragr}

\begin{paragr}
            Let $w_1\in W(L)$ be the element that corresponds to the permutation in  $\SG_{2d}$ defined by
$$w_1(i)=\left\lbrace
              \begin{array}{l}
               d+1-i \text{  if } 1\leq i \leq d\\
                i \text{ if }   d+1\leq i \leq 2d.
              \end{array}
\right.
$$

Let  $w_2=w_0w_1$. Since $w_1^2=1$, one also has $w_2w_1=w_0$. For  $1\leq i\leq d$ we have  $w_2(i)=2i-1$ and $w_2(2d+1-i)=2(d+1-i)$.  

Let $S=w_2^{-1}R w_2$. Then $S\in \fc_2(L)$ and $\xi_L^S=w_2^{-1}\xi w_2$. We can write $\xi_L^S$  as a product of transpositions with  disjoint supports:
$$\xi_L^S=(1 \,\,  d+1)\cdot (2\, \, d+2) \cdots (d \, \, 2d).
$$

Let
$$\nu=\begin{pmatrix}
  I_n & \sqrt{\tau} I_n\\ I_n& - \sqrt{\tau}I_n
\end{pmatrix}
$$
so that  $\nu \iota(\nu)^{-1}=\xi_L^S$.

Let $S_{\nu}=\nu^{-1}Q_L\nu\cap G'$, $L_\nu=\nu^{-1}L\nu\cap G'$,  $M_\nu=\nu^{-1}M \nu\cap G'$and $N_{\nu}=\nu^{-1}N_{Q_L}\nu\cap G'$. Note that $N_\nu$ is the unipotent radical of $S_\nu$ and $L_\nu $ is a Levi factor of $S_\nu$. Moreover $M_\nu$ is reductive and contains $S_\nu$. 

We define for any function $\psi\in \Ac_{Q_L}^0(G)$ which is cuspidal and $\la\in \ago_{L,\CC}^{S,*}$,
\begin{align}\label{eq:JxiLS}
    J(\xi_L^S,\psi,\la)=\int_{  N_{\nu} (\AAA) A_{L_{\nu}}^\infty  L_{\nu}(F) \back G'(\AAA)} \psi_\la(\nu h)\,dh.
  \end{align}
The next proposition shows that the integrals is convergent in some region and admits otherwise an analytic continuation.

\begin{proposition} \label{prop:uneeqfct}(Jacquet-Lapid-Rogawski) Let $\psi\in \Ac_{Q_L}^0(G)$ be a cuspidal automorphic form. For any $\la\in \ago_{L,\CC}^{S,*}$ such that $\Re(\bg \la,\al_0^\vee\bd)$ is large enough for all roots of $A_L$ in $N_P$ the integral in the right-hand side of \eqref{eq:JxiLS} is convergent and satisfies:
      \begin{align}\label{eq:JxiLS=Jxi}
        J(\xi_L^S,\psi,\la)= J(\xi, M(w_2, \la)\psi , w_2\la).
    \end{align}
  \end{proposition}

  \begin{preuve}
    This is a consequence of \cite[theorem 31 and proposition 34]{JLR}.
  \end{preuve}

One can check that
 \begin{align*}
   D_{L,\sigma_0}^{M,+}(\la)=\prod_{\al\in \Delta_{Q_L}^P, w_1 \al<0} (\bg \la, \al^\vee\bd -1).
 \end{align*}
In particular, we can define for $\mu\in \ago_P^{G,*}$
\begin{align*}
  M_{-1}(w_1,\mu)= \lim_{\la \to \Lambda^M_L} D_{L,\sigma_0}^{M,+}(\la)M(w_1,\la+\mu).
\end{align*}

By the functional equation of intertwining operators, we have 
$$ M(w_0, \la+\mu)=M(w_2w_1,\la+\mu)=M(w_2,w_1(\la+\mu)) M(w_1,\la).$$
By multiplying by the factor $D_{L,\sigma_0}^{M,+}(\la)$ and taking the limit $\la \to \Lambda^M_L$ we get
\begin{align*}
  M_{-1}(w_0,\mu)= M(w_2,w_1(\Lambda_L^M+\mu))  M_{-1}(w_1).
\end{align*}
where  $M_{-1}(w_1)=\lim_{\la\to \La_L^M} D_{L,\sigma_0}^{M,+}(\la)M(w_1,\la)$.
Note that $w_1(\La_L^M+\mu)$ belongs to the cone in $\ago_{L,\CC}^{S,*}$ introduced in proposition  \ref{prop:uneeqfct} as soon as  $\Re(\bg \mu,\al_0^\vee\bd)$ is large.  Thus by 
\eqref{eq:JxiLS=Jxi} we have:

\begin{align}
\label{eq:JxiJxiLS}       J(\xi, M_{-1}(w_0,\mu)\varphi, w_0(\Lambda_L^M+\mu)) &=  J(\xi, M(w_2,w_1\Lambda_L^M+\mu)   M_{-1}(w_1)\varphi, w_2 w_1(\Lambda_L^M+\mu))\\
\nonumber &= J(\xi_L^ S,  M_{-1}(w_1)\varphi   , w_1(\Lambda_L^M+\mu)).
    \end{align}

\begin{proposition}\label{prop:etape2} 
  For any $\mu\in \ago_P^{G,*}$ with large  $\Re(\bg \mu,\al_0^\vee\bd)$  we have 
\begin{align*}
  J(\xi, M_{-1}(w_0,\mu)\varphi, w_0(\Lambda_L^M+\mu)) =    \vol(\ago_R^{G}/\ZZ(\Delta_R^\vee))^{-1} J_G(E_{-1}^P(\varphi),\mu).
\end{align*}
\end{proposition}

\begin{preuve}
  Using \eqref{eq:JxiJxiLS}, we see that it amounts to show the equality:
  \begin{align*}
   J(\xi_L^ S,  M_{-1}(w_1)\varphi   , w_1(\Lambda_L^M+\mu))=    \vol(\ago_R^{G}/\ZZ(\Delta_R^\vee))^{-1} J_G(E_{-1}^P(\varphi),\mu).
  \end{align*}
  If $\Re(\bg \mu,\al_0^\vee\bd)$ is large enough, the left-hand side is given by the following convergent integral:
  \begin{align*}
    \int_{  N_{\nu} (\AAA) A_{L_{\nu}}^\infty  L_{\nu}(F) \back G'(\AAA)}   \exp(\bg w_1(\Lambda_L^M+\mu), H_{Q_L}(\nu g)\bd)  M_{-1}(w_1)\varphi_{}(\nu g)\,dg\\
    =    \int_{ M_{\nu}(\AAA) \back G'(\AAA) } 
    \int_{  N_{\nu} (\AAA) A_{L_{\nu}}^\infty  L_{\nu}(F) \back M_\nu(\AAA)} \exp(\bg w_1(\Lambda_L^M+\mu), H_{Q_L}(\nu hg)\bd)      M_{-1}(w_1)\varphi(\nu hg)\,dhdg.
  \end{align*}
  Note that $w_1\mu=\mu$ and $\bg \mu, H_{Q_L}(\nu hg)\bd=\bg \mu, H_{P}(\nu hg)\bd=\bg \mu, H_{P}(\nu g)\bd$. Thus, by proposition  \ref{prop:prod-scal} below, the last line is equal to (note that $w_1\mu=\mu$)

  \begin{align*}
    \vol(\ago_R^{G}/\ZZ(\Delta_R^\vee))^{-1}  \int_{ M_{\nu}(\AAA) \back G'(\AAA) } \exp(\bg \mu, H_{P}(\nu g)\bd)    \int_{ [M_\nu]_0} E_{-1}^P(\varphi,\nu hg)\, dh dg    \\
    =\vol(\ago_R^{G}/\ZZ(\Delta_R^\vee))^{-1}  \int_{ A_{M_\nu}^\infty M_{\nu}(F) \back G'(\AAA) } \exp(\bg \mu, H_{P}(\nu g)\bd) E_{-1}^P(\varphi,\nu g)\,  dg\\
    = \vol(\ago_R^{G}/\ZZ(\Delta_R^\vee))^{-1} J_G(E_{-1}^P(\varphi),\mu).
  \end{align*}
\end{preuve}
\end{paragr}

\begin{paragr}\label{S:preuvefinale}
  
  \begin{proposition}
    \label{prop:prod-scal}
For all $g\in G(\AAA)$, the integral
\begin{align}\label{eq:prod-scal}
 \int_{  N_{\nu} (\AAA) A_{L_{\nu}}^\infty  L_{\nu}(F) \back M_\nu(\AAA)} \exp(\bg w_1\Lambda_L^M, H_{Q_L}(\nu hg)\bd)   M_{-1}(w_1)\varphi(\nu hg)\,dh
\end{align}
is equal to
 \begin{align}\label{eq:prod-scal2}
  \vol(\ago_R^{G}/\ZZ(\Delta_R^\vee))^{-1}     \int_{ [M_\nu]_0} E_{-1}^P(\varphi,\nu hg)\, dh 
    \end{align}
  \end{proposition}

  \begin{preuve}
    We have $\nu M_\nu\nu^{-1}\subset  M$. This embedding  is identified with the embedding $\delta:G_n\to G_n\times G_n$ given by $ m\mapsto \delta(m)=(m, \iota (m))$. The parabolic subgroup  $\nu S_\nu\nu^{-1}$ of  $\nu M_\nu\nu^{-1}$ and its  Levi factor $\nu L_\nu\nu^{-1}$ are identified with the  corresponding subgroups   of $G_n$ respectively denoted by $Q_1$ and $L_1$. Thus the integrals \eqref{eq:prod-scal} and \eqref{eq:prod-scal2} can be written respectively:
\begin{align*}
  \int_{  [G_n]_{Q_1,0}}\exp(\bg w_1\Lambda_L^M, H_{Q_L}(\delta(h)\nu g)\bd)  M_{-1}(w_1,\mu)\varphi(\delta(h)\nu g)\,dh
\end{align*}
and
\begin{align}\label{eq:prod-scal3}
  \vol(\ago_R^{G}/\ZZ(\Delta_R^\vee))^{-1}   \int_{ [G_n]_0} E_{-1}^P(\varphi,\delta(h)\nu g)\, dh.
    \end{align}
We shall show the equality of the two expressions above. From now on we may and we shall replace $\nu g$ by $g$.

Let $T\in \ago_{P_{0,n},\RR}^{G_n}$. On the group $G_n$, we have the Arthur's truncation operator denoted by $\Lambda^T$, see section \ref{sec:disc}. For any function  $\psi$ on $M(\AAA)=G_n(\AAA)\times G_n(\AAA)$,  we denote by $\Lambda^T_1\psi$ the function on  $M(\AAA)$ we get by applying the  operator $\Lambda^T$ on the first variable. When $\psi=E_{-1}^P(\cdot g,\varphi)$, resp.  $\psi=E^P(\cdot g, \varphi,\la)$, the  function $\Lambda^T_1\psi$ evaluated at $\delta(m)$ (for $m\in G_n(\AAA)$) is denoted by $ \Lambda^T_1E_{-1}^P(\delta(m)g,\varphi)$, resp. $\Lambda^T_1E^P)(\delta(m)g,\varphi,\la)$. 

The integral in  \eqref{eq:prod-scal3} (with $\nu g$ replaced by $g$) is the constant term of the polynomial exponential  in the variable $T$ given by 
\begin{align}
   \nonumber \int_{[G_n]_0} \Lambda^T_1E_{-1}^P(\delta(h)g),\varphi) \, dh\\
\label{eq:petite-lim}
= \lim_{\la\to \Lambda_L^M} D_L^M(\la)\cdot  \int_{[G_n]_0  }  (\Lambda^T_1E^P)(\delta(m)g,\varphi,\la)\,dm.
  \end{align}
We  also denote by  $\delta$  the diagonal embedding of $\ago_{P_{0,n}}^{G_n}$ in $\ago_{P_{0,2n}}^{G_{2n}}$ and a dual projection $\delta^*$. For $\la\in \ago_{L,\CC}^*$, the theorem \ref{thm:inversion-calculee}   takes in our situation the following form (this can  also be deduced from \cite[lemma 15.2]{Ar-cours} which gives  an explicit expression of a truncated cuspidal Eisenstein series)
\begin{align*}
  \int_{[G_n]_0  }  (\Lambda^T_1E^P)(\delta(m)g,\varphi,\la)\,dm\\=\sum_{w \in W^{M}(L)}  
 \int_{[G_n]_{Q_1,0}}      \exp(\bg w\la,H_{Q_L}(k(h)g)\bd)   (M(w,\la)\varphi)(\delta(h)g) \, dh \cdot \frac{\exp(\bg w\la, \delta(T)\bd)  }{ \theta_{Q_1}(w\la)}
\end{align*}
where we choose $k\in K_n$ such that $hk^{-1}\in Q_1(\AAA)$ and we set $k(h)=\delta(k)$. Note that thanks to $\delta^*$ we view  $\theta_{Q_1}$  as a polynomial function on $\ago_{P_{0,2n}}^{G_{2n}}$. The  constant term in $T$ in the right-hand side of  \eqref{eq:petite-lim} corresponds in the above sum to the  contributions of the elements $w$ that satisfy $w\Lambda_L^M=0$ that is $w\in W^{G_n}(L_1)w_1$ where  $W^{G_n}(L_1)$ is identified with the diagonal subgroup of $W^M(L)$. By lemma \ref{lem:py-exp-lisse}, this constant term is given by:
\begin{align*}
  \lim_{\la\to \Lambda_L^M} D_L^M(\la)\sum_{w \in W^{G_n}(L_1)} \int_{[G_n]_{Q_1,0}} \exp(\bg ww_1\la, H_{Q_L}( k(h)g)\bd) (M(ww_1,\la)\varphi)(\delta(h)g) \, dh\frac{1  }{ \theta_{Q_1}(ww_1\la)}.
\end{align*}

As we did previously,  we take  $\la=\Lambda_L^M+\eps x$ with  $x=(x_1,\ldots, x_d,x_1',\ldots, x_d')$. If  $x$ is generic, one can take the  limit $\eps\to 0$ term by  term. We identify  $W^{G_n}(L_1)$ with $\SG_d$. Let us consider
\begin{align*}
  \theta_w(x)=\prod_{i=1}^{d-1}   (x_{w^{-1}(d+1-i)}-x_{w^{-1}(d-i)}+ x'_{w^{-1}(i)} -x'_{w^{-1}(i+1)}).
\end{align*}
We have $D_L^M(\la)=D_w^+(x) D_w^-(x)$ where we set  $D_w^{\pm}(x)= \prod_{i=1}^{d-1}D_{w,i}^{\pm}(x)$
 with
\begin{align*}
               D_{w,i}^\pm(x) =\left\lbrace
              \begin{array}{l}
                x_i-x_{i+1} \text{  if  } \pm(w(i)-w(i+1))>0\\
                  x_i'-x_{i+1}' \text{ otherwise}.
              \end{array}\right.
\end{align*}

We have 
\begin{align*}
  \theta_{S}(ww_1\la)=\vol(\ago_{Q_1}^{G_n}/\ZZ(\Delta_{Q_1}^\vee))^{-1} \eps^{d-1} \theta_w(x).
\end{align*}

If $x$ is generic, the above  limit is the product of  $\vol(\ago_{Q_1}^{G_n}/\ZZ(\Delta_{Q_1}^\vee))$ and 
\begin{align*}
 \sum_{w\in W^{G_1}(M_1)}    D_{w}^-(x)  \left(\int_{[G_n]_{Q_1,0}}  \exp(\bg ww_1\Lambda_L^M, H_{Q_L}(  k(h)g)\bd)
 (M_{-1}(ww_1)\varphi)(\delta(h)g) \, dh\right) \theta_w(x)^{-1}
\end{align*}
where $M_{-1}(ww_1)=\lim_{\eps\to 0} \eps^{d-1}   D_{w}^+(x) M(ww_1,\la)$.
We can go further and take  $x_1=x_2=\ldots=x_d=0$. Then all the terms vanish but the one associated to  $w=1$ which gives
\begin{align*}
  \vol(\ago_{Q_1}^{G_n}/\ZZ(\Delta_{Q_1}^\vee))\int_{[G_n]_{Q_1,0}}  \exp(\bg w_1\Lambda_L^M, H_{Q_L}(  k(h)g)\bd)
 (M_{-1}(w_1)\varphi)(\delta(h)g) \, dh.
\end{align*}

This finishes the proof since on one hand we have $\vol(\ago_{Q_1}^{G_n}/\ZZ(\Delta_{Q_1}^\vee))= \vol(\ago_R^{G}/\ZZ(\Delta_R^\vee))$ and on the other hand we have
\begin{align*}
  \bg w_1\Lambda_L^M, H_{Q_L}(  k(h)g)\bd=\bg w_1\Lambda_L^M, H_{Q_L}( \delta( h)g)\bd
\end{align*}
that is $\bg w_1\Lambda_L^M, H_{Q_L}(  \delta(h)g)\bd=0$ for any $h\in Q_1(\AAA)$.
\end{preuve}

\end{paragr}

\section{Relative characters}\label{sec:rel-char}

\subsection{A first example}\label{ssec:EX1}

\begin{paragr}\label{S:Jpair-char} We consider the situation and  notations of  sections  \ref{sec:FR} and \ref{sec:inter-period}. So the group $G$ is $\Res_{E/F}\GL_E(n)$ for some fixed quadratic extension $E/F$. We will also use the notations of section \ref{sec:disc} relative to this group.  Let  $P=MN_P$ be  a standard parabolic subgroup of $G$ and $\pi\in \Pi_{\disc}(M)$.

Let $J$ be a  level and $T$ be a truncation parameter.  Let $\tau\in \hat K_\infty$ and  $\bc_{P,\pi}(\tau,J)$ be an orthonormal basis for the Petersson norm  of the finite dimensional space $\Ac_{P,\pi}(G)^{\tau,J}$. Let $Q$ be a standard parabolic subgroup of $G$ and $w\in \, _QW_{P}$. For  $\la\in \ago_{P,\CC}^{G,*}$, $x\in G(\AAA)$ and $f,f'\in \Sc(G(\AAA))^J$ we set
  \begin{align*}
    \ec_{ (P,\pi,\tau)}^{T,Q}(x, f,f',\la,w)=   \sum_{\varphi\in \bc_{P,\pi}(\tau,J)} E(x, I_{P,\pi}(\la,f)\varphi,\la)    \overline{\Ic^{T,Q}( I_{P,\pi}(-\bar\la,f')\varphi,-\bar\la+\nu_w,w)} .
 \end{align*}

 \begin{remarque}
   The sum is finite and does not depend on the choice of an orthonormal basis.   As in remark \ref{rq:formeB}, we shall view $\ec$ as a sesquilinear form on  $\Sc(G(\AAA))$.
 \end{remarque}

 We also set
  \begin{align*}
    \ec_{ P,\pi}^{T,Q}(x, f,f',\la,w)=   \sum_{\tau\in  \hat K_\infty}    \ec_{ (P,\pi,\tau)}^{T,Q}(x, f,f',\la,w).
  \end{align*}

  \begin{proposition}\label{prop:maj-rel-bilchar1}
     There exist $l,N>0$ and for each standard parabolic subgroup $P$ of $G$ there exists $\ell_P$, a product of non-trivial real linear forms on $\ago_{P}^{G,*}$,      such that for all $q>0$ and all levels $J$ there exist $c>0$ and  a continuous semi-norm $\|\cdot\|_{\Sc}$  on $\Sc(G(\AAA))^J$ such that for all $J$-pairs $(P,\pi)$, all $f,f'\in \Sc(G(\AAA))^J$, all $\la\in \rc_{\pi,c,l}$, all $x\in G(\AAA)^1$, all standard parabolic subgroups $Q$ and all $w\in \, _QW_P$ we have
    \begin{align*}
      \sum_{\tau\in  \hat K_\infty}    \sum_{\varphi\in \bc_{P,\pi}(\tau,J)} |\ell_P(\la) E(x, I_{P,\pi}(\la,f_{})\varphi,\la)    \overline{\Ic^{T,Q}(I_{P,\pi}(-\bar\la,f'_{})\varphi,-\bar\la+\nu_w,w)}| \\
      \leq \frac{  \|x\|_{G}^N\|f\|_{\Sc}\|f'\|_{\Sc}    }{(1+\|\la\|^2)^q(1+\La_{\pi}^2)^{q}}.
    \end{align*}
    Moreover we can take $\ell_P=1$ if we restrict the statement to the elements $w\in W(P;Q)$.
  \end{proposition}

  \begin{remarque}\label{rq:holo}
    As it follows from the proof  we may assume that all the maps $E(x, \varphi,\la)$,  $ \ell_P(\la)  \overline{\Ic^{T,Q}(\varphi,-\bar \la+\nu_w,w)}$ and for  $w\in W(P;Q)$ the map $\overline{\Ic^{T,Q}(\varphi,-\bar \la,w)}$     are holomorphic on  $ \rc_{\pi,c,l}$.
  \end{remarque}

  \begin{preuve} We shall use theorem \ref{thm:maj-Eis} and theorem \ref{thm:maj-EisR} which provides the maps $\ell_P$.     Using Cauchy-Schwartz inequality it suffices to show that the following two expressions satisfy the conclusion of the theorem (with $N=0$ for the second expression)
        \begin{align}\label{eq:CS1}
       \sum_{\tau\in  \hat K_\infty}    \sum_{\varphi\in \bc_{P,\pi}(\tau,J)} |E(x, I_{P,\pi}(\la,f_{})\varphi,\la)  |^2
    \end{align}
    and
    \begin{align}\label{eq:CS2}
             \sum_{\tau\in  \hat K_\infty}    \sum_{\varphi\in \bc_{P,\pi}(\tau,J)} |\ell_P(\la)\Ic^{T,Q}(I_{P,\pi}(-\bar\la,f'_{})\varphi,-\bar\la+\nu_w,w)|^2
    \end{align}
    For \eqref{eq:CS1}, this is the content of theorem \ref{thm:maj-Eis}. To bound the expression \eqref{eq:CS2}, we shall use theorem \ref{thm:maj-EisR} in combination with proposition \ref{prop:cont-LaTm}  and the fact that the map $g\mapsto \|g\|_Q^{-N}$ is integrable for large $N$ on $[G']_{Q'}^1$. Thus there exists $l>0$  such that for all $q>0$ and all levels $J$ there exist $c>0$,  a finite family $(X_i)_{i\in I}$ of elements of $\uc(\ggo_\infty)$  and a continuous semi-norm $\|\cdot\|_{\Sc}$  on $\Sc(G(\AAA))^J$ such that for all $J$-pairs $(P,\pi)$, all $f,f'\in \Sc(G(\AAA))^J$, all $\la\in \rc_{\pi,c,l}$, all $x\in G(\AAA)^1$, all standard parabolic subgroups $Q$ and all $w\in \, _QW_P$ we have that \eqref{eq:CS2} is bounded by
    \begin{align*}
     \frac{  \sum_{i\in I}\|L(X_i) f\|_{\Sc}^2 }{(1+\|\la\|^2)^q(1+\La_{\pi}^2)^{q}}.
    \end{align*}
    The conclusion is clear. The last statement follows from the last statement of theorem \ref{thm:maj-EisR}.
    
  \end{preuve}
  
\end{paragr}

\begin{paragr}\label{S:ecT} With notations as above, we introduce also the linear form
  \begin{align*}
    \ec_{ (P,\pi,\tau)}^{T,Q}(x, f,\la,w)=   \sum_{\varphi\in \bc_{P,\pi}(\tau,J)} E(x, I_{P,\pi}(\la,f)\varphi,\la)    \overline{\Ic^{T,Q}(\varphi,-\bar\la+\nu_w,w)} .
 \end{align*}
and
  \begin{align*}
    \ec_{ P,\pi}^{T,Q}(x, f,\la,w)=   \sum_{\tau\in  \hat K_\infty}    \ec_{ (P,\pi,\tau)}^{T,Q}(x,f,\la,w).
  \end{align*}

  \begin{remarque}
    As before, the definition of $\ec_{ (P,\pi,\tau)}^{T,Q}(x, f,\la,w)$ does not depend on the choice of  $\bc_{P,\pi}(\tau,J)$ and we may view it as a linear form on $S(G(\AAA))$.
  \end{remarque}

If $Q=G$ then $w=1$ and we remove $G$ and $w$ from the notation: we simply write   $\ec_{ P,\pi}^{T}(x, f,\la)$ and $\ec_{ (P,\pi,\tau)}^{T}(x,f,\la)$.

  \begin{proposition}\label{prop:maj-rel-char1}  There exist $l,N>0$ and for each standard parabolic subgroup $P$ of $G$ there exists $\ell_P$, a product of non-trivial real linear forms on $\ago_{P}^{G,*}$,   such that for all $q>0$ and all levels $J$ there exist $c>0$ and  a continuous semi-norm $\|\cdot\|_{\Sc}$  on $\Sc(G(\AAA))^J$ such that for all $J$-pairs $(P,\pi)$, all $f\in \Sc(G(\AAA))^J$, all $\la\in \rc_{\pi,c,l}$, all $x\in G(\AAA)^1$, all standard parabolic subgroups $Q$ and all $w\in \, _QW_P$ we have 
    \begin{align*}
      \sum_{ \tau\in  \hat K_\infty}   |\ell_P(\la)\ec_{ (P,\pi,\tau)}^{T,Q}(x, f,\la,w)|   \leq \frac{  \|x\|_{G}^N\|f\|_{\Sc}  }{(1+\|\la\|^2)^q(1+\La_{\pi}^2)^{q}}
    \end{align*}
    Moreover $\la\mapsto \ell_P(\la) \ec_{ P,\pi}^{T,Q}(x, f,\la,w)$ is holomorphic on $\rc_{\pi,c,l}$.

 Both statements hold with $\ell_P=1$ if we restrict ourselves to elements $w\in W(P;Q)$.
  \end{proposition}

  \begin{preuve}
    As in the proof of theorem \ref{thm:maj-Eis}, for a level $J$ and an integer  $m\geq 1$ large enough, we can find $Z\in \uc(\ggo_\infty)$, $g_1\in \Cc(G(\AAA))$ and $g_2\in C_c^m(G(\AAA))$ such that
  \begin{itemize}
  \item $Z$ is  invariant under $K_\infty$-conjugation ;
  \item $g_1$ and $g_2$ are  invariant under $K$-conjugation and $J$-biinvariant;
  \item  for any $f\in \Sc(G(\AAA))^{J}$  we have:
    \begin{align*}
      f=f*g_1+ f*Z*g_2.
    \end{align*}
  \end{itemize}
  For $i=1,2$ we can observe that $I(\la, g_i)$ preserves the subspace $\Ac_{P,\pi}^{\tau,J}$ for any $\tau\in \hat K_\infty$. It follows that we have:
  \begin{align*}
    \ec_{ (P,\pi,\tau)}^{T,Q}(x, f,\la,w)=   \ec_{ (P,\pi,\tau)}^{T,Q}(x, f,\overline{g_1^\vee},\la,w)+  \ec_{ (P,\pi,\tau)}^{T,Q}(x, f*Z,\overline{g_2^\vee},\la,w).
  \end{align*}
  Then the proposition follows from proposition \ref{prop:maj-rel-bilchar1} and its mild extension to the subspace $C_c^m(G(\AAA))$ as long as $m$ is large enough.
  The holomorphy of  $ \ell_P(\la)\ec_{ (P,\pi,\tau)}^{T,Q}(x, f,\la,w)$ is clear on $\rc_{\pi,c,l}$ since it is a finite sum of product of holomorphic functions (see remark \ref{rq:holo}). The bound shows that the convergence of $\sum_\tau  \ell_P(\la)\ec_{ (P,\pi,\tau)}^{T,Q}(x, f,\la,w)$ is uniform on compact subsets of $\rc_{\pi,c,l}$  thus it is also holomorphic.

  In the same way the last statement follows from the last statement of proposition \ref{prop:maj-rel-bilchar1}.
  \end{preuve}
  
\end{paragr}

\subsection{A second example}\label{ssec:2nd-ex}

\begin{paragr}
  We keep  the notations of §  \ref{S:Jpair-char}. For any  standard parabolic subgroups $R$ and $w\in \, _RW_P$ we set
\begin{align}\label{eq:PcRpitau}
   \pc^R_{(P,\pi,\tau)}(x,f,\la,w)=\sum_{\varphi\in \bc_{P,\pi}(\tau,J)}  E(x,I_P(\la,f)\varphi,\la) \overline{\pc^R(\varphi,-\bar\la+\nu_w,w)}.
      \end{align}
We also set 
\begin{align}\label{eq:PTRmpitau}
&        \pc^{T,m}_{(P,\pi,\tau)}(x,f,\la)=         \sum_{\varphi\in \bc_{P,\pi}(\tau,J)}  E(x,I_P(\la,f)\varphi,\la) \overline{\pc^{T,m}(\varphi,-\bar\la)}      ;                   \\
 &\nonumber   \pc^R_{P,\pi}(x,f,\la,w)=   \sum_{ \tau\in  \hat K_\infty}  \pc^R_{(P,\pi,\tau)}(x,f,\la,w);\\
 & \nonumber   \pc^{T,m}_{P,\pi}(x,f,\la)=   \sum_{ \tau\in  \hat K_\infty}   \pc^{T,m}_{(P,\pi,\tau)}(x,f,\la).
\end{align}

      \begin{proposition}\label{prop:maj-rel-char3} Let $1>\eta>0$.  
    There exist $l, N ,\eps>0$ such that for all $q>0$ and all levels $J$ there exist $c>0$ and  a continuous semi-norm $\|\cdot\|_{\Sc}$  on $\Sc(G(\AAA))^J$ such that for all $J$-pairs $(P,\pi)$, all $f\in \Sc(G(\AAA))^J$, all $\la\in \rc_{\pi,c,l}$, all $x\in G(\AAA)^1$, all standard parabolic subgroups $R$ and all $w\in W(P;R)$ we have 
    \begin{align*}
      \sum_{ \tau\in  \hat K_\infty}   |\pc^R_{ (P,\pi,\tau)}(x, f,\la,w)|   \leq \frac{  \|x\|_{G}^N\|f\|_{\Sc}  }{(1+\|\la\|^2)^q(1+\La_{\pi}^2)^{q}}\\
      \sum_{ \tau\in  \hat K_\infty}   |\ec_{ (P,\pi,\tau)}^{T}(x,f,\la)- \pc^{T,m}_{ (P,\pi,\tau)}(x, f,\la)|   \leq \frac{  \exp(-\eps \| T\| )\|x\|_{G}^N\|f\|_{\Sc}  }{(1+\|\la\|^2)^q(1+\La_{\pi}^2)^{q}}
          \end{align*}
for all $T$ such that $d(T)\geq \eta \|T\|$.
    Moreover the  expressions $\ec_{ P,\pi}^{T}(x, f,\la)$, $ \pc_{ P,\pi}^{T,m}(x, f,\la,w)$ and $\pc^R_{ P,\pi}(x, f,\la,w)$ are  holomorphic on $\rc_{\pi,c,l}$.
      \end{proposition}

      \begin{preuve} We fix  $l$ and $c$ such that the maps in remark \ref{rq:holo}, among them  the Eisenstein series  $E(x,I_P(\la,f)\varphi,\la)$,  are holomorphic on   $\rc_{\pi,c,l}$. Let $R$ be a standard parabolic subgroup and $w\in W(P;R)$. By theorem \ref{thm:singularities}, the possible singularities of  $\overline{\pc^R(\varphi,-\bar \la,w)}$ on $\rc_{\pi,c,l}$ are simple and along the hyperplanes $\bg w \la,\al^\vee\bd=0$ where $\al\in \Delta_{R_w}^R$ and for such a root $\al$, the elementary symmetry $s_\al$ is such that $M(s_\al,0)$ acts by $-1$ on  $\Ac_{R_w,w\pi}(G)$. But by the functional equation of Eisenstein series, we have:
    \begin{align*}
      E(x,\varphi,\la)&=E(x,M(s_\al w,\la)\varphi,s_\al w\la).
    \end{align*}
For $\la\in \rc_{\pi,c,l}$ such that $\bg w \la,\al^\vee\bd=0$ we have $s_\al w\la=w\la$ and we get:
 \begin{align*}
      E(x,\varphi,\la)&=E(x,M(s_\al ,0) M(w,\la)\varphi,s_\al w\la)\\
&=-E(x, M(w,\la)\varphi, w\la)\\
&=- E(x,\varphi, \la)
    \end{align*}
and thus $E(x,\varphi,\la)$ vanishes on the hyperplane $\bg w \la,\al^\vee\bd=0$. Thus the possible singularities cancel with the zeroes of the Eisenstein series. We deduce that $\pc^R_{ (P,\pi,\tau)}(x, f,\la,w)$ is holomorphic. We want to bound it. For this we appeal to the very definition \eqref{eq:PTRc} of the regularized period to write:
\begin{align}\label{eq:NEW expr pcR}
  \pc^R_{ (P,\pi,\tau)}(x, f,\la,w)&=\sum_{ P_0\subset Q\subset R} (-2)^{-\dim(\ago_Q^R)}  \sum_{w ' \in \,_QW^R_{R_w}  w } \ec_{ (P,\pi,\tau)}^{T,Q}(x,f,\la,w') \cdot \frac{\exp(\bg w'(-\la+\nu_{{w'}}^w),T_Q^R\bd)}{\hat{\theta}_Q^R(w'(-\la+\nu_{{w'}}^w))}
\end{align}
where $T$ is an arbitrary truncation parameter. Since we know the holomorphy of $\pc^R_{ (P,\pi,\tau)}(x, f,\la,w)$ we can use the same argument as in the proof of proposition \ref{prop:bound-scalaire}. Indeed it suffices to bound $\pc^R_{ (P,\pi,\tau)}(x, f,\la,w)$ by the derivatives of the product of $\pc^R_{ (P,\pi,\tau)}(x, f,\la,w)$ and  some linear forms. To do this, we start from the right-hand side of \eqref{eq:NEW expr pcR}, we use  proposition \ref{prop:maj-rel-char1} and Cauchy formula.

We have observed in proposition \ref{prop:maj-rel-char1} that $\ec_{ P,\pi}^{T}(x, f,\la)$ is holomorphic on $\rc_{\pi,c,l}$. This is also the case of  $\pc_{ P,\pi}^{T,m}(x, f,\la)$ because it can be identified with the summand of exponent of ``type 1'' of $\ec_{ P,\pi}^{T}(x, f,\la)$ (see lemma \ref{lem:types} and proposition \ref{prop:type1}). Thus the difference $\ec_{ (P,\pi,\tau)}^{T}(x, f,\la)-\pc_{ (P,\pi,\tau)}^{T,m}(x, f,\la)$ is holomorphic and is equal to (by proposition \ref{prop:MS-periods})
\begin{align*}
  \sum_{ P_0\subset Q\subset G} 2^{-\dim(\ago_Q^G)}  \sum_{w  \in \,_QW^G_{P}\setminus W(P;Q)  } \pc_{ (P,\pi,\tau)}^{Q}(x,f,\la,w) \cdot \frac{\exp(\bg w(-\la+\nu_{{w}}),T_Q\bd)}{\theta_Q^G(w(-\la+\nu_{{w}}))}
\end{align*}
By the same techniques as before, the result  follows from proposition \ref{prop:maj-rel-char1} and lemma \ref{lem:types} assertion 2 (of course we can moreover assume that  $c$ satisfies $c<c_1 \eta$ where $c_1$ is the constant that appears in lemma \ref{lem:types} assertion 2).
       \end{preuve}
\end{paragr}

\section{Spectral expansion}\label{sec:Spec-exp}

\subsection{Contribution of a spectral datum as a limit}\label{ssec:asalim}

\begin{paragr}\label{S:Pichi} Let $\Xgo(G)$ be a the set of cuspidal data of $G$, namely the set of equivalence classes of pairs $(M,\pi)$ where $M$ is a standard Levi subgroup of $G$ and $\pi\in \Pi_{\cusp}(M)$. Two data $(M,\pi)$ and $(M',\pi')$  are equivalent if there exists $w\in  W(M,M')$ that sends $\pi$ to $\pi'$. 

Let $P=MN_P$ be a standard parabolic subgroup of $G$. We have the coarse Langlands decomposition 
\begin{align}
  L^2([G]_{P,0})=\hat\oplus_{\chi\in \Xgo(G)} L^2_\chi([G]_{P,0}).
\end{align}
Let $ \Ac_{P,\chi,\disc}(G)$ be the closed subspace of $\Ac_{P,\disc}^0(G)$ generated by the functions whose class belongs to $L^2_\chi([G]_{P,0})$.  By the Langlands construction of the spectral decomposition, see e.g. \cite{MWlivre}, and the description of the discrete spectrum obtained in \cite{MW}, there exists  a finite subset $\Pi_\chi(M)\subset\Pi_{\disc}(M)$ such that we have an isotypical decomposition 
$$ \Ac_{P,\chi,\disc}(G)=\bigoplus_{\pi\in \Pi_\chi(M)} \Ac_{P,\pi}(G).$$
\end{paragr}

\begin{paragr} Let $\chi\in \Xgo(G)$ be a cuspidal datum.  Let $f\in \Sc(\AAA)$ and $K_{f,\chi}^0$ be the kernel associated to the operator given by right convolution by $f$ on $L^2_\chi([G]_{0})$. By \cite[lemma 2.10.1.1]{BCZ}, there exists $N_0>0$ such that for all $N>0$ there exists a semi-norm $\|\cdot\|$ on $\Sc(G(\AAA))$ such that  for all $x,y\in G(\AAA)^1$ we have $|K_{f,\chi}^0(x,y)|\leq \|f\| \|x\|^{N+N_0}  \|y\|^{-N_0} $. In particular the following integral 
  \begin{align}\label{eq:int-Kchi}
J_\chi(x,f)=  \int_{[G']_0} K_{f,\chi}^0(x,y)\,dy.
  \end{align}
  is absolutely convergent for all $x\in G(\AAA)$. 

  Let $T$ be a truncation parameter. To get the spectral expansion of $J_\chi(x,f)$ it will be easier to first consider the spectral expansion of the truncated variant of $J_\chi(x,f)$ namely:
   \begin{align}\label{eq:int-Kchi-T}
    J_\chi^T(x,f)=    \int_{[G']_0} (K_{f,\chi}^0\La^T_m)(x,y)\,dy.
   \end{align}
   It is also absolutely convergent because of the properties of the kernel and the mixed truncation operators, see proposition \ref{prop:cont-LaTm}. The notation $\lim_{T\to +\infty}$ means the limit when $d(T)\to +\infty$. We have:
   
   \begin{proposition}\label{prop:lim-JT}\cite[Proposition 4.3.4.1]{BCZ} For all $x\in G(\AAA)$,
$$     \lim_{T\to +\infty}   J_\chi^T(x,f)=  J_\chi(x,f).$$
   \end{proposition}
   \end{paragr}

   \begin{paragr} We now recall the spectral expansion of  $J_\chi^T(x,f)$ in terms of the relative character $\ec_{P,\pi}^T(x,f,\la)$ defined in § \ref{S:ecT}.

     \begin{proposition}\label{prop:cv-spectral-T} For all $x\in G(\AAA)$,
       \begin{align*}
         J_\chi^T(x,f)=  \sum_{P_0\subset P} |\pc(M)|^{-1} \int_{i\ago_P^{G,*}}  \sum_{\pi \in \Pi_{\chi}(M)}   \ec_{P,\pi}^T(x,f,\la) \, d\la 
       \end{align*}
       where the right-hand side is absolutely convergent.
     \end{proposition}

     \begin{preuve}
       The equality is stated in \cite[Proposition 4.2.3.3]{BCZ}. Without appealing to \cite{BCZ}, we observe that proposition \ref{prop:maj-rel-char1} shows that the right-hand side is absolutely convergent since  the set $\Pi_{\chi}(M)$ is finite.
     \end{preuve}

     \begin{proposition} 
       \label{prop:limT}For all $x\in G(\AAA)$, the expression
       \begin{align*}
               \sum_{P_0\subset P} |\pc(M)|^{-1} \int_{i\ago_P^{G,*}}  \sum_{\pi \in \Pi_{\chi}(M)}   \pc_{P,\pi}^{T,m}(x,f,\la) \, d\la.
       \end{align*}
       is absolutely convergent.  For any $\eta>0$,  the limit when $T\to+\infty$ of the expression above for $T$ such that $d(T)\geq \eta\|T\|$  is $   J_\chi(x,f)$.
     \end{proposition}

     \begin{preuve}
       We fix $\eta>0$ and a level $J$ such that $f\in \Sc(G(\AAA))^J$.   By proposition \ref{prop:maj-rel-char3},  we get $N,\eps>0$  and for any (large) $q$ a continuous semi-norm $\|\cdot\|_{\Sc}$ on $\Sc(G(\AAA))^J$ such that
       \begin{align*}
         \sum_{P_0\subset P} |\pc(M)|^{-1} \int_{i\ago_P^{G,*}}  \sum_{\pi \in \Pi_{\chi}(M)}   |\ec_{ P,\pi}^{T}(x,f,\la)- \pc^{T,m}_{ P,\pi}(x, f,\la)|\, d\la\\
        \leq     \exp(-\eps \| T\| )\|x\|_{G}^N\|f\|_{\Sc}   \sum_{P_0\subset P} |\pc(M)|^{-1} \left(\sum_{\pi \in \Pi_{\chi}(M)} (1+\La_{\pi}^2)^{-q}\right) \int_{i\ago_P^{G,*}}(1+\|\la\|^2)^{-q}\, d\la
       \end{align*}
       for all $T$ such that $d(T)\geq \eta\|T\|$.

       The two statements of the theorem are then easily deduced since the set $\Pi_{\chi}(M)$ is finite and since we can choose $q$ large enough so that all the integrals converge: the  absolute convergence follows from that of proposition \ref{prop:cv-spectral-T} and the computation of the limit follows from proposition \ref{prop:lim-JT}.

     \end{preuve}
\end{paragr}

\subsection{Computation of the limit} \label{ssec:comput-limit}

\begin{paragr} We continue with the notations of the previous section.  Let $\tau\in \hat K_\infty$. Let $x\in G(\AAA)$. Using corollary \ref{cor:main-part}, we get the following  expression for the relative character $\pc^{T,m}_{ (P,\pi,\tau)}(x, f,\la)$ in terms of the intertwining periods:
  \begin{align}
 \nonumber   \pc^{T,m}_{(P,\pi,\tau)}(x,f,\la)&=         \sum_{\varphi\in \bc_{P,\pi}(\tau,J)}  E(x,I_P(\la,f)\varphi,\la) \overline{\pc^{T,m}(\varphi,-\bar\la)}\\
    \label{eq:exp-JQtau}             &=\sum_{Q\in \fc_2(M)}  2^{-\dim(\ago_Q^G)}    \Jc_{(Q,\pi,\tau)}(x,f,\la) \frac{\exp(-\bg \la,T_Q\bd)}{\theta_Q^G(-\la)},
      \end{align}
    where for $Q\in \fc_2(M)$ we introduce the relative character
 \begin{align}\label{eq:JcQpitau}
   \Jc_{(Q,\pi,\tau)}(x,f,\la)= \sum_{\varphi\in \bc_{P,\pi}(\tau,J)}  E(x,I_P(\la,f)\varphi,\la) \overline{ J_Q(\varphi,-\bar \la)} .
 \end{align}

 Let $Q\in \fc_2(M)$. There exist $R$ a standard parabolic subgroup and  $w\in W(P;R)$ such that $Q=w^{-1}Rw$. Then by theorem \ref{thm:period-entrelac} we have
 \begin{align}\label{eq:JcQ=pc}
   \Jc_{(Q,\pi,\tau)}(x,f,\la)= \pc^{R}_{(P,\pi,\tau)}(\varphi,\la,w).
 \end{align}
 Thus $\Jc_{(Q,\pi,\tau)}(x,f,\la)$ inherits all the properties of $\pc^{R}_{(Q,\pi,\tau)}(\varphi,\la,w)$ given in subsection \ref{ssec:2nd-ex}. We can set:
 \begin{align}\label{eq:JcQpi}
   \Jc_{Q,\pi}(x,f,\la)=\sum_{\tau\in \hat K_\infty}   \Jc_{(Q,\pi,\tau)}(x,f,\la).
 \end{align}
 Then we have
 \begin{align}  \label{eq:exp-JQ}
      \pc^{T,m}_{P,\pi}(x,f,\la)&=\sum_{Q\in \fc_2(M)}  2^{-\dim(\ago_Q^G)}    \Jc_{Q,\pi}(x,f,\la) \frac{\exp(-\bg \la,T_Q\bd)}{\theta_Q^G(-\la)}
    \end{align}
\end{paragr}

\begin{paragr}[A $(G,M)$-family of relative characters.] ---

  \begin{proposition}\label{prop:GMrelchar}
    Let $L\in \lc_2(M)$. For any $f\in \Sc(G(\AAA))$, $x\in G(\AAA)$ and $\la_0\in i\ago_{M}^{L,*}$, the family $(\Jc_{Q,\pi}(x,f,\la_0+ \la))_{Q\in \pc(L)}$ is a $(G,L)$-family of Schwartz functions of the variable $\la\in i\ago_{L}^{G,*}$.
  \end{proposition}

  \begin{preuve}
We may assume $x\in G(\AAA)^1$.    We fix a level $J$ such that  $f\in \Sc(G(\AAA))^J$. According to \eqref{eq:JcQ=pc} and proposition \ref{prop:maj-rel-char3}, there exists $l>0$ such that for all $q>0$ there exists $c>0$ and $C>0$ such that we have:
    \begin{align*}
      \sum_{ \tau\in  \hat K_\infty}   |\Jc_{(Q,\pi,\tau)}(x,f,\la)|   \leq \frac{  C  }{(1+\|\la\|^2)^q}
\end{align*}
for all $\la\in \rc_{\pi,c,l}$ and $Q\in \fc_2(M)$. The function $\Jc_{Q,\pi}(x,f,\la)$ is holomorphic on  $\rc_{\pi,c,l}$ thus smooth on $i\ago_P^{G,*}$. By the equality above it is also rapidly decreasing. Using Cauchy formula, all its real derivatives are also rapidly decreasing on the imaginary axis. Thus $\Jc_{Q,\pi}(x,f,\la)$ is a Schwartz function.
On the other hand we can assume that   $\pc^{T,m}_{P,\pi}(x,f,\la)$ is also holomorphic on $\rc_{\pi,c,l}$, still by proposition \ref{prop:maj-rel-char3}.  To conclude we apply   proposition \ref{prop:GMfam} to the smooth expression given by the right-hand side of \eqref{eq:exp-JQ}.
  \end{preuve}
\end{paragr}

\begin{paragr} \label{S:relch-JL}As a consequence of the definition of a $(G,L)$-family we have:

  \begin{corollaire}
    Let $L\in \lc_2(M)$ and $\la\in i\ago_M^{L,*}$. Then $\Jc_{Q,\pi}(x,f,\la)$ does not depend on the choice of $Q\in \pc(L)$.
  \end{corollaire}

As a consequence, we set
\begin{align}\label{eq:relch-JL}
  \Jc_{L,\pi}(x,f,\la)=\Jc_{Q,\pi}(x,f,\la),\ \ \la\in i\ago_M^{L,*}, Q\in \pc(L).
\end{align}

\begin{proposition}\label{prop:relch-JL}
  There exists $N>0$ such that for all $q>0$ there exists a continuous semi-norm $\|\cdot\|_{\Sc}$ on $\Sc(G(\AAA))$ such that for all $f\in \Sc(G(\AAA))$, $x\in G(\AAA)^1$, for all standard parabolic subgroup $P=MN_P$, all $\pi\in \Pi_{\disc}(M)$ and all $L\in \lc_2(M)$ we have:
  \begin{align*}
    \int_{i\ago_M^{L,*}} |\Jc_{L,\pi}(x,f,\la)|\, d\la \leq \frac{\|x\|_G^N \|f\|_{\Sc}}{(1+\La_\pi^2)^q}
  \end{align*}
\end{proposition}

\begin{preuve}
  According to proposition \ref{prop:maj-rel-char3}, the equality \eqref{eq:JcQ=pc} and the definition \eqref{eq:relch-JL}, we see that for all $q>0$ and all levels $J$ there exists  a continuous semi-norm $\|\cdot\|_{\Sc}$ on $\Sc(G(\AAA))^J$  such that for all $f\in \Sc(G(\AAA))^J$, $x\in G(\AAA)^1$ and  pairs $ (P,\pi)$, all $L\in \lc_2(M)$ and all $\la\in i\ago_M^{L,*}$
  \begin{align*}
    |\Jc_{L,\pi}(x,f,\la)|\leq \frac{\|x\|_G^N \|f\|_{\Sc}}{  (1+\|\la\|^2)^q(1+\La_\pi^2)^q}
  \end{align*}
The results follows.
\end{preuve}
\end{paragr}

\begin{paragr}

  \begin{theoreme}\label{thm:Jchi}
    For all $f\in \Sc(G(\AAA))$ and $x\in G(\AAA)$ we have
    \begin{align*}
      J_\chi(x,f)= \sum_{P_0\subset P} |\pc(M)|^{-1}   \sum_{L\in \lc_2(M)}  2^{-\dim(\ago_L^G)} \sum_{\pi \in \Pi_{\chi}(M)} \int_{i\ago_M^{L,*}} \Jc_{L,\pi}(x,f,\la)\, d\la
    \end{align*}
where the right-hand side is absolutely convergent and $M$ stands for the standard Levi factor of $P$.
      \end{theoreme}

      \begin{preuve}
        Since the sums are finite, the absolute convergence is that of the inner integral which follows from proposition \ref{prop:relch-JL}. By proposition \ref{prop:limT},  the   expression $J_\chi(x,f)$ is the limit when $T\to+\infty$  of the sum over $P_0\subset P=MN_P$ and $\pi \in \Pi_{\chi}(M)$ of the constant  $|\pc(M)|^{-1}  $ times the expression
        $$\int_{i\ago_P^{G,*}}    \pc_{P,\pi}^{T,m}(x,f,\la) \, d\la.$$
        Let fix $(P,\pi)$. Using the expression \eqref{eq:exp-JQ} of $\pc^{T,m}_{P,\pi}$, we can fix also $L\in \lc_2(M)$ and compute for the limit when $T\to+\infty$ of
        \begin{align*}
          2^{-\dim(\ago_L^G)} \int_{i\ago_{P}^{G,*}} \sum_{Q\in \pc(L)}     \Jc_{Q,\pi}(x,f,\la) \frac{\exp(-\bg \la,T_Q\bd)}{\theta_Q^G(-\la)}\, d\la.
        \end{align*}
       Since    $(\Jc_{(Q,\pi,\tau)}(x,f,\la))_{Q\in \pc(L)}$ is a Schwartz $(G,L)$-family by proposition \ref{prop:GMrelchar}, not only the integrand is also a Schwartz function, see \cite[corollary 3]{Linner}, but also, see \cite[lemma 8]{Linner}, the limit exists and is equal to: 
        \begin{align*}
          2^{-\dim(\ago_L^G)}       \int_{i\ago_{P}^{L,*}}   \Jc_{L,\pi}(x,f,\la)\, d\la.
        \end{align*}
      \end{preuve}

\end{paragr}

\subsection{Spectral decomposition of the Flicker-Rallis period of the automorphic kernel}\label{ssec:spec-FR}

\begin{paragr}
  Let $f\in \Sc(G(\AAA))$. The automorphic kernel associated to $f$ is the kernel of the operator given by right convolution by $f$ on $L^2([G])$ namely
  \begin{align*}
    K_f(x,y)=\sum_{\gamma \in G(F)} f(x^{-1}\gamma y), \ \ \ x,y\in G(\AAA).
  \end{align*}
We set for $x\in G(\AAA)$
\begin{align}\label{eq:Jf}
  J^G(x,f)=\int_{[G']} K_f(x,y)\, dy.
\end{align}

\begin{lemme}\label{lem:Jf}
  \begin{enumerate}
  \item   The integral that defines $J^G(x,f)$ is absolutely convergent and the map  $f\mapsto J(x,f)$ is a continuous distribution on $\Sc(G(\AAA))$.
  \item We have also:
    \begin{align}\label{eq:Jf0}
      J^G(x,f)=\frac12 \sum_{\chi\in \Xgo(G)} J_\chi(x,f)   
    \end{align}
where the right-hand side is absolutely convergent.
\end{enumerate}
\end{lemme}

\begin{preuve}
  The assertion 1 comes from \cite[lemma 2.10.1.1]{BCZ}. Let $K_f^0(x,y)$ be the  kernel of the operator given by right convolution by $f$ on $L^2([G]_0)$. Then we have:
  \begin{align*}
    K_f^0(x,y)=\int_{A_G^\infty}   K_f(ax,y)\,da.
  \end{align*}
  Taking into account the discrepancy between the measures on $A_G^\infty$ and $A_{G'}^\infty$, we get:
  \begin{align*}
     J^G(x,f)=\frac12 \int_{[G']_0} K_f^0(x,y)\, dy.
  \end{align*}
However we have the (coarse) spectral decomposition 
$$K_f^0(x,y)=\sum_{\chi\in \Xgo(G)} K_{f,\chi}^0(x,y).
$$
Integrating term by term over $y\in [G']_0$ which is possible by \cite[lemma 2.10.1.1]{BCZ} we get assertion 2.
\end{preuve}
\end{paragr}

\begin{paragr} We can now state and prove our main theorem.

  \begin{theoreme}
    \label{thm:spec-kernel}
 For any $f\in \Sc(G(\AAA))$ and $x\in G(\AAA)$ we have
 \begin{align*}
      J(x,f)= \sum_{P_0\subset P} |\pc(M)|^{-1}   \sum_{L\in \lc_2(M)}  2^{-\dim(\ago_L)} \sum_{\pi \in \Pi_{\disc}(M)} \int_{i\ago_M^{L,*}} \Jc_{L,\pi}(x,f,\la)\, d\la
    \end{align*}
where the right-hand side is absolutely convergent.
  \end{theoreme}

  \begin{preuve}
Using lemma \ref{lem:Jf} and  then theorem \ref{thm:Jchi} we get:
\begin{align*}
   J(x,f)&= \frac12 \sum_{\chi\in \Xgo(G)} J_\chi(x,f)\\
&= \sum_{\chi\in \Xgo(G)} \sum_{P_0\subset P} |\pc(M)|^{-1}   \sum_{L\in \lc_2(M)}  2^{-\dim(\ago_L)} \sum_{\pi \in \Pi_{\chi}(M)} \int_{i\ago_M^{L,*}} \Jc_{L,\pi}(x,f,\la)\, d\la.
\end{align*}
To conclude we have to show that the expression above is absolutely convergent. Let $J\subset K_f$ be a level such that $f$ is $J$-biinvariant. Then the terms attached to $\pi\in \Pi_{\disc}(M)=\bigcup_{\chi\in \Xgo(G)}\Pi_{\chi}(M)$ vanish unless $\pi\in \Pi_{\disc}(M)^J$. Thus by proposition \ref{prop:relch-JL} we are reduced to the following statement: for all standard Levi subgroup $M$ and for large $q>0$ we have
\begin{align*}
  \sum_{\pi \in \Pi_{\disc}(M)^J} \frac1{(1+\La_\pi^2)^q} <\infty.
\end{align*}
This is due to Müller \cite[line (6.17) p. 711 and below]{Muller02}.
  \end{preuve}
  \end{paragr}

\bibliography{biblio}

\begin{thebibliography}{BPLZZ21}

\bibitem[Art78]{ar1}
J.~Arthur.
\newblock A trace formula for reductive groups {I}. {T}erms associated to
  classes in {$G(\mathbb{Q})$}.
\newblock {\em Duke Math. J.}, 45:911--952, 1978.

\bibitem[Art80]{ar2}
J.~Arthur.
\newblock A trace formula for reductive groups {II}.
\newblock {\em Comp. Math.}, 40:87--121, 1980.

\bibitem[Art81]{ar-inv}
J.~Arthur.
\newblock The trace formula in invariant form.
\newblock {\em Ann. of Math. (2)}, 114(1):1--74, 1981.

\bibitem[Art82]{ar-truncated}
J.~Arthur.
\newblock On the inner product of truncated {E}isenstein series.
\newblock {\em Duke Math. J.}, 49(1):35--70, 1982.

\bibitem[Art83]{Ar-PW}
J.~Arthur.
\newblock A {P}aley-{W}iener theorem for real reductive groups.
\newblock {\em Acta Math.}, 150(1-2):1--89, 1983.

\bibitem[Art05]{Ar-cours}
J.~Arthur.
\newblock An introduction to the trace formula.
\newblock In {\em Harmonic analysis, the trace formula, and {S}himura
  varieties}, volume~4 of {\em Clay Math. Proc.}, pages 1--263. Amer. Math.
  Soc., Providence, RI, 2005.

\bibitem[Ber88]{Ber}
J.~Bernstein.
\newblock On the support of {P}lancherel measure.
\newblock {\em J. Geom. Phys.}, 5(4):663--710 (1989), 1988.

\bibitem[BK14]{BK}
J.~Bernstein and B.~Kr\"{o}tz.
\newblock Smooth {F}r\'{e}chet globalizations of {H}arish-{C}handra modules.
\newblock {\em Israel J. Math.}, 199(1):45--111, 2014.

\bibitem[BL24]{BL-mero}
J.~Bernstein and E.~Lapid.
\newblock On the meromorphic continuation of {E}isenstein series.
\newblock {\em J. Amer. Math. Soc.}, 37(1):187--234, 2024.

\bibitem[Bou68]{Bki-Lie}
N.~Bourbaki.
\newblock {\em \'{E}l\'{e}ments de math\'{e}matique. {F}asc. {XXXIV}. {G}roupes
  et alg\`ebres de {L}ie. {C}hapitre {IV}: {G}roupes de {C}oxeter et syst\`emes
  de {T}its. {C}hapitre {V}: {G}roupes engendr\'{e}s par des r\'{e}flexions.
  {C}hapitre {VI}: syst\`emes de racines}.
\newblock Actualit\'{e}s Scientifiques et Industrielles [Current Scientific and
  Industrial Topics], No. 1337. Hermann, Paris, 1968.

\bibitem[BPC23]{BPC}
R.~Beuzart-Plessis and P.-H. Chaudouard.
\newblock The global {G}an-{G}ross-{P}rasad conjecture for unitary groups.
  {II}. {F}rom {E}isenstein series to {B}essel periods, Forum Pi, to appear,
  preprint 2023.

\bibitem[BPCZ22]{BCZ}
R.~Beuzart-Plessis, P.-H. Chaudouard, and M.~Zydor.
\newblock The global {G}an-{G}ross-{P}rasad conjecture for unitary groups: the
  endoscopic case.
\newblock {\em Pub. Math. IHES}, pages 1--154, 2022.
\newblock en ligne sur https://www.springer.com/journal/10240.

\bibitem[BPLZZ21]{BPLZZ}
R.~Beuzart-Plessis, Y.~Liu, W.~Zhang, and X.~Zhu.
\newblock Isolation of cuspidal spectrum, with application to the
  {G}an-{G}ross-{P}rasad conjecture.
\newblock {\em Ann. of Math. (2)}, 194(2):519--584, 2021.

\bibitem[FL17]{FL-Iran}
T.~Finis and E.~Lapid.
\newblock On the analytic properties of intertwining operators {I}: global
  normalizing factors.
\newblock {\em Bull. Iranian Math. Soc.}, 43(4):235--277, 2017.

\bibitem[Fli88]{Flicker}
Y.~Flicker.
\newblock Twisted tensors and {E}uler products.
\newblock {\em Bull. Soc. Math. France}, 116(3):295--313, 1988.

\bibitem[FLM12]{FLM-degree}
T.~Finis, E.~Lapid, and W.~M\"{u}ller.
\newblock On the degrees of matrix coefficients of intertwining operators.
\newblock {\em Pacific J. Math.}, 260(2):433--456, 2012.

\bibitem[FLM15]{FLM-limit-mult}
T.~Finis, E.~Lapid, and W.~M\"{u}ller.
\newblock Limit multiplicities for principal congruence subgroups of {${\rm
  GL}(n)$} and {${\rm SL}(n)$}.
\newblock {\em J. Inst. Math. Jussieu}, 14(3):589--638, 2015.

\bibitem[Jac97]{Jac-Edin}
H.~Jacquet.
\newblock Automorphic spectrum of symmetric spaces.
\newblock In {\em Representation theory and automorphic forms ({E}dinburgh,
  1996)}, volume~61 of {\em Proc. Sympos. Pure Math.}, pages 443--455. Amer.
  Math. Soc., Providence, RI, 1997.

\bibitem[JLR99]{JLR}
H.~Jacquet, E.~Lapid, and J.~Rogawski.
\newblock Periods of automorphic forms.
\newblock {\em J. Amer. Math. Soc.}, 12(1):173--240, 1999.

\bibitem[JPSS83]{JPSS}
H.~Jacquet, I.~I. Piatetskii-Shapiro, and J.~A. Shalika.
\newblock Rankin-{S}elberg convolutions.
\newblock {\em Amer. J. Math.}, 105(2):367--464, 1983.

\bibitem[Kra95]{Kraski}
W.~Kra\'skiewicz.
\newblock Reduced decompositions in {W}eyl groups.
\newblock {\em European J. Combin.}, 16(3):293--313, 1995.

\bibitem[Lap06]{LapFRTF}
E.~Lapid.
\newblock On the fine spectral expansion of {J}acquet's relative trace formula.
\newblock {\em J. Inst. Math. Jussieu}, 5(2):263--308, 2006.

\bibitem[Lap08]{Lap-remark}
E.~Lapid.
\newblock A remark on {E}isenstein series.
\newblock In {\em Eisenstein series and applications}, volume 258 of {\em
  Progr. Math.}, pages 239--249. Birkh\"{a}user Boston, Boston, MA, 2008.

\bibitem[Lap11a]{Linner}
E.~Lapid.
\newblock On {A}rthur's asymptotic inner product formula of truncated
  {E}isenstein series.
\newblock In {\em On certain {$L$}-functions}, volume~13 of {\em Clay Math.
  Proc.}, pages 309--331. Amer. Math. Soc., Providence, RI, 2011.

\bibitem[Lap11b]{Lap-asym}
E.~Lapid.
\newblock On {A}rthur's asymptotic inner product formula of truncated
  {E}isenstein series.
\newblock In {\em On certain {$L$}-functions}, volume~13 of {\em Clay Math.
  Proc.}, pages 309--331. Amer. Math. Soc., Providence, RI, 2011.

\bibitem[Lap13]{LapHC}
E.~Lapid.
\newblock On the {H}arish-{C}handra {S}chwartz space of {$G(F)\backslash G(\Bbb
  A)$}.
\newblock In {\em Automorphic representations and {$L$}-functions}, volume~22
  of {\em Tata Inst. Fundam. Res. Stud. Math.}, pages 335--377. Tata Inst.
  Fund. Res., Mumbai, 2013.
\newblock With an appendix by Farrell Brumley.

\bibitem[LR03]{LR}
E.~Lapid and J.~Rogawski.
\newblock Periods of {E}isenstein series: the {G}alois case.
\newblock {\em Duke Math. J.}, 120(1):153--226, 2003.

\bibitem[LRS99]{LRS}
W.~Luo, Z.~Rudnick, and P.~Sarnak.
\newblock On the generalized {R}amanujan conjecture for {${\rm GL}(n)$}.
\newblock In {\em Automorphic forms, automorphic representations, and
  arithmetic ({F}ort {W}orth, {TX}, 1996)}, volume 66, Part 2 of {\em Proc.
  Sympos. Pure Math.}, pages 301--310. Amer. Math. Soc., Providence, RI, 1999.

\bibitem[LW13]{LabWal}
J.-P. Labesse and J.-L. Waldspurger.
\newblock {\em La formule des traces tordue d'apr\`es le {F}riday {M}orning
  {S}eminar}, volume~31 of {\em CRM Monograph Series}.
\newblock American Mathematical Society, Providence, RI, 2013.
\newblock With a foreword by Robert Langlands [dual English/French text].

\bibitem[Mic22]{Michel}
P.~Michel.
\newblock Recent progresses on the subconvexity problem.
\newblock {\em Ast\'{e}risque}, (438, S\'{e}minaire Bourbaki. Vol. 2021/2022.
  Expos\'{e}s 1181--1196):353--401, 2022.

\bibitem[MS04]{MuSpeh}
W.~M\"{u}ller and B.~Speh.
\newblock Absolute convergence of the spectral side of the {A}rthur trace
  formula for {${\rm GL}_n$}.
\newblock {\em Geom. Funct. Anal.}, 14(1):58--93, 2004.
\newblock With an appendix by E. M. Lapid.

\bibitem[MW89]{MW}
C.~M{\oe}glin and J.-L. Waldspurger.
\newblock Le spectre r\'{e}siduel de {${\rm GL}(n)$}.
\newblock {\em Ann. Sci. \'{E}cole Norm. Sup. (4)}, 22(4):605--674, 1989.

\bibitem[MW94]{MWlivre}
C.~M{\oe}glin and J.-L. Waldspurger.
\newblock {\em D\'{e}composition spectrale et s\'{e}ries d'{E}isenstein},
  volume 113 of {\em Progress in Mathematics}.
\newblock Birkh\"{a}user Verlag, Basel, 1994.
\newblock Une paraphrase de l'\'{E}criture.

\bibitem[Mül98]{Muller98}
W.~Müller.
\newblock The trace class conjecture in the theory of automorphic forms. {II}.
\newblock {\em Geom. Funct. Anal.}, 8(2):315--355, 1998.

\bibitem[Mül02]{Muller02}
W.~Müller.
\newblock On the spectral side of the {A}rthur trace formula.
\newblock {\em Geom. Funct. Anal.}, 12(4):669--722, 2002.

\bibitem[Mül07]{Muller07}
W.~Müller.
\newblock Weyl's law for the cuspidal spectrum of {${\rm SL}_n$}.
\newblock {\em Ann. of Math. (2)}, 165(1):275--333, 2007.

\bibitem[Off06]{Offen-symplectic}
O.~Offen.
\newblock On symplectic periods of the discrete spectrum of {${\rm GL}_{2n}$}.
\newblock {\em Israel J. Math.}, 154:253--298, 2006.

\bibitem[Ram15]{Ram1}
D.~Ramakrishnan.
\newblock A mild {T}chebotarev theorem for {${\rm GL}(n)$}.
\newblock {\em J. Number Theory}, 146:519--533, 2015.

\bibitem[{Ram}18]{Ram2}
D.~{Ramakrishnan}.
\newblock {A Theorem on GL(n) \`a la Tchebotarev}.
\newblock {\em arXiv e-prints}, page arXiv:1806.08429, Jun 2018.

\bibitem[RS96]{RuSar}
Z.~Rudnick and P.~Sarnak.
\newblock Zeros of principal {$L$}-functions and random matrix theory.
\newblock volume~81, pages 269--322. 1996.
\newblock A celebration of John F. Nash, Jr.

\bibitem[Wal03]{Walds-Plan}
J.-L. Waldspurger.
\newblock La formule de {P}lancherel pour les groupes {$p$}-adiques (d'apr\`es
  {H}arish-{C}handra).
\newblock {\em J. Inst. Math. Jussieu}, 2(2):235--333, 2003.

\bibitem[Yam14]{Ysymp}
S.~Yamana.
\newblock Symplectic periods of the continuous spectrum of {${\rm GL}(2n)$}.
\newblock {\em Ann. Inst. Fourier (Grenoble)}, 64(4):1561--1580, 2014.

\bibitem[Yam15]{Yquad}
S.~Yamana.
\newblock Periods of residual automorphic forms.
\newblock {\em J. Funct. Anal.}, 268(5):1078--1104, 2015.

\bibitem[Zha14]{Zhang2}
W.~Zhang.
\newblock Automorphic period and the central value of {R}ankin-{S}elberg
  {$L$}-function.
\newblock {\em J. Amer. Math. Soc.}, 27:541--612, 2014.

\bibitem[Zyd20]{Z3}
M.~Zydor.
\newblock Les formules des traces relatives de {J}acquet--{R}allis
  grossi\`eres.
\newblock {\em J. Reine Angew. Math.}, 762:195--259, 2020.

\end{thebibliography}
\bibliographystyle{alpha}

\begin{flushleft}
Pierre-Henri Chaudouard \\
\medskip
Université Paris Cité \\
IMJ-PRG \\
Bâtiment Sophie Germain\\
8 place Aurélie Nemours\\
F-75013 Paris CEDEX 13 \\
France\\
\medskip
Institut universitaire de France (IUF)\\
\medskip
email:\\
Pierre-Henri.Chaudouard@imj-prg.fr \\
\end{flushleft}

\end{document}